\documentclass[a4paper,twoside,12pt]{book}

\usepackage{ifthen}

\newboolean{smallbook}
\newboolean{PSdiagram}
\newboolean{hyper}
\newboolean{debug}

\setboolean{smallbook}{false} 
\setboolean{PSdiagram}{false} 
\setboolean{hyper}{false} 
\setboolean{debug}{true}

\usepackage{mathrsfs}
\usepackage{amssymb}
\usepackage{amsmath}
\usepackage{dsfont}
\usepackage{stmaryrd}
\usepackage{exscale}
\usepackage{extraipa}
\usepackage{tocloft}
\usepackage[dvips]{graphicx}
\usepackage{makeidx}

\usepackage[latin9]{inputenc}
\usepackage{ae}

\usepackage{a4}
\usepackage{fancyhdr}
\usepackage[compact,nobottomtitles,rigidchapters]{titlesec}

\usepackage{theorem}
\usepackage{diagrams}

\makeindex

\ifthenelse{\boolean{hyper}}
    {\usepackage[pdftex]{hyperref}}
    {}

\ifthenelse{\boolean{smallbook}}
   {\setboolean{PSdiagram}{true}}
   {}

\ifthenelse{\boolean{PSdiagram}}
  {\setboolean{debug}{false}}
  {}

\ifthenelse{\boolean{debug}}
   { \usepackage{srcltx}
   }
   {}

\titleformat{\part}[block]{\Huge\bf}{\thepart. }{0pt}{}

\titleformat{\chapter}[block]{\LARGE\bf}{\thechapter. }{0.5em}{}

\titlespacing*{\chapter}{0in}{-2ex}{2in}
\titlespacing*{\section}{0in}{*4}{*2}
\titlespacing*{\subsection}{0in}{*4}{*2}

\let\oldsection\section
\let\oldsubsection\subsection

\newcommand\mysection[1]{\renewcommand{\bottomtitlespace}{0.42\textheight}\oldsection{#1}}
\newcommand\mysubsection[1]{\renewcommand{\bottomtitlespace}{0.33\textheight}\oldsubsection{#1}}

\let\section\mysection
\let\subsection\mysubsection

\fancypagestyle{plain}{
  \fancyhf{}
}

\ifthenelse{\boolean{smallbook}}
  {
   \setlength{\topmargin}{0cm}
   \setlength{\oddsidemargin}{-0.5cm} %
   \setlength{\evensidemargin}{3.8cm} 
  }
  {}

\renewcommand\textrm[1]{{\rm #1}}

\newlength{\listindent}
\newcommand\romanenum{\renewcommand{\labelenumi}{\textup{(}\roman{enumi}\textup{)}}}

\newcommand\enumroman[1]{{\rm(\/\textit{\roman{#1}}\/)}}
\newcommand\enumarabic[1]{\arabic{#1}.}
\newcounter{enumcounter}

\newenvironment{enumsrm}
   {\begin{list}{\enumroman{enumiv}}{\usecounter{enumiv}
                      \setlength{\topsep}{0.1ex plus0.1ex minus0.1ex}
                      \setlength{\itemsep}{0.1ex plus0.1ex minus0.1ex}
                      \setlength{\parsep}{0.1ex plus0.1ex minus0.1ex}
                      \setlength{\itemindent}{0in}
                      \setlength{\labelwidth}{2.5em}
                      \setlength{\leftmargin}{2\listindent}
                      \setlength{\partopsep}{\medskipamount} } }
   {\end{list}}

\ifthenelse{\boolean{PSdiagram}} %
 {
  \diagramstyle[height=0.3in,width=0.4in,labelstyle=\scriptstyle,dpi=1200,nohug,PostScript=dvips]
 }{
  \diagramstyle[height=0.3in,width=0.4in,labelstyle=\scriptstyle,dpi=600,nohug,noPostScript]
 }

\newarrowtail{<=}\Leftarrow\Rightarrow\Uparrow\Downarrow

\newarrow{Corresponds} <--->
\newarrow{Equal} =====
\newarrow{Equiv} {<=}==={=>}
\newarrow{Dash} {}{dash}{}{dash}{}
\newarrow{Dashinto} C{dash}{}{dash}>
\newarrow{Dashonto} {}{dash}{}{dash}{>>}
\newarrow{Dotmapsto} {mapsto}...>

\newarrow{CTo} ----{->}
\newarrow{CMapsto} |---{->}
\newarrow{CInto} C---{->}
\newarrow{COnto} ----{->>}

\newcommand\diagramiso{\hbox{\lower1.5ex\hbox{$_\sim$}}}

\theoremstyle{change} %
\theorembodyfont{\it}

\newtheorem{defi}{Definition.}[subsection]
\newtheorem{defi*}{Definition.}

\newtheorem{lem}[defi]{Lemma.}
\newtheorem{prop}[defi]{Proposition.}
\newtheorem{theorem}[defi]{Theorem.}
\newtheorem{cor}[defi]{Corollary.}
\newtheorem{remark}[defi]{Remark.}

\newtheorem{example}[defi]{Example.}
\DeclareMathAlphabet{\mathpzc}{OT1}{pzc}{m}{it}
                                 
\newenvironment{interlude}{\begin{trivlist}\setlength\topsep{\theorempreskipamount}\item}{\end{trivlist}}
\newenvironment{smallinterlude}{\begin{trivlist}\setlength\topsep{\theorempreskipamount}\item\small}{\normalsize\end{trivlist}}

\newenvironment{remarks*}{\begin{interlude}\pagebreak[2]\noindent\textbf{Bemerkungen.\ }\nopagebreak}{\end{interlude}}

\newenvironment{examples*}{\begin{interlude}\pagebreak[2]\noindent\textbf{Beispiele.\ }\nopagebreak}{\end{interlude}}

\newenvironment{proof}{\begin{interlude}\noindent\textit{Proof.}}{\qed\end{interlude}}

\newenvironment{proof*}{\begin{interlude}\noindent\noindent\textit{Beweis.}}{\end{interlude}}

\newcommand\qedsymbol{$\square$}

\newcommand\qed{\nopagebreak\penalty10000\hbox{}\nobreak\penalty10000\nopagebreak\penalty10000\hfill\hbox{\qedsymbol}\par}

\def\phi{\varphi}
\def\theta{\vartheta}
\def\epsilon{\varepsilon}
\def\sh{^\#}

\DeclareMathOperator{\id}{id} %
\DeclareMathOperator{\im}{im} %
\DeclareMathOperator{\coker}{coker} %
\DeclareMathOperator{\preim}{pre\kern0.13em im} %
\DeclareMathOperator{\preker}{pre\kern0.13em ker} %
\DeclareMathOperator{\precoker}{pre\kern0.13em coker} %
\DeclareMathOperator{\rad}{rad} %

\newcommand\set[2][auto]{
     \ifthenelse{\equal{#1}{auto}}{\left\lbrace}{\csname #1\endcsname\lbrace} #2 \ifthenelse{\equal{#1}{auto}}{\right\rbrace}{\csname #1\endcsname\rbrace} }

\let\oplus\varoplus

\newcommand{\indlim}{\operatorname*{\underrightarrow{\rm lim}}} %

\newcommand\isoto{\mbox{$\hspace{7.5pt}\raise 3pt\hbox{$\sim$}\hspace{-17pt}\longrightarrow\hspace{3pt}$}\linebreak[0]}
\newcommand\isoot{\mbox{$\hspace{8.5pt}\raise 3pt\hbox{$\sim$}\hspace{-18pt}\longleftarrow\hspace{3pt}$}\linebreak[0]}

\newcommand{\FF}{{\rm F}} %

\newenvironment{packed_enum}{
\begin{enumerate}
  \setlength{\itemsep}{1pt}
  \setlength{\parskip}{0pt}
  \setlength{\parsep}{0pt}
}{\end{enumerate}}

\title{\textsc{\Large Uniformly rigid spaces and Néron models of formally finite type}}

\author{Christian Kappen}
\date{30. 6. 2005}

\def\romanenum{\renewcommand{\labelenumi}{\textup{(}\roman{enumi}\textup{)}}}

\def\sh{\textup{sh}}

\def\b{\mathfrak{b}}

\def\Res{\textup{Res}}

\def\red{\textup{red}}
\def\Bl{\textup{Bl}}

\def\Q{\mathbb{Q}}

\def\Spec{\textup{Spec }}
\def\Spf{\textup{Spf }}
\def\Max{\textup{Max }}

\def\ker{\textup{ker }}
\def\im{\textup{im }}
\def\coker{\textup{coker }}

\def\Hom{\textup{Hom}}

\def\A{\mathbb{A}}
\def\O{\mathcal{O}}

\def\J{\mathcal{J}}
\def\fI{\mathfrak{I}}

\def\Z{\mathbb{Z}}
\def\N{\mathbb{N}}

\def\G{\mathbb{G}}

\def\m{\mathfrak{m}}
\def\a{\mathfrak{a}}

\def\q{\mathfrak{q}}
\def\n{\mathfrak{n}}

\def\plim{\varprojlim}

\def\P{\mathbb{P}}
\def\p{\mathfrak{p}}

\def\id{\textup{id}}
\def\sep{{\textup{\footnotesize sep}}}
\def\Gal{\textup{Gal}}

\def\Lie{\textup{Lie}}

\def\length{\textup{length}}

\def\an{\textup{an}}
\def\rig{\textup{rig}}
\def\Rig{\textup{Rig}}
\def\FF{\textup{FF}}
\def\sRig{\textup{uRig}}
\def\uRig{\textup{uRig}}
\def\D{\mathbb{D}}
\def\srig{\textup{urig}}
\def\urig{\textup{urig}}

\def\Sp{\textup{Sp\,}}
\def\sAff{\textup{sAff}}

\def\sSp{\textup{sSp\,}}
\def\sp{\textup{sp}}
\def\sI{\mathcal{I}}
\def\sJ{\mathcal{J}}
\def\sK{\mathcal{K}}

\def\sT{\mathcal{T}}

\def\B{\mathbb{B}}
\def\fU{\mathfrak{U}}
\def\fV{\mathfrak{V}}
\def\fW{\mathfrak{W}}

\def\fG{\mathfrak{G}}
\def\fH{\mathfrak{H}}
\def\fX{\mathfrak{X}}
\def\fY{\mathfrak{Y}}
\def\fZ{\mathfrak{Z}}
\def\fS{\mathfrak{S}}
\def\fT{\mathfrak{T}}
\newcommand{\ul}[1]{\underline{#1}}

\def\cl{\textup{cl}}
\def\rad{\textup{rad}}
\def\F{\mathbb{F}}

\def\sup{\textup{sup}}

\def\sF{\mathcal{F}}
\def\sG{\mathcal{G}}
\def\Jac{\textup{Jac}\,}
\def\Gauss{\textup{Gauss}\,}
\def\M{\mathcal{M}}

\def\fC{\mathfrak{C}}
\def\parent{\textup{par}}
\def\children{\textup{ch}}
\def\leaves{\textup{lv}}
\def\subtree{\textup{subt}}

\def\Frac{\textup{Frac}}

\def\inner{\textup{inn}}
\def\sr{\textup{ur}}
\def\ur{\textup{ur}}
\def\r{\textup{r}}
\def\fr{\mathfrak{r}}

\def\Zar{\textup{Zar}}

\def\sa{\textup{sa}}
\def\aux{\textup{aux}}

\def\fD{\mathfrak{D}}
\def\comp{\textup{comp}}
\def\Spa{\textup{Spa}}
\def\ad{\textup{ad}}
\def\sCon{\textup{sCon}}
\def\Tori{\textup{Tori}}
\def\End{\textup{End}}
\def\sG{\mathscr{G}}

\def\bA{\textup{bA}}
\def\Ind{\textup{Ind}}
\def\FS{\textup{FS}}

\begin{document}

\frontmatter

\begin{titlepage}
\begin{center}
\thispagestyle{empty}
\bigskip
{\Large \textbf{Mathematik}}

\vspace{3.5cm}

\huge{\textbf{Uniformly rigid spaces and\\ Néron models of formally finite type}}

\vspace{3cm}

\large{Inaugural-Dissertation\\
zur Erlangung des Doktorgrades\\
der Naturwissenschaften im Fachbereich\\
Mathematik und Informatik\\
der Mathematisch-Naturwissenschaftlichen Fakultät\\
der Westfälischen Wilhelms-Universität Münster
}

\vspace{2.5cm}
\normalsize
vorgelegt von\\
Christian Kappen\\
aus Münster\\
-- 2009 --
\end{center}
\end{titlepage}
\normalsize
\newpage
\thispagestyle{empty}

\vspace*{15cm}
\linespread{1.6}
\begin{table}[b]
\begin{tabular}{ll}
Dekan:&Prof.\ Dr.\ Dr.\ h.c.\ Joachim Cuntz\\
Erster Gutachter:&Prof.\ Dr.\ Siegfried Bosch\\
Zweiter Gutachter:&Prof.\ Dr.\ Peter Schneider\\
Tag der mündlichen Prüfung:&\rule{5cm}{0.03cm}\\
Tag der Promotion:&\rule{5cm}{0.03cm}
\end{tabular}
\end{table}

\linespread{1.0}

%


\chapter{Introduction}


\textbf{\large Motivation}
\vspace{0.5cm}

Néron models are important objects in arithmetic algebraic geometry. Let $K$ be the fraction field of a discrete valuation ring $R$. If $X$ is a smooth and separated $K$-scheme, then the \emph{Néron lft-model} $\ul{X}$ of $X$ is a smooth and separated $R$-model of $X$ with the universal property that for every smooth $R$-scheme $\ul{Z}$, any morphism from the generic fiber of $\ul{Z}$ to $X$ extends uniquely to a morphism from $\ul{Z}$ to $\ul{X}$. A Néron lft-model is called a Néron model if it is quasi-compact. Every abelian $K$-variety $A$ admits a Néron model, cf.\ \cite{BLR} 1.3/2. If $A$ has semi-abelian reduction in the sense of \cite{BLR} 7.4, then the identity component of its Néron model is preserved under base change, cf.\ \cite{BLR} 7.4/4. The situation is more complicated in the unstable case. Chai has introduced the so-called \emph{base change conductor}, a rational number which provides a measure for the defect of stability of $A$ over $K$. It measures how the volume form of the Néron model of $A$ changes under a stabilizing field extension.

Let us assume that $K$ is complete with perfect residue field of positive characteristic $p$. In their article \cite{ChaiYudeSh}, Chai, Yu and de Shalit prove that the base change conductor of an algebraic $K$-torus $T$ is an isogeny invariant, and they state an explicit formula for the base change conductor of $T$ in terms of the $\Q$-rational Galois representation attached to $T$. In \cite{Chai}, Chai shows that the base change conductor of an abelian $K$-variety $A$ with potentially \emph{multiplicative} reduction is likewise an isogeny invariant, which fails to hold for general abelian $K$-varieties. Motivated by this result, he conjectures that there exists a formula for the base change conductor of $A$ in terms of the Galois action on the $\Q_p$-rational character group of the formal torus that is obtained from the stable Néron model of $A$ via formal completion along its identity section. 

In his unpublished preprint \cite{ChaiPre}, Chai states such a formula, and he describes an ingenious method for deriving it, by means of an analytic localization process, from the formula for algebraic tori. His arguments depend crucially on a notion of \emph{formal Néron models of formally finite type} which he introduces in \cite{ChaiPre}, Section 5. Chai does not develop the general theory of these formal Néron models; he only indicates in an ad hoc way what is needed to deduce the formula for the base change conductor, using nontrivial input such as the Néron-Popescu Desingularization Theorem (cf.\ \cite{Sw}). Chai's methods prove to be unsuitable in non-affine situations, which is essentially due to the existence of unbounded global functions on the generic fibers of formal $R$-schemes of formally finite type. 

In this thesis, we lay the foundation for the general theory of Néron models of formally finite type, thereby providing a
conceptual approach to Chai's applications. In doing so, we modify Chai's original definition of formal Néron models of formally finite type, so as to avoid unbounded morphisms of generic fibers. Our methods do not depend on the Néron-Popescu Theorem.

Let us describe the content of this thesis in more detail.

\vspace{0.7cm}
\textbf{\large Formal Néron models}
\vspace{0.5cm}

Let $R$ be a complete discrete valuation ring with fraction field $K$ and residue field $k$, and let $\pi\in R$ be a uniformizer. 

In their article \cite{BS}, Bosch and Schlöter have introduced and studied formal Néron models of smooth rigid $K$-spaces. The \emph{formal Néron model} $\fX$ of a smooth rigid $K$-space $X$ is the universal smooth formal $R$-scheme of locally \emph{topologically finite} (tf) type equipped with a morphism $\phi$ from its generic fiber $\fX^\rig$ to $X$. In contrast with the algebraic situation, $\fX$ is not necessarily a model of $X$: The universal morphism $\phi$ needs not be an isomorphism. In all interesting cases, $\phi$ is an open immersion, or at least $\phi$ is injective and an open immersion on every quasi-compact subspace of its domain. The formal Néron model of $X$ exists in many important cases; for example, it exists when $X$ is a smooth rigid $K$-group whose group of unramified points is bounded or when $X$ is the analytification of a smooth $K$-group scheme admitting a Néron lft-model, cf.\ \cite{BS} 1.2, 6.2 and 6.6.

Formal $R$-schemes of locally tf type are locally the formal spectra of quotients of \emph{strictly convergent power series} rings in finitely many variables over $R$
\[
R\langle T_1,\ldots,T_m\rangle
\]
by ideals of strictly convergent power series; the topology on these $R$-algebras is the $\pi$-adic one. Formal $R$-schemes of locally tf type have been investigated by Raynaud and others also in the case where the valuation on $R$ is not discrete. In our setting, $R$ is discretely valued, and all formal $R$-schemes of interest to us will be locally noetherian; we refer to \cite{EGAIn} Section 10 for the basic theory of locally noetherian formal schemes. As we have already mentioned above, Chai has introduced an alternative notion of formal Néron models for rigid spaces by allowing formal $R$-schemes of locally \emph{formally finite} (ff) type; these are locally the formal spectra of quotients of \emph{mixed power series} rings in finitely many variables over $R$
\[
R[[S_1,\ldots,S_n]]\langle T_1,\ldots,T_m\rangle
\]
by ideals of mixed power series, where the topology on $R[[S_1,\ldots,S_n]]$ is defined by the ideal $(\pi,S_1,\ldots,S_n)$; we refer to Section \ref{locfftypesec} of this thesis for a more detailed discussion of these rings. The rigid-analytic generic fibers of formal $R$-schemes of locally ff type have been defined and studied by Rapoport, Berthelot and de Jong in \cite{RZ} Section 5.5, \cite{Berth} Section 0.2 and \cite{dJ} Section 7 respectively; for a synopsis we refer to Section \ref{berthsec} of this treatise. One commonly says that the rigid generic fiber $\fX^\rig$ of a formal $R$-scheme of locally ff type $\fX$ is obtained from $\fX$ by means of \emph{Berthelot's construction}. An instructive example is provided by the rigid open unit disc, which is the generic fiber $\D^1_K=\mathbb{D}^\textup{rig}$ of the affine formal $R$-scheme $\mathbb{D}\mathrel{\mathop:}=\textup{Spf } R[[S]]$ and which is not quasi-compact. 

In \cite{ChaiPre} 5.2, Chai suggests the following definition: The \emph{formal N\'eron model of locally ff type} of a smooth rigid $K$-space $X$ is the universal smooth formal $R$-scheme of locally ff type equipped with a morphism from its rigid generic fiber $\fX^\rig$ to $X$. We refer to Section \ref{smoothnesssec} for a discussion of the notion of smoothness for morphisms of ff type. Chai's definition provides a workable notion in affine situations:

Using Theorem \cite{dJ}  7.4.1 concerning bounded functions on rigid generic fibers of normal formal $R$-schemes of locally ff type, one easily shows that every \emph{affine} smooth formal $R$-scheme of ff type is the formal Néron model of locally ff type of its rigid generic fiber. The connected components of the Néron lft-model $\ul{\mathscr{T}}$ of a split $K$-torus $\mathscr{T}$ are \emph{affine}; therefore an elementary calculation shows that the $\pi$-adic completion of $\ul{\mathscr{T}}$ is the formal Néron model of locally ff type of $\mathscr{T}^\an$: While being of locally tf type over $R$, it satisfies the Néron extension property for all smooth formal $R$-schemes of locally ff type.  Using Néron-Popescu desingularization (cf.\ \cite{Sw}) in conjunction with the universal property of algebraic Néron models, Chai establishes the existence of formal Néron models of locally ff type for rigid spaces which, after some finite Galois base field extension, are generic fibers of completions of algebraic Néron lft-models along \emph{affine} closed formal subschemes of their special fibers, cf.\ \cite{ChaiPre} Section 5. 

These methods do not apply in \emph{non-affine} situations; in particular, it is unclear under which general conditions a formal Néron model in the sense of \cite{BS} satisfies the stronger mapping property of the formal Néron model of locally ff type. The main problem lies in the fact that Raynaud's theory of formal models, as it is comprehensively documented in \cite{BL1} and \cite{FRG}, does not carry over to the ff type setting. For instance, an \emph{unbounded} function on the open rigid unit disc induces a morphism to the rigid projective line which does not extend to a morphism $\fD\rightarrow\P^{1,\wedge}_R$, where $\fD$ is an $R$-model of ff type for $\D^1_K$ and where $\P^{1,\wedge}_R$ denotes the completion of $\P^1_R$ along its special fiber.


We overcome this fundamental problem by introducing a new category of non-archimedean analytic $K$-spaces: The category of \textit{uniformly rigid $K$-spaces}. Using this new category, we can define a notion of formal Néron models of formally finite type which is workable also in non-affine situations. 

A uniformly rigid $K$-space may be considered as a rigid $K$-space equipped with a \emph{more rigid} G-topology and with a sheaf of functions that are \emph{bounded} on rather big admissible open subsets; these subsets may for example be defined by non-strict or \emph{strict} inequalities. It is sometimes instructive to think of a uniformly rigid $K$-space as a rigid $K$-space equipped with additional rigidifying structure. Morphisms of uniformly rigid $K$-spaces may be viewed as morphisms of underlying rigid $K$-spaces satisfying certain \emph{boundedness conditions}; under suitable separatedness assumptions, they admit descriptions in terms of formal $R$-models of ff type, cf.\ Corollary \ref{graphpropercor}.

A formal $R$-scheme of locally ff type induces a uniformly rigid structure on its rigid generic fiber: The Berthelot generic fiber functor $\rig$ admits a factorization
\[
\begin{diagram}
\textup{FS}_R&&\rTo^{\rig}&&\Rig_K&\\
&\rdTo<\srig&&\ruTo>\r&&\\
&&\sRig_K&&&\quad,
\end{diagram}
\]
where $\FS_R$, $\Rig_K$ and $\sRig_K$ denote the categories of formal $R$-schemes of locally ff type, of rigid $K$-spaces and of uniformly rigid $K$-spaces respectively. Berthelot's construction is effectively carried out by the functor $\r$, while the uniformizer $\pi$ of $R$ is already a unit in the structural sheaves of the objects in $\uRig_K$. We postulate that the generic fiber of a formal $R$-scheme of locally ff type should naturally be viewed as a \emph{uniformly} rigid $K$-space. While yielding a locally G-ringed $K$-space, the functor $\urig$ preserves more integral structure than Berthlot's functor $\rig$.

Raynaud's theory of formal models of locally tf type induces a canonical uniformly rigid structure on every quasi-paracompact and quasi-separated rigid $K$-space. For example, the rigid open unit disc $\D^1_K$ carries at least two canonical uniformly rigid structures: one induced by a quasi-paracompact $R$-model $\fD$ of locally tf type in the sense of Raynaud and another one induced by the affine $R$-model $\D=\Spf R[[T]]$ of ff type. The uniformly rigid $K$-space $\D^\urig$ is quasi-compact in its uniformly rigid G-topology, while $\fD^\urig$, just like the rigid open unit disc itself, is not quasi-compact. If $X$ is a quasi-paracompact and quasi-separated rigid $K$-space, we let $X^\ur$ denote the uniformly rigid $K$-space associated to a quasi-paracompact formal $R$-model of locally tf type for $X$. By \cite{BL1} 4.1, $\ur$ defines a functor, cf. also Section \ref{compfuncsec}.


A uniformly rigid $K$-space will be called \emph{smooth} if its underlying rigid $K$-space is smooth. We define the \emph{formal Néron model} of a smooth uniformly rigid $K$-space $X$ to be the universal smooth formal $R$-scheme of locally ff type $\fX$ equipped with a morphism of uniformly rigid $K$-spaces $\fX^\srig\rightarrow X$. 

Formal Néron models of uniformly rigid spaces can be used to derive Chai's formula for the base change conductor, and they are flexible enough to deal with non-affine situations.

\vspace{0.7cm}
\textbf{\large Uniformly rigid spaces}
\vspace{0.5cm}

In the first chapter of this thesis, we systematically develop the theory of uniformly rigid $K$-spaces. According to Definition \ref{srigspacedefi}, a \emph{uniformly rigid $K$-space} is a saturated locally G-ringed $K$-space admitting an admissible \emph{semi-affinoid} covering: A uniformly rigid $K$-space is composed of semi-affinoid $K$-spaces just like a rigid $K$-space is composed of affinoid $K$-spaces. Let us recall that a locally G-ringed $K$-space is called \emph{saturated} if admissibility of subsets and coverings can be checked locally, cf.\ p. 339 of \cite{BGR}. 

Semi-affinoid $K$-spaces are saturated locally G-ringed $K$-spaces which correspond to \emph{semi-affinoid $K$-algebras}. According to Definition \ref{semaffalgdefi}, a $K$-algebra is called semi-affinoid if it admits an $R$-lattice of ff type, that is, if it is obtained from an $R$-algebra of ff type by inverting $\pi$. Equivalently, the semi-affinoid $K$-algebras are the quotients of the \emph{free} semi-affinoid $K$-algebras
\[
R[[S_1,\ldots,S_m]]\langle T_1,\ldots,T_n\rangle\otimes_RK
\]
by ideals of semi-affinoid functions. The free semi-affinoid $K$-algebra above is easily identified with the ring of \emph{bounded} functions on the rigid $K$-space $\D^m_K\times\B^n_K$, where $\D^m_K$ is the $m$-dimensional open rigid unit disc and where $\B^n_K$ denotes the $n$-dimensional closed rigid unit disc.

One does not dispose of a general Noether normalization theorem for semi-affinoid $K$-algebras. Nonetheless, we can develop the commutative algebra of semi-affinoid $K$-algebras without too much effort, cf.\ Section \ref{semaffalg}: We first recall that the residue fields in the maximal ideals of a semi-affinoid $K$-algebra $A$ are finite extensions of $K$, and we describe the points of $\Max A$ in terms of rig-points of $R$-lattices of ff type for $A$. We then establish the universal property of the free semi-affinoid $K$-algebra above, cf. Proposition \ref{freesemiaffinoidprop}. To state it, we need the notion of topological quasi-nilpotency for semi-affinoid functions, cf.\ Definition \ref{pbtqndefi}: An element $f$ in $A$ is called topologically quasi-nilpotent if $|f(x)|<1$ for all $x\in \Max A$. 
We also show that semi-affinoid $K$-algebras are Jacobson rings and that the category of semi-affinoid $K$-algebras has amalgamated sums.

The maximum principle fails for general semi-affinoid $K$-algebras. For example, the coordinate function $S\in R[[S]]\otimes_RK$ is topologically quasi-nilpotent, while the supremum of $|S(x)|$ over all $x\in\Max (R[[S]]\otimes_RK)$ equals $1$.

A semi-affinoid $K$-algebra has a unique $K$-Banach algebra structure, which we discuss in Section \ref{semaffbanachsec}. For example, a $K$-Banach algebra norm on $A$  is obtained by presenting $A$ as a quotient of a free semi-affinoid $K$-algebra and by considering the residue norm of the Gauss norm. The maximum of a semi-affinoid function $f\in A$ is attained on the compact Berkovich spectrum $\M(A)$ of $A$. One must beware of the fact that the formation of $\M(A)$ is in general incompatible with complete localization of an $R$-lattice of ff type $\ul{A}$ for $A$, the Banach topology on $A$ being induced by the \emph{$\pi$-adic} topology on $\ul{A}$. In general, there is no ring topology on $A$ such that $\ul{A}\subseteq A$ is an open subring carrying the subspace topology: In order to consider $A$ as an f-adic ring in the sense of Huber, cf.\ \cite{Huberbuch}	 p.37, such that $\ul{A}$ is a subring of definition of $A$, one has to replace the natural topology on $\ul{A}$ by the $\pi$-adic one, cf.\ Example \ref{hubernosalvationex}.

If $A$ is a semi-affinoid $K$-algebra, there is a universal rigid $K$-space $X_A$ together with a $K$-homomorphism $\phi_A\colon A\rightarrow\Gamma(X_A,\O_{X_A})$. It is obtained from an $R$-lattice of ff type $\ul{A}$ of $A$ by means of Berthelot's construction: $X_A=(\Spf\ul{A})^\rig$. The homomorphism $\phi_A$ yields an identification of $\Max A$ with the set of physical points of $X_A$, and it induces isomorphisms of completed stalks, cf.\ Section \ref{assocrigspacesubsec}.

In Section \ref{gtopsec}, we define a G-topology on $\Max A$ which is coarser than the rigid-analytic G-topology induced via $\phi_A$. We first define an auxiliary G-topology $\sT_\aux$ and then obtain the uniformly rigid $G$-topology $\sT_\urig$ by means of the enhancement techniques described in \cite{BGR} 9.1.2. The $\sT_\aux$-admissible subsets of $\Max A$ are the \emph{semi-affinoid subdomains} of $\Max A$. These are essentially \emph{nested} rational subdomains defined by strict and non-strict inequalities. Technically speaking, semi-affinoid subdomains are defined in terms of compositions of admissible formal blowups, open immersions and completion morphisms starting out from an $R$-lattice of ff type. A formal blowup is called \emph{admissible} if it is defined by an ideal that is $\pi$-adically open but not necessarily open, cf.\ Section \ref{blowupssec}. Admissible formal blowups have been defined and studied independently by Nicaise, Temkin and the author, cf.\ \cite{N} and \cite{T}. 
They are stable under composition, but an admissible formal blowup defined on some open part $\fU$ of a formal $R$-scheme of ff type $\fX$ needs not extend to an admissible formal blowup on $\fX$, as it is shown by Example \ref{thegabberex}.

If $U\subseteq\Max A$ is a semi-affinoid subdomain, there exists a unique homomorphism $\psi_U\colon A\rightarrow A_U$ of semi-affinoid $K$-algebras which is universal with the property that the induced map $\Max\psi_U\colon\Max A_U\rightarrow \Max A$ factorizes through $U$, cf.\ Proposition \ref{subdomunivprop}. One easily sees that $\Max\psi_U$ is injective onto $U$. In particular, a semi-affinoid subdomain $U$ in $\Max A$ is equipped with an essentially unique ring $A_U$ of semi-affinoid functions. Semi-affinoid subdomains behave reasonably well. For instance, if $U$ is a semi-affinoid subdomain in $\Max A$ and if $V$ is a subset of $U$, then $V$ is a semi-affinoid subdomain in $U$ if and only if it is a semi-affinoid subdomain in $\Max A$. Let us point out that semi-affinoid subdomains are not defined in terms of their universal property, for lack of a Gerritzen-Grauert theorem in the uniformly rigid setting, cf.\ the discussion at the end of Section \ref{semaffsubdomsec}.

A semi-affinoid subdomain is called \emph{retrocompact} if its description involves no strict inequalities or, what amounts to the same, no completion morphisms. A \emph{retrocompact covering} of $\Max A$ is a finite covering by retrocompact semi-affinoid subdomains. A covering of $\Max A$ by semi-affinoid subdomains is called $\sT_\aux$-admissible if it has a leaflike refinement, where a \emph{leaflike covering} is a retrocompact covering that is obtained from an $R$-lattice of ff type by iterations of Zariski coverings and admissible formal blowups. We prove that every retrocompact covering has a leaflike refinement, cf.\ Proposition \ref{retroprop}. 
An explicit description of the saturated uniformly rigid G-topology $\sT_\urig$ obtained from $\sT_\aux$ is given in Section \ref{gtopsec}. We prove that the uniformly rigid G-topology is coarser than the rigid G-topology induced via $\phi_A$, while it is finer than the Zariski topology defined by the non-vanishing loci of functions in $A$.

The universal property satisfied by semi-affinoid subdomains yields a natural presheaf of semi-affinoid $K$-algebras on $\sT_\aux$. We prove that it is actually a sheaf for $\sT_\aux$; more generally, we establish the uniformly rigid version of Tate's Acyclicity Theorem, using methods from \cite{EGAIII1} and arguing via induction along a rooted tree describing the complexity of a given leaflike covering. In particular, we obtain a sheaf of $K$-algebras on $\sT_\urig$. Let $\sSp A$ denote the resulting G-ringed $K$-space; by Corollary \ref{srigpropcor}, its G-ringed $K$-space structure is compatible with restrictions to semi-affinoid subdomains. The universal $K$-homomorphism $\phi_A\colon A\rightarrow \Gamma(X_A,\O_{X_A})$ extends to a universal comparison morphism $X_A\rightarrow\sSp A$ of G-ringed $K$-spaces which is a bijection on physical points and which induces isomorphisms of stalks. The locally G-ringed $K$-space $\sSp A$ is called the semi-affinoid $K$-space associated to $A$, and the resulting functor from the dual of the category of semi-affinoid $K$-algebras to the category of locally G-ringed $K$-spaces is fully faithful, cf.\ Proposition \ref{ffprop}.

Uniformly rigid $K$-spaces can be obtained from semi-affinoid $K$-spaces by means of standard glueing procedures. In particular, the category of uniformly rigid $K$-spaces admits fibered products. The comparison morphisms $X_A\rightarrow\sSp A$ globalize: If $X$ is a uniformly rigid $K$-space, there exists a rigid $K$-space $X^\r$ together with a \emph{universal} comparison morphism of locally G-ringed $K$-spaces $X^\r\rightarrow X$ which is a bijection on physical points and which induces isomorphisms of stalks. The functor $\Spf\ul{A}\mapsto\sSp(\ul{A}\otimes_RK)$ globalizes as well, yielding the uniformly rigid generic fiber functor $\urig\colon\FF_S\rightarrow\uRig_K$. If $X$ is a semi-affinoid $K$-space and if $U\subseteq X$ is an admissible open semi-affinoid subspace, we do not know whether $U$ is necessarily a semi-affinoid subdomain; however, semi-affinoid subspaces of semi-affinoid $K$-spaces admit a description akin to one provided by the Gerritzen-Grauert theorem in the classical rigid setting, cf.\ Lemma \ref{semaffsubspacecharlem}.

In Section \ref{basefieldchangesec}, we study the base field extension functor for potentially non-finite extension fields and for uniformly rigid $K$-spaces satisfying a local finiteness condition. In Section \ref{cohmodsec}, we study coherent modules on uniformly rigid $K$-spaces; in particular, we single out a class of  coherent modules for which Kiehl's patching theorem \cite{KiehlAB} 1.2 holds. We present two different approaches. On the one hand, we study schematic closures of submodules in formal models; on the other hand, we use Berthelot's construction in conjunction with Kiehl's patching theorem for rigid spaces. In particular, we prove that a coherent \emph{ideal} on a semi-affinoid $K$-space is associated to its ideal of global sections. We thus dispose of a well-behaved notion of closed uniformly rigid subspaces. Using our results on schematic closures of submodules, we can establish the existence of schematic closures of closed uniformly rigid subspaces in given formal $R$-models of locally ff type. In Section \ref{galdescsec}, we discuss finite Galois descent for uniformly rigid spaces. Since the uniformly rigid G-topology depends on data over the valuation ring $R$, it seems difficult here to go beyond the semi-affinoid case.

Semi-affinoid $K$-algebras have also been investigated by Lipshitz and Robinson in \cite{LiRo} I, under the name of \emph{quasi-affinoid} $K$-algebras. The approach in \cite{LiRo} includes the situation where $R$ is not discrete and where the machinery of locally noetherian formal geometry is not available. In this more general setting, considerable effort is required to establish the universal properties of the free semi-affinoid $K$-algebras, cf.\ \cite{LiRo} I.5.2.3. In part III of \cite{LiRo}, Lipshitz and Robinson define a G-topology on the maximal spectrum of a semi-affinoid $K$-algebra, and they prove the corresponding acyclicity theorem. The resulting locally G-ringed $K$-spaces are not very different from our semi-affinoid $K$-spaces equipped with their auxiliary G-topology $\sT_\aux$. Lipshitz and Robinson do not consider the induced saturated G-topology, and they do not study the global aspects of the theory. Moreover, the definition of the G-topology in \cite{LiRo} is much less explicit than our definition, cf.\ \cite{LiRo} III.2.3.2, which is why the proof of the acyclicity theorem in \cite{LiRo} depends on the nontrivial quantifier elimination theorem \cite{LiRo} II Theorem 4.2. Our approach avoids quantifier elimination.

\vspace{0.7cm}
\textbf{\large Formal Néron models of uniformly rigid spaces}
\vspace{0.5cm}

In the second part of this thesis, we lay the foundations for the general theory of formal Néron models for uniformly rigid $K$-spaces. The major part of Chapter \ref{fnmchap} culminates in Theorem \ref{maincompthm}, showing that in many important cases, the formal Néron model (of tf type) of a rigid $K$-group in the sense of \cite{BS} already satisfies the stronger universal property of the formal Néron model (of ff type) of the associated uniformly rigid $K$-group:

\textbf{\ref{maincompthm} Theorem.} \textit{Let $G$ be a quasi-paracompact rigid $K$-group admitting a quasi-paracompact formal N\'eron model $\mathfrak{G}$ of locally tf type in the sense of \textup{\cite{BS}} such that the universal morphism $\mathfrak{G}^\textup{rig}\rightarrow G$ is a retrocompact open immersion and such that $\mathfrak{G}^\textup{rig}$ contains all formally unramified points of $G$. Then $\mathfrak{G}$ is the formal N\'eron model of the uniformly rigid $K$-space $G^\ur$.}

An open immersion is called retrocompact if it is a quasi-compact morphism. A formally unramified point of $G$ is a $K'$-valued point of $G\hat{\otimes}_KK'$, where $K'$ is the fraction field of a not necessarily finite cdvr-extension $R'$ of $R$. For example, the conditions of Theorem \ref{maincompthm} are fulfilled when $G$ is quasi-compact or when $G$ is the analytification of a smooth algebraic $K$-group scheme $\mathscr{G}$  admitting a N\'eron lft-model $\underline{\mathscr{G}}$, where $\mathscr{G}$ is commutative or where $\underline{\mathscr{G}}$ is of finite type over $R$. If $G$ is quasi-compact, the extension property for formally unramified points is satisfied by \cite{Wegel} Theorem 4; in the algebraic case, the extension property follows from \cite{BLR} 10.1/3.

On the basis of Theorem \ref{maincompthm}, further examples of formal Néron models of uniformly rigid $K$-spaces can be constructed by means of completion, descent and product decomposition techniques, as we explain in Section \ref{constrtechressec}. The resulting formal Néron models are in general not of locally tf type, and their universal morphisms are in general not surjective. General existence criteria and construction procedures for formal Néron models of uniformly rigid $K$-spaces seem to be out of reach.

The \emph{proof} of Theorem \ref{maincompthm} involves certain compactifications of affine formal $R$-schemes of locally ff type which we call \emph{envelopes}: If $\mathfrak{X}$ is an affine formal $R$-scheme of ff type and if $\ul{A}$ denotes the topological ring of functions on $\mathfrak{X}$, then the envelope $\mathfrak{X}_\pi$ of $\mathfrak{X}$ is defined to be $\textup{Spf } \ul{A}^\pi$, where $\ul{A}^\pi$ denotes the ring $\ul{A}$ equipped with its $\pi$-adic topology. Envelopes can also be defined for formal $R$-schemes which are proper over affine formal $R$-schemes of ff type, cf.\ Section \ref{envsec}; these envelopes are unique if they exist, but their existence is not always guaranteed. If $\mathfrak{X}$ is smooth over $R$, new formally unramified points may appear on $\fX_\pi$. For example, the Gauss point on the Berkovich spectrum $\M(R[[S]]\otimes_RK)$ extends to a formally unramified point of $\Spf (R[[S]]^\pi)$ which does not lie in $\Spf R[[S]]$, cf.\ Proposition \ref{integralredshilovprop}. Generic fibers of envelopes exist naturally in Huber's category of adic spaces, cf. Section \ref{envsec}.

The formation of envelopes does not commute with complete localization. None\-theless, the concept of envelopes becomes useful if it is put into practice within the framework of uniformly rigid $K$-spaces: If $\fY$ is an affine flat formal $R$-scheme of ff type, if $\fX$ is a flat separated formal $R$-scheme of locally tf type and if $\phi\colon\fY^\urig\rightarrow \fX^\urig$ is a morphism of \emph{uniformly rigid} generic fibers, then $\phi$ extends to a morphism from a (not necessarily algebraic) proper modification of $\fY_\pi$ to $\fX$, cf.\ Proposition \ref{graphenvexprop}. 

The proof of this statement relies on Grothendieck's Formal Existence Theorem \cite{EGAIII1} 5.1.4 and on our careful study of the category of uniformly rigid $K$-spaces. We first consider the situation where $\fX$ is allowed to be separated and of locally ff type over $R$ and where $\fY$ is any flat formal $R$-scheme of locally ff type. Under these assumptions, we establish the fact that the schematic closure $\ul{\Gamma}_\phi$ of the graph of $\phi$ in $\fY\times\fX$ is proper over $R$, cf.\ Theorem \ref{graphproperthm}. To do so, we use schematic images under (not necessarily algebraic) proper formal morphisms, which we study in Section \ref{schemimsec}, using the formal proper flat base change Theorem \ref{mainflatbasechangethm} for direct images of coherent sheaves. We also need to invoke the fact that a morphism of flat formal $R$-schemes of locally ff type is adic if and only if the associated morphism of rigid generic fibers, obtained via Berthelot's construction, is retrocompact, cf.\ Proposition \ref{rqcadicprop}. In the case where $\fX$ is of locally tf type, we apply Grothendieck's Formal Existence Theorem to the reductions modulo powers of $\pi$ of $\ul{\Gamma}_\phi$.

If the affine formal $R$-scheme of ff type $\fY$ is \emph{smooth}, cf.\ Section \ref{smoothnesssec}, the proper modification of $\fY_\pi$ provided by Proposition \ref{graphenvexprop} is an isomorphism generically on $\fY_\pi$, cf.\ Corollary \ref{genisocor}. This is a uniformly rigid analog of \cite{BS} Lemma 2.2; cf.\ also Lemma \ref{genisolem} where we recall the well-known fact that admissible blowups of locally noetherian $R$-schemes with reduced special fibers are isomorphisms on $R$-dense open subsets. To prove Corollary \ref{genisocor}, we essentially apply Lemma \ref{genisolem}, using induction along a rooted tree describing the complexity of a suitable leaflike covering of $\fY^\urig$. Technical problems arise due to the fact that the formation of envelopes does not commute with complete localization. Some of these difficulties are dealt with in Section \ref{envlocsec}, where we consider a proper formal scheme $\fY'$ over an affine formal $R$-scheme of ff type $\fY$ admitting an envelope $\fY'_\pi\rightarrow\fY_\pi$ and where we show that we do not lose any points on $\fY'_\pi$ if we consider a suitable affine open covering of $\fY'$ and then form the envelopes of the given affine parts.

Building upon these results, we can prove Weil-type extension theorems for morphsims to uniformly rigid $K$-groups, cf.\ Theorem \ref{mainweilextthm}. Our proofs require some further preparations: In Section \ref{ancontsec}, we provide an ff type version of Bartenwerfer's and Lütkebohmert's Hartogs continuation theorem, cf.\ \cite{Ltke} Satz 6, and in Section \ref{tubessec} we prove a uniformly rigid version of Kiehl's tube lemma, cf.\ \cite{KiehldeRh} Satz 1.6, by adopting methods of Kisin, cf.\ \cite{Kisin} Section 2. To finally achieve the proof of Theorem \ref{mainweilextthm}, we must argue locally on generic fibers of envelopes, which seems to be possible only within Huber's category of adic spaces.

In Section \ref{nmsec}, we derive Theorem \ref{maincompthm} from Theorem \ref{mainweilextthm}, and we explain how it can be used to construct new interesting examples of uniformly rigid $K$-spaces admitting formal Néron models. Assuming that the residue field $k$ of $R$ has positive characteristic, we define strongly concordant uniformly rigid $K$-groups, transcribing definitions of Chai, and we establish the existence of their formal Néron models. Finally, in Section \ref{chaismethod} we assume that $k$ is perfect of positive characteristic, and we explain how Chai's strategy for computing the base change conductor of abelian varieties with potentially multiplicative reduction can be carried out in the uniformly rigid setting. Following our more conceptual approach, we avoid referring to the Néron-Popescu Theorem, thereby simplifying Chai's original arguments.

\vspace{0.7cm}
\textbf{\large Acknowledgements}
\vspace{0.5cm}

Numerous people and institutions have supported my work on this thesis. First of all, I would like to express my heartfelt gratitude to my thesis advisor Siegfried Bosch -- for bringing Chai's work on the base change conductor to my attention, for many very helpful discussions and for providing excellent working conditions. Above all, I would like to thank Professor Bosch for initially introducing me to the beauty of modern mathematics. The mathematical knowledge which he has passed on to me has profoundly enriched my life.

I would like to thank Ahmed Abbes, Brian Conrad, Ching-Li Chai, Christopher Deninger, Ofer Gabber, Alexander Gurevich, Urs Hartl, Simon Thomas Hüsken, Kiran Kedlaya, Klaus Loerke, Werner Lütkebohmert, Johannes Nicaise, Jérôme Poineau, Björn Selander, Matthias Strauch and Christian Wahle for helpful mathematical discussions. Moreover, I would like to thank Urs Hartl, Peter Schneider and their coworkers for offering interesting and challenging seminars, thereby adding to the inspiring and pleasant ambience in the department. In addition, I would like to thank my friends and colleagues Simon Thomas Hüsken, Katharina Küper and Christian Wahle for encouragement, friendship and support.

My work was financially supported by the \emph{Hamburger Stiftung für internationale Forschungs- und Studienvorhaben}, by the \emph{German National Academic Foundation} and by the Graduiertenkolleg \emph{Analytic Topology and Metageometry}. I would like to extend my gratitude to these institutions. Moreover, I would like to thank the \emph{Massachusetts Institute of Technology}, where this work was initiated during the academic year 2005/06, for its hospitality.

I would like to take the opportunity to thank Professor Peter Boschan for supervising my undergraduate studies in theoretical physics.

Last but not least, I thank my parents for their unrestricted love and support.
\newpage



\renewcommand{\contentsname}{Contents}

\cleardoublepage
\pagestyle{fancy}
\def\MakeUppercase#1{#1}
\vspace*{0.5cm}
\tableofcontents

\mainmatter
\pagestyle{fancy}

\romanenum


\chapter{Uniformly rigid spaces}\label{semrigspchap}	

	\section{Complements on formal geometry, Part I}\label{semrigformalsec}
	
		\subsection{Morphisms of locally ff type}\label{locfftypesec}

We begin by reviewing the notion of a morphism of locally formally finite type of locally noetherian formal schemes, cf.\ also \cite{Berk2} Section 1 in \cite{Yek} Section 1.

Let $R$ be a noetherian adic ring with an ideal of definition $\fr$. If $T_1,\ldots,T_n$ is a finite system of variables, then
\[
R\langle T_1,\ldots,T_n\rangle
\]
denotes the ring of strictly convergent power series in the $T_i$ with coefficients in $R$, cf.\ \cite{EGAIn} 0.7.5. If $S$ is a variable, one easily checks that $R[[S]]$ is a noetherian adic ring for the $(\fr,S)$-adic topology and that this topology does not depend on the choice of $\fr$. Let us now consider a finite system of variables $S_1,\ldots,S_m$. We equip $R[[S_1,\ldots,S_m]]$ with the $(\fr,S_1,\ldots,S_m)$-adic topology; then
\[
R[[S_1,\ldots,S_m]]\langle T_1,\ldots,T_n\rangle
\]
is a noetherian adic ring, carrying the $(\fr,S_1,\ldots,S_m)$-adic topology, such that the $S_i$ are topologically nilpotent.
One verifies without difficulty that it is characterized by the following universal property: If $A$ is a noetherian adic topological $R$-algebra and if $s_1,\ldots,s_m$, $t_1,\ldots,t_n$ are elements in $A$ where the $s_i$ are topologically nilpotent, there exists a unique continuous $R$-homomorphism 
\[
R[[S_1,\ldots,S_m]]\langle T_1,\ldots,T_n\rangle\rightarrow A
\]
sending $S_i$ to $s_i$ and $T_j$ to $t_j$ for all $1\leq i\leq m$, $1\leq j\leq n$. In particular, there is a canonical topological isomorphism $R[[S]]\langle T\rangle\cong R\langle T\rangle[[S]]$; for example, if $R=k$ is a discrete field, we obtain $k[[S]]\langle T\rangle\cong k[T][[S]]$.

The following statement is proven by means of standard lifting arguments in topolo\-gical commutative algebra:

\begin{lem}\label{basicfftlem}
Let $R$ be a noetherian adic ring, and let $A$ be a noetherian adic topological $R$-algebra. Then the following are equivalent:
\begin{packed_enum}
\item There exist finite systems of variables $S_1,\ldots,S_m$, $T_1,\ldots,T_n$ such that $A$ is a topological quotient of $R[[S_1,\ldots,S_m]]\langle T_1,\ldots,T_n\rangle$.
\item For some or any pair of ideals of definition $\fr\subseteq R$, $\a\subseteq A$ such that $\fr A\subseteq\a$, the discrete quotient $A/\a$ is an $R/\fr$-algebra of finite type.
\end{packed_enum}
\end{lem}

\begin{defi}
Let $R$ be a noetherian adic ring. A noetherian adic topological $R$-algebra $A$ is said to be of \emph{formally finite}\index{algebra! topological, of ff type} or \emph{ff} type if it satisfies the equivalent conditions of Lemma \ref{basicfftlem}.
\end{defi}

It is clear from Lemma \ref{basicfftlem} that the category of topological $R$-algebras of ff type is stable under complete localization, under the formation of quotients, under the formation of complete tensor products and under formal completion. Moreover, the property of being of ff type can be checked locally on formal spectra, and any continuous homomorphism of topological $R$-algebras of ff type is itself of ff type.

\begin{defi}\label{fftypemorphdefi}
If $\phi\colon\fY\rightarrow\fX$ is a morphism of locally noetherian formal schemes, then $\phi$ is said to be of locally \emph{formally finite}\index{morphism!of locally ff type} or locally \emph{ff} type if the underlying morphism of reduced subschemes of definition is of locally finite type. We say that $\phi$ is of \emph{formally finite} or \emph{ff} type if it is of locally ff type and quasi-compact.
\end{defi}

The above remarks also apply to morphisms of locally ff type.

\begin{defi}\label{fftalgdefi}
Let $R$ be a noetherian adic ring. An $R$-algebra $A$ is said to be of \emph{formally finite} \index{algebra!of ff type} or \emph{ff} type if there exists a ring topology on $A$ such that $A$ is a topological $R$-algebra of ff type.
\end{defi}

An $R$-algebra $A$ is of ff type if and only if there exist finite systems of variables $S_1,\ldots,S_m$, $T_1,\ldots,T_n$ such that $A$ is a quotient of $R[[S_1,\ldots,S_m]]\langle T_1,\ldots,T_n\rangle$. If $A$ is a (topological) $R$-algebra of ff type, a (strict) epimorphism 
\[
R[[S_1,\ldots,S_m]]\langle T_1,\ldots,T_n\rangle\rightarrow A
\]
is called a formally finite presentation\index{presentation!formally finite} of $A$, and the image of the tuple 
\[
((S_1,\ldots,S_m),(T_1,\ldots,T_n))
\]
in $A$ is called a formal generating system of $A$. \index{generating system!formal}

\subsubsection{Continuity of homomorphisms}

\begin{lem}\label{jacobsontoplem}
Let $R$ be a noetherian adic ring with the property that the quotient of $R$ modulo an ideal of definition is a Jacobson ring. Let moreover $A$ be a topological $R$-algebra of ff type. Then the biggest ideal of definition of $A$ coincides with the Jacobson radical of $A$.
\end{lem}
\begin{proof}
Let $\a$ and $\fr$ be ideals of definition of $A$ and $R$ respectively. Then $A/\a$ is of finite type over $R/\fr$ and thereby a Jacobson ring. Hence, the radical of $\a$ is the intersection of all maximal ideals above $\a$. Since $A$ is $\a$-adically complete, every maximal ideal of $A$ contains $\a$. It follows that the radical of $\a$, which is the biggest ideal of definition of $A$, coincides with the Jacobson radical of $A$, as claimed.
\end{proof}

\begin{lem}\label{jacobsonequivlem}
Let $R$ be a noetherian adic ring with the property that the quotient of $R$ modulo its biggest ideal of definition is of finite type over a field. Then the category of $R$-algebras of ff type is naturally isomorphic to the category of topological $R$-algebras of ff type.
\end{lem}
\begin{proof}
We construct an inverse to the forgetful functor. Let us first choose a field $k$ together with a ring homomorphism $k\rightarrow R/\fr$ of finite type, where $\fr$ is the biggest ideal of definition of $R$; this is possible by our assumption on $R$. If $A$ is an $R$-algebra of ff type, we choose any formally finite presentation $\alpha$ of $A$; then $\alpha$ defines a topology on $A$ such that $A$ becomes a topological $R$-algebra of ff type. By Lemma \ref{jacobsontoplem}, this topology is the Jacobson-adic topology and, hence, does not depend on the choice of $\alpha$. Indeed, $R/\fr$ is a Jacobson ring since it is of finite type over $k$. It remains to see that if $\phi\colon A\rightarrow B$ is a homomorphism of $R$-algebras, then $\Jac A\subseteq\phi^{-1}(\Jac B)$. In other words, we must show that if $f$ is an element in $A$ that vanishes in every closed point of $\Spec A$, then the $\phi$-image of $f$ vanishes in every closed point of $\Spec B$. Thus, it suffices to see that for any maximal ideal $\m$ in $B$, the $\phi$-pullback $\p\mathrel{\mathop:}=\phi^{-1}(\m)$ is maximal in $A$. Let $\b$ denote the biggest ideal of definition of $B$. Since $B$ is $\b$-adically complete, $\m$ contains $\b$, and hence the field $B/\m B=(B/\b)/\m(B/\b)$ is of finite type over $R/\fr$ which is of finite type over $k$. Thus, we have inclusions $k\subseteq A/\p\subseteq B/\m$, where $k\subseteq B/\m$ is a finite field extension; it follows that the domain $A/\p$ is a field and, hence, that $\p$ is maximal in $A$, as desired.
\end{proof}

The conditions of Lemma \ref{jacobsonequivlem} are satisfied if $R$ is a field or a complete discrete valuation ring. In the following, the statement of Lemma \ref{jacobsonequivlem} will be used without further reference. For example, if $R$ is as in the statement of Lemma \ref{jacobsonequivlem} and if $A\rightarrow B_i$, $i=1,2$, are $R$-homomorphisms of $R$-algebras of ff type, then the completed tensor product $B_1\hat{\otimes}_RB_2$ is well-defined, and it is the amalgamated sum of $B_1$ and $B_2$ over $A$ in the category of $R$-algebras of ff type.
				
		\subsection{Excellence and normality}\label{exnormalsec}		
			Let $R$ be a field or a complete discrete valuation ring with fraction field $K$. In this section, we briefly gather some results on excellence and normality for formal $R$-schemes of locally ff type.

$R$-algebras of ff type are excellent. This follows from \cite{V1} Proposition 7 in the field case or in the equal characteristic case and from \cite{V2} Theorem 9 when $R$ has mixed characteristic. The field case is contained in the case where $R$ is a complete discrete valuation ring of equal characteristic whenever no $R$-flatness assumption is involved.

Let $\fX$ be a formal $R$-scheme of locally ff type, and let $x\in\fX$ be a point. We do not know whether the stalk $\O_{\fX,x}$ of $\fX$ in $x$ or its maximal-adic completion $\hat{\O}_{\fX,x}$ are excellent. Let $\textbf{P}$ denote the property of a noetherian ring of being reduced, normal, regular, Gorenstein, Cohen-Macaulay, a complete intersection, $R_i$ or $S_i$. Let us assume that $\fX=\Spf A$ is affine, and let $x\in \fX$ correspond to an open prime ideal $\p\subseteq A$. By \cite{ConradIrr} Lemma 1.2.1 and its proof, the conditions $\textbf{P}(A_\p)$, $\textbf{P}(\O_{\fX,x})$ and $\textbf{P}(\hat{\O}_{\fX,x})$ are equivalent. It is thus unambiguous to ask whether $\fX$ satisfies $\textbf{P}$ in a given point. By \cite{ConradIrr} Lemma 1.2.1, the non-\textbf{P}-locus in $\fX$ is a Zariski-closed subset of $\fX$. 

From now on, let $R$ be a complete discrete valuation ring with fraction field $K$. By \cite{ConradIrr} Section 1.1, affinoid $K$-algebras and stalks of rigid $K$-spaces are excellent, so the analogous statements hold for rigid $K$-spaces.

If $A$ is an excellent ring, the normalization $\tilde{A}$ of $A$ is defined to be the integral closure of $A$ in the total fraction ring of its biggest reduced quotient $A/\rad(A)$ of $A$. The ring $\tilde{A}$ is normal, and the natural ring homomorphism $A\rightarrow\tilde{A}$ is finite. If $A$ is an affinoid $K$-algebra or an $R$-algebra of ff type, the formation of $\tilde{A}$ commutes with rigid or formal localization, cf.\ \cite{ConradIrr} 1.2.3. It follows that the homomorphisms $A\rightarrow\tilde{A}$ glue to finite morphisms $\tilde{X}\rightarrow X$ and $\tilde{\fX}\rightarrow\fX$, for every rigid $K$-space $X$ and every formal $R$-scheme of locally ff type $\fX$, such that $\tilde{X}$ and $\tilde{\fX}$ are normal. A rigid $K$-space $Y$ over $X$ is called a \emph{normalization}\index{normalization} of $X$ if it is $X$-isomorphic to $\tilde{X}$, and similarly in the formal case. By \cite{ConradIrr} Theorem 2.1.3, normalization morphisms commute with Berthelot's generic fiber functor which we discuss in Section \ref{berthsec}:

\begin{theorem}
Let $\fX$ be a formal $R$-scheme of locally ff type, and let $\phi\colon\tilde{\fX}\rightarrow\fX$ be a normalization of $\fX$. Then $\phi^\rig$ is a normalization of $\fX^\rig$.
\end{theorem}

We say that a formal $R$-scheme $\fX$ of locally ff type is \emph{generically normal}\index{generically normal} if $\fX^\rig$ is normal. In Corollary \ref{normisblowupcor}, we will show that normalizations of generically normal flat formal $R$-schemes of ff type are finite admissible formal blowups. To do so, we will need the following trivial observation:

\begin{lem}\label{flatgennormalnormalizationlem}
Let $\fX$ be a generically normal flat formal $R$-scheme of locally ff type, and let $\tilde{\fX}$ be a normalization of $\fX$. Then $\tilde{\fX}$ is $R$-flat.
\end{lem}
\begin{proof}
We may assume that $\fX$ is affine, $\fX=\Spf A$. Then $A$ is reduced, being contained in the $K$-algebra $A\otimes_RK$ which is normal and, hence, reduced. Now $\tilde{\fX}=\Spf\tilde{A}$, where $\tilde{A}$ is the integral closure of $A$ in its total ring of fractions. Since $A$ has no $\pi$-torsion, the same holds for its total ring of fractions and, hence, for $\tilde{A}$.
\end{proof}

		\subsection{Formal blowups}\label{blowupssec}
			Let $\fX$ be a locally noetherian formal scheme, and let $\sI\subseteq\O_\fX$ be a coherent ideal. The formal blowup of $\fX$ in $\sI$ has been defined and studied in \cite{BL1} in the case where $\sI$ is \emph{open}. It has been independently observed by Nicaise, Temkin, the author and presumably others that the basic definitions and many properties of formal blowups carry over to the situation where $\sI$ is not open.

\begin{defi}\index{blowup!formal}
Let $\fX$ be a locally noetherian formal scheme, and let $\sI\subseteq\O_\fX$ be a coherent ideal. The formal blowup of $\sI$ on $\fX$ is a morphism of locally noetherian formal schemes $b_\sI\colon\fX_\sI\rightarrow\fX$ satisfying the following universal property: $b_\sI^*\sI$ is invertible, and for any morphism of locally noetherian formal schemes $\phi\colon\fY\rightarrow\fX$ with the property that $\phi^*\sI$ is invertible, there exists a unique morphism $\phi_\sI\colon\fY\rightarrow\fX_\sI$ such that $\phi=b_\sI\circ\phi_\sI$.
\end{defi}

It is easily seen that if $\fX=\Spf A$ is affine and if $\sI$ is associated to $I\subseteq A$, then the formal blowup of $\sI$ on $\fX$ is obtained from the algebraic blowup of $\Spec A$ in $I$ via completion with respect to an ideal of definition of $A$. Algebraic blowups commuting with flat base change, this construction globalizes, and we obtain the following result, cf.\ \cite{N} 2.16:

\begin{prop}\label{basicblowupprop}
Let $\fX$ be a locally noetherian formal scheme, and let $\sI\subseteq\O_\fX$ be a coherent ideal. Then the formal blowup $b_\sI\colon\fX_\sI\rightarrow\fX$ exists. It is proper, and it is algebraic over any affine open part of $\fX$. Moreover, if $\phi\colon\fY\rightarrow\fX$ is a \emph{flat} morphism of locally noetherian formal schemes, then 
\[
\phi^*b_\sI\colon\fX_\sI\times_\fX\fY\rightarrow\fY
\]
is the formal blowup of $\fY$ in $\phi^*\sI$; this applies in particular to the case where $\phi$ is a completion morphism. 
\end{prop}

In \cite{N} 2.16 (4), the flat base change morphism $\phi$ is assumed to be adic, but one easily sees that this assumption is unnecessary. We refer to Section \ref{faithflatsec} for a discussion of flatness in locally noetherian formal geometry.

If $\fX=\Spf A$ is affine and if $\sI$ is generated by sections $f_1,\ldots,f_r\in A$, then $\fX_\sI$ is covered by the affine open parts
\[
\Spf A\langle T_1,\ldots,T_r\rangle/(f_i T_j-f_j)/(f_i\textup{-torsion})\quad,\quad 1\leq i\leq r\quad.
\]
Indeed, this follows from the description of algebraic blowups and the fact that completion morphisms of locally noetherian formal schemes are flat, cf.\ \cite{N} 2.18 and \cite{T} 2.1.4.

If $\phi\colon\fY\rightarrow\fX$ is a not necessarily flat morphism of locally noetherian formal schemes, the strict transform of $\phi$ under $b_\sI$ is defined to be the formal blowup of $\fY$ in $\phi^*\sI\subseteq\O_\fY$, together with the induced morphism $\phi_\sI\colon\fY_{\phi^*\sI}\rightarrow\fX_\sI$. If $\phi$ is a closed immersion, $\phi_\sI$ is a closed immersion as well. Indeed, this statement can be checked locally on $\fX$, where it follows from the corresponding algebraic fact and the observation that closed immersions are preserved under formal completion.

\subsubsection{Admissible formal blowups}

Let $R$ be a complete discrete valuation ring with fraction field $K$ and residue field $k$, and let $\pi\in R$ be a uniformizer. This convention will be effective throughout the rest of this thesis.

\begin{defi}\index{blowup!admissible formal}
Let $\fX$ be a formal $R$-scheme of locally ff type. A formal blowup of $\fX$ is called \emph{admissible} if it can be defined by a $\pi$-adically open ideal.
\end{defi}

Let us emphasize that in the locally ff type situation, a $\pi$-adically open ideal needs not be open. Similarly, we say that a blowup of an $R$-scheme is \emph{admissible}\index{blowup!admissible} if it is defined by a $\pi$-adically open quasi-coherent ideal.

If $\fX=\Spf A$ is a \emph{flat} affine formal $R$-scheme of ff type and if $\sI\subseteq\O_\fX$ is a $\pi$-adically open coherent ideal, generated by sections $f_1,\ldots,f_r\in A$, then $\fX_\sI$ is covered by the affine open parts
\[
\Spf A\langle T_1,\ldots,T_r\rangle/(f_i T_j-f_j)/(\pi\textup{-torsion})\quad,\quad 1\leq i\leq r\quad.
\]
In particular, $\fX_\sI$ is $R$-flat. The proof of this statement is identical to the proof in the tf type situation which is given in \cite{BL1} 2.2, cf.\ \cite{N} 2.18. 

Admissible formal blowups induce isomorphisms of rigid generic fibers. Again, we refer to Section \ref{berthsec} for a synopsis of Berthelot's construction.

\begin{prop}\label{rigisoprop}
Let $\fX$ be a formal $R$-scheme of locally ff type, and let $b$ be an admissible formal blowup of $\fX$. Then the rigid-analytic generic fiber $b^\rig$ of $b$ is an isomorphism.
\end{prop}

This statement is easily reduced to the corresponding one in the tf type setup, cf.\ \cite{N} 2.19.

\begin{lem}\label{admblowupcartlem}
Let $\phi\colon\fY\rightarrow \fX$ be a morphism of flat formal $R$-schemes of locally ff type, and let $b_\sI\colon\fX_\sI\rightarrow\fX$ be the admissible formal blowup of $\fX$ in a $\pi$-adically open coherent ideal $\sI\subseteq\O_\fX$.
Let us write $\sJ\mathrel{\mathop:}=\phi^*\sI$, $\psi\mathrel{\mathop:}=\phi_\sI$. Then the associated strict transform diagram
\[
\begin{diagram}
\fY_{\sJ}&\rTo^{\psi}&\fX_\sI\\
\dTo<{b_\sJ}&&\dTo>{b_\sI}\\
\fY&\rTo^{\phi}&\fX
\end{diagram}
\]
is \emph{cartesian} in the category of \emph{flat} formal $R$-schemes of locally ff type. In particular, if $\phi$ is an admissible formal blowup, then $\psi$ is an admissible formal blowup as well.
\end{lem}	
\begin{proof}
By our above discussion, $\fY_\sJ$ is $R$-flat. Let $\fZ$ be a flat formal $R$-scheme of locally ff type, and let $\sigma\colon\fZ\rightarrow \fY$, $\tau\colon\fZ\rightarrow \fX_\sI$ be $R$-morphisms such that the outer square in the diagram
\[
\begin{diagram}
\fZ&&&&\\
&\rdTo(4,2)^{\tau}\rdTo(2,4)_{\sigma}&&&\\
&&\fY_\sJ&\rTo_{\psi}&\fX_\sI\\
&&\dTo>{b_\sJ}&&\dTo>{b_\sI}\\
&&\fY&\rTo^{\phi}&\fX
\end{diagram}
\]
commutes. Then $\sI\O_\fZ$ is invertible, and by the universal property of $b_\sJ$ there exists a unique morphism $\mu\colon\fZ\rightarrow \fY_\sJ$ such that the resulting lower triangle in the following diagram commutes:
\[
\begin{diagram}
\fZ&&&&\\
&\rdTo(4,2)^{\tau}\rdTo(2,4)_{\sigma}\rdDashto_{\mu}&&&\\
&&\fY_\sJ&\rTo_{\psi}&\fX_\sI\\
&&\dTo>{b_\sJ}&&\dTo>{b_\sI}\\
&&\fY&\rTo^{\phi}&\fX
\end{diagram}
\]
It remains to show that the upper triangle commutes as well. By Corollary \ref{rigflatfaithfulcor}, the functor $\rig$ is faithful on the category of flat formal $R$-schemes of locally ff type; hence we may check commutativity after passing to rigid generic fibers. However, by Proposition \ref{rigisoprop} the morphism $b_\sI^\rig$ is an isomorphism, and so we obtain the desired result.
\end{proof}

We introduce the following notation:

\begin{defi}\index{fibered product!R@$R$-flat}
Let $\fX\rightarrow\fZ$ and $\fY\rightarrow\fZ$ be morphisms of flat formal $R$-schemes of locally ff type. We let
\[
\fX\times'_\fZ\fY
\]
denote the resulting fibered product in the category of \emph{flat} formal $R$-schemes of locally ff type; it is called the \emph{$R$-flat fibered product} of $\fX$ and $\fY$ over $\fZ$.
\end{defi}

We easily see that $R$-flat fibered products exist. In fact, $\fX\times'_\fZ\fY\subseteq\fX\times_\fZ\fY$ is the closed formal subscheme that is defined by the coherent ideal of $\pi$-torsion. 

If $\sI$ and $\sJ$ are $\pi$-adically open coherent ideals in $\O_\fX$, then the above proposition shows that $\fX_\sI\times'_\fX\fX_\sJ$ is the admissible formal blowup of $\sI\sJ$ on $\fX$. Indeed, one uses the universal property of formal blowups together with the fact that a product of two ideals is invertible if and only if both factors are invertible.

\subsubsection{Formal dilatations}

\begin{defi}\index{formal dilatation}
Let $\fX$ be a flat formal $R$-scheme of locally ff type, and let $\sI\subseteq\O_\fX$ be a coherent ideal containing $\pi$. The formal dilatation of $\fX$ in $\sI$ is the restriction of the admissible formal blowup $b_\sI$ of $\sI$ on $\fX$ to the open part of $\fX_\sI$ where $b_\sI^*\sI$ is generated by $\pi$.
\end{defi}

In particular, a formal dilatation is an affine adic morphism that is obtained, locally, from an algebraic dilatation via formal completion. Of course, the coherent ideals $\sI$ in $\O_\fX$ containing $\pi$ correspond to the closed formal subschemes of $\fX$ that are contained in the special fiber $\fX_k$ of $\fX$, which is the closed formal subscheme defined by $\pi$. Formal dilatations are characterized by the expected universal property:

\begin{prop}\label{admformaldilprop}
Let $\fX$ be a formal $R$-scheme of locally ff type, let $\sI\subseteq\O_\fX$ be a coherent ideal containing $\pi$, and let $d_\sI\colon\fX'_\sI\rightarrow\fX$ be the formal dilatation of $\sI$ on $\fX$. If $\fZ$ is a flat formal $R$-scheme of locally ff type and if $\phi\colon\fZ\rightarrow\fX$ is any morphism whose special fiber factorizes through $V(\sI)$, then $\phi$ factorizes uniquely through $d_\sI$.
\end{prop}

For the proof, we refer to \cite{N} 2.22. It follows in particular that the formal dilatation of a formal $R$-group scheme of locally ff type along a closed formal $k$-subgroup scheme of its special fiber is a formal $R$-group scheme of ff type, again. Of course, the rigid generic fiber of a flat formal $R$-scheme of locally ff type may shrink under formal dilatations.

\subsubsection{Compositions of formal blowups}

A composition of formal blowups of noetherian formal schemes in open ideals is a formal blowup in an open ideal again, cf.\ \cite{BL1} 2.5. In virtue of Grothendieck's Formal Existence Theorem, the openness assumption can be dropped. To make a precise statement, let us introduce the following terminology: Let $\fX$ be a locally noetherian formal scheme, and let $\fT\subseteq\fX$ be a closed formal subscheme. A formal blowup of $\fX$ is called $\fT$-supported if it can be defined by an ideal whose associated closed formal subscheme lies, locally, in an infinitesimal neighborhood of $\fT$ in $\fX$. In particular, if $\fX$ is of locally ff type over $R$, then a formal blowup of $\fX$ is admissible if and only if it is supported in the special fiber of $\fX$. For the proof of the following important result, we refer to \cite{T} 2.1.6:

\begin{prop}\label{compstableprop}
Let $\fX$ be a noetherian formal scheme, let $\fT\subseteq\fX$ be a closed formal subscheme, let $\phi\colon\fX'\rightarrow\fX$ be a $\fT$-supported formal blowup of $\fX$, and let $\psi\colon\fX''\rightarrow\fX'$ be a $\phi^{-1}(\fT)$-supported formal blowup of $\fX'$; then $\phi\circ\psi$ is a $\fT$-supported formal blowup of $\fX$.
\end{prop}

Here we wrote $\phi^{-1}(\fT)$ to denote the schematic preimage $\fT\times_\fX\fX'$ of $\fT$ in $\fX'$. In particular, if $\fX$ is of ff type over $R$, a composition of admissible formal blowups as above is an admissible formal blowup, again.

\subsubsection{Continuation of formal blowups}

If $\fX$ is a quasi-paracompact locally noetherian formal scheme and if $\fU$ is a quasi-compact open formal subscheme of $\fX$, then any formal blowup of $\fU$ that is definable by an \emph{open} coherent ideal $\sI\subseteq\O_\fU$ extends to a formal blowup of $\fX$. It is even possible to extend $\sI$ to an open coherent ideal in $\O_\fX$, cf.\ \cite{FRG} 2.6/13. In fact, after reducing modulo a sufficiently high power of an ideal of definition, the problem becomes algebraic in nature, and it suffices to apply \cite{EGAIn} 6.9.6. These extension results fail to hold if we drop the openness assumption on $\sI$. 
Before giving a counterexample, we prove two preparatory lemmas. We begin with an important statement which will also be used later in this thesis:

\begin{lem}\label{genisolem}
Let $X$ be a flat locally noetherian $R$-scheme, let $\phi\colon X'\rightarrow X$ be an admissible blowup, and let $\eta\in X_k$ be a generic point. On any quasi-compact open neighborhood of $\eta$ in $X$ whose special fiber is integral, $\phi$ is defined by a $\pi$-adically open coherent ideal whose vanishing locus does not contain $\eta$. In particular, if $X_k$ is reduced in $\eta$, then $\phi$ is an isomorphism around an open neighborhood of $\eta$ in $X$.
\end{lem}
\begin{proof}
We may assume that $X$ is quasi-compact and that $X_k$ is integral. Let $\sI\subseteq\O_X$ be a $\pi$-adically open ideal defining $\phi$. If $\sI\subseteq\pi \O_\fX$, we may uniquely divide $\sI$ by $\pi$ without changing $\phi$ and without changing the fact that $\sI$ is $\pi$-adically open. If the resulting coherent ideal still lies in $\pi \O_X$, we repeat the procedure. This process terminates after finitely many steps. Indeed, if it did not terminate, then $\pi$ would be a unit in $\O_X$, $\sI$ being $\pi$-adically open. However, $X_k$ contains $\eta$ and, hence, is nonempty. We may therefore assume that $\sI$ is not contained in $\pi \O_X$. Since $X_k$ is integral, this means precisely that $\eta\notin V(\sI)$. Since $\phi$ is an isomorphism on the open complement of $V(\sI)$, the claim has been shown.
\end{proof}

\begin{lem}\label{codimtwoisolocuslem}
Let $X$ be a locally noetherian scheme, let $\sI\subseteq\O_X$ be a coherent ideal such that $V(\sI)$ has codimension $\geq 2$ in $X$, and let $b_\sI$ denote the blowup of $X$ in $\sI$. Then for any open subscheme $U\subseteq X$,  the restriction 
\[
b_\sI^{-1}(U)\rightarrow U
\]
of $b_\sI$ is an isomorphism if and only if $U$ is disjoint to $V(\sI)$.
\end{lem}
\begin{proof}
We know that $b_\sI$ is an isomorphism over the complement of $V(\sI)$. By \cite{EGAIV1} 0.14.2.3, the codimension of $V(\sI)$ in $X$ can at most increase under localization, so we may assume that $X$ is affine and that $V(\sI)$ is nonempty.  We must show that $b_\sI$ is not an isomorphism. If $b_\sI$ was an isomorphism, then $\sI$ would be locally principal. By \cite{EGAIV1} 0.14.2.3 and by Krull's Principal Ideal Theorem, this would contradict the fact that $V(\sI)$ has codimension $\geq 2$ in $X$.
\end{proof}

\begin{example}\label{thegabberex}
Let $\fX$ be the completion of $\Spec R[X,Y,Z]$ along the closed subscheme $V(\pi, XY, Z)$, and let $\fU\subseteq\fX$ be the open complement of $V(X,Y)$. Then $\fU$ has two connected components $\fU\cap D(X)$ and $\fU\cap D(Y)$. We define a closed formal subscheme $\fV\subseteq\fU$ by setting
\begin{eqnarray*}
\fV\cap D(X)&\mathrel{\mathop:}=&V(\pi,Y,Z)\\
\fV\cap D(Y)&\mathrel{\mathop:}=&V(\pi,Z)\quad.
\end{eqnarray*}
Then $\fV$ is contained in the special fiber of $\fU$; let $\sI\subseteq\O_\fU$ denote the corresponding coherent ideal. The admissible formal blowup $b_\sI$ of $\fU$ in $\sI$ does \emph{not} extend to an admissible formal blowup of $\fX$.
\end{example}
\begin{proof}
Let us first observe that $\fX$ is a connected affine smooth formal $R$-scheme of ff type, cf.\ Section \ref{smoothnesssec}. Viewing the special fiber $\fX_k$ as a smooth formal $\Spf k[[Z]]$-scheme of ff type, it makes sense to call $\fX_{kk}\mathrel{\mathop:}=V(\pi,Z)$ the \emph{very special fiber} of $\fX$. We thus use a subscript $kk$ to indicate reduction modulo $\pi$ and $Z$. The biggest subscheme of definition of $\fX$ is a transversal intersection of two lines,
\[
\fX_\red\,=\,V(\pi,XY,Z)\,=\,\Spec k[X,Y]/(XY)\quad,
\]
so we see that $\fU$ has indeed two affine connected components $\fU\cap D(X)$ and $\fU\cap D(Y)$; in particular, $\fV$ is well-defined.

Let us assume that there exists an admissible formal blowup $\phi$ of $\fX$ extending $b_\sI$. Let $A$ denote the ring of global functions on $\fX$. By Lemma \ref{genisolem}, $\phi$ is defined by a $\pi$-adically open ideal $J\subseteq A$ whose image in the domain $A_k$ is nonzero. Let $B_1$ and $B_2$ denote the rings of global functions in $\fU\cap D(X)$ and $\fU\cap D(Y)$ respectively. The restriction homomorphisms $A_k\rightarrow B_{i,k}$ are injective since they are flat homomorphisms of domains. Hence, $J B_{i,k}$ is nonzero for $i=1,2$, and so it follows that $J B_i$ has height $\geq 2$ for $i=1,2$. By Lemma \ref{codimtwoisolocuslem}, it follows that $\rad(J B_{1,k})=(\pi,Y,Z)$ and $\rad(JB_{2,k})=(\pi,Z)$. Hence, $J B_{2,kk}=0$. Since $A_kk\rightarrow B_{2,kk}$ is injective, being a flat homomorphism of domains, it follows that $J A_{kk}=0$. However, this implies $JB_{1,kk}=0$, which contradicts the fact that up to nilpotent elements, $J B_{1,kk}$ agrees with $Y B_{1,kk}$ which has height one.
\end{proof}


\subsubsection{Generic fibers of coherent modules}\label{cohmodgenfibsec}


\begin{defi}\label{semafffuncmodeldef}\index{generic fiber!of a coherent module}
Let $\fX$ be a formal $R$-scheme of locally ff type, and let $\sF$ be a coherent $\O_\fX$-module. The sheaf 
\[
\sF_K\,\mathrel{\mathop:}=\,\sF\otimes_RK
\]
is called the generic fiber of $\sF$. The generic fiber of $\O_\fX$ is called the sheaf of semi-affinoid functions on $\fX$.
\end{defi}

Here we view $R$ and $K$ as constant sheaves on $\fX$, and the tensor product is formed in the category of sheaves on the underlying topological space of $\fX$. Since $K$ is $R$-flat and since $\fX$ is locally noetherian, the sheaf tensor product essentially coincides with the presheaf tensor product: For any coherent $\O_\fX$-module $\sF$ and any quasi-compact open subset $\fU\subseteq\fX$,
\[
\sF_K(\fU)\,=\,\sF(\fU)\otimes_RK\quad.
\]
If $\fX$ is $R$-flat, the structural sheaf $\O_\fX$ is a subsheaf of $\O_{\fX,K}$. Let us emphasize that the $\O_\fX$-modules $\sF_K$ are \emph{not} coherent.

The following Proposition is an immediate generalization of \cite{LtkeFRG} 2.1 to the ff type situation:

\begin{prop}\label{blowupisoonsemiafffuncprop}
Let $\fX$ be a flat formal $R$-scheme of locally ff type, let $\phi\colon\fX'\rightarrow\fX$ be an admissible formal blowup, let $\sF$ be a coherent $\O_\fX$-module, and let $\sG$ be a coherent $\O_{\fX'}$-module. Then the canonical morphisms
\begin{eqnarray*}
\sF_K&\rightarrow&(\phi_*\phi^*\sF)_K\\
(\phi^*\phi_*\sG)_K&\rightarrow&\sG_K
\end{eqnarray*}
are isomorphisms. Moreover, $(R^n\phi_*\sG)_K=0$ for all $n\geq 1$.
\end{prop}

The \emph{proof} is literally the same as the proof of \cite{LtkeFRG} 2.1.

\subsubsection{Finite admissible formal blowups}

We discuss admissible formal blowups of flat formal $R$-schemes of locally ff type that are \emph{finite}\index{blowup!finite admissible formal} in the sense of \cite{EGAIII1} 4.8.2. We begin by proving an analog of \cite{BL1} 4.5; the proof is not different from the proof in the tf type case. From now on, we will occasionally use the following notation: We will write an underscored latin character to denote an $R$-algebra of ff type, and we will write the same character without underscore to denote its \emph{semi-affinoid generic fiber}, that is, its localization with respect to $\pi$.

\begin{prop}\label{affineblowupprop}
Let $\ul{A}$ be a flat $R$-algebra of ff type with generic fiber $A\mathrel{\mathop:}=\ul{A}\otimes_RK$, and let $\ul{A}'$ be a subring of $A$ containing $\ul{A}$ such that the inclusion $\ul{A}\subseteq\ul{A}'$ is finite. Then the induced finite morphism $b\colon\Spf \ul{A}'\rightarrow \Spf \ul{A}$ is an admissible formal blowup. More precisely speaking, let $f_1,\ldots,f_n\in\ul{A}'$ be elements that generate $\ul{A}'$ as a module over $\ul{A}$, and let $c\in R$ be a nonzero element such that $cf_1,\ldots,cf_n\in\ul{A}$; then $b$ is the admissible formal blowup of $(c,cf_1,\ldots,cf_n)\subseteq \ul{A}$.

\end{prop} 
\begin{proof}
Let $I\subseteq\ul{A}$ denote the ideal $(c,cf_1,\ldots,cf_n)$; then $I$ is $\pi$-adically open. We verify that the natural morphism $\Spf \ul{A}'\rightarrow\Spf \ul{A}$ satisfies the universal property of the formal blowup of $\Spf \ul{A}$ in $I$: The ideal $I\ul{A}'$ is generated by $c$. Since $c$ is not a zero-divisor in $\ul{A}'\subseteq \ul{A}_K$,  it follows that $I\ul{A}'$ is invertible. Let $\ul{D}$ be a flat $R$-algebra of ff type, and let $\ul{\phi}\colon\ul{A}\rightarrow \ul{D}$ be an $R$-homomorphism such that $I\ul{D}$ is invertible; we must show that $\ul{\phi}$ extends uniquely to an $R$-homomorphism $\ul{\phi}'\colon\ul{A'}\rightarrow\ul{D}$. We let $g_i\in D\mathrel{\mathop:}=\ul{D}\otimes_RK$ denote the image of $f_i$ under the induced $K$-homomorphism $\phi\colon A\rightarrow D$. Let us first assume that $I\ul{D}$ is generated by $c$. Then $cg_i\in c\ul{D}$ for all $i$, and hence $g_i\in \ul{D}$ for all $i$ since $c$ is not a zero-divisor in $\ul{D}$. Thus, $\phi$ maps $\ul{A}'$ to $\ul{D}$, and we obtain the desired extension $\ul{\phi}'\colon\ul{A}'\rightarrow \ul{D}$ of $\ul{\phi}$. Let us now consider the general situation. We may work locally on $\Spf \ul{D}$ and thereby assume that $I\ul{D}$ is freely generated by $cg_i$ for some index $i$. In that case, $c\in cg_i \ul{D}$, and hence $g_i\in D$ has an inverse $g_i^{-1}\in \ul{D}$. Now $g_i$ is integral over $\ul{D}$; indeed, the homomorphism $\ul{\phi}$ takes integral equations for $f_i$ over $\ul{A}$ to integral equations for $g_i$ over $\ul{D}$. Multiplying an integral equation for $g_i$ over $\ul{D}$ with a suitable power of $g_i^{-1}$, we see that $g_i\in \ul{D}$. Hence, $g_i$ is a unit in $\ul{D}$, and we are back in the situation where $I\ul{D}$ is generated by $c$.
\end{proof}

We say that a finite morphism of flat formal $R$-schemes of locally ff type is \emph{generically an isomorphism} if it satisfies the equivalent conditions of the following Lemma. Again, we refer to Section \ref{berthsec} for a discussion of Berthelot's construction.

\begin{lem}\label{genisocharlem}
Let $\phi\colon\fY\rightarrow\fX$ be a finite morphism of formal $R$-schemes of locally ff type. Then the following are equivalent:
\begin{packed_enum}
\item The rigid generic fiber $\phi^\rig$ of $\phi$ is an isomorphism.
\item The associated morphism of sheaves of semi-affinoid functions
\[
\phi^\sharp_K\colon\O_{\fX,K}\rightarrow\phi_*\O_{\fY,K}
\]
is an isomorphism.
\end{packed_enum}
\end{lem}
\begin{proof}
Let us assume that ($i$) holds. To prove statement ($ii$), we may assume that $\fX$ and, hence, $\fY$ are affine. Let us write $\fX=\Spf\ul{A}$, $\fY=\Spf\ul{B}$, and let $\ul{\sigma}\colon\ul{A}\rightarrow\ul{B}$ denote the finite $R$-homomorphism corresponding to $\phi$. Let $\sigma\colon A\rightarrow B$ denote the finite $K$-homomorphism obtained from $\ul{\sigma}$ via $\cdot\otimes_RK$; we must show that $\sigma$ is an isomorphism. For any maximal ideal $\m\subseteq A$, we let $\hat{A}_\m$ denote the $\m$-adic completion of $A$. By  ($i$) and by \cite{dJ} 7.1.9, $\sigma\otimes_A \hat{A}_\m$ is an isomorphism for all maximal ideals $\m\subseteq A$, cf. our discussion before Corollary \ref{rigflatfaithfulcor}. Since the natural morphism $A\rightarrow\bigoplus_{\m\in\Max A} \hat{A}_\m$ is faithfully flat, we conclude that $\sigma$ is an isomorphism, as desired.

Conversely, let us assume that ($ii$) holds. We may again assume that $\fX$ and, hence, $\fY$ are affine. The result then follows from \cite{dJ} 7.1.9, 7.2.2 together with the fact that a morphism of affinoid rigid $K$-spaces that is a bijection on physical points and that induces isomorphisms of complete stalks is an isomorphism, cf.\ \cite{BGR} 7.3.3/5 and 8.2.1/4.
\end{proof}


\begin{cor}\label{fingenisoisblowupcor}
Let $\phi$ be a finite morphism of flat formal $R$-schemes of ff type that is generically an isomorphism. Then $\phi$ is a finite admissible formal blowup.
\end{cor}
\begin{proof}
Let $\fY$ and $\fX$ denote the domain and the target of $\phi$ respectively. Since $\fX$ and $\fY$ are $R$-flat and since $\phi$ is finite and generically an isomorphism, we may, by means of $\phi^\sharp$, consider $\phi_*\O_\fY$ as a coherent $\O_\fX$ module containing $\O_\fX$ as a submodule. Since $\fX$ is quasi-compact, Proposition \ref{affineblowupprop} shows that there exists an element $c$ in $R$ such that the $\O_\fX$-module homomorphism
\[
\phi_*\O_{\fY} \rightarrow \phi_*\O_{\fY}\quad,\quad f\mapsto c\cdot f
\] 
has image in $\O_\fX$. This image is a coherent ideal $\sI\subseteq\O_\fX$, and by Proposition \ref{affineblowupprop}, $\phi$ coincides with the admissible formal blowup of $\fX$ in $\sI$.
\end{proof}

\begin{cor}\label{normisblowupcor}\index{normalization!by finite blowups}
Let $\fX$ be a generically normal flat formal $R$-scheme of ff type. Then the normalization of $\fX$ is a finite admissible formal blowup.
\end{cor}
\begin{proof}
Let $\phi\colon\tilde{\fX}\rightarrow\fX$ denote the normalization of $\fX$. By Corollary \ref{flatgennormalnormalizationlem}, $\tilde{\fX}$ is $R$-flat, and $\phi$ is a finite morphism which is generically an isomorphism. Hence, the statement follows from Corollary \ref{fingenisoisblowupcor}.
\end{proof}



We conclude our discussion of finite admissible formal blowups by establishing results which we will need in order to introduce the uniformly rigid G-topology of a semi-affinoid $K$-space. We begin with a formal version of Stein factorization:

\begin{lem}\label{standardblowupfactlem}
Let $\fX$ be a flat affine formal $R$-scheme of ff type, and let $\phi\colon\fX'\rightarrow\fX$ be an admissible formal blowup. Then $\phi$ admits a unique factorization
\[
\phi\colon\fX'\overset{\sigma}{\longrightarrow}\fX''\overset{\tau}{\longrightarrow} \fX\quad,
\]
where $\sigma$ induces an isomorphism on rings of global sections and where $\tau$ is a finite admissible formal blowup.
\end{lem}
\begin{proof}
Let us write $\fX=\Spf A$, and let $A'$ denote the $R$-flat topological ring of global sections of $\fX'$. By \cite{EGAIII1} 3.4.2, $A'$ is a finite $A$-algebra, and its topology is induced by the topology on $A$. We set $\fX''\mathrel{\mathop:}=\Spf A'$. By \cite{EGAIn} 10.4.6, morphisms to affine formal schemes are given by continuous homomorphisms of rings of global sections; thus we obtain existence and uniqueness of the above factorization. Since $\phi$ is an admissible formal blowup, Proposition \ref{blowupisoonsemiafffuncprop} shows that $\tau$ is generically an isomorphism. By Proposition \ref{affineblowupprop}, it follows that $\tau$ is a finite admissible formal blowup, as desired.
\end{proof}

\begin{cor}\label{affadmblowupisfinitecor}
Let $\fX$ be a flat affine formal $R$-scheme of ff type, and let $\phi\colon\fX'\rightarrow\fX$ be an admissible formal blowup. If $\fX'$ is affine, then $\phi$ is finite.
\end{cor}
\begin{proof}
Let us consider the factorization of $\phi$ that is obtained from Lemma \ref{standardblowupfactlem}. If $\fX'$ is affine, then $\sigma$ is a morphism of affine formal $R$-schemes that induces an isomorphism of rings of global sections and, hence, an isomorphism. It follows that $\phi=\tau$ is a finite admissible formal blowup.
\end{proof}

\begin{cor}\label{affinefineblowupcor}
Let us consider a commutative diagram of flat formal $R$-schemes of ff type
\[
\begin{diagram}
\fZ'&\rTo^{\psi'}&\fX'&\\
\dTo<{\phi'}&&\dTo>{\phi}&\\
\fZ&\rTo^{\psi}&\fX&\quad,
\end{diagram}
\]
where $\fX$, $\fX'$ and $\fZ$ are \emph{affine} and where $\phi'$ is an admissible formal blowup. Let
\[
\phi'\colon\fZ'\overset{\sigma}{\longrightarrow}\fZ''\overset{\tau}{\longrightarrow} \fZ
\]
be the factorization provided by Lemma \ref{standardblowupfactlem}. There exists a unique morphism 
\[
\psi''\colon\fZ''\rightarrow \fX'
\]
such that the diagram
\[
\begin{diagram}
\fZ'&\rTo^{\psi'}&\fX'\\
\dTo<{\sigma}&&\dEqual\\
\fZ''&\rDashto^{\psi''}&\fX'\\
\dTo<{\tau}&&\dTo>{\phi}\\
\fZ&\rTo^{\psi}&\fX
\end{diagram}
\]
commutes.
\end{cor}
\begin{proof}
Indeed, by \cite{EGAIn} 10.4.6, a morphism to a an affine formal scheme corresponds to a continuous homomorphism of rings of global sections, so the statement is obvious.
\end{proof}

		\subsection{Berthelot's construction}\label{berthsec}

In this section, we briefly explain Berthelot's generic fiber construction which we have already referred to above. Let $\fX$ be a formal $R$-scheme of locally ff type, let $\fX_0\subseteq\fX$ denote the smallest subscheme of definition of $\fX$, and for $n\in\N$ let $\fX_n$ denote the $n$-th infinitesimal neighborhood of $\fX_0$ in $\fX$; then $\fX_n$ is a scheme of finite type over $R/\pi^{n+1}R$. Let $\sI\subseteq\O_\fX$ be the biggest ideal of definition of $\fX$ such that $\fX_0$ is defined by $\sI$ and such that $\fX_n$ is defined by $\sI^{n+1}$; then the special fiber $\fX_n\cap\fX_k$ of $\fX_n$ is defined by $(\pi,\fI^{n+1})$. For each $n\in\N$, we let $\fX_{(n)}$ denote the formal dilatation of $\fX$ in $\fX_n\cap\fX_\pi$. The following Lemma is a global version of \cite{dJ} 7.1.2:

\begin{lem}\label{berthconstrbasiclem}
For each $n\in\N$, the formal $R$-scheme $\fX_{(n)}$ is of locally tf type, and it is identified with the formal dilatation of $\fX_{(n+1)}$ along the strict transform of $\fX_n\cap\fX_k$.
\end{lem}
\begin{proof}
To prove the first statement, we may assume that $\fX$ is affine, $\fX=\Spf A$. Then $\fX_{(n)}$ is affine as well, and if $I\subseteq A$ denotes the ideal of global sections of $\sI$, then $\fX_{(n)}$ is identified with the formal spectrum of the $I$-adic completion of the subring $A[I^{n+1}/\pi]$ of $A\otimes_RK$. The $I$-adic topology on $A[I^{n+1}/\pi]$ coinciding with the $\pi$-adic topology, the formal $R$-scheme of ff type $\fX_{(n)}$ is adic over $R$ and, hence, of tf type. The second statement follows immediately from the universal properties of formal dilatations, cf.\ Proposition \ref{admformaldilprop}.
\end{proof}

In particular, Raynaud's generic fiber functor $\rig$ induces open immersions 
\[
\fX_{(n)}^\rig\hookrightarrow\fX_{(n+1)}^\rig\quad.
\]
One defines
\[
\fX^\rig\,\mathrel{\mathop:}=\,\bigcup_{n\in\N}\fX_{(n)}^\rig\quad;
\]
the rigid $K$-space $\fX^\rig$ is called the \emph{rigid generic fiber} or the \emph{Berthelot generic fiber} of $\fX$. Since morphisms of formal $R$-schemes of locally ff type restrict to morphisms of subschemes of definition, it is clear from our above discussion that $\rig$ defines a functor. It generalizes Raynaud's $\rig$-functor: Indeed, to verify this claim it suffices to consider the situation where $\fX=\Spf A$ is affine of tf type over $R$; if $I$ denotes the biggest ideal of definition of $A$, then $I^{n_0}\subseteq\pi A$ for some $n_0\in\N$, the $\pi$-adic topology on $A$ coinciding with the $I$-adic one. Hence, $\fX_{(n)}=\Spf A/(\textup{$\pi$-torsion})$ for all $n\geq n_0$.

The functorial points of the rigid generic fiber of $\fX$ with values in quasi-compact quasi-separated rigid $K$-spaces are intrinsically characterized as follows, cf.\ also \cite{dJ} 7.1.7:

\begin{lem}\label{funcpoitslem}
Let $\fY$ be a flat formal $R$-scheme of tf type; then every morphism $\phi\colon\fY^\rig\rightarrow\fX^\rig$ extends to a morphism $\ul{\phi}'\colon\fY'\rightarrow\fX$, where $\fY'$ is a suitable admissible formal blowup of $\fY$. Moreover, if $\ul{\phi}''\colon\fY''\rightarrow\fX$ is another such extension, there exists an admissible formal blowup $\fY'''$ of $\fY$ dominating both $\fY'$ and $\fY''$ such that $\ul{\phi}'$ and $\ul{\phi}''$ coincide on $\fY'''$.
\end{lem}
\begin{proof}
Since $\fY^\rig$ is quasi-compact, $\phi$ factorizes over $\fX_{(n)}^\rig$ for $n$ big enough. By \cite{BL1} 4.1, there exists an admissible formal blowup $\fY'$ of $\fY$ together with an extension $\ul{\tilde{\phi}}'\colon\fY'\rightarrow\fX_{(n)}$ of $\phi$; the morphism $\ul{\phi}'$ is thus obtained by composing $\ul{\tilde{\phi}}'$ with the dilatation morphism $\fX_{(n)}\rightarrow\fX$. Let $\ul{\phi}''\colon\fY''\rightarrow\fX$ be any extension of $\phi$, where $\fY''$ is an admissible formal blowup of $\fY$. The strict transform of $\fY''$ under the formal dilatation $\fX_{(n)}\rightarrow\fX$ is an admissible formal blowup, since it induces an isomorphism of rigid generic fibers and since the specialization map $\sp_{\fY''}$ is surjective. After replacing $\fY''$ by this strict transform, we may therefore assume that $\ul{\phi}''$ factorizes over $\fX_{(n)}$. The existence of $\fY'''$ and $\ul{\phi}'''$ now again follows from \cite{BL1} 4.1.
\end{proof}

Lemma \ref{funcpoitslem} allows us to define a functorial specialization map $\sp_\fX$ from $\fX^\rig$ to the set of closed physical points of $\fX$, by considering points of $\fX^\rig$ with values in finite field extensions of $K$. 
If $\fX$ is $R$-flat, \cite{dJ} 7.2.4 and 7.2.5 show that $\sp_\fX$ is surjective.

The naturalness of Berthelot's construction also becomes apparent within the category of adic spaces: Viewing the category of locally noetherian formal schemes and the category of rigid $K$-spaces as full subcategories of the category of adic spaces, the Berthelot generic fiber $\fX^\rig$ is literally the fibered product of $\fX$ over $\Spf R$ with $\Sp K$. This follows by inspection of the proof of \cite{Huberbuch} 1.2.2, cf.\ also \cite{RZ} Remark 5.11.

Let $\fX=\Spf A$ be an affine formal $R$-scheme of ff type. By \cite{dJ} 7.1.9, the natural $K$-homomorphism $A\otimes_RK\rightarrow\Gamma(\fX^\rig,\O_{\fX^\rig})$ induces a bijection from $\Max A\otimes_RK$ onto the set of physical points of $\fX^\rig$, and it induces isomorphisms of completed stalks. Let us note the following useful corollary of this statement:

\begin{cor}\label{rigflatfaithfulcor}
The functor $\rig$ is faithful on the category of \emph{flat} formal $R$-schemes of locally ff type.
\end{cor}
\begin{proof}
Let $\fX$ and $\fY$ be flat formal $R$-schemes of locally ff type, and let $\phi,\psi\colon\fY\rightarrow\fX$ be morphisms whose rigid generic fibers $\phi^\rig$ and $\psi^\rig$ coincide. Since $\sp_\fY$ is surjective, $\phi$ and $\psi$ must coincide as maps of closed points. Since $\fY$ and $\fX$ are Jacobson spaces, $\phi$ and $\psi$ must coincide on all physical points. We may thus assume that $\fX$ and $\fY$ are affine, $\fX=\Spf A$ and $\fY=\Spf B$. Then $\phi$ and $\psi$ correspond to continuous $R$-homomorphisms from $A$ to $B$. Since $B$ is $R$-flat, it suffices to see that the induced $K$-homomorphisms from $A\otimes_RK$ to $B\otimes_RK$ coincide. This now follows from $\phi^\rig=\psi^\rig$ by using \cite{dJ} 7.1.9 quoted above together with Krull's Intersection Theorem.
\end{proof}

Let us note that in the proof of Corollary \ref{rigflatfaithfulcor}, we only used $R$-flatness of $\fY$.

	\section{Semi-affinoid algebras}\label{semaffalg}

		We now dispose of the requisite concepts from formal geometry in order to define semi-affinoid $K$-algebras and to study the main aspects of their commutative algebra. As before, $K$ denotes the fraction field of a complete discrete valuation ring $R$ with residue field $k$, and $\pi\in R$ is a uniformizer.		

		\subsection{Lattices of ff type}
			If $A$ is a $K$-algebra, an \emph{$R$-lattice}\index{lattice@$R$-lattice!of a $K$-algebra} for $A$ is an $R$-subalgebra $\underline{A}\subseteq A$ of $A$ such that the natural morphism $\underline{A}\otimes_RK\rightarrow A$ is an isomorphism. For example, $A$ is an $R$-lattice of $A$. Clearly $A$ is a domain if and only if it admits an $R$-lattice that is a domain, and then any $R$-lattice of $A$ is a domain.

\begin{defi}\label{semaffalgdefi}
Let $A$ be a $K$-algebra.
\begin{enumerate}
\item $A$ is called semi-affinoid if it admits an $R$-lattice of ff type.\index{algebra!semi-affinoid}
\item If $\underline{A}$ is an $R$-lattice of ff type for $A$ and if $(\underline{s},\underline{t})$ is a formal generating system of $\underline{A}$, we say that $(\underline{s},\underline{t})$ is a semi-affinoid generating system of $A$.\index{generating system!semi-affinoid}
\item If $\underline{A}$ is an $R$-lattice of ff type for $A$ and if $\ul{\phi}$ is a formally finite presentation of $\ul{A}$, we say that $\ul{\phi}\otimes_RK$ is a semi-affinoid presentation of $A$.\index{presentation!semi-affinoid}
\item The category of semi-affinoid $K$-algebras is the full subcategory of the category of $K$-algebras defined by the property that its objects are the semi-affinoid $K$-algebras.\index{category!of semi-affinoid $K$-algebras}
\end{enumerate}
\end{defi}

\begin{lem}\label{semiaffquotientlem}
Semi-affinoid $K$-algebras are noetherian, and any quotient of a semi-affinoid $K$-algebra, naturally considered as a $K$-algebra, is semi-affinoid again.
\end{lem}
\begin{proof}
Let $A$ be a semi-affinoid $K$-algebra, and let $\ul{A}$ be an $R$-lattice of ff type for $A$; then $\ul{A}$ is noetherian. Since $A$ is a localization of $\ul{A}$, it follows that $A$ is noetherian. Let $I\subseteq A$ be an ideal; we set $\ul{I}\mathrel{\mathop:}=I\cap\ul{A}$. Then $\ul{I}$ is an ideal in $\ul{A}$, and the natural homomorphism $\ul{A}/\ul{I}\rightarrow A/I$ is injective. Moreover, $I=\ul{I}\cdot A$. Indeed, since $I\subseteq A$ is finitely generated and since $\pi$ is a unit in $A$, there exists a generating system of $I$ that is contained in $\ul{A}$. Hence $(\ul{A}/\ul{I})\otimes_RK=A/I$, and we see that $\ul{A}/\ul{I}$ is an $R$-lattice of $A/I$. Since $\ul{A}/\ul{I}$ is of ff type over $R$, it follows that the $K$-algebra $A/I$ is semi-affinoid.
\end{proof}

			\subsubsection{Some remarks on Noether normalization}\label{noethersec}
				The commutative algebra of affinoid $K$-algebras is studied my means of the Noether Normalization Theorem, cf.\ \cite{BGR} 6.1.2/1. There is no general analog of this theorem in the semi-affinoid setup, as it is shown by the following example, cf.\ \cite{LiRo} I.2.3.5:

\begin{example}\label{nonoetherex}
Let us consider the semi-affinoid $K$-algebra
\[
A\mathrel{\mathop:}=R[[X]]\langle Y\rangle/(XY)\otimes_RK\quad.
\]
There exists no finite $K$-monomorphism 
\[
\phi\colon R[[S_1,\ldots,S_m]]\langle T_1,\ldots,T_n\rangle\otimes_RK\rightarrow A\quad.
\]
\end{example}
\begin{proof}
It is easily verified that $A$ is one-dimensional. Hence, if $\phi$ exists, then either $m=1$ and $n=0$, or $m=0$ and $n=1$. Let us assume that $\phi$ exists with $m=1$. Composing $\phi$ with the projection sending $X$ to zero, we obtain a finite $K$-homomorphism $R[[S]]\otimes_RK\rightarrow R\langle Y\rangle\otimes_RK$. By \cite{BGR} 3.8.1/7, by Remark \ref{pbfuncrem} and by Proposition \ref{normallatticeprop} proven below, it restricts to an integral $R$-homomorphism $R[[S]]\rightarrow R\langle Y\rangle$ which, then again, reduces modulo $\pi$ to an integral $k$-homomorphism $k[[S]]\rightarrow k[Y]$. Since $k[[S]]$ and $k[Y]$ are domains of equal dimension, this $k$-homomorphism must be injective, which is impossible, the fraction field of $k[[S]]$ having infinite transcendence degree over $k$. Conversely, let us assume that $n=1$. Arguing along the same lines, we obtain an integral $k$-monomorphism $k[T]\rightarrow k[[X]]$. However, since the fraction field of $k[[X]]$ has infinite transcendence degree over $k$, it cannot be algebraic over $k(T)$, so we again arrive at a contradiction.
\end{proof}

In \cite{LiRo} I.2, Weierstrass preparation techniques and partial results on Noether normalization are discussed, but we will not need them. 

In order to study the commutative algebra of semi-affinoid $K$-algebras, one may try to reduce to a situation where Noether normalization is available. One strategy consists in reducing to a affinoid situation, by means of Berthelot's construction. Alternatively, one may try to reduce to the case where $A$ admits a \emph{local}\index{lattice@$R$-lattice!local} lattice of ff type, by means of formal completion of a given ff type lattice  along a maximal ideal. Indeed, in the presence of a local lattice, Noether normalization works:

\begin{lem}\label{localnoethnormlem}
Let $A$ be an integral flat local $R$-algebra of ff type. There exists a finite $R$-monomorphism
\[
\phi\colon R[[S_1,\ldots,S_{d-1}]]\hookrightarrow A\quad,
\]
where $d=\dim A$.
\end{lem}
\begin{proof}
Let us choose a system of parameters $s_0,\ldots,s_{d-1}$ of $A$ such that $s_0=\pi$. Since $A$ is maximal-adically complete, there exists a unique continuous $R$-homomorphism
\[
\phi\colon R[[S_1,\ldots,S_{d-1}]]\rightarrow A
\]
sending $S_i$ to $s_i$ for $1\leq i\leq d-1$. Let us write $\ul{S}$ for the system $S_1,\ldots,S_{d-1}$, and let $\fr$ denote the maximal ideal of $R[[\ul{S}]]$; then $A/\fr A$ is $k$-finite since the quotient of the noetherian ring $A$ modulo its maximal ideal $\m$ is $k$-finite and since $\fr A$ is $\m$-primary. Since $R[[\ul{S}]]$ is $\fr$-adically complete and since $A$ is $\fr$-adically separated, it follows that $\phi$ is finite, cf.\ \cite{Eis} Ex. 7.2. The dimensions of $R[[\ul{S}]]$ and $A$ being equal, $\phi$ must be injective, as desired.
\end{proof}

\begin{cor}\label{localnoethnormcor}
Let $A$ be an integral semi-affinoid $K$-algebra admitting a local $R$-lattice of ff type, and let $d$ denote the dimension of $A$. There exists a finite $K$-algebra monomorphism
\[
\phi\colon R[[S_1,\ldots,S_d]]\otimes_RK\hookrightarrow A\quad.
\]
\end{cor}

			\subsubsection{Rig-points and the specialization map}\label{specmapsec}

If $A$ is a semi-affinoid $K$-algebra and if $\ul{A}\subseteq A$ is an $R$-lattice of ff type for $A$, there exists a natural specialization map $\sp_{\ul{A}}$ from the set of maximal ideals in $A$ onto the set of maximal ideals in $\ul{A}$. It is defined as follows, compare \cite{dJ} 7.1.9 and 7.1.10: If $\m\subseteq A$ is a maximal ideal, we set
\[
\sp_{\ul{A}}(\m)\,\mathrel{\mathop:}=\,\sqrt{(\m\cap\ul{A})+\pi\ul{A}}\quad.
\]
\index{specialization map}

\begin{lem}\label{specmapwelldeflem}
The map $\sp_{\ul{A}}$ is well-defined.
\end{lem}
\begin{proof}
Let $\m\subseteq A$ be a maximal ideal; we must check that the ideal $\sp_{\ul{A}}(\m)\subseteq\ul{A}$ is maximal. Let us write $\p\mathrel{\mathop:}=\m\cap\ul{A}$; then $(\ul{A}/\p)_\pi=A/\m$ is a field, and by the Artin-Tate Theorem \cite{EGAIV1} 0.16.3.3 it follows that $\ul{A}/\p$ is a semi-local ring of dimension $\leq 1$. Moreover, $\ul{A}/\p$ is of ff type over $R$ and, hence, $\pi$-adically complete. Since $\ul{A}/\p\subseteq A/\m$ is $R$-flat and since $(\ul{A}/\p)_\pi$ is local, it thus follows from Hensel's Lemma that $(\ul{A}/\p)/\pi(\ul{A}/\p)$ is local as well, cf.\ \cite{Bourbaki} III.4.6 Proposition 8. Since $\p A=\m$, the class of $\pi$ in $\ul{A}/\p$ is nonzero, and so the local noetherian ring $(\ul{A}/\p)/\pi(\ul{A}/\p)$ is zero-dimensional. Thus its quotient modulo its nilradical is a field, and it follows that the radical of $\p+\pi\ul{A}$ is maximal in $\ul{A}$, as desired.
\end{proof}

The proof of Lemma \ref{specmapwelldeflem} shows that if $\m\subseteq A$ is a maximal ideal and if $\p$ denotes $\m\cap\ul{A}$, then $\ul{A}/\p$ is a one-dimensional local domain of ff type over $R$ whose fraction field $K'\mathrel{\mathop:}=(\ul{A}/\p)_\pi=A/\m$ is an extension of $K$. 

\begin{lem}\label{finiteresiduefieldlem}
The field extension $K'/K$ is finite.
\end{lem}
\begin{proof}
It suffices to show that $\ul{A}/\p$ is finite over $R$. Since $R$ is $\pi$-adically complete and since $\ul{A}/\p$ is $\pi$-adically separated, it suffices to show that $\ul{A}/(\p+\pi\ul{A})$ is finite over $k$, cf. \cite{Eis} Ex. 7.2. The ring $\ul{A}/(\p+\pi\ul{A})$ being noetherian, its nilradical is nilpotent; hence it suffices to show that the quotient of $\ul{A}$ modulo the maximal ideal $\sqrt{\p+\pi\ul{A}}$ is $k$-finite. Since $\ul{A}$ is of ff type over $R$, since maximal ideals are open and since field extensions of finite type are finite, the desired statement follows.
\end{proof}

We immediately obtain the following:

\begin{cor}\label{saffmaxpbismalcor}
Let $\phi\colon A\rightarrow B$ be a $K$-homomorphism of semi-affinoid $K$-algebras, and let $\m\subseteq B$ be a maximal ideal; then the $\phi$-preimage of $\m$ in $A$ is maximal as well.
\end{cor}
\begin{proof}
Let us write $\p\mathrel{\mathop:}=\phi^{-1}(\m)$; then $\phi$ induces inclusions of domains $K\subseteq A/\p \subseteq B/\m$, where $K$ is a field and where $B/\m$ is a field which, by Lemma \ref{finiteresiduefieldlem}, is finite over $K$. If follows that $A/\p$ is a field.
\end{proof}

Moreover, it follows that the absolute value on $K$ uniquely extends to an absolute value on $K'$. Let $R'$ denote the integral closure of $\ul{A}/\p$ in $K'$; we easily see from the above discussion that $R'/R$ is a finite local extension of discrete valuation rings corresponding to the extension of discretely valued fields $K'/K$.

The above discussion also shows that the maximal ideals in $A$ correspond bijectively to the $\rig$-points of $\ul{A}$, which are equivalence classes of $R$-homomorphisms
\[
r\colon\ul{A}\rightarrow R'\quad,
\]
where $R'/R$ is a finite local extension of discrete valuation rings and where two such homomorphisms are called equivalent if they are dominated by a third one. Our discussion also implies that if $\m\subseteq A$ is a maximal ideal with associated rig-point $r_\m\colon\ul{A}\rightarrow R'$, then $\sp_{\ul{A}}(\m)$ is the $r_\m$-preimage of the valuation ideal in $R'$. 
It easily follows that $\sp_{\ul{A}}$ is functorial in $\ul{A}$:

\begin{cor}\label{specmapfunctorialcor}
Let $\ul{\phi}\colon\ul{A}\rightarrow\ul{B}$ be an $R$-homomorphism of $R$-lattices of ff type of semi-affinoid $K$-algebras $A$ and $B$, let $\phi\colon A\rightarrow B$ denote the generic fiber of $\ul{\phi}$, and let $\m\subseteq B$ be a maximal ideal. Then 
\[
\sp_{\ul{A}}(\phi^{-1}(\m))\,=\,\ul{\phi}^{-1}(\sp_{\ul{B}}(\m))\quad.
\]
\end{cor}

\begin{lem}\label{spsurjlem}
The map $\sp_{\ul{A}}$ is surjective onto the maximal ideals of $\ul{A}$.
\end{lem}
\begin{proof}
Let $\n\subseteq\ul{A}$ be a maximal ideal, and let $\ul{\phi}\colon\ul{A}\rightarrow\ul{B}$ denote the $\n$-adic completion of $\ul{A}$; then $\ul{\phi}$ is flat and of ff type. Let $\phi\colon A\rightarrow B$ denote the $K$-homomorphism obtained from $\ul{\phi}$ by inverting $\pi$. Since $\ul{B}$ is $R$-flat, the semi-affinoid $K$-algebra $B$ is nonzero; let $\m'\subseteq B$ be a maximal ideal, and let $\m\subseteq A$ denote its $\phi$-preimage. 
Since $\n$ is the $\ul{\phi}$-preimage of the unique maximal ideal of $\ul{B}$, Corollary \ref{specmapfunctorialcor} implies that $\sp_{\ul{A}}(\m)=\n$. Hence, $\sp_{\ul{A}}$ is surjective.
\end{proof}

More generally, these arguments show that if $\ul{A}$ is an $R$-algebra of ff type, if $I\subseteq\ul{A}$ is any ideal and if $\ul{A}|_I$ denotes the completion of $\ul{A}$ along $I$, then the induced map $\Max(\ul{A}|_I\otimes_RK)\rightarrow\Max(\ul{A}\otimes_RK)$ is injective onto $\sp_{\ul{A}}^{-1}(V(I))$.
				
			\subsubsection{Power-bounded functions}\label{topqnilpotsec}
				Let $A$ be a semi-affinoid $K$-algebra, let $f\in A$ be an element in $A$, and let $x\in \Max A$ be a maximal ideal in $A$. By Lemma \ref{finiteresiduefieldlem} and the subsequent discussion, the absolute value $|f(x)|$ of $f$ in $x$ is well-defined. From now on, the elements of a semi-affinoid $K$-algebra will be called \emph{$K$-semi-affinoid functions} or simply semi-affinoid functions in $A$,\index{function!semi-affinoid} and we view them as functions on $\Max A$. If 
\[
\phi^*\colon A\rightarrow B
\]
is a $K$-homomorphism of semi-affinoid $K$-algebras and if $\phi\colon\Max B\rightarrow\Max A$ is the associated map of maximal spectra, then clearly $|(\phi^*f)(x)|=|f(\phi(x))|$ for all $f\in A$, $x\in\Max B$.

\begin{defi}\label{pbtqndefi}
Let $A$ be a semi-affinoid $K$-algebra, and let $f\in A$ be a semi-affinoid function. Then
\begin{packed_enum}
\item $f$ is called \emph{power-bounded} if $|f(x)|\leq 1$ for all $x\in\Max A$, and\index{function!power-bounded}
\item $f$ is called \emph{topologically quasi-nilpotent} if $|f(x)|<1$ for all $x\in\Max A$.\index{function!topologically quasi-nilpotent}
\end{packed_enum}
We let $\mathring{A}$ and $\check{A}$ denote the subsets of power-bounded and topologically quasi-nilpotent elements in $A$ respectively.
\end{defi}
For example, the coordinate function $S\in R[[S]]\otimes_RK$ is topologically quasi-nilpotent, while
\[
\sup_{x\in\Max (R[[S]]\otimes_RK)}|S(x)|=1\quad.
\]

Let us note some immediate observations:
\begin{remark}\label{pbfuncrem}Let $A$ be a semi-affinoid $K$-algebra. 
\begin{packed_enum}
\item $\mathring{A}\subseteq A$ is an $R$-subalgebra, and $\check{A}\subseteq\mathring{A}$ is an ideal.
\item If $\phi\colon A\rightarrow B$ is a homomorphism of semi-affinoid $K$-algebras, then $\phi(\mathring{A})\subseteq\mathring{B}$, and $\phi(\check{A})\subseteq\check{B}$.
\item $\mathring{A}$ contains the nilradical of $A$. 
\end{packed_enum}
\end{remark}

\begin{lem}\label{topquasinilpotinlatticelem}
Let $A$ be a semi-affinoid $K$-algebra, and let $\underline{A}\subseteq A$ be an $R$-lattice of ff type. Then $\ul{A}\subseteq\mathring{A}$, and $\check{A}\cap\ul{A}$ is the biggest ideal of definition of $\ul{A}$.
\end{lem}
\begin{proof}
By Lemma \ref{jacobsontoplem}, the biggest ideal of definition of $\ul{A}$ equals the Jacobson radical of $\ul{A}$. Let $f$ be an element in $\ul{A}$. By Lemma \ref{spsurjlem}, the specialization map $\sp_{\ul{A}}$ is surjective onto the maximal ideals of $\ul{A}$, so we must show that for every $x\in\Max A$ with specialization $\n\subseteq\ul{A}$, we have 
\begin{packed_enum}
\item $|f(x)|\leq 1$ and
\item $|f(x)|<1$ if and only if $f\in\n$.
\end{packed_enum}
Let $x$ be an element of $\Max A$ with specialization $\n$, let $r_x\colon\underline{A}\rightarrow R'$ be a representative of the rig-point given by $x$, and let $\fr'$ denote the maximal ideal of $R'$. Then $|f(x)|=|r_x(f)|\leq 1$, and $|f(x)|<1$ if and only if $r_x(f)\in\fr'$. Since  $r_x^{-1}(\fr')=\n$, the latter holds if and only if $f\in\n$, as desired.
\end{proof}

\begin{prop}\label{normallatticeprop}\index{lattice@$R$-lattice!normal}
Let $A$ be a semi-affinoid $K$-algebra, and let $\ul{A}$ be a normal $R$-lattice of ff type for $A$. Then $\ul{A}=\mathring{A}$.
\end{prop}
\begin{proof}
One could give a direct proof using Lemma \ref{localnoethnormcor}, but for the sake of brevity, we invoke Berthelot's construction, cf. Section \ref{berthsec}, and use results of de Jong: By \cite{dJ} 7.1.9, we may view $A$ as a subring of the ring of global functions on $(\Spf\ul{A})^\rig$, and by \cite{dJ} 7.4.1, $\ul{A}$ coincides with the ring of power-bounded global functions on $(\Spf\ul{A})^\rig$ under this identification. Hence, every power-bounded function in $A$ is contained in $\ul{A}$, as desired.
\end{proof}

\begin{prop}\label{powerboundedintegralprop}
Let $A$ be a semi-affinoid $K$-algebra, and let $\ul{A}\subseteq A$ be an $R$-lattice of ff type. The inclusion of $R$-subalgebras $\ul{A}\subseteq\mathring{A}$ is integral, and it is finite when $A$ is reduced.
\end{prop}
\begin{proof}
Let $\ul{\phi}\colon\ul{A}\rightarrow\ul{\tilde{A}}$ denote the normalization of $\ul{A}$, and let $\phi\colon A\rightarrow\tilde{A}$ denote the homomorphism of semi-affinoid $K$-algebras that is obtained from $\ul{\phi}$ by inverting $\pi$. Since $\ul{\phi}$ factorizes over an injective $R$-homomorphism 
\[
\ul{A}/\rad(\ul{A})\hookrightarrow\ul{\tilde{A}}\quad,
\]
since $K$ is $R$-flat and since $\rad(A)=\rad(\ul{A}) A$, it follows that $\phi$ factorizes over an injective $K$-homomorphism 
\[
A/\rad(A)\hookrightarrow\tilde{A}\quad.
\]
By Proposition \ref{normallatticeprop},
\[
(\tilde{A})^\circ=\ul{\tilde{A}}\quad.
\]
Let us consider a power-bounded semi-affinoid function $f\in\mathring{A}$; then $\phi(f)\in\ul{\tilde{A}}$. Since $\ul{\phi}$ is finite, there exists an integral equation $P(T)\in\ul{A}[T]$ for $\phi(f)$ over $\ul{A}$. By the factorization of $\phi$ mentioned above, we conclude that $P(f)\in A$ is nilpotent. Let $s\in\N$ be an integer such that $P(f)^s=0$; then $P(T)^s\in\ul{A}[T]$ is an integral equation for $f$ over $\ul{A}$.

If $A$ is reduced, $\phi$ is injective, and hence $\mathring{A}$ is an $\ul{A}$-submodule of the finite $\ul{A}$-module $\ul{\tilde{A}}$. Since $\ul{A}$ is noetherian, it follows that $\mathring{A}$ is a finite $\ul{A}$-module. 
\end{proof}

If $A$ is not reduced, then $\mathring{A}$ cannot be  finite over any $R$-lattice of ff type for $A$. Indeed, let $f\in A$ be a nonzero nilpotent element in $A$; then $f$ is infinitely $\pi$-divisible in $\mathring{A}$. If $\mathring{A}$ was finite over an $R$-lattice of ff type, then $\mathring{A}$ would itself be of ff type over $R$ and, hence, $\pi$-adically separated; we would thus arrive at a contradiction.

\begin{cor}\label{latticeexcor}
Let $A$ be a semi-affinoid $K$-algebra.
\begin{packed_enum}
\item If $A$ is reduced, then $\mathring{A}$ is the biggest $R$-lattice of ff type for $A$.
\item Let $\ul{A}$ be an $R$-lattice of ff type for $A$, and let $M\subseteq\mathring{A}$ be a finite set of power-bounded semi-affinoid functions in $A$. Then the subring $\ul{A}[M]$ of $A$ generated by $\ul{A}$ and the elements in $M$ is an $R$-lattice of ff type for $A$.
\end{packed_enum}
\end{cor}
\begin{proof}
Statement ($i$) is clear by Proposition \ref{powerboundedintegralprop}; let us show statement ($ii$). By Proposition \ref{powerboundedintegralprop}, the finitely many elements of $M$ are integral over $\ul{A}$, and hence $\ul{A}[M]$ is finite over $\ul{A}$. It follows that $\ul{A}[M]$ is of ff type over $R$.
\end{proof}

If $\ul{A}$ is an $R$-lattice of ff type for a semi-affinoid $K$-algebra $A$ and if $\ul{A}'\subseteq A$ is any $R$-algebra of ff type, then the sum $\ul{A}+\ul{A}'$ is an $R$-lattice of ff type for $A$ containing both $\ul{A}$ and $\ul{A}'$. Building upon this easy observation, we can establish the universal properties of the \emph{free} semi-affinoid $K$-algebras:\index{algebra!free semi-affinoid}

\begin{prop}\label{freesemiaffinoidprop}
Let $A$ be a semi-affinoid $K$-algebra, let $f_1,\ldots,f_m\in A$ be topologically quasi-nilpotent, and let $g_1,\ldots,g_n\in A$ be power-bounded. There exists a unique $K$-homomorphism
\[
\phi\colon R[[S_1,\ldots,S_m]]\langle T_1,\ldots,T_n \rangle\otimes_RK\rightarrow A
\]
sending $S_i$ to $s_i$ and $T_j$ to $t_j$, for $1\leq i\leq m$ and $1\leq j\leq n$.
\end{prop}

\begin{proof}
Let us write $\ul{S}$ and $\ul{T}$ to denote the systems of the $S_i$ and the $T_j$. By Corollary \ref{latticeexcor}, $A$ admits an $R$-lattice of ff type $\ul{A}$ containing the $f_i$ and the $g_j$. By Lemma \ref{topquasinilpotinlatticelem}, the $f_i$ are topologically quasi-nilpotent in $\ul{A}$; hence there exists a \emph{unique} $R$-homomorphism $\ul{\phi}\colon R[[\ul{S}]]\langle\ul{T}\rangle\rightarrow\ul{A}$ sending $S_i$ to $f_i$ and $T_j$ to $g_j$ for all $i$ and $j$, so $\phi\mathrel{\mathop:}=\ul{\phi}\otimes_RK$ is a $K$-homomorphism with the desired properties. It remains to show that these properties determine $\phi$ uniquely. Let $\phi'\colon R[[\ul{S}]]\langle\ul{T}\rangle\otimes_RK\rightarrow A$ be any $K$-homomorphism sending $S_i$ to $f_i$ and $T_j$ to $g_j$ for all $i$ and $j$, and let $\ul{A}'$ denote the $\phi'$-image of $R[[\ul{S}]]\langle\ul{T}\rangle$. Then $\ul{A}'\subseteq A$ is an $R$-subalgebra of ff type, so after replacing $\ul{A}$ by $\ul{A}+\ul{A}'$ we may assume that $\phi'$ restricts to an $R$-homomorphism $\ul{\phi}'\colon R[[\ul{S}]]\langle\ul{T}\rangle\rightarrow\ul{A}$. But then $\ul{\phi}'=\ul{\phi}$ by uniqueness of $\ul{\phi}$, and hence $\phi'=\phi$, as desired.
\end{proof}

\begin{cor}\label{freesemaffcor}
Let $A$ be a semi-affinoid $K$-algebra.
\begin{enumerate}
\item Let $\ul{A}_1$, $\ul{A}_2$ be $R$-lattices of ff type for $A$. If $\ul{A}_2$ contains a formal generating system of $\ul{A}_1$, then $\ul{A}_1$ is contained in $\ul{A}_2$.
\item Let $\ul{A}_1\subseteq\ul{A}_2$ be an inclusion of $R$-lattices of ff type for $A$. Then $\ul{A}_2$ is finite over $\ul{A}_1$. Moreover, the corresponding $R$-morphism $\Spf\ul{A}_2\rightarrow\Spf\ul{A}_1$ is a finite admissible formal blowup.
\item Let $\phi\colon A\rightarrow B$ be a homomorphism of semi-affinoid $K$-algebras, and let $\ul{A}\subseteq A$, $\ul{B}\subseteq B$ be $R$-lattices of ff type such that there exists a formal generating system of $\ul{A}$ mapping to $\ul{B}$ via $\phi$. Then $\phi(\ul{A})\subseteq\ul{B}$.
\item Let $\phi\colon A\rightarrow B$ be a homomorphism of semi-affinoid $K$-algebras, and let $\ul{A}$ be an $R$-lattice of ff type for $A$. There exists an $R$-lattice of ff type $\ul{B}$ for $B$ such that $\phi(\ul{A})\subseteq\ul{B}$. Moreover, if $\ul{B}'$ is any $R$-lattice of ff type for $B$, we can choose $\ul{B}$ such that $\ul{B}'\subseteq\ul{B}$. 
\end{enumerate}
\end{cor}
\begin{proof}
Let $(\ul{f},\ul{g})$ be a formal generating system of $\ul{A}_1$ that is contained in $\ul{A}_2$, and let $\ul{\phi}$ be the associated formally finite presentation of $\ul{A}_1$. By Proposition \ref{freesemiaffinoidprop} and its proof, $\ul{\phi}$ factorizes through $\ul{A}_2$. Since $\ul{A}_1=\im\ul{\phi}$, it follows that $\ul{A}_1$ is contained in $\ul{A}_2$. We have thus shown statement ($i$).

Let us prove ($ii$). Let $\ul{A}_1\subseteq\ul{A}_2$ be an inclusion of $R$-lattices of ff type for $A$, and let $M\subseteq\ul{A}_2$ be a finite set whose elements are the components of a formal generating system for $\ul{A}_2$ over $R$. Then $\ul{A}_1[M]\subseteq\ul{A}_2$ is an $R$-lattice of ff type for $A$, and by Proposition \ref{powerboundedintegralprop} it is finite over $\ul{A}_1$. By part ($i$) of this corollary, we know that $\ul{A}_2$ is contained in $\ul{A}_1[M]$; hence we obtain $\ul{A}_2=\ul{A}_1[M]$. Thus $\ul{A}_2$ is finite over $\ul{A}_1$. By Proposition \ref{affineblowupprop}, the induced morphism of formal spectra is a finite admissible formal blowup.


Let us prove part $(iii)$. We choose a formal generating system  $(\ul{f},\ul{g})$ of $\ul{A}$ such that the components of $\phi(\ul{f})$ and $\phi(\ul{g})$ are contained in $\ul{B}$. The components of $\phi(\ul{f})$ are nilpotent in $\ul{B}$ since they are topologically quasi-nilpotent in $B$. Composing $\phi$ with the semi-affinoid presentation of $A$ corresponding to the formal generating system $(\ul{f},\ul{g})$ of $\ul{A}$, we see from Proposition \ref{freesemiaffinoidprop} and its proof that $\phi(\ul{A})\subseteq\ul{B}$.

Let us prove statement ($iv$). Let $\ul{A}\subseteq A$ be an $R$-lattice of ff type, and let $(\ul{f},\ul{g})$ be a formal generating system of $\ul{A}$. Then the components of $\phi(\ul{f})$ are topologically quasi-nilpotent, and the components of $\phi(\ul{g})$ are power-bounded in $B$. According to Corollary \ref{latticeexcor}, there exists an $R$-lattice $\ul{B}$ of ff type for $B$ containing $\ul{B}'$ and the components of $\phi(\ul{f})$ and $\phi(\ul{g})$, and it follows from part ($iii$) of this corollary that $\phi(\ul{A})\subseteq\ul{B}$.
\end{proof}

It follows that if $\ul{A}_1$ and $\ul{A}_2$ are $R$-lattices of ff type for a semi-affinoid $K$-algebra $A$, there exists an $n\in\N$ such that $\pi^n\ul{A}_2\subseteq\ul{A}_1$. Indeed, after replacing $\ul{A}_2$ by $\ul{A}_1+\ul{A}_2$ we may assume that $\ul{A}_1\subseteq\ul{A}_2$, and then the statement follows from Corollary \ref{freesemaffcor} \textup{(}$ii$\textup{)}.

\begin{cor}\label{afflattfintypecor}
Let $A$ be an \emph{affinoid} $K$-algebra, and let $\ul{A}\subseteq A$ be an $R$-lattice of ff type. Then $\ul{A}$ is of tf type over $R$.
\end{cor}
\begin{proof}
Let $\ul{A}'$ be an $R$-lattice of tf type for $A$, and let $\ul{A}''$ be an $R$-lattice of ff type for $A$ containing both $\ul{A}$ and $\ul{A}'$. By Corollary \ref{freesemaffcor} ($ii$), $\ul{A}''$ is finite over $\ul{A}'$ and, hence, an $R$-algebra of tf type. After replacing $\ul{A}'$ by $\ul{A}''$, we may thus assume that $\ul{A}\subseteq\ul{A}'$. Again by Corollary \ref{freesemaffcor} ($ii$), this inclusion is finite. Let $J'$ denote the Jacobson radical of $\ul{A}'$; then by \cite{MatsumuraCRT} \S 9 Lemma 2, $J'\cap\ul{A}$ is the Jacobson radical $J$ of $\ul{A}$. By Lemma \ref{jacobsontoplem}, these Jacobson radicals agree with the biggest ideals of definition. Reducing modulo Jacobson radicals, we obtain a finite inclusion $\ul{A}/J\subseteq\ul{A}'/J$. Since $\ul{A}'/J$ is finitely generated over $k$, it follows from the Artin-Tate Lemma that $\ul{A}/J$ is finitely generated over $k$ as well. Hence $\ul{A}$ is of tf type over $R$, as claimed.
\end{proof}

If $A$ is a semi-affinoid $K$-algebra and if $\ul{B}\subseteq A$ is an $R$-subalgebra of ff type such that $\ul{B}$ contains the components of a semi-affinoid generating system for $A$, then $\ul{B}$ needs of course not be an $R$-lattice for $A$. For example, $R\langle T\rangle\subseteq R[[T]]\otimes_RK$ is not an $R$-lattice.

		\subsection{Amalgamated sums}\label{amalgsumssec}
			\begin{prop}\label{amalgsumsprop}
Let $A$, $B_1$ and $B_2$ be semi-affinoid $K$-algebras, and let 
\begin{eqnarray*}
&\phi_1\colon A\rightarrow B_1&\\
&\phi_2\colon A\rightarrow B_2&
\end{eqnarray*}
be $K$-homomorphisms. Then there exists a semi-affinoid $K$-algebra $B_1\hat{\otimes}_AB_2$ together with a commutative diagram
\[
\begin{diagram}
A&\rTo^{\phi_1}&B_1\\
\dTo<{\phi_2}&&\dTo>{\sigma_1}\\
B_2&\rTo^{\sigma_2}&B_1\hat{\otimes}_AB_2
\end{diagram}
\]
that is cocartesian in the category of semi-affinoid $K$-algebras. More precisely speaking, if $\ul{A}\subseteq A$, $\ul{B}_1\subseteq B_1$ and $\ul{B}_2\subseteq B_2$ are $R$-lattices of ff type for $A$, $B_1$ and $B_2$ respectively such that $\phi_i(\ul{A})\subseteq\ul{B}_i$ for $i=1,2$, we may set 
\begin{eqnarray*}
B_1\hat{\otimes}_AB_2&\mathrel{\mathop:}=&(\ul{B}_1\hat{\otimes}_{\ul{A}}\ul{B}_2)\otimes_RK\\
\sigma_i&\mathrel{\mathop:}=&\ul{\sigma}_i\otimes_RK\quad,
\end{eqnarray*}
where we view $\ul{B}_i$ as an $\ul{A}$-algebra via the $R$-homomorphism induced by $\phi_i$ and where $\ul{\sigma}_i\colon\ul{B}_i\rightarrow \ul{B}_1\hat{\otimes}_{\ul{A}}\ul{B}_2$ is the natural coprojection, for $i=1,2$.
\end{prop}
\begin{proof}
By Corollary \ref{freesemaffcor} ($iv$), we may choose $R$-lattices $\ul{A}$, $\ul{B}_1$ and $\ul{B}_2$ as in the statement of the proposition. Let $C$ be a semi-affinoid $K$-algebra, and for $i=1,2$ let 
\[
\tau_i\colon B_i\rightarrow C
\]
be a $K$-homomorphism such that $\tau_1\circ\phi_1=\tau_2\circ\phi_2$. By Corollary \ref{freesemaffcor} ($iv$), there exists an $R$-lattice $\ul{C}$ of ff type for $C$ such that $\tau_i(\ul{B}_i)\subseteq\ul{C}$ for $i=1,2$. For $i=1,2$, we let $\ul{\tau}_i\colon\ul{B}_i\rightarrow\ul{C}$ denote the induced $R$-homomorphism. Then $\ul{\tau}_1\circ\ul{\phi}_1=\ul{\tau}_2\circ\ul{\phi}_2$ since the same holds after inverting $\pi$ and since $\pi$ is not a zero divisor in $\ul{A}$. By the universal property of the complete tensor product in the category of $R$-algebras of ff type, there exist a unique $R$-homomorphism 
\[
\ul{\tau}\colon\ul{B}_1\hat{\otimes}_{\ul{A}}\ul{B}_2\rightarrow\ul{C}
\]
such that $\ul{\tau}_i=\ul{\tau}\circ\ul{\sigma}_i$ for $i=1,2$. Setting $\tau\mathrel{\mathop:}=\ul{\tau}\otimes_RK$, we obtain $\tau_i=\tau\circ\sigma_i$ for $i=1,2$. We must show that $\tau$ is uniquely determined. Let 
\[
\tau'\colon(\ul{B}_1\hat{\otimes}_{\ul{A}}\ul{B}_2)\otimes_RK\rightarrow C
\]
be any $K$-homomorphism such that $\tau_i=\tau'\circ\sigma_i$ for $i=1,2$. By Corollary \ref{freesemaffcor} ($iv$), there exists an $R$-lattice $\ul{C}'$ of ff type for $C$ containing $\ul{C}$ such that $\tau'$ restricts to an $R$-morphism $\ul{\tau}'\colon\ul{B}_1\hat{\otimes}_{\ul{A}}\ul{B}_2\rightarrow\ul{C}'$; then $\tau'=\ul{\tau}'\otimes_RK$. It suffices to show that $\ul{\tau}'$ coincides with $\ul{\tau}$ composed with the inclusion $\ul{\iota}:\ul{C}\subseteq\ul{C}'$. For $i=1,2$, the compositions $\ul{\tau}'\circ\ul{\sigma}_i$ and $\ul{\iota}\circ\ul{\tau}\circ\ul{\sigma}_i$ coincide after inverting $\pi$, hence they coincide because $\pi$ is not a zero divisor in $\ul{B}_i$, for $i=1,2$. The universal property of $(\ul{B}_1\hat{\otimes}_{\ul{A}}\ul{B}_2,\ul{\sigma}_1,\ul{\sigma}_2)$ implies that $\ul{\tau}'=\ul{\iota}\circ\ul{\tau}$, as desired.
\end{proof}

		\subsection{The Jacobson property}
			We prove that semi-affinoid $K$-algebras are Jacobson rings; in other words, we show that they satisfy the statement of Hilbert's Nullstellensatz. The corresponding fact for affinoid $K$-algebras is proved by means of Noether normalization techniques, which are not available to us in the semi-affinoid setting, cf.\ Example \ref{nonoetherex}. We reduce  to the Jacobson property of affinoid $K$-algebras by means of Berthelot's construction.

\begin{prop}\label{hilbertprop}
Semi-affinoid $K$-algebras are Jacobson rings.
\end{prop}
\begin{proof}
By Lemma \ref{semiaffquotientlem}, any quotient of a semi-affinoid $K$-algebra is again semi-affinoid as a $K$-algebra; hence it suffices to show that if $A$ is a semi-affinoid $K$-algebra and if $f\in A$ is a semi-affinoid function such that $f(x)=0$ for all $x\in\Max A$, then $f$ is nilpotent. We may divide $A$ by its nilradical and thereby assume that $A$ is reduced. Let $\ul{A}$ be an $R$-lattice of ff type for $A$, and let $X=(\Spf\ul{A})^\rig$ denote the rigid-analytic generic fiber of $\Spf \ul{A}$. Since $A$ is excellent and since rigid $K$-spaces are excellent, \cite{dJ} Lemma 7.1.9 shows that the space $X$ is reduced and that we may view $A$ as a subring of $\Gamma(X,\O_X)$ such that the value of $f$ in a point $x\in X$ agrees with the value of $f$ in the corresponding maximal ideal of $A$. Since $f(x)=0$ for all $x\in X$, we see that $f=0$ as a function on $X$ and, hence, in $A$.
\end{proof}
			
		\subsection{Examples of semi-affinoid $K$-algebras}
			The results in this section are of a complementary  and illustrative nature; they are not essential for our discussion of uniformly rigid spaces.

\subsubsection{Relatively free semi-affinoid algebras}

Let $A$ be a semi-affinoid $K$-algebra. 

\begin{prop}\label{relcloseddiscprop}
There exists a semi-affinoid $K$-algebra $A\langle T\rangle^\sa$ together with a $K$-algebra homomorphism 
\[
\phi\colon A\rightarrow A\langle T\rangle^\sa
\]
and a semi-affinoid function $T\in A\langle T\rangle^\sa$ such that the following universal property is satisfied: For any semi-affinoid $K$-algebra $B$, any $K$-homomorphism $\psi\colon A\rightarrow B$ and any power-bounded function $f\in B$, there exists a unique $K$-homomorphism $\psi'\colon A\langle T\rangle^\sa\rightarrow B$ such that $\psi'(T)=f$ and such that the diagram
\[
\begin{diagram}
A&\rTo^\phi&A\langle T\rangle^\sa\\
&\rdTo<\psi&\dTo>{\psi'}\\
&&B
\end{diagram}
\]
commutes. If $\ul{A}$ is an $R$-lattice of $A$, we may set $A\langle T\rangle^\sa\mathrel{\mathop:}=\ul{A}\langle T\rangle\otimes_RK$.
\end{prop}
\begin{proof}
Let $\ul{A}$ be any $R$-lattice of ff type for $A$, and let $\ul{\phi}\colon\ul{A}\rightarrow\ul{A}\langle T\rangle$ be the natural $R$-homomorphism. We set $\phi\mathrel{\mathop:}=\ul{\phi}\otimes_RK$, and we claim that $\phi$ has the asserted universal property. Let $B$ and $\psi$ be as in the statement of the proposition. We choose an $R$-lattice $\ul{B}$ of ff type for $B$ such that $\psi$ restricts to an $R$-homomorphism $\ul{\psi}\colon\ul{A}\rightarrow\ul{B}$; this is possible by Corollary \ref{freesemaffcor} ($iv$). By Corollary \ref{latticeexcor}, we may assume that $f\in\ul{B}$. There exists a unique $R$-homomorphism $\ul{\psi}'\colon\ul{A}\langle T\rangle\rightarrow\ul{B}$ such that the diagram
\[
\begin{diagram}
\ul{A}&\rTo^{\ul{\phi}}&\ul{A}\langle T\rangle\\
&\rdTo<{\ul{\psi}}&\dTo>{\ul{\psi}'}\\
&&\ul{B}
\end{diagram}
\]
commutes, and setting $\psi'\mathrel{\mathop:}=\ul{\psi}'\otimes_RK$, we obtain a $K$-homomorphism $\psi'$ such that the diagram in the statement of the proposition commutes. It remains to see that $\psi'$ is uniquely determined by this property. Let $\psi''$ be another such morphism. Enlarging $\ul{B}$ such that it contains the $\psi''$-image of a formal generating system of $\ul{A}\langle T\rangle$, Proposition \ref{freesemiaffinoidprop} allows us to assume that $\psi''$ restricts to an $R$-homomorphism $\ul{\psi}''\colon\ul{A}\langle T\rangle\rightarrow\ul{B}$. Now $\ul{\psi}'$ and $\ul{\psi}''$ must coincide, which implies that $\psi'=\psi''$
\end{proof}

Arguing along the same lines, one proves the following:

\begin{prop}\label{relopendiscprop}
There exists a semi-affinoid $K$-algebra $A[[T]]^\sa$ together with a $K$-algebra homomorphism 
\[
\phi\colon A\rightarrow A[[T]]^\sa
\]
and a semi-affinoid function $T\in A[[T]]^\sa$ such that the following universal property is satisfied: For any semi-affinoid $K$-algebra $B$, any $K$-homomorphism $\psi\colon A\rightarrow B$ and any topologically quasi-nilpotent function $f\in B$, there exists a unique $K$-homomorphism $\psi'\colon A[[T]]^\sa\rightarrow B$ such that $\psi'(T)=f$ and such that the diagram
\[
\begin{diagram}
A&\rTo^\phi&A[[T]]^\sa\\
&\rdTo<\psi&\dTo>{\psi'}\\
&&B
\end{diagram}
\]
commutes. If $\ul{A}$ is an $R$-lattice of $A$, we may set $A[[T]]^\sa\mathrel{\mathop:}=\ul{A}[[T]]\otimes_RK$.
\end{prop}

By iteration, we obtain semi-affinoid $K$-algebras $A[[\ul{S}]]\langle \ul{T}\rangle^\sa$ with the expected universal properties, for pairs of finite systems of variables $\ul{S}$, $\ul{T}$.

\subsubsection{Special subdomains}

Let us define and characterize the semi-affinoid $K$-algebras that will be associated to certain special semi-affinoid subdomains.

\begin{prop}\label{ratsubdoms}
Let $A$ be a semi-affinoid $K$-algebra, and let $f_1,\ldots,f_m$, $g_1,\ldots,g_n$, $h\in A$ be semi-affinoid functions without a common zero. There exists a homomorphism of semi-affinoid $K$-algebras 
\[
\phi\colon A\rightarrow A\left\llbracket \frac{f_1}{h},\ldots,\frac{f_m}{h}\right\rrbracket\left\langle \frac{g_1}{h},\ldots,\frac{g_n}{h}\right\rangle^\sa
\]
such that the following holds:
\begin{packed_enum}
\item $\phi(h)$ is invertible,
\item for all $1\leq i\leq m$,  $\phi(f_i)/\phi(h)$ is topologically quasi-nilpotent,
\item for all $1\leq j\leq n$,  $\phi(g_j)/\phi(h)$ is power-bounded and
\item for any semi-affinoid $K$-algebra $B$ and any $K$-homomorphism $\psi\colon A\rightarrow B$ where $\psi(h)$ is invertible, where the $\psi(f_i)/\psi(h)$ are topologically quasi-nilpotent and where the $\psi(g_j)/\psi(h)$ are power-bounded, there exists a unique $K$-homomorphism 
\[
\psi'\colon A\left\llbracket \frac{f_1}{h},\ldots,\frac{f_m}{h}\right\rrbracket\left\langle \frac{g_1}{h},\ldots,\frac{g_n}{h}\right\rangle^\sa\rightarrow B
\]
such that $\psi=\psi'\circ\phi$.
\end{packed_enum}
 More precisely speaking, this universal property is satisfied by 
 \[
 A\left\llbracket S_1,\ldots,S_m\right\rrbracket\langle T_1,\ldots,T_n\rangle^\sa/(hS_i-f_i,hT_j-g_j)\quad.
 \]
If $\ul{A}$ is an $R$-lattice of ff type for $A$ and if $n$ is a natural number such that $\pi^nh$, $\pi^nf_i$ and $\pi^ng_j$ are contained in $\ul{A}$ for all $i$ and $j$, then 
\[
\ul{A}\left\llbracket S_1,\ldots,S_m\right\rrbracket\langle T_1,\ldots,T_n\rangle/(\pi^nhS_i-\pi^nf_i,\pi^nhT_j-\pi^ng_j)/(\textup{$\pi$-torsion})\quad.
\]
is an $R$-lattice of ff type for the semi-affinoid $K$-algebra that we are characterizing.
\end{prop}
\begin{proof}
Let 
\[
\phi\colon A\rightarrow A\left\llbracket S_1,\ldots,S_m\right\rrbracket\langle T_1,\ldots,T_n\rangle^\sa/(hS_i-f_i,hT_j-g_j)
\]
be the natural $K$-homomorphism. Since the system $h,f_i,g_j$ has no common zero on $\Max A$, the relations $gS_i=f_i$, $hT_j=g_j$ show that $h$ has no zero on the maximal spectrum of the target of $\phi$; by the Hilbert Nullstellensatz \ref{hilbertprop}, $\phi(h)$ is invertible. It now follows from Proposition \ref{relcloseddiscprop} and Proposition \ref{relopendiscprop} that $\phi$ has the asserted universal property. Since $\ul{A}\left\llbracket \ul{S}\right\rrbracket\langle \ul{T}\rangle$ is an $R$-lattice of ff type for $A\left\llbracket \ul{S}\right\rrbracket\langle \ul{T}\rangle^\sa$, we see that 
\begin{eqnarray*}
&&\ul{A}\left\llbracket \ul{S}\right\rrbracket\langle \ul{T}\rangle/(\pi^nhS_i-\pi^nf_i,\pi^nhT_j-\pi^ng_j)\otimes_RK\\
&=&A\left\llbracket \ul{S}\right\rrbracket\langle \ul{T}\rangle^\sa/(hS_i-f_i,hT_j-g_j)\quad,
\end{eqnarray*}
which proves the last statement of the proposition.
\end{proof}

Of course, the lattice that we describe in Proposition \ref{ratsubdoms} does not depend on the choice of $n$, since we divide by $\pi$-torsion.

			
		\subsection{$K$-Banach norms on semi-affinoid $K$-algebras}\label{semaffbanachsec}

As we will discuss in this section, a semi-affinoid $K$-algebra $A$ carries a unique $K$-Banach algebra structure. The resulting topology induces the $\pi$-adic topology on any $R$-lattice $\ul{A}$ of ff type for $A$. The Berkovich spectrum of $A$ may be viewed as a compactification of the rigid space obtained from $A$ via Berthelot's construction; its construction does not commute with complete localization of $\ul{A}$ and, hence, it does not globalize. The Berkovich compactification point of view that we just mentioned is related to $\pi$-envelopes and their adic generic fibers which we will discuss later. Again, this section is of a rather complementary and illustrative nature, and its reading is not mandatory for the rest of this thesis.

\subsubsection{Norms associated to $R$-lattices of ff type}

Let $A$ be a semi-affinoid $K$-algebra, and let $\ul{A}$ be an $R$-lattice of ff type for $A$. If $f\in A$ is a nonzero semi-affinoid function, we set
\[
\nu_{\ul{A}}(f)\mathrel{\mathop:}=\max\{n\in\Z\,;\,\pi^{-n} f\in\ul{A}\}\quad.
\]
Moreover, we set $|f|_{\ul{A}}\mathrel{\mathop:}=|\pi|^{\nu_{\ul{A}}(f)}\in \mathbb{R}$ and $|0|_{\ul{A}}\mathrel{\mathop:}=0\in\mathbb{R}$. For any nonzero $f\in\N$, the integer $\nu_{\ul{A}}(f)$ is well-defined. Indeed, for $n\in\Z$ small enough,  $\pi^{-n}f$ must be contained in $\ul{A}$ because $\ul{A}$ is an $R$-lattice in $A$. In particular, the set of integers $n$ with $\pi^{-n}f\in\ul{A}$ is nonempty. If there was no largest such $n$, there would exist a nonzero element in $\ul{A}$ that is infinitely $\pi$-divisible, but $\ul{A}$ is $\pi$-adically separated. 

\begin{lem}The map $|\cdot|_{\ul{A}}\colon A\rightarrow\mathbb{R}$ is a $K$-Banach algebra norm on $A$, and the induced topology on $\ul{A}$ is the $\pi$-adic topology.
\end{lem}
\begin{proof}
The map $|\cdot|_{\ul{A}}$ is readily seen to enjoy the following properties:
\begin{packed_enum}
\item $|1|_{\ul{A}}=1$, and for $f\in A$, $|f|_{\ul{A}}=0$ if and only if $f=0$
\item $|\alpha f|_{\ul{A}}=|\alpha|\cdot|f|_{\ul{A}}$ for any $\alpha\in K$, $f\in A$; furthermore $|A|_{\ul{A}}=|K|$
\item $\ul{A}=\{f\in A\,;\,|f|_{\ul{A}}\leq 1\}$
\item $|f+g|_{\ul{A}}\leq |f|_{\ul{A}}+|g|_{\ul{A}}$ for all $f,g\in A$
\item $|f\cdot g|_{\ul{A}}\leq|f|_{\ul{A}}\cdot|f|_{\ul{A}}$ for all $f,g\in A$
\end{packed_enum}
Indeed, properties ($i$)-($iii$) are immediate from the definition of $|\cdot|_{\ul{A}}$, and properties ($iv$) and ($v$) follow via rescaling from ($ii$), ($iii$) and the fact that $\ul{A}$ is closed under addition and multiplication. We conclude that $|\cdot|_{\ul{A}}$ is a non-archimedean $K$-algebra norm on $A$. 
A fundamental basis of open neighborhoods of zero for the induced topology is given by the additive subgroups $\pi^{n}\ul{A}\subseteq A$, for $n\in\Z$. Hence, $A$ is complete and separated with respect to $|\cdot|_{\ul{A}}$, the $R$-lattice $\ul{A}$ being $\pi$-adically complete and separated.
\end{proof}

We say that $|\cdot|_{\ul{A}}$ is the lattice norm associated to $\ul{A}$. Of course, an $R$-lattice of ff type $\ul{A}$ is uniquely determined by its associated lattice norm since an element $f\in A$ lies in $\ul{A}$ if and only if $|f|_{\ul{A}}\leq 1$. Thus, the norms associated with $R$-lattices of ff type provide another way of thinking about $R$-lattices of ff type themselves.

We easily verify that $\ul{A}/\pi\ul{A}$ is a domain if and only if  $|\cdot|_{\ul{A}}$ is multiplicative. Moreover, we see that if $\phi\colon A\rightarrow B$ is a homomorphism of semi-affinoid $K$-algebras and if $\ul{A}$ and $\ul{B}$ are $R$-lattices of ff type for $A$ and $B$ respectively such that $\phi(\ul{A})\subseteq\ul{B}$, then $\phi$ is \emph{contractive} with respect to $|\cdot|_{\ul{A}}$ and $|\cdot|_{\ul{B}}$. 

The following Lemma allows us to reduce some questions on lattice norms to the case of a \emph{local} $R$-lattice of ff type. 

\begin{lem}\label{isometryballslem}
Let $A$ be a semi-affinoid $K$-algebra, let $\ul{A}$ be an $R$-lattice of ff type for $A$, and let $\p_1,\ldots,\p_n$ be the associated prime ideals of $\ul{A}/\pi\ul{A}$. For each $1\leq i\leq n$, let $\m_i\subseteq\ul{A}$ be a maximal ideal lying above $\p_i$, let $\ul{B}_i$ be the $\m_i$-adic completion of $\ul{A}$, and let $B_i$ denote the semi-affinoid generic fiber of $\ul{B}_i$. Then the natural $K$-homomorphism
\[
\phi\colon A\rightarrow\bigoplus_{1\leq i\leq n}B_i
\]
is an isometry for $|\cdot|_{\ul{A}}$ and for the lattice norm associated to $\bigoplus_{1\leq i\leq n}\ul{B}_i$.
\end{lem}
\begin{proof}
Let $f\in A$ be a semi-affinoid function such that $|f|_{\ul{A}}=1$; we must show that $|\phi(f)|_{\ul{B}}=1$. Clearly $\phi$ is contractive for the considered norms, so let us assume that $\phi(f)$ has norm $<1$. Then the $i$-th component $\phi(f)_i$ of $\phi(f)$ is $\pi$-divisible in $\ul{B}_i$ for each $i$. It follows that the class $\bar{f}$ of $f$ modulo $\pi$ vanishes in the maximal-adic completion of $(\ul{A}/\pi\ul{A})_{\m_i}$ for every $i$ and, hence, in $(\ul{A}/\pi\ul{A})_{\m_i}$ itself. In particular, $\bar{f}$ vanishes in the local rings $(\ul{A}/\pi\ul{A})_{\p_i}$, and it follows that $\bar{f}$ vanishes in an open subset of $\Spec\ul{A}/\pi\ul{A}$ containing all associated points, which implies that $\bar{f}=0$. Hence $|f|_{\ul{A}}<1$, which contradicts our assumption on $f$.
\end{proof}


\subsubsection{The Gauss norm and residue norms}

Let $\ul{S}$ and $\ul{T}$ be finite systems of variables. The semi-affinoid $K$-algebra 
\[
R[[\ul{S}]]\langle\ul{T}\rangle\otimes_RK
\]
admits the distinguished $R$-lattice of ff type $R[[\ul{S}]]\langle\ul{T}\rangle$. 

\begin{defi}\label{gaussnormdef}
Let $f=\sum_{i,j}a_{ij}\ul{S}^i\ul{T}^j$ be an element of $R[[\ul{S}]]\langle\ul{T}\rangle\otimes_RK$, where $i$ and $j$ are multi-indices.  We set
\[
|f|_\Gauss\,\mathrel{\mathop:}=\,\max\{|a_{ij}|\,;\,i,j\}\quad,
\]
and we call $|\cdot|_\Gauss$ the \emph{Gauss norm}\index{Gauss norm} on $R[[\ul{S}]]\langle\ul{T}\rangle\otimes_RK$.
\end{defi}

The maximum exists since the coefficients of $f$ are bounded and since the valuation on $R$ is discrete. It is immediate from the definitions that the Gauss norm coincides with the distinguished lattice norm $|\cdot|_{R[[\ul{S}]]\langle\ul{T}\rangle}$. In particular, $|\cdot|_\Gauss$ is multiplicative.

\begin{lem}Let $A$ be a semi-affinoid $K$-algebra, let $\ul{A}$ be an $R$-lattice of ff type for $A$, and let 
\[
\ul{\phi}\colon R[[\ul{S}]]\langle\ul{T}\rangle\rightarrow\ul{A}
\]
be a formal presentation of $\ul{A}$, where $\ul{S}$ and $\ul{T}$ are suitable finite systems of variables. Then $|\cdot|_{\ul{A}}$ is the residue norm of $|\cdot|_\Gauss$ under $\phi\mathrel{\mathop:}=\ul{\phi}\otimes_RK$.
\end{lem}
\begin{proof}
Let $|\cdot|_\phi$ denote the residue norm attached to $\phi$; then $|f|_\phi$ is the infimum of the Gauss norms of the elements in the $\phi$-preimage of $f$. Since $|\cdot|_\Gauss$ takes values in $|K|$ and since the valuation on $K$ is discrete, $|\cdot|_\phi$ takes values in $|K|$ as well. Hence, it suffices to show that for any $f\in A$, we have $|f|_\phi\leq 1$ if and only if $|f|_{\ul{A}}\leq 1$. But $|f|_\phi\leq 1$ if and only if $f$ has a $\phi$-preimage in $R[[\ul{S}]]\langle\ul{T}\rangle$, which is the case if and only if $f\in\ul{A}$, that is, if and only if $|f|_{\ul{A}}\leq 1$.
\end{proof}

\subsubsection{Equivalence of $K$-Banach algebra norms}

We show that all $K$-Banach algebra norms on a semi-affinoid $K$-algebra are equivalent. Let $A$ be a semi-affinoid $K$-algebra, and let $\ul{A}_1$ and $\ul{A}_2$ be $R$-lattices of ff type for $A$. It is obvious that the norms $|\cdot|_{\ul{A}}$ and $|\cdot|_{\ul{B}}$ are equivalent. Indeed, we may assume that $\ul{A}_1\subseteq\ul{A}_2$, and we have seen that there exists an $n\in\N$ such that $\pi^n\ul{A}_2\subseteq\ul{A}_1$. 

\begin{prop}\label{banachequivprop}
Let $A$ be a semi-affinoid $K$-algebra, let $B$ be any Banach $K$-algebra, and let $\phi\colon B\rightarrow A$ by any $K$-homomorphism. Then $\phi$ is continuous with respect to any $K$-Banach algebra topology on $A$. In particular, all $K$-Banach algebra norms on $A$ are equivalent.
\end{prop}
\begin{proof}
According to [BGR] Proposition 3.7.5/2, it suffices to see that $A$ is noetherian, that $A/\m^{n+1}$ is $K$-finite for any maximal ideal $\m\subseteq A$ and any $n\in\N$, and that
\[
\bigcap_{\m\subseteq A\textup{ maximal}}\quad\bigcap_{n\in\N}\m^{n+1}\,=\,0\quad.
\]
We have seen in Lemma \ref{semiaffquotientlem} that semi-affinoid $K$-algebras are noetherian. Let $\m\subseteq A$ be a maximal ideal. By Lemma \ref{finiteresiduefieldlem}, $A/\m$ is $K$-finite. Let us argue by induction on $n\in\N$ and assume that $A/\m^{n+1}$ is $K$-finite for a given $n$. Since $A$ is noetherian, the $A/\m$-module $\m^{n+1}/\m^{n+2}$ is finitely generated, and the short exact sequence
\[
0\rightarrow \m^{n+1}/\m^{n+2}\rightarrow A/\m^{n+2}\rightarrow A/\m^{n+1}\rightarrow 0
\]
shows that $A/\m^{n+2}$ is $K$-finite, as desired. Let finally $f\in A$ be an element that is contained in the intersection of all $\m^{n+1}$. By Krull's Intersection Theorem, for every maximal ideal $\m\subseteq A$ there exists an $m\in\m$ such that $(1-m)f=0$. It follows that the annihilator ideal of $f$ cannot be contained in any maximal ideal of $A$, and we conclude that $f=0$.
\end{proof}

Let us make one important remark: We have seen that a semi-affinoid $K$-algebra $A$ carries a unique $K$-Banach topology and that this topology induces the $\pi$-adic topology on any $R$-lattice of ff type $\ul{A}$ for $A$. Let us fix such an $R$-lattice of ff type $\ul{A}$ for $A$. There exists a unique \emph{group} topology on $(A,+)$ such that $\ul{A}\subseteq A$ is open and such that the induced topology on $\ul{A}$ is the Jacobson-adic topology: Let $\a\subseteq\ul{A}$ be an ideal of definition; then a subset $U\subseteq A$ is open in this group topology if for every $x\in U$ there exists an $n\in\N$ such that $x+\a^n\subseteq U$. However, this topology on $A$ needs not be a ring topology:

\begin{example}\label{hubernosalvationex}
Multiplication by $\pi^{-1}$ in $R[[S]]\otimes_RK$ is not continuous for the group topology induced by the $(\pi,S)$-adic topology on $R[[S]]$.
\end{example}
\begin{proof}
Let $m$ denote the $R[[S]]\otimes_RK$-module endomorphism of $R[[S]]\otimes_RK$ that is given by $\pi^{-1}$-multiplication. If $m$ was continuous, then $m^{-1}(R[[S]])$ would contain $(\pi,S)^n$ for some $n\in\N$; however, there exists no $n\in\N$ such that $\pi^{-1}S^n$ is contained in $R[[S]]$.
\end{proof}

Thus, it is in general not possible to consider the pair $(A,\ul{A})$ as an f-adic ring in the sense of \cite{Huberbuch} p. 37 unless we replace the topology on $\ul{A}$ by the $\pi$-adic topology. Here lies the reason for the fact that the theory of uniformly rigid spaces is \emph{not} essentially contained in Huber's theory of adic spaces.

\subsubsection{The supremum semi-norm}

Let $A$ be a semi-affinoid $K$-algebra, provided with its Banach algebra topology. The map
\[
|\cdot|_\sup\colon A\rightarrow\mathbb{R}\quad,\quad f\mapsto\sup_{x\in\Max A}|f(x)|\quad,
\]
is called the supremum semi-norm on $A$. 

If $\p_1,\ldots,\p_n$ are the minimal prime ideals of $A$ and if for $f\in A$ we let $f_i$ denote the image of $f$ in $A/\p_i$, then $|f|_\sup\,=\,\max_{1\leq i\leq n}|f_i|_\sup$, the set $\Max A$ being the union of the $\Max A_i$.

\begin{lem}\label{supnormimmediatelem}
The map $|\cdot|_\sup$ is well-defined, and it is a power-multiplicative $K$-algebra semi-norm. Moreover, $|\cdot|_\sup$ is a norm if and only if $A$ is reduced.
\end{lem}
\begin{proof}
Let $\ul{A}$ be an $R$-lattice of ff type for $A$. Any element $f\in A$ can be written as a product $f=\pi^{-n}g$ with $n\in\N$ and $g\in\ul{A}$. Since $\ul{A}\subseteq\mathring{A}$, we obtain $|f(x)|\leq |\pi|^{-n}$ for every $x\in\Max A$, which implies that $|f|_\sup$ is a well-defined real number. The last statement is obvious from the fact that semi-affinoid $K$-algebras are Jacobson rings, cf.\ Proposition \ref{hilbertprop}.
\end{proof}

In the following Lemma, we gather some standard facts on $|\cdot|_\sup$ which are proved in the usual way. Power-boundedness and topological nilpotency are understood with respect to the $K$-Banach topology.

\begin{lem}\label{supseminormlem}
Let $|\cdot|$ be any $K$-Banach algebra norm on $A$, and let $f\in A$ be a semi-affinoid function.
\begin{packed_enum}
\item The inequality $|f|_\sup\leq|f|$ holds. 
\item The function $f$ is power-bounded if and only if $|f|_\sup\leq 1$.
\item The function $f$ is topologically nilpotent if and only if $|f|_\sup< 1$.
\item We have identities
\[
|f|_\sup\,=\,\lim_{n\in\N}\sqrt[n]{|f|^n}\,=\,\inf_{n\in\N}\sqrt[n]{|f|^n}\quad.
\]
\end{packed_enum}
\end{lem}
\begin{proof}
Statement ($i$) holds for general $K$-Banach algebras by [BGR] 3.8.2/2. 

Let us prove statement ($ii$). If $f$ is power-bounded, there exists a real number $r$ such that $|f^n|\leq r$ for all $n\in\N$. Since $|\cdot|_\sup$ is power-multiplicative and since $|f^n|_\sup\leq |f^n|$ for all $n$, we obtain $|f|_\sup\leq\sqrt[n]{r}$, which converges to $1$ when $n$ tends to infinity. Conversely, let us assume that $|f|_\sup\leq 1$; then $f\in\mathring{A}$, and by Corollary \ref{latticeexcor} there exists an $R$-lattice $\ul{A}$ of ff type for $A$ containing $f$. Then $f^n\in\ul{A}$ for all $n\in\N$, and it follows that $|f^n|_{\ul{A}}\leq 1$ for all $n\in\N$; hence $f$ is power-bounded.

We prove statement ($iii$). If $f$ is topologically nilpotent, the sequence $(|f^n|)_{n\in\N}$ is a zero sequence, and in particular we may choose $n$ big enough such that $|f^n|<1$; then $|f|_\sup\leq\sqrt[n]{|f^n|}<1$ because $|f|^n_\sup=|f^n|_\sup\leq |f^n|$. Conversely, let us assume that $|f|_\sup<1$. Since $|\cdot|_\sup$ is power-multiplicative, there exists an integer $n$ such that $|f^n|_\sup\leq|\pi|$. Then $\pi^{-1}f^n\in\mathring{A}$, and by Corollary \ref{latticeexcor} there exists an $R$-lattice $\ul{A}$ of ff type for $A$ containing $\pi^{-1}f^n$. Now $f^n$ is $\pi$-divisible in $\ul{A}$ and, hence, topologically nilpotent for the $\pi$-adic topology on $\ul{A}$. It follows that $f$ is topologically nilpotent for the Banach topology on $A$.

The second equality in statement ($iv$) is established for general semi-normed rings in [BGR] at the beginning of section 1.3.2. Let us prove the first equality. The inequality '$\leq$' follows from part ($i$) of this lemma. Let us assume that the strict inequality '$<$' holds. Let us write 
\[
|\cdot|'\mathrel{\mathop:}=\inf_{n\in\N}\sqrt[n]{|\cdot|}\quad;
\]
then by [BGR] 1.3.2/1, $|\cdot|'$ is a power-multiplicative semi-norm on $A$. Since $\sqrt{|K^*|}$ is dense in the positive real numbers, there exists an $r\in\sqrt{|K^*|}$ such that 
\[
|f|_\sup<r<|f|'\quad.
\]
Let $m$ be a natural number such that $r^m\in |K^*|$, and let $c$ be an element in $K$ such that $|c|=r$; then 
\[
|f^m/c|_\sup<1<|f^m/c|'\quad.
\]
Hence, $|(f^m/c)^n|'$ tends to infinity when $n$ goes to infinity. But $|\cdot|'\leq |\cdot|$ since $|\cdot|$ is sub-multiplicative, and it follows that $f^m/c$ is not power-bounded. However, this contradicts the inequality $|f^m/c|_\sup<1$. 
\end{proof}

Lemma \ref{supseminormlem} ($ii$) shows that we are using the adjective 'power-bounded' in an unambiguous way: A semi-affinoid function $f\in A$ is power-bounded for the unique $K$-Banach topology on $A$ if and only if $f\in\mathring{A}$, that is, if and only if $f$ is power-bounded in the sense of Definition \ref{pbtqndefi}.
Moreover, Lemma \ref{supseminormlem} ($iv$) shows that if $A$ admits a power-multiplicative $K$-Banach norm $|\cdot|$, then $|\cdot|=|\cdot|_\sup$. In particular, the Gauss norm on a free semi-affinoid $K$-algebra coincides with the supremum semi-norm.

Let us remark that the maximum principle does \emph{not} hold for general semi-affinoid $K$-algebras. For examle, we have already observed that if $S$ is a variable, then the power-bounded semi-affinoid function $S\in R[[S]]\otimes_RK$ is topologically quasi nilpotent without being topologically nilpotent.

The following result is also shown in \cite{LiRo} 5.1.11:

\begin{prop}\label{supnormvaluegroupprop}
The supremum semi-norm on a semi-affinoid $K$-algebra $A$ takes values in $\sqrt{|K|}$.
\end{prop}
\begin{proof}
If there exists a semi-affinoid $K$-algebra $B$ together with a $K$-ho\-mo\-mor\-phism $\phi\colon A\rightarrow B$ such that $|B|_\sup\subseteq\sqrt{|K|}$ and such that $\phi$ is an isometry for some choice of $K$-Banach norms on $A$ and $B$, then $|A|_\sup\subseteq\sqrt{|K|}$ by Lemma \ref{supseminormlem} ($iv$). Hence, by Lemma \ref{isometryballslem} we may assume that $A$ admits a local lattice $\ul{A}$. By the remark following Lemma \ref{supnormimmediatelem}, we may moreover assume that $A$ is integral. By Corollary \ref{localnoethnormcor}, there exists a finite monomorphism $R[[\ul{S}]]\otimes_RK\hookrightarrow A$ for some finite system of variables $\ul{S}$. Since the supremum semi-norm on $R[[\ul{S}]]\otimes_RK$ coincides with the Gauss norm which takes values in $|K|$, we conclude from [BGR] 3.8.1/7 that $|A|_\sup\subseteq\sqrt{|K|}$, as desired.
\end{proof}

\begin{prop}If $A$ is reduced, then $A$ is a Banach function algebra; that is, $|\cdot|_\sup$ is a $K$-Banach norm on $A$.
\end{prop}
\begin{proof}
By Lemma \ref{supnormimmediatelem}, $|\cdot|_\sup$ is a $K$-algebra norm on $A$. By Proposition \ref{powerboundedintegralprop}, $\mathring{A}$ is an $R$-lattice of ff type for $A$; we must verify that  $|\cdot|_\sup$ is equivalent to $|\cdot|_{\mathring{A}}$. By Lemma  \ref{supseminormlem}, we know that $|\cdot|_\sup\leq|\cdot|_{\mathring{A}}$, so it suffices to show that 
\[
|f|_{\mathring{A}}\leq |\pi|^{-1}\cdot|f|_\sup
\]
for all $f\in A$. We may assume that $|f|_{\mathring{A}}=1$, which means that $f\in\mathring{A}$ and $\pi^{-1}f\notin\mathring{A}$. The latter condition is equivalent to $|\pi^{-1}f|_\sup=|\pi|^{-1}\cdot|f|_\sup>1$, so we obtain the desired inequality.
\end{proof}

\subsubsection{The Berkovich spectrum}

Let $A$ be a semi-affinoid $K$-algebra, equipped with its unique $K$-Banach algebra topology, and let $\M(A)$ denote the Berkovich spectrum of $A$, that is, the set of all real-valued bounded multiplicative semi-norms on $A$, provided with the weakest topology such that all evaluation maps are continuous, cf.\ \cite{Berk1} 1.2. We recall that any bounded multiplicative semi-norm on a $K$-Banach algebra is automatically a $K$-algebra semi-norm. By \cite{Berk1} Theorem 1.2.1, for $A\neq 0$ the space $\M(A)$ is a nonempty compact Hausdorff space. By Lemma \ref{supseminormlem} and \cite{Berk1} Theorem 1.3.1, we know that
\[
|f|_\sup\,=\,\max_{x\in\M(A)}|f(x)|\quad,\quad(*)
\]
so while the maximum of a semi-affinoid function might not be attained on any point $x\in\Max A$, it \emph{is} attained in $\M(A)$. We may thus consider $\M(A)$ as a \emph{compactification} of $\Max A$ and, hence, of the rigid space that is associated to $A$ via Berthelot's construction, cf. Section \ref{assocrigspacesubsec}. For example, 
\[
|\cdot|_\Gauss\in\M(R[[\ul{S}]]\langle\ul{T}\rangle\otimes_RK)
\]
realizes the maximum for any semi-affinoid function in $R[[\ul{S}]]\langle\ul{T}\rangle\otimes_RK$.

An analytic field extension of $K$ is a field $\sK$ over $K$, together with an absolute value extending the absolute value on $K$ such that $\sK$ is complete with respect to this absolute value. If $\sK$ is an analytic extension of $K$, we may consider its valuation subring $\mathring{\sK}\subseteq\sK$, the valuation ideal $\check{\sK}\subseteq\mathring{\sK}$ and the residue field $\tilde{\sK}=\mathring{\sK}/\check{\sK}$. The elements of $\M(A)$ may be considered as equivalence classes of $K$-algebra homomorphisms from $A$ to analytic extensions of $K$, cf.\ \cite{Berk1} Remark 1.2.2 ($ii$).

If $\ul{A}$ is an $R$-lattice of ff type for $A$ and if $x$ in an element in $\M(A)$, then a representative $\chi_x\colon A\rightarrow\sK$ of $x$ restricts to an $R$-homomorphism 
\[
\chi_{x,\ul{A}}\colon\ul{A}\rightarrow\mathring{\sK}
\]
by $(*)$, since $|f|_\sup\leq 1$ for all $f\in\ul{A}$. It thus induces a $k$-homomorphism 
\[
\tilde{\chi}_{x,\ul{A}}\colon \ul{A}/\pi \ul{A}\rightarrow\tilde{\sK}\quad.
\]
The kernel of $\tilde{\chi}_{x,\ul{A}}$ is a prime ideal which is independent of the choice of the representative $\chi_x$ of $x$; hence we obtain a \emph{specialization map}
\[
\sp_{\ul{A}}\colon\M(A)\rightarrow\Spec \ul{A}/\pi\ul{A}
\]
extending the specialization map that we have studied in Section \ref{specmapsec}.
\begin{lem}
For any choice of $\ul{A}$, the extended specialization map $\sp_{\ul{A}}$ is surjective.
\end{lem}
\begin{proof}
Let $\p$ be any prime ideal in $\ul{A}$ containing $\pi$, and let $\q\subseteq$ be a minimal prime ideal of $\ul{A}$ that is contained in $\p$. By \cite{EGAII} 7.1.7, there exists a discrete valuation subring $R'$ of $Q(\ul{A}/\q)$ dominating $(\ul{A}/\q)_\p$; its completion yields a point in $\M(A)$ specializing to $\p$. Hence, $\sp_{\ul{A}}$ is surjective.
\end{proof}

In particular, $\sp_{\ul{A}}$ needs not factor through the set of \emph{open} prime ideals in $\ul{A}/\pi\ul{A}$.

Let us recall that a closed subset $\Gamma\subseteq\M(A)$ is called a \emph{boundary} of $\M(A)$ if every function $f\in A$ attains its maximum on $\Gamma$. If there exists a unique minimal boundary, it is called the \emph{Shilov boundary} of $\M(A)$, and it is denoted by $\Gamma(A)$. We briefly discuss the Shilov boundary of a semi-affinoid $K$-algebra $A$ admitting an $R$-lattice of ff type $\ul{A}$ whose reduction modulo $\pi$ is \emph{integral}.

\begin{prop}\label{integralredshilovprop}
Let $A$ be a semi-affinoid $K$-algebra with an $R$-lattice of ff type $\ul{A}$ having the property that $\ul{A}/\pi\ul{A}$ is a domain. Then $|\cdot|_\sup=|\cdot|_{\ul{A}}$, and
\[
\Gamma(A)\,=\,\{|\cdot|_\sup\}\,=\,\sp_{\ul{A}}^{-1}(\{\eta\})\quad,
\]
where $\eta\in\Spec\ul{A}/\pi\ul{A}$ is the generic point.
\end{prop}
\begin{proof}
Since $\ul{A}/\pi\ul{A}$ is integral, the lattice norm $|\cdot|_{\ul{A}}$ is multiplicative; hence $|\cdot|_{\ul{A}}\in\M(A)$. By Lemma \ref{supseminormlem}, multiplicativity of $|\cdot|_{\ul{A}}$ implies $|\cdot|_\sup=|\cdot|_{\ul{A}}$. Since $\M(A)$ is Hausdorff, the one-point set $\{|\cdot|_\sup\}$ is closed and, hence, a minimal boundary for $\M(A)$. Let $\Gamma$ be another boundary of $\M(A)$; we claim that $|\cdot|_\sup\in\Gamma$. Indeed, let us assume that $|\cdot|_\sup\notin\Gamma$. Since $\Gamma$ is closed, there exists an $f\in A$ and a positive real number $\epsilon$ such that the set
\[
\{x\in\M(A)\,;\,|f|_\sup-|f|_x<\varepsilon\}
\]
is disjoint to $\Gamma$. However, for every $f\in A$ there exists an $x\in\Gamma$ such that $|f|_x=|f|_\sup$ since $\Gamma$ is a boundary. We thus arrive at a contradiction, and it follows that $\{|\cdot|_\sup\}$ is the \emph{unique} minimal boundary of $\M(A)$; that is, $\Gamma(A)=\{|\cdot|_\sup\}$. Since $|\cdot|_\sup=|\cdot|_{\ul{A}}$, the set of elements $f\in \ul{A}$ satisfying $|f|_\sup<1$ is precisely the prime ideal $\pi\ul{A}$ in $\ul{A}$; hence $|\cdot|_\sup$ specializes to $\pi\ul{A}$. The local ring $\ul{A}_{\pi\ul{A}}$ is a one-dimensional principal ideal domain, hence a discrete valuation ring; thus the $\sp_{\ul{A}}$-fiber of $\eta$ contains only one point.
\end{proof}

We finally illustrate that the Berkovich spectrum $\M(A)$ of a semi-affinoid $K$-algebra $A$ does not behave well with respect to complete localization of $R$-lattices of ff type: If $\ul{A}\rightarrow\ul{A}_{\{f\}}$ is a complete localization with respect to the Jacobson-adic topology on $\ul{A}$, then the induced homomorphism of semi-affinoid $K$-algebras needs not be a complete localization of Banach algebras, and the induced map of valuation spectra needs not be injective. It is for this reason that the theory of uniformly rigid $K$-spaces is not contained in Berkovich's theory of $K$-analytic spaces.

\begin{example}
Let us consider the flat $R$-algebra of ff type
\[
\ul{A}\,\mathrel{\mathop:}=\,R\langle X,Y\rangle[[Z]]/(Z-XY)\quad.
\]
The natural map 
\[
\M(\ul{A}_{\{X-Y\}}\otimes_RK)\rightarrow\M(\ul{A}\otimes_RK)
\]
is \emph{not} injective.
\end{example}
\begin{proof}
By the Jacobian criterion, we see that $\ul{A}$ is $R$-smooth and, hence, $R$-flat, cf.\ Section \ref{smoothnesssec}. It follows that the $\pi$-reductions of $\ul{A}$ and $\ul{A}_{\{X-Y\}}$ are $k$-smooth and, hence, locally integral. The reduction of $\ul{A}$ modulo its biggest ideal of definition is $k[X,Y]/(XY)$, so we see that while $\ul{A}$ is integral, $\ul{A}_{\{X-Y\}}$ decomposes as a direct sum of two nontrivial domains. Since complete localization is flat and since flat homomorphisms of domains are dominant, we see from Proposition \ref{integralredshilovprop} that $\M(\ul{A}_{\{X-Y\}}\otimes_RK)$ has two 'Gauss points' that both map to the one 'Gauss point' of $\M(\ul{A}\otimes_RK)$.
\end{proof}

We thus conclude our discussion of semi-affinoid $K$-algebras.

	\section{Semi-affinoid and uniformly rigid spaces}\label{semaffspacessec}
		We first give an ad hoc definition of the category of semi-affinoid $K$-spaces and later explain how to view it as a full subcategory of the category of locally G-ringed $K$-spaces. The category of \emph{affinoid} $K$-spaces is defined as the dual of the category of \emph{affinoid} $K$-algebras, cf.\ \cite{BGR} 7.1.4/1. By analogy, we define:

\begin{defi}The category $\sAff_K$ of \emph{semi-affinoid $K$-spaces} \index{category!of semi-affinoid $K$-spaces} is the dual of the category of semi-affinoid $K$-algebras. 
\end{defi}\index{semi-affinoid space}

If $A$ is a semi-affinoid $K$-algebra, the corresponding semi-affinoid $K$-space is denoted by $\sSp A$. If $\phi$ is a morphism of semi-affinoid $K$-spaces, we write $\phi^*$ for the corresponding morphism of semi-affinoid $K$-algebras. If $X=\sSp A$ is a semi-affinoid $K$-space, we write $\Gamma(X,\O_X)\mathrel{\mathop:}= A$. A semi-affinoid $K$-space $X$ is called \emph{reduced}\index{semi-affinoid space!reduced} if its ring of global semi-affinoid functions $\Gamma(X,\O_X)$ is reduced.



Let us note that the category of semi-affinoid $K$-spaces has fibered products.\index{fibered product!of semi-affinoid spaces} Indeed, this is a formal consequence of Proposition \ref{amalgsumsprop}. Moreover, Proposition \ref{amalgsumsprop} shows that cartesian squares of flat affine formal $R$-schemes in the category of formal $R$-schemes of ff type or in the category of flat formal $R$-schemes of ff type induce cartesian squares of associated semi-affinoid spaces.

By definition, the category of semi-affinoid $K$-spaces contains the category of affinoid $K$-spaces as a full subcategory. Nonetheless, we will need to distinguish between the affinoid $K$-space $\Sp A$ and the semi-affinoid $K$-space $\sSp A$ that is associated with an affinoid $K$-algebra $A$ when we consider these spaces as locally G-ringed spaces: The rigid-analytic G-topology on the maximal spectrum of an affinoid K-algebra will be different from the uniformly rigid G-topology that we are going to define.

Let $A$ be a semi-affinoid $K$-algebra. We call 
\[
|\sSp A|\,\mathrel{\mathop:}=\,\Max A
\]
the set of physical points of $\sSp A$.\index{point!physical} By Corollary \ref{saffmaxpbismalcor}, a morphism
\[
\phi\colon\sSp B\rightarrow\sSp A
\]
of semi-affinoid $K$-spaces induces a map 
\[
|\phi|\colon|\sSp B|\rightarrow |\sSp A|\quad.
\]
We will sometimes drop the vertical lines from the notation if no confusion can result. By Lemma \ref{finiteresiduefieldlem}, for any $f\in A$ and any $x\in\sSp A$, the absolute value $|f(x)|$ is well-defined.



If $X=\sSp A$ is a semi-affinoid $K$-space, a flat affine \emph{$R$-model}\index{model!of a semi-affinoid space} of $X$ is a flat affine formal $R$-scheme $\fX$ of ff type together with an identification of its ring of global sections with an $R$-lattice of ff type for $A$. We say that $X$ is the semi-affinoid generic fiber of $\fX$, and we write $X=\fX^\srig$.\index{generic fiber!semi-affinoid} The same notation is used to denote semi-affinoid generic fibers of morphisms.

			\subsection{The rigid space associated to a semi-affinoid space}\label{assocrigspacesubsec}
				Let $X=\sSp A$ be a semi-affinoid $K$-space, and let $\fX$ be a flat affine $R$-model of ff type for $X$. Berthelot's construction yields a rigid $K$-space $X^\r\mathrel{\mathop:}=\fX^\rig$ together with a $K$-homomorphism 
\[
\phi\colon A\rightarrow\Gamma(X^\r,\O_{X^\r})\quad,
\]
cf.\ \cite{dJ} 7.1.8. By our discussion in Section \ref{specmapsec} and by \cite{dJ} 7.1.9, the homomorphism $\phi$ induces a bijection $|X^\r|\rightarrow|X|$ and local homomorphisms $A_\m\rightarrow\O_{X^\r,x}$ which are isomorphisms on maximal-adic completions, where $x$ is a point of $X^\r$ and where $\m\in\Max A$ is the image of $x$ under the above bijection. The pair $(X^\r,\phi)$ is \emph{universal} in the following sense: 

\begin{prop}\label{semaffspaceassocprop}
Let $Y$ be a rigid $K$-space, and let $\psi\colon A\rightarrow\Gamma(Y,\O_Y)$ be a $K$-algebra homomorphism. There exists a unique morphism of rigid $K$-spaces $\sigma\colon Y\rightarrow X^\r$ such that $\psi=\Gamma(\sigma^\sharp)\circ\phi$. 
\end{prop}
\begin{proof}
We may assume that $Y$ is affinoid, $Y=\Sp B$. Uniqueness of $\sigma$ follows from the above-mentioned fact that $\phi$ induces a bijection of points and isomorphisms of completed stalks. Let $\ul{A}\subseteq A$ be the $R$-lattice of ff type corresponding to $\fX$. By Corollary \ref{freesemaffcor} ($iv$) and Corollary \ref{afflattfintypecor}, $\psi$ restricts to an $R$-homomorphism $\ul{\psi}\colon\ul{A}\rightarrow\ul{B}$, where $\ul{B}$ is a suitable $R$-lattice of tf type for $B$; now $\sigma\mathrel{\mathop:}=(\Spf\ul{\psi})^\rig$ has the required properties.
\end{proof}

We say that $X^\r$ is the \emph{rigid space associated to} $X$ via Berthelot's construction.\index{rigid space!of a semi-affinoid space} 
If $\ul{\tau}\colon\fY\rightarrow\fX$ is a morphism of affine flat formal $R$-schemes of ff type, we easily see that the unique morphism $(\ul{\tau}^\srig)^r$ provided by Proposition \ref{semaffspaceassocprop} is given by $\ul{\tau}^\rig$.

			\subsection{Closed immersions of semi-affinoid spaces}
				\begin{defi}\label{closedimsemaffdefin}\index{closed immersion!of semi-affinoid spaces}
A morphism $\phi\colon Y\rightarrow X$ of semi-affinoid $K$-spaces is called a \emph{closed immersion} if the underlying homomorphism $\phi^*$ of semi-affinoid $K$-algebras is surjective.
\end{defi}

Let $A$ be a semi-affinoid $K$-algebra. By Lemma \ref{semiaffquotientlem}, every ideal $I\subseteq A$ gives rise to a closed immersion $\sSp A/I\hookrightarrow\sSp A$ of semi-affinoid $K$-spaces. It follows that the isomorphism classes of closed immersions of semi-affinoid $K$-spaces with target $\sSp A$ correspond to the ideals in $A$. Clearly a closed immersion $\sSp A/I\rightarrow\sSp A$ is injective, and its image consists of the maximal ideals in $A$ containing $I$.

\begin{lem}\label{closedimbclem}
Let $\phi\colon Y\hookrightarrow X$ be a closed immersion of semi-affinoid $K$-spaces, let $\psi\colon Z\rightarrow X$ be any morphism of semi-affinoid $K$-spaces, and let us consider the cartesian diagram
\[
\begin{diagram}
Y\times_XZ&\rTo^{\phi'}&Z\\
\dTo<{\psi'}&&\dTo>\psi\\
Y&\rTo^\phi&X
\end{diagram}
\]
in the category of semi-affinoid $K$-spaces. Then $\phi'$ is a closed immersion as well, and its image in $Z$ is the $\psi$-preimage of $\phi(Y)$. Moreover, if $\phi^*\colon A\rightarrow A/I$ and $\psi^*\colon A\rightarrow C$ are the homomorphisms of semi-affinoid $K$-algebras corresponding to $\phi$ and $\psi$, then $A/I\hat{\otimes}_AC$ with its coprojections is naturally isomorphic to $C/IC$.
\end{lem}
\begin{proof}
It suffices to verify that the natural diagram
\[
\begin{diagram}
A&\rTo&A/I\\
\dTo<{\psi^*}&&\dTo>{\psi^*\textup{mod}\,I}\\
C&\rTo&C/IC
\end{diagram}
\]
is cocartesian in the category of semi-affinoid $K$-algebras. By Lemma \ref{semiaffquotientlem}, $C/IC$ is a semi-affinoid $K$-algebra; hence the desired statement follows from the fact that the above diagram is already cocartesian in the category of all $K$-algebras.
\end{proof}
				
			\subsection{Semi-affinoid pre-subdomains}				
				In analogy with \cite{BGR} 7.2.2, we define semi-affinoid pre-subdomains of semi-affinoid $K$-spaces in terms of a universal property:

\begin{defi}\label{semiaffsubdomdefi}\index{subdomain!semi-affinoid pre-}
Let $X$ be a semi-affinoid $K$-space. A subset $U\subseteq |X|$ is called a semi-affinoid pre-subdomain of $X$ if there exists a morphism $\phi\colon Y\rightarrow X$ of semi-affinoid $K$-spaces with $\im|\phi|\subseteq U$ such that for any morphism $\psi\colon Z\rightarrow X$ of semi-affinoid $K$-spaces with $\im|\psi|\subseteq U$, there exists a unique morphism $\eta\colon Z\rightarrow Y$ of semi-affinoid $K$-spaces such that the diagram
\[
\begin{diagram}
Z&&\rTo^{\psi}&&X\\
&\rdDashto<{\eta}&&\ruTo>{\phi}&\\
&&Y&&
\end{diagram}
\]
commutes. We say that $\phi$ represents all semi-affinoid morphisms to $X$ with image in $U$. 
\end{defi}

Of course, this universal property determines $\phi$ up to unique isomorphism. For trivial reasons, $|X|$ is a semi-affinoid pre-subdomain of $X$, the associated morphism $\phi$ being given by the identity on $X$. Similarly, the empty subset is a semi-affinoid pre-subdomain, the universal morphism $\phi$ being associated to the zero homomorphism.

\begin{prop}\label{presubdomstalkscompprop}
Let $X$ be a semi-affinoid $K$-space, let $U\subseteq |X|$ be a semi-affinoid pre-subdomain, and let $\phi\colon Y\rightarrow X$ be a morphism of semi-affinoid $K$-spaces representing all semi-affinoid morphisms to $X$ with image in $U$. Then the map $|\phi|$ is injective, and $\im|\phi|=U$. Moreover, if $X=\sSp A$, $Y=\sSp B$ and if $y$ is a point in $Y$, then viewing $B$ as an $A$-algebra via $\phi^*$, the following holds: $\m_{y}=\m_{\phi(y)}B$, and for all $n\in\N$, the induced homomorphism 
\[
A/\m_{\phi(y)}^n\rightarrow B/\m_y^n
\]
is bijective. 
\end{prop}
\begin{proof}
The proof of the corresponding statement for affinoid subdomains in affinoid $K$-spaces given in BGR 7.2.2/1 works in the present situation. Let us recall the argument. For each $n\in\N_{\geq 1}$ and for each $x\in U$, let us consider the canonical commutative diagram
\[
\begin{diagram}
A&\rTo^{\phi^*}&B&\\
\dTo&\ldDashto&\dTo&\\
A/\m_x^n&\rTo& B/\m_x^nB&\quad.
\end{diagram}
\]
We observe that $A/\m_x^n$ is a semi-affinoid $K$-algebra with the property that the semi-affinoid $K$-space $\sSp A/\m_x^n$ has only one physical point. The natural morphism $\sSp A/\m_x^n\rightarrow X$ has image $x\in U$; by the universal property of $\phi$, there exists a unique $K$-homomorphism $B\rightarrow A/\m_x^n$ such that the resulting left upper triangle commutes. Now the resulting lower right triangle commutes after composition with $\phi^*$; since $B/\m_x^nB$ is a semi-affinoid $K$-algebra and since the image of $|\sSp B/\m_x^nB|$ in $|\sSp A|$ must be contained in $U$, it follows from the universal property of $\phi$ that the lower right triangle commutes as well. It follows that the lower horizontal map is surjective. Furthermore, the diagonal map is surjective, and its kernel contains $\m_x^n B$; hence the lower horizontal map is in fact bijective. Taking $n=1$, we see that $\m_x B$ is a maximal ideal in $B$. Hence, the fiber of $|\phi|$ over $x$ contains exactly one point. Varying $x$, we see that $|\phi|$ is injective with image $U$. For fixed $x\in U$, let $y\in Y$ be the unique point over $x$; then $\m_y=\m_x B$ as we have just seen, and hence $\m_y^n=\m_x^n B$ for all integers $n\geq 1$. Since the lower horizontal map is bijective for all $n$, we obtain the final statement of the proposition.
\end{proof}

If $X$ is a semi-affinoid $K$-space, the universal property in Definition \ref{semiaffsubdomdefi} allows us to define a presheaf $\O_{X}$ on the category of semi-affinoid pre-subdomains in $X$ with values in the category of semi-affinoid $K$-algebras. We need to make choices in order to define $\O_X$, but the result is independent of choices up to unique isomorphism. 

\begin{defi}\label{semaffpresheafdefi}
Let $X$ be a semi-affinoid $K$-space. The presheaf $\O_X$ on the category of semi-affinoid pre-subdomains in $X$ is called the presheaf of semi-affinoid functions on $X$.
\end{defi}

Let us already mention that $\O_X$ is smaller than the presheaf of rigid-analytic functions on $X^\r$; for example, by Lemma \ref{supnormimmediatelem} all global sections of $\O_X$ are \emph{bounded}. 

If $U\subseteq|X|$ is a semi-affinoid pre-subdomain, we may view $U$ as the set of physical points on an essentially unique semi-affinoid $K$-space; thereby it makes sense to consider the set of semi-affinoid pre-subdomains of $U$. By definition, the semi-affinoid pre-subdomains of $U$ correspond to the semi-affinoid pre-subdomains of $X$ that are contained in $U$. In particular, the notion of semi-affinoid pre-subdomain is transitive.

Let us point out that so far, we have not yet established the existence of any non-trivial semi-affinoid pre-subdomain.

			\subsection{Semi-affinoid subdomains}\label{semaffsubdomsec}				
				In the following, we will focus on semi-affinoid pre-subdomains of a specific type, which we call semi-affinoid subdomains. We will not pursue the study of general semi-affinoid pre-subdomains any further. For example, we do not know whether all semi-affinoid pre-subdomains induce admissible open subsets under the functor $\r$ which is induced by Berthelot's construction.

\begin{defi}\label{semaffsubdomdefdef}
Let $\ul{\phi}\colon\fY\rightarrow\fX$ be a morphism of flat formal $R$-schemes of ff type. Then $\ul{\phi}$ is called \emph{simple}\index{morphism!simple} if it is 
\begin{packed_enum}
\item an open immersion,
\item a completion morphism or
\item an admissible formal blowup.
\end{packed_enum}
We say that  $\ul{\phi}$ \emph{defines a semi-affinoid subdomain}\index{morphism!defining a subdomain} if $\fX$ and $\fY$ are affine and if $\ul{\phi}$ is a composition of finitely many simple morphisms.
\end{defi}

Simple morphisms induce open immersions of rigid generic fibers: Indeed, by Proposition \ref{rigisoprop}, admissible formal blowups induce isomorphisms of generic fibers; open immersions induce retrocompact open immersions, cf. \cite{dJ} 7.2.2 and 7.2.4 (d), and completion morphisms induce possibly non-retrocompact open immersions, cf. \cite{dJ} 7.2.5.

Let us emphasize that if $\ul{\phi}\colon\fY\rightarrow\fX$ is a morphism of flat formal $R$-schemes of ff type that defines a semi-affinoid subdomain and if
\[
\fY\rightarrow\fY_{n}\rightarrow\cdots\rightarrow\fY_1\rightarrow\fX
\]
is a factorization of $\ul{\phi}$ into simple morphisms, then by definition $\fY$ and $\fX$ are affine, but the flat formal $R$-schemes $\fY_{i}$ need not be affine.

If $\ul{\psi}\colon\fZ\rightarrow\fY$ and $\ul{\phi}\colon\fY\rightarrow\fX$ are morphisms of affine flat formal $R$-schemes defining semi-affinoid subdomains, then the composition $\ul{\phi}\circ\ul{\psi}$ defines a semi-affinoid subdomain as well. This is immediate from Definition \ref{semaffsubdomdefdef}. Moreover, it is clear that if $\ul{\phi}\colon\fY\rightarrow\fX$ defines a semi-affinoid subdomain, then $\ul{\phi}^\srig$ is \emph{flat} in the sense that the underlying $K$-algebra homomorphism is flat: Indeed, flatness can be checked in completions along maximal ideals and, hence, after passing to Berthelot generic fibers.

\begin{defi}Let $X$ be a semi-affinoid $K$-space. 
\begin{enumerate}
\item A subset $U\subseteq|X|$ is called a \emph{semi-affinoid subdomain}\index{subdomain!semi-affinoid} if there exist a flat affine formal $R$-model $\fX$ of $X$, an affine flat formal $R$-scheme $\fY$ and a morphism $\ul{\phi}\colon\fY\rightarrow\fX$ defining a semi-affinoid subdomain such that $U=\im|\ul{\phi}^\srig|$. We say that $\ul{\phi}$ defines $U$ as a semi-affinoid subdomain of $X$.
\item If we can choose $\ul{\phi}$ to be simple, then $U$ is called a \emph{simple} semi-affinoid subdomain.\index{subdomain!simple} 
\end{enumerate}
\end{defi}

\begin{example}
Let $X$ be a semi-affionid $K$-space, and let us consider semi-affinoid functions $f_1,\ldots,f_m$, $g_1,\ldots,g_n$, $h$ on $X$ without a common zero. Then the subset
\[
\{x\in X\,;\,\forall i,j\,:\,|f_i(x)|<|h(x)|,|g_j(x)|\leq|h(x)|\}
\]
of $X$ is a semi-affinoid subdomain. 
\end{example}
\begin{proof}
Let $\ul{A}$ be an $R$-lattice of ff type of the semi-affinoid $K$-algebra $A$ of functions on $X$, and let $n\in \N$ be a natural number such that all $\pi^nh,\pi^nf_i,\pi^ng_j$ are contained in $\ul{A}$. By Proposition \ref{ratsubdoms}, the above set is the semi-affinoid subdomain defined by the morphism of affine formal $R$-schemes that is associated to the natural $R$-homomorphism
\[
\ul{A}\rightarrow\ul{A}\left\llbracket S_1,\ldots,S_m\right\rrbracket\langle T_1,\ldots,T_n\rangle/(\pi^nhS_i-\pi^nf_i,\pi^nhT_j-\pi^ng_j)/(\textup{$\pi$-torsion})\quad,
\]
which is obtained as follows:
\begin{packed_enum}
\item Let us write $\fX=\Spf\ul{A}$, and let $\fX'$ be the admissible formal blowup of $\fX$ in the ideal generated by $(\pi^nh,\pi^nf_i,\pi^ng_j)$; this ideal is $\pi$-adically open since $h$, the $f_i$ and the $g_j$ generate the unit ideal in $A$.
\item Let $\fX_h\subseteq\fX'$ be the affine open part where $h$ is invertible.
\item Let $\fX_h^\wedge$ be the completion of $\fX_h$ along the ideal generated by the $f_i/h$.
\end{packed_enum}
Then $\fX_h^\wedge\rightarrow\fX$ corresponds to the homomorphism above; this is easily seen by checking universal properties.
\end{proof}

\begin{prop}\label{subdomunivprop}Let $X$ be a semi-affinoid $K$-space, and let $U\subseteq |X|$ be a semi-affinoid subdomain. Then $U$ is a semi-affinoid pre-subdomain in $X$, and if $\ul{\phi}\colon\fY\rightarrow\fX$ is an $R$-morphism defining $U$ as a semi-affinoid subdomain in $X$, then $\ul{\phi}^\srig$ represents all semi-affinoid morphisms to $X$ with image in $U$.
\end{prop}
\begin{proof}
Let us write $Y\mathrel{\mathop:}=\fY^\srig$, $\phi\mathrel{\mathop:}=\ul{\phi}^\srig$. We need to show that for any morphism $\psi\colon Z\rightarrow X$ of semi-affinoid $K$-spaces with $\im|\psi|\subseteq U$, there exists a unique morphism $\xi\colon Z\rightarrow Y$ such that the diagram
\[
\begin{diagram}
&&Y\\
&\ruDashto<\xi&\dTo>\phi\\
Z&\rTo^\psi&X
\end{diagram}
\]
commutes.

By Proposition \ref{presubdomstalkscompprop}, $|\phi|$ is an injection, and $\phi^*$ induces isomorphisms of completed stalks; hence uniqueness of $\xi$ follows from Krull's Intersection Theorem. 

Let us show existence. By Corollary \ref{freesemaffcor} ($iv$), there exists an affine flat formal $R$-model $\fZ$ of $Z$ such that $\psi$ extends to an $R$-morphism $\ul{\psi}\colon\fZ\rightarrow \fX$. It suffices to establish the following claim: There exist a \emph{finite} admissible formal blowup $\ul{\beta}\colon\fZ'\rightarrow \fZ$ and a morphism $\ul{\xi}'\colon\fZ'\rightarrow \fY$ such that the diagram
\[
\begin{diagram}
\fZ'&\rDashto^{\ul{\xi}}&\fY\\
\dDashto<{\ul{\beta}}&&\dTo>{\ul{\phi}}\\
\fZ&\rTo^{\ul{\psi}}&\fX
\end{diagram}
\]
of affine flat formal $R$-schemes commutes. To prove this claim, let us consider a factorization
\[
\fY\overset{\ul{\phi}_n}{\longrightarrow}\fY_{n}\overset{\ul{\phi}_{n-1}}{\longrightarrow}\cdots\overset{\ul{\phi}_1}{\longrightarrow}\fY_1\overset{\ul{\phi}_0}{\longrightarrow}\fX
\]
of $\ul{\phi}$ into simple morphisms $\ul{\phi}_i$.  Let us write $\fY_{n+1}\mathrel{\mathop:}=\fY$, $\fY_0\mathrel{\mathop:}=\fX$. We first show that for each $0\leq i\leq n+1$, there exist an admissible formal blowup $\ul{\beta}_i\colon\fZ_i\rightarrow\fZ$ of $\fZ$ and a morphism $\ul{\xi}_i\colon\fZ_i\rightarrow\fY_i$ such that the diagram
\[
\begin{diagram}
\fZ_i&\rDashto^{\ul{\xi}_i}&\fY_i\\
\dDashto<{\ul{\beta}_i}&&\dTo>{\ul{\phi}_0\circ\cdots\circ\ul{\phi}_{i-1}}\\
\fZ&\rTo^{\ul{\psi}}&\fX
\end{diagram}
\]
commutes. We proceed by induction on $i$. For $i=0$ there is nothing to show. Let us assume that the statement holds for 
some $0\leq i\leq n$, and let us show that it holds for $i+1$ as well. If $\ul{\phi}_i$ is an open immersion or a completion morphism, we may choose $\fZ_{i+1}=\fZ_i$, $\ul{\beta}_{i+1}=\ul{\beta}_i$, and we let $\ul{\xi}_{i+1}$ denote the unique factorization of $\ul{\xi}_i$ that exists by surjectivity of the specialization map. 
If $\ul{\phi}_i$ is an admissible formal blowup, we let $\ul{\alpha}_{i+1}\colon\fZ_{i+1}\rightarrow\fZ_i$ denote the strict transform of $\fZ_i$ under $\ul{\phi}_i$ via $\ul{\xi}_i$; then $\ul{\xi}_{i+1}$ exists by the universal property of $\ul{\alpha}_{i+1}$, and by Proposition \ref{compstableprop} the composition $\ul{\beta}_{i+1}\mathrel{\mathop:}=\ul{\beta}_i\circ\ul{\alpha}_{i+1}$ is an admissible formal blowup. The claim now follows from Corollary \ref{affinefineblowupcor}, applied to $\ul{\beta}_{n+1}$.
\end{proof}




Among all semi-affinoid subdomains of a given semi-affinoid $K$-space, we are particularly interested in those that are defined by a morphism $\ul{\phi}\colon\fY\rightarrow\fX$ that admits a factorization into simple morphisms that are not completion morphisms.

\begin{defi}A simple morphism of flat formal $R$-schemes of ff type that is not a completion morphism is called a \emph{retrocompact} simple morphism.
\end{defi}

In other words, the retrocompact simple morphisms of flat formal $R$-schemes of ff type are the open immersions and the admissible formal blowups. In particular, retrocompact morphisms are adic, and under Berthelot's construction they induce retrocompact open immersions. This explains our terminology.

\begin{defi}
Let $X$ be a semi-affinoid $K$-space. A semi-affinoid subdomain $U\subseteq X$ is called
\begin{packed_enum}\index{subdomain!retrocompact}
\item \emph{retrocompact} if there exists an $R$-morphism $\ul{\phi}\colon\fY\rightarrow\fX$ defining $U$ that is a composition of finitely many retrocompact simple morphisms. 
\item \emph{formally open} if it is defined by an open immersion $\ul{\phi}$.\index{subdomain!formally open}
\end{packed_enum}
Such a morphism $\ul{\phi}$ is said to represent $U$ as a \emph{retrocompact} / \emph{formally open} semi-affinoid subdomain in $X$.
\end{defi}

The retrocompact semi-affinoid subdomains in affinoid $K$-spaces are affinoid:

\begin{lem}\label{veryspecialinafflem}
Let $A$ be an affinoid $K$-algebra. A retrocompact semi-affinoid subdomain $U\subseteq|\Sp A|=|\sSp A|$ in $\sSp A$ is an affinoid subdomain in $\Sp A$.
\end{lem}
\begin{proof}
Let $\ul{\phi}\colon\fY\rightarrow \fX$ be a retrocompact $R$-morphism defining $U$ in $X$. By Corollary \ref{afflattfintypecor}, $\fX$ is of tf type over $R$. Since $\ul{\phi}$ is a retrocompact morphism, it is adic, and it follows that $\fY$ is of tf type over $R$ as well. Hence, $\phi\mathrel{\mathop:}=\ul{\phi}\otimes_RK$ is a morphism of affinoid $K$-spaces. By  Proposition \ref{subdomunivprop}, $\phi$ represents all semi-affinoid maps with image in $U$; in particular it represents all affinoid maps with image in $U$. Hence, $U$ is an affinoid subdomain in $\Sp A$.
\end{proof}

Conversely, it is clear that for any affinoid $K$-algebra $A$, the rational subdomains in $\Sp A$ define semi-affinoid subdomains in $\sSp A$. Let $U\subseteq\Sp A$ be a general affinoid subdomain in $\Sp A$. By the Theorem of Gerritzen and Grauert (\cite{BGR} 7.3.5/1), $U$ is a finite union of rational subdomains. Let $\fX$ be any affine flat formal $R$-model of tf type for $\Sp A$. By \cite{BL1} Lemma 4.4, there exist an admissible formal blowup $\fX'\rightarrow\fX$ of $\fX$ and an open formal subscheme $\fU\subseteq\fX'$ such that $U=\fU^\rig$. However, we do not know whether $\fU$ is affine, so we do not know whether a general affinoid subdomain $U$ in $\Sp A$ is a (retrocompact) semi-affinoid subdomain or even a semi-affinoid pre-subdomain in $\sSp A$. Nonetheless, we will see that affinoid subdomains in $\Sp A$ are admissible open in the uniformly rigid G-topology on $\sSp A$, cf. Corollary \ref{finretroadmcor}.

There exists a semi-affinoid analogon of the Gerritzen-Grauert theorem, saying that every semi-affinoid pre-subdomain is a finite union of semi-affinoid subdomains, cf.\ \cite{LiRo} I.6.2. However, a finite union of semi-affinoid subdomains needs not be admissible in the uniformly rigid G-topology that we will define; for example, the uniformly rigid closed unit disc $\sSp K\langle T\rangle$ is \emph{non-admissibly} covered by the uniformly rigid open unit disc $\sSp (R[[T]]\otimes_RK)$ and the uniformly rigid unit circle $\sSp K\langle T, T^{-1}\rangle$, cf.\ Example \ref{nonadmdisccovex}. Hence, the Gerritzen-Grauert result just mentioned does not make general semi-affinoid pre-subdomains more useful. We shall thus focus on semi-affinoid pre-subdomains that are semi-affinoid subdomains. However, let us already mention that semi-affinoid pre-subdomains which are semi-affinoid sub\emph{spaces} in the sense of Section \ref{srigspacessec} \emph{do} admit a Gerritzen-Grauert type description, cf.\ Lemma \ref{semaffsubspacecharlem}.


			\subsection{Preimages of semi-affinoid subdomains}\label{preimsemaffsec}
				\begin{lem}\label{presubdompreimlem}
Let $X$ be a semi-affinoid $K$-space, let $U\subseteq|X|$ be a semi-affinoid pre-subdomain, and let $\psi\colon Z\rightarrow X$ be any morphism of semi-affinoid $K$-spaces. Then $\psi^{-1}(U)$ is a semi-affinoid pre-subdomain in $Z$. If $\phi\colon Y\rightarrow X$ is a morphism of semi-affinoid $K$-spaces representing all semi-affinoid morphisms to $X$ with image in $U$ and if
\[
\begin{diagram}
Y\times_XZ&\rTo^{\phi'}&Z\\
\dTo<{\psi'}&&\dTo>\psi\\
Y&\rTo^\phi&X
\end{diagram}
\]
is the induced cartesian square in the category of semi-affinoid $K$-spaces, then $\phi'$ represents all semi-affinoid morphisms to $Z$ with image in $\psi^{-1}(U)$.
\end{lem}
\begin{proof}
It suffices to verify the second assertion, which is an immediate consequence of the universal property of fibered products in the category of semi-affinoid $K$-spaces.
\end{proof}

In fact, if $U$ is a semi-affinoid subdomain in $X$, then $\psi^{-1}(U)$ is a semi-affinonid subdomain in $Z$. To prove this, we need to use the fact that simple morphisms are preserved by base change in the category of \emph{flat} formal $R$-schemes of ff type: Let
\[
\begin{diagram}
\fY'&\rTo^{\ul{\phi}'}&\fZ&\\
\dTo<{\ul{\psi}'}&&\dTo>{\ul{\psi}}&\\
\fY&\rTo^{\ul{\phi}}&\fX&
\end{diagram}
\]
be a diagram of flat formal $R$-schemes of ff type that is cartesian in the category of flat formal $R$-schemes of ff type. If $\ul{\phi}$ is a simple morphism, then $\ul{\phi}'$ is a simple morphism as well. Indeed, if $\ul{\phi}$ is an open immersion or a completion morphism, then $\ul{\phi}$ is flat. By Proposition \ref{flatnessbasechangeprop}, flatness is preserved under base change with respect to morphisms of ff type; hence the above square is cartesian in the category of all formal schemes. Open immersions are clearly preserved by any base change, and the same visibly holds for completion morphisms; hence we see that $\ul{\phi}'$ is an open immersion or a completion morphism. If $\ul{\phi}$ is an admissible formal blowup, then $\ul{\phi}'$ is an admissible formal blowup by Lemma \ref{admblowupcartlem}. 

\begin{prop}\label{subdomdefibcstableprop}
Let 
\[
\begin{diagram}
\fY'&\rTo^{\ul{\phi}'}&\fZ&\\
\dTo<{\ul{\psi}'}&&\dTo>{\ul{\psi}}&\\
\fY&\rTo^{\ul{\phi}}&\fX&
\end{diagram}
\]
be a diagram of morphisms of affine flat formal $R$-schemes of ff type that is cartesian in the category of flat formal $R$-schemes of ff type. If $\ul{\phi}$ defines a semi-affinoid subdomain, then $\ul{\phi}'$ defines a semi-affinoid subdomain as well. Moreover, if $\phi$, $\phi'$ and $\psi$ denote the semi-affinoid generic fibers of $\ul{\phi}$, $\ul{\phi}'$ and $\ul{\psi}$ respectively, then
\[
\im\phi'\,=\,\psi^{-1}(\im\phi)\quad.
\]
\end{prop}
\begin{proof}
By construction of the fibered product in the category of flat formal $R$-schemes of ff type, $\fY'$ is affine. Hence, to prove that $\ul{\phi}'$ defines a semi-affinoid subdomain, it suffices to see that $\ul{\phi}'$ admits a factorization into simple morphisms. Let 
\[
\fY\overset{\ul{\phi}_n}{\longrightarrow}\fY_{n}\overset{\ul{\phi}_{n-1}}{\longrightarrow}\cdots\overset{\ul{\phi}_1}{\longrightarrow}\fY_1\overset{\ul{\phi}_0}{\longrightarrow}\fX
\]
be a factorization of $\ul{\phi}$ into simple morphisms; we claim that it induces a factorization of $\ul{\phi}'$ into simple morphisms via base change in the category of flat formal $R$-schemes of ff type. Indeed, this follows immediately from the discussion preceding this Lemma, together with the fact that compositions of cartesian diagrams are cartesian.

The statement on the generic images follows for example from Lemma \ref{presubdompreimlem} combined with the fact that cartesian diagrams of flat affine formal $R$-schemes induce cartesian diagrams of semi-affinoid $K$-spaces.
\end{proof}

\begin{cor}\label{semaffpreimcor}
Let $\phi\colon Y\rightarrow X$ be a morphism of semi-affinoid $K$-spaces, and let $U\subseteq X$ be a semi-affinoid subdomain in $X$. Then $\phi^{-1}(U)$ is a semi-affinoid subdomain in $Y$.
\end{cor}
\begin{proof}
Let $\ul{\psi}\colon\fZ\rightarrow\fX$ be an $R$-morphism representing $U$ as a semi-affinoid subdomain in $X$. By Corollary \ref{freesemaffcor} ($iv$), there exists an affine flat formal $R$-model of ff type $\fY$ for $Y$ together with a model $\ul{\phi}\colon\fY\rightarrow\fX$ of $\phi$. By Proposition \ref{subdomdefibcstableprop}, the projection $\fY\times'_\fX\fZ\rightarrow\fY$ represents $\phi^{-1}(U)$ as a semi-affinoid subdomain in $Y$.
\end{proof}


Semi-affinoid subdomains in $X$ can be represented with respect to any given flat affine formal $R$-model of ff type for $X$:

\begin{cor}\label{subdomdefgiventargetcor}
Let $X$ be a semi-affinoid $K$-space, let $U\subseteq |X|$ be a semi-affinoid subdomain, and let $\fX$ be an affine flat formal $R$-model of ff type for $X$. There exists an $R$-morphism
\[
\ul{\psi}\colon\fY\rightarrow\fX
\]
representing $U$ as a semi-affinoid subdomain in $X$.
\end{cor}
\begin{proof}
Let us write $\fX_1\mathrel{\mathop:}=\fX$, and let $\ul{\phi}\colon\fY_2\rightarrow\fX_2$ be an $R$-morphism representing $U$ as a semi-affinoid subdomain in $X$. By Corollary \ref{freesemaffcor} ($ii$), there exists a flat affine formal $R$-scheme $\fX_{12}$ together with finite admissible formal blowups 
\[
\ul{\beta}_i\colon\fX_{12}\rightarrow\fX_i\quad,\quad i=1,2\quad.
\]
By Proposition \ref{subdomdefibcstableprop} and Lemma \ref{admblowupcartlem}, the strict transform 
\[
\ul{\phi}'\colon\fY_{12}\rightarrow\fX_{12}
\]
of $\ul{\phi}$ under $\ul{\beta}_2$ defines $U$ as a semi-affinoid subdomain in $X$ as well. Hence,
\[
\ul{\psi}\,\mathrel{\mathop:}=\,\ul{\beta}_1\circ\ul{\phi}'\,:\,\fY_{12}\rightarrow\fX
\]
also represents $U$ as a semi-affinoid subdomain in $X$.
\end{proof}

\begin{cor}\label{domainmodcor}
Let $X$ be a semi-affinoid $K$-space, let $U$ be a semi-affinoid subdomain in $X$, and let $\ul{\phi}\colon\fY\rightarrow\fX$ be a morphism of affine flat formal $R$-schemes of ff type that is a model of the natural morphism of semi-affinoid $K$-spaces $\phi\colon U\rightarrow X$. There exists a finite admissible formal blowup $\ul{\beta}\colon\fY'\rightarrow\fY$ such that $\ul{\phi}\circ\ul{\beta}$ defines $U$ as a semi-affinoid subdomain in $X$.
\end{cor}
\begin{proof}
By Corollary \ref{subdomdefgiventargetcor}, there exists an $R$-morphism $\ul{\psi}\colon\fY''\rightarrow\fX$ defining $U$ as a semi-affinoid subdomain in $X$. By Corollary \ref{freesemaffcor} ($ii$), there exists an affine flat formal $R$-scheme $\fY'$ together with finite admissible formal blowups $\ul{\beta}\colon\fY'\rightarrow\fY$ and $\ul{\gamma}\colon\fY'\rightarrow\fY''$. Now $\ul{\psi}\circ\ul{\gamma}$ defines $U$ as a semi-affinoid subdomain in $X$, and $\ul{\psi}\circ\ul{\gamma}=\ul{\phi}\circ\ul{\beta}$ since both morphisms agree on semi-affinoid generic fibers.
\end{proof}


\begin{cor}\label{semaffsubdomgoodnotioncor}
Let $X$ be a semi-affinoid $K$-space, let $U\subseteq |X|$ be a semi-affinoid subdomain with its canonical structure of semi-affinoid $K$-space, and let $V\subseteq U$ be a subset. Then the following are equivalent:
\begin{packed_enum}
\item $V$ is a semi-affinoid subdomain in $U$.
\item $V$ is a semi-affinoid subdomain in $X$.
\end{packed_enum}
\end{cor}
\begin{proof}
The implication ($ii$)$\Rightarrow$($i$) follows from Corollary \ref{semaffpreimcor}. To show the converse implication ($i$)$\Rightarrow$($ii$), let us assume that $V$ is a semi-affinoid subdomain in $U$, and let 
\[
\ul{\phi}_U\colon\fY\rightarrow\fX
\]
be an $R$-morphisms representing $U$ as a semi-affinoid subdomain in $X$. 
By Corollary \ref{subdomdefgiventargetcor}, there exists an $R$-morphism
\[
\ul{\phi}_V\colon\fZ\rightarrow\fY
\]
representing $V$ as a semi-affinoid subdomain on $U$. The composition $\ul{\phi}_U\circ\ul{\phi}_V$ then defines $V$ as a semi-affinoid subdomain in $X$.
\end{proof}

\begin{cor}\label{semiaffintcor}
The set of semi-affinoid subdomains in a semi-affinoid $K$-space is stable under the formation of intersections.
\end{cor}
\begin{proof}
Let $X$ be a semi-affinoid $K$-space, and let $U,V\subseteq |X|$ be semi-affinoid subdomains in $X$. By Corollary \ref{semaffpreimcor}, $U\cap V$ is a semi-affinoid subdomain in $U$; by Corollary \ref{semaffsubdomgoodnotioncor}, it follows that $U\cap V$ is a semi-affinoid subdomain in $X$.
\end{proof}

All results in this section remain valid if we consider retrocompact simple morphisms and retrocompact semi-affinoid subdomains,
the proofs remaining unchanged. For example, if $X$ is a semi-affinoid $K$-space and if $U$ is a retrocompact semi-affinoid subdomain in $X$, then $V\subseteq U$ is a retrocompact semi-affinoid subdomain in $X$ if and only if it is a retrocompact semi-affinoid subdomain in $U$. Of course, if $V$ is a retrocompact semi-affinoid subdomain of a general semi-affinoid subdomain in $X$, then $V$ needs not be retrocompact in $X$.

			\subsection{G-topologies on semi-affinoid spaces}\label{gtopsec}
				In this section, we define a saturated G-topology on the category of semi-affinoid $K$-spaces.

Let us first recall the notion of a G-topology on a pair $(\mathfrak{C},|\cdot|)$, where $\mathfrak{C}$ is a category and where $|\cdot|\colon\mathfrak{C}\rightarrow\text{Sets}$ is a functor from $\mathfrak{C}$ to the category of sets; cf \cite{BGR} 9.1.2. Let GTop denote the category of G-topological spaces.

\begin{defi}\label{gtopdefi}
Let $\mathfrak{C}$ be a category, and let $|\cdot|\colon\mathfrak{C}\rightarrow\text{Sets}$ be a functor. A G-topology $\sT$ on the pair $(\mathfrak{C},|\cdot|)$ is a factorization
\[
\begin{diagram}
\mathfrak{C}&&\rTo^{|\cdot|}&&\text{Sets}&\\
&\rdTo<{|\cdot|_\sT}&&\ruTo>F&&\\
&&\textup{GTop}&&&\quad,
\end{diagram}
\]
where $F$ denotes the forgetful functor from \textup{GTop} to Sets.
\end{defi}


The category $\sAff_K$ of semi-affinoid $K$-spaces gives rise to a pair $(\sAff_K,|\cdot|)$ as above, where $|\cdot|$ is the functor sending a semi-affinoid $K$-space to the set of its physical points. 
We first define an auxiliary G-topology $\sT_\textup{aux}$ on $(\sAff_K,|\cdot|)$ and then obtain the uniformly rigid G-topology $\sT_\srig$ on $(\sAff_K,|\cdot|)$ by means of the general saturation procedure described in \cite{BGR} 9.1.2.

Let $X$ be a semi-affinoid $K$-space. The $\sT_\textup{aux}$-admissible subsets of $X$ will be given by the semi-affinoid subdomains of $X$. We will say that a covering of a semi-affinoid subdomain in $X$ by semi-affinoid subdomains is $\sT_\textup{aux}$-admissible if it admits a leaflike refinement. In order to define leaflike coverings, we first need to define treelike coverings.

\subsubsection{Treelike and leaflike coverings}

First, we introduce the notion of a simple covering.

\begin{defi}\label{simplecovdefi}
Let $X$ be a semi-affinoid $K$-space. A finite family $(X_i)_{i\in I}$ of semi-affinoid subdomains of $X$ is called a \emph{simple} covering\index{covering!simple} of $X$ if there exist
\begin{packed_enum}
\item an affine flat formal $R$-model of ff type $\fX$ for $X$,
\item an admissible formal blowup $\ul{\beta}\colon\fX'\rightarrow\fX$ and
\item an affine open covering $(\fX_i)_{i\in I}$ of $\fX'$
\end{packed_enum}
such that for each $i\in I$, the restriction $\ul{\beta}|_{\fX_i}\colon\fX_i\rightarrow\fX$ defines $X_i$ as a semi-affinoid subdomain in $X$. The formal data given by \textup{(}i\textup{)}-\textup{(}iii\textup{)} is called a \emph{formal presentation}\index{presentation!of a simple covering} of the simple covering $(X_i)_{i\in I}$.
\end{defi}

In particular, if $(X_i)_{i\in I}$ is a simple covering of $X$, then $X_i$ is a retrocompact simple semi-affinoid subdomain of $X$ for all $i\in I$. The simple coverings generalize the rational coverings of affinoid rigid geometry. Let us observe that simple coverings are stable under pullbacks:

\begin{lem}\label{simplecovintlem}
Let $X$ be a semi-affinoid $K$-space, let $(X_i)_{i\in I}$ be a simple covering of $X$, and let $\phi\colon Y\rightarrow X$ be a morphism of semi-affinoid $K$-spaces. Then $(\phi^{-1}(X_i))_{i\in I}$ is a simple covering of $Y$.
\end{lem}
\begin{proof}
Let $\fX$, $\ul{\beta}\colon\fX'\rightarrow\fX$ and $(\fX_i)_{i\in I}$ be a formal presentation of $(X_i)_{i\in I}$ as in Definition \ref{simplecovdefi}. Let $\ul{\phi}\colon\fY\rightarrow\fX$ be any $R$-morphism of affine flat formal $R$-schemes of ff type that is a model of $\phi$. Let $\ul{\beta}'\colon\fY'\rightarrow\fY$ be the induced admissible formal blowup of $\fY$, and let $(\fY_i)_{i\in I}$ be the induced admissible open covering of $\fY$. By Proposition \ref{subdomdefibcstableprop}, this formal data presents $(\phi^{-1}(X_i))_{i\in I}$ as a simple covering of $Y$.
\end{proof}

Before discussing treelike coverings, let us introduce some terminology concerning rooted trees. We recall that a \emph{tree}\index{tree} is a connected graph without circles. A \emph{rooted} tree is a tree with a distinguished vertex, called the \emph{root} of the rooted tree. If $I$ is a rooted tree, then the set of vertices of $I$ is equipped with a natural partial ordering $\prec$, where $i\prec j$ for vertices $i$ and $j$ of $I$ if and only if the unique path from the root of $I$ to $j$ passes through $i$. The \emph{parent} $\parent(i)$ of a non-root vertex $i$ is the vertex connected to $i$ on the path to the root, and a \emph{child} of a vertex $i$ is a vertex of which $i$ is the parent. A \emph{sibling} of non-root vertex $i$ is a non-root vertex whose parent coincides with the parent of $i$. If $i$ is a vertex of a rooted tree, we let $\textup{ch}(i)$ denote the set of the children of $i$. A \emph{leaf} of a rooted tree is a vertex without children, and a vertex that is not a leaf is called an \emph{inner} vertex; we let $\leaves(I)$ denote the set of leaves of $I$, and we let $\inner(I)$ denote the set of inner vertices of $I$. A \emph{linear} rooted tree is a rooted tree where each vertex has at most one child, and an \emph{almost linear} tree is a rooted tree where each vertex has at most one child except possibly for the root. 
For each vertex $i$ of a finite rooted tree $I$, we let $\subtree(i)$ denote the rooted tree with root $i$ given by all vertices $j$ of $I$ with the property that $i$ lies on the path of $j$ to the root of $I$. We let $v(I)$ denote the number of vertices of $I$, and we call $v(I)$ the \emph{volume} of $I$. 
For each vertex $i\in I$ we let $l(i)$ denote the \emph{length} of the unique path from $i$ to the root, that is, the number of edges of this path. 

If $I$ is a set, the structure of a rooted tree on $I$ is an identification of $I$ with the set of vertices of a rooted tree.

\begin{defi}\label{treelikecovdefi}
Let $X$ be a semi-affionid $K$-space, and let $(X_i)_{i\in I}$ be a finite family of semi-affinoid subdomains of $X$. A rooted tree structure on $I$ is said to be \emph{suitable for $(X_i)_{i\in I}$} or simply \emph{suitable} if the following holds:
\begin{packed_enum}
\item If $r\in I$ is the root, then $X_r=X$.
\item For all inner vertices $i$ of $I$, $(X_j)_{j\in\textup{ch}(i)}$ is a simple covering of $X_i$.
\end{packed_enum}
A finite covering of $X$ by semi-affinoid subdomains is called \emph{treelike}\index{covering!treelike} if its index set admits a suitable rooted tree structure.
\end{defi}

\begin{lem}\label{treelikecovintlem}
Let $X$ be a semi-affinoid $K$-space, let $(X_i)_{i\in I}$ be a treelike covering of $X$, and let $\phi\colon Y\rightarrow X$ be a morphism of semi-affinoid $K$-spaces. Then any rooted tree structure on $I$ that is suitable for $(X_i)_{i\in I}$ is suitable for $(\phi^{-1}(X_i))_{i\in I}$. In particular, this covering of $Y$ is treelike.
\end{lem}
\begin{proof}
Let us fix a rooted tree structure on $I$ that is suitable for the treelike covering $(X_i)_{i\in I}$ of $X$. Let $r\in I$ be the root; then $\phi^{-1}(X_r)=\phi^{-1}(X)=Y$. Moreover, let $i\in I$ be an inner vertex. Since $(X_j)_{j\in\textup{ch}(i)}$ is a simple covering of $X_i$, we conclude from Lemma \ref{simplecovintlem} that $(\phi^{-1}(X_j))_{j\in\textup{ch}(i)}$ is a simple covering of $\phi^{-1}(X_i)$. Hence, the chosen rooted tree structure on $I$ is suitable for the covering $(\phi^{-1}(X_i))_{i\in I}$ of $Y$.
\end{proof}

\begin{defi}\label{leaflikecovdefi}
Let $X$ be a semi-affinoid $K$-space. A finite family $(X_i)_{i\in I}$ of semi-affinoid subdomains of $X$ is called a \emph{leaflike} covering\index{covering!leaflike} of $X$ if it can be enlarged to a treelike covering $(X_i)_{i\in J}$, $I\subseteq J$, where $J$ admits a suitable rooted tree structure such that $I$ is identified with the set of leaves of $J$.
\end{defi}

Of course, simple coverings are leaflike coverings. Moreover, any constituent of a treelike or leaflike covering of $X$ is a retrocompact semi-affinoid subdomain of $X$. Intuitively speaking, leaflike coverings are iterated rational coverings.

\begin{lem}\label{leaflikecovintlem}
Let $X$ be a semi-affinoid $K$-space, let $(X_i)_{i\in I}$ be a leaflike covering of $X$, and let $\phi\colon Y\rightarrow X$ be a morphism of semi-affinoid $K$-spaces. Then $(\phi^{-1}(X_i))_{i\in I}$ is a leaflike covering of $Y$.
\end{lem}
\begin{proof}
Let $(X_i)_{i\in J}$ be an extension of $(X_i)_{i\in I}$ to a treelike covering of $X$ satisfying the condition that there exists a suitable rooted tree structure on $J$ such that $I\subseteq J$ is identified with the set of leaves. By Lemma \ref{treelikecovintlem}, this rooted tree structure on $J$ is suitable for $(\phi^{-1}(X_i))_{i\in J}$; hence $(\phi^{-1}(X_i))_{i\in I}$ is a leaflike covering of $Y$.
\end{proof}

\begin{lem}\label{leaflikestablelem}
Let $U\subseteq X$ be a semi-affinoid subdomain, let $(U_i)_{i\in I}$ be a leaflike covering of $U$, and for each $i\in I$ let $(V_{ij})_{j\in J_i}$ be a leaflike covering of $U_i$. Then $(V_{ij})_{i\in I,j\in J_i}$ is a leaflike covering of $U$.
\end{lem}
\begin{proof}
Let us choose a treelike covering $(U_i)_{i\in I'}$  of $U$ extending $(U_i)_{i\in I}$ together with a suitable rooted tree structure on $I'$ such that $I\subseteq I'$ is the set of leaves. Similarly, for each $i\in I$ we choose a treelike covering $(V_{ij})_{j\in J_i'}$ extending $(V_{ij})_{j\in J_i}$ together with a suitable rooted tree structure on $J_i'$ such that $J_i\subseteq J_i'$ is identified with the set of leaves for all $i\in I$. For each $i\in I$, we glue the rooted tree $J_i'$ to the rooted tree $I'$ by identifying the root of $J_i'$ with the leaf $i$ of $I'$. We obtain a rooted tree $J'$ whose set of leaves is identified with the disjoint union of the sets $J_i$, $i\in I$. For each $i\in I$, $U_i=V_{i r_i}$, where $r_i$ is the root of $J_i'$; hence we obtain a covering $(V_j)_{j\in J'}$ such that the given rooted tree structure on $J'$ is suitable for $(V_j)_{j\in J'}$; indeed, this can be checked locally on the rooted tree $J'$. We conclude that the composite covering $(V_{ij})_{i\in I,j\in J_i}$ of $U$ is leaflike.
\end{proof}


\subsubsection{Integral models of treelike coverings}





A treelike covering is built up from simple coverings, and simple coverings are defined in terms of formal presentations. We now show that for any treelike covering, the formal presentations of the underlying simple coverings can be chosen in a coherent way. To prove this result, we need to know that simple coverings of a semi-affinoid $K$-space $X$ admit formal presentations with respect to any given affine flat formal $R$-model of ff type for $X$:

\begin{lem}\label{simplecovadaptedtomodellem}
Let $X$ be a semi-affinoid $K$-space, let $(X_i)_{i\in I}$ be a simple covering of $X$, and let $\fX$ be an affine flat formal $R$-model of ff type for $X$. Then there exist
\begin{packed_enum}
\item an admissible formal blowup $\ul{\beta}\colon\fX'\rightarrow\fX$ and
\item an affine open covering $(\fX_i)_{i\in I}$ of $\fX'$
\end{packed_enum}
such that for each $i\in I$, the restriction $\ul{\beta}|_{\fX_i}\colon\fX_i\rightarrow\fX$ defines $X_i$ as a semi-affinoid subdomain in $X$.
\end{lem}
\begin{proof}
By assumption, there exist an affine flat formal $R$-model of ff type $\fY$ for $X$, an admissible formal blowup $\ul{\gamma}\colon\fY'\rightarrow\fY$ and an affine open covering $(\fY_i)_{i\in I}$ of $\fY'$ such that for each $i\in I$, the restriction $\ul{\gamma}|_{\fY_i}\colon\fY_i\rightarrow\fY$ defines $X_i$ as a semi-affinoid subdomain in $X$. By Proposition \ref{affineblowupprop}, there exists an affine flat formal $R$-scheme $\fZ$ together with finite admissible formal blowups $\ul{\beta}_1\colon\fZ\rightarrow\fX$ and $\ul{\beta}_2\colon\fZ\rightarrow\fY$.
Let $\ul{\delta}\colon\fX'\rightarrow\fZ$ denote the strict transform of $\ul{\gamma}$ under $\ul{\beta}_2$; then $\ul{\delta}$ is an admissible formal blowup by Lemma \ref{admblowupcartlem}, and hence the composition $\ul{\beta}\mathrel{\mathop:}=\ul{\beta}_1\circ\ul{\delta}\colon\fX'\rightarrow\fX$ is an admissible formal blowup by Proposition \ref{compstableprop}. For $i\in I$, let $\fX_i\subseteq\fX'$ be the preimage of $\fY_i$ under the projection from $\fX'=\fY'\times'_\fY\fZ$ to $\fY'$. Since this projective is finite, $\fX_i$ is affine. We claim that $\ul{\beta}|_{\fX_i}\colon\fX_i\rightarrow\fX$ represents $X_i$ as  a semi-affinoid subdomain in $X$ for all $i\in I$. And indeed, $\ul{\beta}|_{\fX_i}$ is the strict transform of $\ul{\gamma}|_{\fY_i}$ under $\ul{\beta_2}$ composed with $\ul{\beta}_1$, so the statement follows from Proposition \ref{subdomdefibcstableprop}.
\end{proof}

\begin{defi}\label{formalprecovdefi}
A \emph{treelike formal covering}\index{covering!treelike formal} is the following data:
\begin{packed_enum}
\item a rooted tree $I$,
\item a family $(\fX_i)_{i\in I}$ of affine flat formal $R$-schemes of ff type,
\item for each inner vertex $i\in I$ an admissible formal blowup $\ul{\beta}_i\colon\fX_i'\rightarrow\fX_i$ and
\item for each inner vertex $i\in I$ and for each child $j$ of $i$ an open immersion 
\[
\ul{\phi}_j\colon\fX_j\hookrightarrow\fX_i'
\]
\end{packed_enum}
such that for each inner vertex $i\in I$,
\[
\fX_i'\,=\,\bigcup_{j\in\children(i)}\ul{\phi}_j(\fX_j)\quad.
\]
We sometimes drop the morphisms $\ul{\phi}_i$ and $\ul{\beta}_i$ from the notation and simply say that $(\fX_i)_{i\in I}$ is a treelike formal covering.
\end{defi}

\begin{lem}\label{simpadlem}
Let $X$ be a semi-affinoid $K$-space, and let $\fX$ be an affine flat formal $R$-model of ff type for $X$. Every treelike covering of $X$ together with any suitable rooted tree structure is induced by a treelike formal covering with root $\fX$.
\end{lem}
\begin{proof}
Let $(X_i)_{i\in I}$ be a treelike covering of $X$. We fix a suitable rooted tree structure on $I$, and we let $r$ denote the root of $I$. We set $\fX_r\mathrel{\mathop:}=\fX$; then $(X_i)_{i\in\children(r)}$ is a simple covering of $X$, and by Lemma \ref{simplecovadaptedtomodellem} it follows that there exist an admissible formal blowup $\ul{\beta}\colon\fX'\rightarrow\fX$ and an affine open covering $(\fX_i)_{i\in\children(r)}$ of $\fX'$ such that for each $i\in\children(r)$, the restriction of $\ul{\beta}$ to $\fX_i$ defines $X_i$ as a semi-affinoid subdomain of $X$. We may now apply the same reasoning to the children of $r$ and conclude by induction on the volume of $I$.
\end{proof}

\subsubsection{Retrocompact coverings}

\begin{defi}\index{covering!retrocompact}
A \emph{retrocompact covering} of a semi-affinoid $K$-space $X$ is a finite covering of $X$ by retrocompact semi-affinoid subdomains of $X$.
\end{defi}

We will prove that retrocompact coverings admit leaflike refinements.





\begin{lem}\label{buildtreelem}
Let $I$ be an almost linear rooted tree, and let us consider the following formal data parametrized by $I$:
\begin{packed_enum}
\item a family $(\fX_i)_{i\in I}$ of affine flat formal $R$-schemes of ff type,
\item for each inner vertex $i\in I$ an admissible formal blowup $\ul{\beta}_i\colon\fX_i'\rightarrow\fX_i$ and
\item for each inner vertex $i\in I$ and each child $j$ of $i$ an open immersion 
\[
\ul{\phi}_j\colon\fX_j\hookrightarrow\fX_i'
\]
\end{packed_enum}
such that
\[
\fX_r^\srig\,=\,\bigcup_{i\in\,\leaves(I)}\fX_i^\srig\quad,
\]
where $r$ denotes the root of $I$ and where we view $\fX_i^\srig$ as a subset of $\fX_r^\srig$ in the natural way. Then $I$ together with this data extends to a treelike formal covering $(\fX_i)_{i\in J}$ such that $(\fX_i^\srig)_{i\in\leaves(J)}$ is a refinement of $(\fX_i^\srig)_{i\in\leaves(I)}$.
\end{lem}
\begin{proof}
We may assume that $r$ is an inner vertex of $I$ and that $r$ has more than one child. We argue by induction on the volume $v(I)$ of $I$. By the surjectivity of the specialization map in the flat case,
\[
\fX_r'\,=\,\bigcup_{i\in\children(r)}\ul{\phi}_i(\fX_i)\quad,
\]
$\fX_r^\srig$ being covered by the $\fX_i^\srig$. Let $i$ be any child of $r$, and let us consider the linear subtree $\subtree(i)$ of $I$. For any $j\in\children(r)$ different from $i$, we glue a copy of $\subtree(j)$ to the root $i\in\subtree(i)$, and we provide it with the formal data obtained from the formal data attached to $\subtree(j)$ via $R$-flat pullback under $\ul{\beta}_i\circ\ul{\phi}_i$. We thereby replace $\subtree(i)$ by an almost linear rooted tree $I_i$ with formal data attached to it such that $v(I_i)=v(I)-1$, such that the leaves of $I_i$ generically refine the leaves of $I$ and such that $I_i$ with its formal data satisfies the assumptions of this Lemma. By our induction hypothesis, $I_i$ with its formal data can be extended to a treelike formal covering whose leaves generically form a refinement of the leaves of $I_i$. Since this holds for all $i\in\children(r)$, we conclude that $I$ with its formal data can be extended to a treelike formal covering whose leaves generically form a refinement of the leaves of $I$.
\end{proof}

\begin{prop}\label{retroprop}
Any retrocompact covering of a semi-affinoid $K$-space admits a leaflike refinement.
\end{prop}
\begin{proof}
Let $X$ be a semi-affinoid $K$-space, and let $(X_i)_{i\in I}$ be a retrocompact covering of $X$. For each $i\in I$, we choose a morphism $\ul{\phi}_i$ defining $X_i$ as a retrocompact semi-affinoid subdomain in $X$. We factor each $\ul{\phi}_i$ as a composition of an alternating sequence of open immersions and admissible formal blowups. After replacing $(X_i)_{i\in I}$ by a finite refinement, we may assume that the domains of the open immersions occurring in these factorizations are affine. By Corollary \ref{subdomdefgiventargetcor} and in consideration of our remark at the end of Section \ref{preimsemaffsec}, we may assume that the targets of the $\ul{\phi}_i$ coincide with a given affine flat formal $R$-model of ff type $\fX$ for $X$. Thus, after possibly adding trivial open immersions and trivial admissible formal blowups, the $\ul{\phi}_i$ and their factorizations yield a finite almost linear rooted tree with formal data attached to it as in Lemma \ref{buildtreelem}. By Lemma \ref{buildtreelem}, this rooted tree with its formal data can be extended to a treelike formal covering $(\fX_i)_{i\in J}$ whose leaves generically refine $(X_i)_{i\in I}$. Passing to semi-affinoid generic fibers, we thus obtain a treelike covering of $X$ whose associated leaflike covering refines $(X_i)_{i\in I}$.
\end{proof}

\subsubsection{An auxiliary G-topology}\label{auxgtopsec}

We now define an auxiliary G-topology $\sT_\textup{aux}$ on the category of semi-affinoid $K$-spaces, equipped with its its physical points functor, as follows: 

\begin{defi}
Let $X$ be a semi-affinoid $K$-space.
\begin{enumerate}
\item The $\sT_\textup{aux}$-admissible subsets of $|X|$ are the semi-affinoid subdomains.
\item A covering of a semi-affinoid subdomain in $X$ by semi-affinoid subdomains is called $\sT_\textup{aux}$-admissible if and only if it has a leaflike refinement.
\end{enumerate}
\end{defi}

By Proposition \ref{retroprop}, retrocompact coverings of semi-affinoid $K$-spaces are $\sT_\textup{aux}$-admissible, so one may consider $\sT_\textup{aux}$ as a semi-affinoid analog of the weak rigid G-topology on an affinoid $K$-space, cf.\ \cite{BGR} 9.1.4. 

Let us note some properties of $\sT_\textup{aux}$ which are immediate from the definition:
\begin{packed_enum}
\item Every $\sT_\textup{aux}$-admissible covering of a semi-affinoid subdomain of $X$ has a refinement consisting of \emph{retrocompact} semi-affinoid subdomains. 
\item A covering of a semi-affinoid subdomain of $X$ by semi-affinoid subdomains is $\sT_\textup{aux}$-admissible if and only if it has a $\sT_\textup{aux}$-admissible refinement. 
\end{packed_enum}

\begin{prop}\label{auxgrothtopprop}
$\sT_\textup{aux}$ defines a G-topology on $(\sAff_K,|\cdot|)$. 
\end{prop}
\begin{proof}
Let $X$ be a semi-affinoid $K$-space. To see that $\sT_\aux$ defines a G-topology on $X$, we must verify the axioms listed in \cite{BGR} 9.1.1/1:
\begin{packed_enum}
\item By Corollary \ref{semiaffintcor}, the set of semi-affinoid subdomains in $X$ is stable under the formation of \emph{intersections}.
\item If $U\subseteq X$ is a semi-affinoid subdomain, then the \emph{trivial covering} of $U$ is clearly simple and, hence, $\sT_\textup{aux}$-admissible.
\item Let $V\subseteq U$ be semi-affinoid subdomains in $X$, and let $(U_i)_{i\in I}$ be a $\sT_\textup{aux}$-admissible covering of $U$. Then $(U_i\cap V)_{i\in I}$ is a $\sT_\textup{aux}$-admissible covering of $V$. Indeed, we may replace $(U_i)_{i\in I}$ be a refinement and thereby assume that $(U_i)_{i\in I}$ is leaflike; then the claim follows from Lemma \ref{leaflikecovintlem}. 
\item Let $U\subseteq X$ be a semi-affinoid subdomain, let $(U_i)_{i\in I}$ be a $\sT_\textup{aux}$-admissible covering of $U$, and for each $i\in I$ let $(V_{ij})_{j\in J_i}$ be a $\sT_\textup{aux}$-admissible covering of $U_i$. Then $(V_{ij})_{i\in I,j\in J_i}$ is a $\sT_\textup{aux}$-admissible covering of $U$. Indeed, by part ($iii$) of this proof we may pass to refinements and thereby assume that the given coverings are leaflike. Now the claim follows from Lemma \ref{leaflikestablelem}. 
\end{packed_enum}
We finally have to show that the map $|\phi|\colon|Y|\rightarrow|X|$ induced by a morphism of semi-affinoid $K$-spaces $\phi\colon Y\rightarrow X$ is continuous for $\sT_\textup{aux}$. 
\begin{packed_enum}
\item If $U$ is a semi-affinoid subdomain in $X$, then Corollary \ref{semaffpreimcor} shows that $\phi^{-1}(U)$ is a semi-affinoid subdomain in $Y$. 
\item Let $U$ be a semi-affinoid subdomain in $X$, and let $(U_i)_{i\in I}$ be a $\sT_\textup{aux}$-ad\-mis\-si\-ble covering of $U$. To show that the induced covering $(\phi^{-1}(U_i))_{i\in I}$ of $\phi^{-1}(U)$ is $\sT_\textup{aux}$-admissible, we may replace $(U_i)_{i\in I}$ be a refinement and thereby assume that $(U_i)_{i\in I}$ is leaflike. The desired statement then follows from Lemma \ref{leaflikecovintlem}. 
\end{packed_enum}
\end{proof}

Let us note that $\sT_\textup{aux}$ has the additional property that for each semi-affinoid $K$-space $X$, the subsets $|X|$ and $\emptyset$ are admissible.

\subsubsection{The uniformly rigid G-topology}

Let $(\fC,|\cdot|)$ be a pair as in Definition \ref{gtopdefi}. The set of G-topologies on $(\fC,|\cdot|)$ is partially ordered by the relation of being \emph{slightly finer}; cf.\ \cite{BGR} 9.1.2/1 and the subsequent comments. By \cite{BGR} 9.1.2/2, for any G-topology $\sT$ on $(\fC,|\cdot|)$ there exists a unique finest G-topology $\sT'$ on $(\fC,|\cdot|)$ that is slightly finer than $\sT$.

\begin{defi}Let $\sT_\srig$ denote the unique finest G-topology on $(\sAff_K,|\cdot|)$ that is slightly finer than $\sT_\textup{aux}$. It is called the \emph{uniformly rigid} G-topology.
\end{defi}

The following proposition is contained in the statement of \cite{BGR} 9.1.2/2:

\begin{prop}\label{saturatedprop}
Let $X$ be a semi-affinoid $K$-space. The uniformly rigid G-topology on $X$ satisfies conditions \textup{(G$_0$) -- (G$_2$)} in \cite{BGR} 9.1.2; that is:
\begin{enumerate}
\item The subsets $\emptyset$ and $X$ of $X$ are $\sT_\srig$-admissible.
\item Let $U\subseteq X$ be $\sT_\srig$-admissible, let $(U_i)_{i\in I}$ be a $\sT_\srig$-admissible covering of $U$, and let $V\subseteq U$ be any subset. If $V\cap U_i$ is $\sT_\srig$-admissible for each $i$, then $V$ is $\sT_\srig$-admissible.
\item Let $U\subseteq X$ be $\sT_\srig$-admissible, and let $(U_i)_{i\in I}$ be a covering of $U$ by $\sT_\srig$-admissible subsets. If $(U_i)_{i\in I}$ admits an $\sT_\srig$-admissible refinement, then $(U_i)_{i\in I}$ is $\sT_\srig$-admissible.
\end{enumerate}
\end{prop}

In other words, admissibility of subsets can be checked locally with respect to admissible coverings, and admissibility of a covering can be checked after refinement. A G-topology satisfying these properties is called \emph{saturated}.

Let us introduce the following terminology: A \emph{semi-affinoid subdomain covering}\index{covering!semi-affinoid subdomain} of a subset $U$ of a semi-affinoid $K$-space $X$ is a not necessarily finite covering of $U$ by semi-affinoid subdomains in $X$. As a corollary of [BGR] 9.1.2/3, we obtain the following explicit description of the uniformly rigid G-topology on a semi-affinoid $K$-space:

\begin{prop}\label{explicitprop}
Let $X$ be a semi-affinoid $K$-space.
\begin{enumerate}
\item A subset $U\subseteq X$ is $\sT_\srig$-admissible if and only if it admits a semi-affinoid subdomain covering $(U_i)_{i\in I}$ such that for any morphism $\phi\colon Y\rightarrow X$ of semi-affinoid $K$-spaces with $\phi(Y)\subseteq U$, the induced semi-affinoid subdomain covering of $Y$ has a leaflike refinement.
\item A covering $(U_i)_{i\in I}$ of a $\sT_\srig$-admissible subset $U$ in $X$ by $\sT_\srig$-admissible subsets is $\sT_\srig$-admissible if and only if for any morphism $\phi\colon Y\rightarrow X$ of semi-affinoid $K$-spaces with $\phi(Y)\subseteq U$, the induced covering of $Y$ has a leaflike refinement.
\end{enumerate}
\end{prop}

Let us recall from \cite{BGR}  p. 337 that a G-topological space is called \emph{quasi-compact} if every admissible covering admits a finite admissible refinement. It follows from Proposition \ref{explicitprop} ($ii$) that semi-affinoid subdomains of semi-affinoid $K$-spaces are quasi-compact.

If $X$ is a quasi-compact G-topological space satisfying condition \textup{(G$_2$)} of \cite{BGR} 9.1.2 and if $(X_i)_{i\in I}$ is an admissible covering of $X$, then $(X_i)_{i\in I}$ has a finite admissible \emph{subcovering}. 
Indeed, let $(Y_j)_{j\in J}$, $\phi\colon J\rightarrow I$ be a finite admissible refinement of $(X_i)_{i\in I}$; then the finite subcovering $(X_{\phi(j)})_{j\in J}$ of $(X_i)_{i\in I}$ has an admissible refinement and, hence, is admissible, in virtue of the \textup{(G$_2$)} condition. By Proposition \ref{saturatedprop}, semi-affinoid $K$-spaces with their uniformly rigid G-topology satisfy \textup{(G$_2$)}; hence every admissible covering of a semi-affinoid $K$-space has a finite admissible subcovering.

For example, it follows that any covering of the semi-affinoid open unit disc $\sSp (R[[S]]\otimes_RK)$ by nested closed semi-affinoid subdiscs is 
\emph{not} $\sT_\srig$-admissible since it has no finite subcovering. More generally, the uniformly rigid open unit disc admits no $\sT_\srig$-admissible covering by semi-affinoid subdomains whose rings of functions are \emph{affinoid}. Indeed, there would exist a finite subcovering, and by the maximum principle for affinoid $K$-algebras the coordinate function on the semi-affinoid open unit disc would attain its maximum, which is not the case.

Let us recall from \cite{BGR} p. 338  that a collection of subsets of a G-topological space is called a \emph{basis} for the G-topology if every admissible open subset admits an admissible covering by elements of this collection. By Proposition \ref{explicitprop} ($i$) and ($ii$), the semi-affinoid subdomains of a semi-affinoid $K$-space form a basis for the uniformly rigid G-topology.

The semi-affinoid open unit disc can be identified with a semi-affinoid subdomain of the semi-affinoid closed unit disc $\sSp K\langle T\rangle$. As we have seen, it admits no $\sT_\srig$-admissible covering by semi-affinoid subdomains whose rings of functions are affinoid. By Lemma \ref{veryspecialinafflem}, we conclude that the \emph{retrocompact} semi-affinoid subdomains of the closed unit disc do \emph{not} form a basis for the uniformly rigid G-topology. In particular, the uniformly rigid G-topology on the semi-affinoid closed unit disc does \emph{not} agree with the rigid-analytic G-topology that is obtained via Berthelot's construction.

Even though the retrocompact semi-affinoid subdomains in a semi-affinoid $K$-space $X$ do not need to form a basis for the uniformly rigid G-topology on $X$, they are general enough to detect all local properties of $X$: 

\begin{lem}\label{retrolocallem}
Let $X$ be a semi-affinoid $K$-space. For every point $x\in X$ and any $\sT_\srig$-admissible subset $U\subseteq X$ containing $x$, there exists a retrocompact semi-affinoid subdomain $V\subseteq X$ whose ring of functions is $K$-affinoid such that $x\in V\subseteq U$.
\end{lem}
\begin{proof}
We may assume that $U$ is a semi-affinoid subdomain in $X$. Let us consider a morphism $\ul{\phi}$ of affine flat formal $R$-schemes defining $U$ as a semi-affinoid subdomain in $X$. Let $\fX$ denote the domain of $\ul{\phi}$. We choose a finite factorization of $\ul{\phi}$ into simple morphisms. We replace all completion morphisms in our factorization by suitable formal dilatations, as follows: Let $\ul{\psi}\colon\fY|_V\rightarrow\fY$ be the first completion morphism occurring in the chosen factorization of $\ul{\phi}$, and let $\sI$ be the biggest ideal of definition of $\fY|_V$. Since $x$ lies in the rigid generic fiber of $\fY|_V$, our discussion in Section \ref{berthsec} shows that there exists an $n\in\N$ such that $x$ is contained in the generic fiber of the formal dilatation $(\fY|_V)_{(n)}$ of $\fY|_V$ in the coherent ideal generated by $\pi$ and $\sI^{n+1}$. By Lemma \ref{berthconstrbasiclem}, $(\fY|_V)_{(n)}$ is of tf type over $R$. On the other hand, let $V$ be the reduced subscheme of the smallest subscheme of definition of $\fY$ such that $\ul{\phi}$ is the completion of $\fY$ along $V$, and let $\sJ\subseteq\O_\fY$ denote the open coherent sheaf corresponding to $V$; then $\sI=\sJ\O_{\fY|_V}$, and we readily see that $(\fY|_V)_{(n)}$ satisfies the universal property of the formal dilatation of $\fY$ in the coherent ideal generated by $\pi$ and $\J^{n+1}$. We replace $\ul{\psi}$ by the formal dilatation $(\fY|_V)_{(n)}\rightarrow\fY$, and we replace the remaining part of the factorization of $\ul{\phi}$ by its strict transform under the formal dilatation $(\fY|_V)_{(n)}\rightarrow\fY|_V$. After iterating this procedure, we obtain a morphism $\ul{\phi}'$ defining a retrocompact semi-affinoid subdomain of $X$ that contains $x$, together with a factorization of $\ul{\phi}'$ into retrocompact simple morphisms. We may clearly assume from the outset that the first step in our factorization of $\ul{\phi}$ is the identity on $\fX$, considered as the completion of $\fX$ along its smallest subscheme of definition; then all formal $R$-schemes  occurring in the resulting factorization of $\ul{\phi}'$ are of tf type over $R$, except possibly for the domain $\fX$ of $\ul{\phi}'$. The statement is now clear.
\end{proof}

Proposition \ref{explicitprop} implies that the uniformly rigid G-topology localizes well, which will make it possible to globalize the theory:

\begin{cor}\label{srigpropcor}
Let $X$ be a semi-affinoid $K$-space. For any semi-affinoid subdomain $U$ of $X$, the uniformly rigid G-topology on $X$ restricts to the uniformly rigid G-topology on $U$.
\end{cor}
\begin{proof}
By Corollary \ref{semaffsubdomgoodnotioncor}, the semi-affinoid subdomains in $U$ are the semi-affinoid subdomains in $X$ contained in $U$, and by Proposition \ref{subdomunivprop} the semi-affinoid morphisms to $X$ with image in $U$ correspond to the semi-affinoid morphisms to $U$. Hence, the statement follows from Proposition \ref{explicitprop} ($i$) and ($ii$).
\end{proof}

Moreover, Proposition \ref{explicitprop} shows that finite unions of retrocompact semi-affinoid subdomains are $\sT_\srig$-admissible, even when they are not semi-affinoid subdomains:

\begin{cor}\label{finretroadmcor}
Let $X$ be a semi-affinoid $K$-space, and let $U\subseteq X$ be a finite union of retrocompact semi-affinoid subdomains in $X$. Then $U$ is a $\sT_\srig$-admissible subset of $X$, and any finite covering of $U$ by retrocompact semi-affinoid subdomains in $X$ is $\sT_\srig$-admissible.
\end{cor}
\begin{proof}
Let $(U_i)_{i\in I}$ be a finite family of retrocompact semi-affinoid subdomains of $X$ such that $U$ is the union of the $U_i$. Let $Y$ be any semi-affinoid $K$-space, and let $\phi\colon Y\rightarrow X$ be any semi-affinoid morphism whose image is lies in $U$. Then $(\phi^{-1}(U_i))_{i\in I}$ is a retrocompact covering of $Y$; by Propostion \ref{retroprop}, it admits a leaflike refinement. By Proposition \ref{explicitprop} ($i$), we conclude that $U$ is a $\sT_\srig$-admissible subset of $X$, and by Proposition \ref{explicitprop} ($ii$) we see that the covering $(U_i)_{i\in I}$ of $U$ is $\sT_\srig$-admissible.
\end{proof}

In particular, if $A$ is an \emph{affinoid} $K$-algebra, then the Theorem of Gerritzen and Grauert (\cite{BGR} 7.3.5/1) shows that every affinoid subdomain in $\Sp A$ corresponds to a $\sT_\urig$-admissible subset of $\sSp A$.

\subsubsection{Comparison with the rigid G-topology}

Let us recall from \cite{BGR} p. 342 that an \emph{affinoid} covering of an affinoid $K$-space is a \emph{finite} covering by affinoid subdomains. We recall the following known fact from rigid geometry:

\begin{lem}\label{easyrigidlem}
Let $X$ be an affinoid $K$-space, and let $U\subseteq X$ be a subset admitting a covering $(U_i)_{i\in I}$ by admissible open subsets $U_i\subseteq X$ such that for any affinoid $K$-space $Y$ and any morphism $\phi\colon Y\rightarrow X$ with image in $U$, the induced covering $(\phi^{-1}(U_i))_{i\in I}$ of $Y$ has an affinoid refinement. Then $U\subseteq X$ is admissible.
\end{lem}
\begin{proof}
For each $i\in I$, we choose an affinoid covering $(U_{ij})_{j\in J_i}$ of $U_i$ such that for any affinoid $K$-space $Z$ and any morphism $\psi\colon Z\rightarrow X$ with image in $U_i$, the induced covering $(\psi^{-1}(U_{ij}))_{j\in J_i}$ of $Z$ has an affinoid refinement. This is possible since the $U_i$ are admissible. Then $(U_{ij})_{i\in I,j\in J_i}$ is a covering of $U$. Let $\phi\colon Y\rightarrow X$ be any morphism of affinoid $K$-spaces with image in $U$;  it suffices to show that the induced covering $(\phi^{-1}(U_{ij}))_{i\in I,j\in J_i}$ of $Y$ has an affinoid refinement. Let $(Z_k)_{k\in K}$ be an affinoid refinement of $(\phi^{-1}(U_i)_{i\in I}$; it exists by assumption. Then for each $k$, the restriction of $\phi$ to $Z_k$ has image in some $U_{i(k)}$, and there exists an affinoid covering of $Z_k$ refining the covering of $Z_k$ induced from $(U_{i(k),j})_{j\in J_{i(k)}}$. Collecting these affinoid refinements, we obtain the desired affinoid refinement of $(\phi^{-1}(U_{ij}))_{i\in I,j\in J_i}$.
\end{proof}


As we have seen in Section \ref{assocrigspacesubsec}, Berthelot's construction functorially associates a rigid $K$-space $X^\r$ to each semi-affionid $K$-space $X$, and it provides a functorial bijection between the sets of points of $X$ and $X^\r$. We thus obtain a rigid-analytic G-topology $\sT_\rig$ on the category of semi-affinoid $K$-spaces.

\begin{prop}\label{srigrigcompprop}
$\sT_\rig$ is finer than $\sT_\srig$. 
\end{prop}
\begin{proof}
Let $X$ be a semi-affinoid $K$-space. If $U\subseteq X$ is a semi-affinoid subdomain, then $U$ is $\sT_\rig$-admissible. Indeed, this is clear from the fact that open immersions, admissible formal blowups and completion morphisms of flat formal $R$-schemes induce open immersions of rigid-analytic generic fibers. Let us now consider a general $\sT_\srig$-admissible subset $U\subseteq X$. We may check $\sT_\rig$-admissibility of $U$ locally with respect to a $\sT_\rig$-admissible covering of $X$; hence we may replace $X$ by the semi-affinoid generic fiber of a suitable formal dilatation, cf. our description of Berthelot's construction in Section \ref{berthsec}, and thereby assume that $X$ is the semi-affinoid space associated to an \emph{affinoid} $K$-algebra. Let $(U_i)_{i\in I}$ be a semi-affinoid subdomain covering of $U$ such that condition ($i$) of Proposition \ref{explicitprop} is satisfied. Let $Y$ be an affinoid $K$-space, and let $\phi\colon Y\rightarrow X^\r$ be any morphism of rigid spaces that factorizes through $U$. Since $X$ is the uniformly rigid space associated to an affinoid $K$-algebra, $\phi$ can also be viewed as a morphism of semi-affinoid $K$-spaces. Hence $(\phi^{-1}(U_i))_{i\in I}$ has a leaflike refinement. By Lemma \ref{veryspecialinafflem}, this refinement is affinoid. It now follows from Lemma \ref{easyrigidlem} that $U\subseteq X$ is $\sT_\rig$-admissible.

Let us now consider a $\sT_\srig$-admissible subset $U$ in $X$ together with a $\sT_\srig$-admissible covering $(U_i)_{i\in I}$ of $U$ by $\sT_\srig$-admissible subsets $U_i$. We have seen that $U$ and the $U_i$ are $\sT_\rig$-admissible; we claim that the covering $(U_i)_{i\in I}$ is $\sT_\rig$-admissible as well. Let us choose a $\sT_\rig$-admissible covering of $X$ by affinoid subspaces $X_j$, $j\in\N$, obtained via suitable formal dilatations according to our discussion of Berthelot's construction in Section \ref{berthsec}. Let us assume that for each $j\in\N$, the induced covering of $(U_i\cap X_j)_{i\in I}$ of $U\cap X_j$ is $\sT_\rig$-admissible. Then $(U_i\cap X_j)_{i\in I,j\in \N}$ is $\sT_\rig$-admissible by transitivity of admissibility for coverings; hence $(U_i)_{i\in I}$ has a $\sT_\rig$-admissible refinement and, hence, is $\sT_\rig$-admissible. We have thus reduced to the case where $X$ is the semi-affinoid space associated to an affinoid $K$-algebra. Let $Y$ be an affinoid $K$-space, and let $\phi\colon Y\rightarrow X^\r$ be a morphism of affinoid $K$-spaces. We may also view $\phi$ as a morphism of semi-affinoid $K$-spaces. Since $(U_i)_{i\in I}$ is $\sT_\srig$-admissible, we see by Corollary \ref{explicitprop} ($ii$) that $(\phi^{-1}(U_i))_{i\in I}$ has a leaflike and, hence, affinoid refinement. It follows that $(U_i)_{i\in I}$ is $\sT_\rig$-admissible.
\end{proof}

The G-topologies $\sT_\rig$ and $\sT_\urig$ differ mainly on the level of admissible coverings, not on the level of admissible subsets:

\begin{lem}\label{rigurigadmsubsetcomparelem}
Let $X$ be a semi-affinoid $K$-space, and let $U\subseteq X$ be a subset which is $\sT_\rig$-admissible and quasi-compact with respect to $\sT_\rig$. Then $U$ is $\sT_\urig$-admissible.
\end{lem}
\begin{proof}
In view of Berthelot's construction which we recalled in Section \ref{berthsec}, $U$ is contained is a semi-affinoid subdomain $W$ of $X$ whose $K$-algebra of global functions is affinoid. By the Gerritzen-Grauert Theorem \cite{BGR} 7.3.5/1, $U$ is a finite union of retrocompact semi-affinoid subdomains of $W$. By Corollary \ref{finretroadmcor}, $U$ is $\sT_\urig$-admissible in $W$ and, hence, in $X$.
\end{proof}

\subsubsection{Comparison with the Zariski topology}


Let $X$ be a semi-affinoid $K$-space, and let $A$ denote its ring of functions; then $|X|=\Max A$, and hence $|X|$ carries a Zariski topology $\sT_\Zar$. Clearly $\sT_\Zar$ defines a G-topology on $(\sAff_K,|\cdot|)$. We emphasize that $\sT_\Zar$ is defined by functions in $A$, not by global functions on the rigid space $X^\r$ associated to $X$.

\begin{prop}\label{zarsrigprop}
$\sT_\srig$ is finer than $\sT_\Zar$.
\end{prop}
\begin{proof}
Let $U\subseteq X$ be a Zariski-open subset, and let $f_1,\ldots,f_n\in A$ be semi-affinoid functions such that $U$ is the union of the Zariski-open subsets $D(f_i)=\{x\in\Max A\,;\,f_i(x)\neq 0\}$. Let $Y$ be a semi-affinoid $K$-space, and let $\phi\colon Y\rightarrow X$ be a semi-affinoid morphism whose image in contained in $U$.
For each $i$, the preimage $\phi^{-1}(D(f_i))$ is the set of points $y\in Y$ where $\phi^*f_i\neq 0$. Since $Y$ is covered by the $\phi^{-1}(D(f_i))$, the $\phi^*f_i$ generate the unit ideal in $B$. That is, there exist elements $b_1,\ldots,b_n$ in $B$ such that $b_1\phi^*f_1+\ldots +b_n\phi^*f_n=1$. Let us set $\gamma\mathrel{\mathop:}=(\max_i|b_i|_\sup)^{-1}$; this number is well-defined since the $b_i$ are bounded functions on $Y$ without a common zero. By the strict triangle inequality, $\max_i|\phi^*f_i(y)|\geq\gamma$ for all $y\in Y$. For each $i$, let $Y_i\subseteq Y$ denote the set of points $y\in Y$ where $|\phi^*f_i(y)|\geq\gamma$; then $(Y_i)_{1\leq i\leq n}$ is a retrocompact covering of $Y$ refining $(\phi^{-1}(D(f_i)))_{1\leq i\leq n}$. By Proposition \ref{retroprop}, retrocompact coverings are $\sT_\srig$-admissible; hence $U\subseteq X$ is $\sT_\srig$-admissible. If $(U_j)_{j\in J}$ is a Zariski-covering of $U$, we may pass to a refinement and assume that the $U_j=D(g_j)\subseteq X$ for some semi-affinoid function $g_j$ on $X$; we can then argue along the same lines to prove that $(U_j)_{j\in J}$ is a $\sT_\srig$-admissible covering of $U$.
\end{proof}

Our proof shows that if $f_1,\ldots,f_n$ are semi-affinoid functions on $X$ and if 
\[
U\,=\,\bigcup_{i=1}^nD(f_i)
\]
is the associated Zariski-open subset of $X$, then setting 
\[
U_{\geq\varepsilon}\,=\,\bigcup_{i=1}^n\,\{x\in X\,;\,|f_i(x)|\geq\varepsilon\}
\]
for $\varepsilon\in\sqrt{|K^*|}$, the resulting covering $(U_{\geq\varepsilon})_{\varepsilon}$ of $U$ by finite unions of retrocompact semi-affinoid subdomains of $X$ is $\sT_\srig$-admissible.

			\subsection{The semi-affinoid Acyclicity Theorem}\label{semafffuncsec}
				Let $X$ be a semi-affinoid $K$-space. In Definition \ref{semaffpresheafdefi}, we have introduced a presheaf $\O_X$ on the category of semi-affinoid pre-subdomains of $X$ with values in the category of semi-affinoid $K$-algebras. In this section we show that $\O_X$ is a \emph{sheaf} for $\sT_\textup{aux}$ and, hence, extends uniquely to a sheaf on $\sT_\srig$. More generally, we show that every $\O_X$-module associated to a finite module over the ring of global functions on $X$ is \emph{acyclic} for any $\sT_\srig$-admissible covering of $X$.

Let $\sF$ be a presheaf of $\O_X$-modules on the category of semi-affinoid subdomains of $X$. Let us recall from \cite{BGR} p. 324 that a semi-affinoid subdomain covering $(X_i)_{i\in I}$ of $X$ is called $\sF$-\emph{acyclic} if the associated augmented \v{C}ech complex is acyclic. The covering $(X_i)_{i\in I}$ is called \emph{universally} $\sF$-acyclic if $(X_i\cap U)_{i\in I}$ is $\sF|_U$-acyclic for any semi-affinoid subdomain $U\subseteq X$.

We now establish the semi-affinoid version of Tate's acyclicity theorem.
Adopting methods from \cite{LtkeFRG}, we will derive it from results in formal geometry.

\begin{theorem}\label{acyclicitytheorem}
Let $X$ be a semi-affinoid $K$-space, and let $(X_i)_{i\in I}$ be a $\sT_\aux$-admissible covering of $X$. Then $\O_X$ is $(X_i)_{i\in I}$-acyclic.
\end{theorem}
\begin{proof}
Let us first consider the case where the covering $(X_i)_{i\in I}$ is simple. We choose a formal presentation $(\fX,\ul{\beta}\colon\fX'\rightarrow\fX,(\fX_i)_{i\in I})$ of $(X_i)_{i\in I}$. By Proposition \ref{blowupisoonsemiafffuncprop}, the natural homomorphism $\O_{\fX,K}\rightarrow(\ul{\beta}_*\O_{\fX'})_K$ is an isomorphism, so it induces an isomorphism of rings of global sections. It follows that $\ul{\beta}$ induces a natural identification of augmented \v{C}ech complexes
\[
C^\bullet_{\textup{aug}}((X_i)_{i\in I},\O_X)\cong C^\bullet_{\textup{aug}}((\fX_i)_{i\in I},\O_{\fX',K})\quad.
\]
We have to show that the complex on the right hand side is acyclic. Since $\O_{\fX',K}$ is a sheaf on $\fX'$, it suffices to show that
\[
\check{H}^q((\fX_i)_{i\in I},\O_{\fX',K})\,=\,0
\]
for all $q\geq 1$. Since $I$ is finite, we have an identification
\[
\check{H}^q((\fX_i)_{i\in I},\O_{\fX',K})\,=\,\check{H}^q((\fX_i)_{i\in I},\O_{\fX'})\otimes_RK\quad.
\]
By the Comparison Theorem \cite{EGAIII1} 4.1.5, 4.1.7 and by the Vanishing Theorem \cite{EGAIII1} 1.3.1, the higher cohomology groups of a coherent sheaf on an affine noetherian formal scheme vanish. Since the $\fX_i$ are affine, Leray's theorem implies that
\[
\check{H}^q((\fX_i)_{i\in I},\O_{\fX'})\,=\,H^q(\fX',\O_{\fX'})\quad.
\]
By \cite{EGAIII1} 1.4.11, $H^q(\fX',\O_{\fX'})=\Gamma(\fX,R^q\ul{\beta}_*\O_{\fX'})$, and by Proposition \ref{blowupisoonsemiafffuncprop} this module is $\pi$-torsion. We have thus finished the proof in the case where $(X_i)_{i\in I}$ is a simple covering.

Any $\sT_\aux$-admissible covering of $X$ has a leaflike refinement; by \cite{BGR} 8.1.4/3, it thus suffices to show that the leaflike coverings of $X$ are universally $\O_X$-acyclic. By Lemma \ref{leaflikecovintlem}, the restriction of a leaflike covering to a semi-affinoid subdomain is a leaflike covering; hence it suffices to show that any leaflike covering of $X$ is $\O_X$-acyclic. We may thus assume that $(X_i)_{i\in I}$ is leaflike. 

Let $(X_j)_{j\in J}$ be a treelike covering of $X$ extending $(X_i)_{i\in I}$, and let us choose a suitable rooted tree structure on $J$ such that $I\subseteq J$ is identified with the set of leaves of $J$. We argue by induction on the volume of $J$. If $J$ has only one vertex, the covering $(X_i)_{i\in I}$ is trivial and, hence, $\O_X$-acyclic. Let us assume that $J$ has more than one vertex. Let $\iota\in I$ be a leaf of $J$ such that the length $l(\iota)$ of the path from $\iota$ to the root is maximal in $\{l(i)\,;\,i\in I\}$. Let $\iota'\mathrel{\mathop:}=\parent(\iota)$ denote the parent of $\iota$. By maximality of $l(\iota)$, all siblings $i\in\children(\iota')$ of $\iota$ are leaves of $J$. Let $J'\mathrel{\mathop:}=J\setminus\children(\iota')$ be the rooted subtree of $J$ that is obtained by removing the siblings of $\iota$ (including $\iota$ itself). Then
\begin{packed_enum}
\item the set of leaves of $J'$ is $I'\mathrel{\mathop:}=(I\setminus\children(\iota'))\cup\{\iota'\}$,
\item $(X_j)_{j\in J'}$ is a treelike covering of $X$, and
\item $v(J')<v(J)$.
\end{packed_enum}
By our induction hypothesis, the covering $(X_i)_{i\in I'}$ is $\O_X$-acyclic. Since $(X_i)_{i\in I}$ is a refinement of $(X_i)_{i\in I'}$, we know from \cite{BGR} 8.1.4/3 that it suffices to prove that for any $r\geq 0$ and any tuple $(i_0,\ldots,i_r)\in(I')^{r+1}$, the covering $(X_i\cap X_{i_0\cdots i_r})_{i\in I}$ of $X_{i_0\cdots i_r}$ is $\O_X$-acyclic, where $X_{i_0\cdots i_r}$ denotes the intersection $X_{i_0}\cap\ldots\cap X_{i_r}$. Let us assume that there exists some $0\leq s\leq r$ such that $i_s\neq\iota'$. Then $i_s\in I$. Since $X_{i_0\cdots i_r}\subseteq X_{i_s}$, we see that the trivial covering of $X_{i_0\cdots i_r}$ refines $(X_i\cap X_{i_0\cdots i_r})_{i\in I}$. Since trivial coverings restrict to trivial coverings and since trivial coverings are acyclic, we deduce from \cite{BGR} 8.1.4/3 that $(X_i\cap X_{i_0\cdots i_r})_{i\in I}$ is acyclic. It remains to consider the case where all $i_s$, $0\leq s\leq r$, coincide with $\iota'$. That is, it remains to see that the covering $(X_i\cap X_{\iota'})_{i\in I}$ of $X_{\iota'}$ is $\O_X$-acyclic. It admits the simple covering $(X_i)_{i\in\children(\iota')}$ as a refinement. Since simple coverings restrict to simple coverings and since simple coverings are $\O_X$-acyclic by what we have shown so far, we conclude by \cite{BGR} 8.1.4/3 that $(X_i\cap X_{\iota'})_{i\in I}$ is $\O_X$-acyclic, as desired.
\end{proof}

In particular, $\O_X$ is a \emph{sheaf} for $\sT_\aux$. By \cite{BGR} 9.2.3/1, it extends uniquely to a sheaf for $\sT_\urig$ which we again denote by $\O_X$. Now that we dispose of the sheaf $\O_X$, we can easily discuss a fundamental example of a \emph{non-admissible} finite covering of a semi-affinoid $K$-space by semi-affinoid subdomains:

\begin{example}\label{nonadmdisccovex}
The natural covering of the semi-affinoid closed unit disc $\sSp K\langle T\rangle$ by the semi-affinoid open unit disc $\sSp (R[[T]]\otimes_RK)$ and the semi-affinoid unit circle $\sSp K\langle T,T^{-1}\rangle$ is not $\sT_\srig$-admissible and, hence, not $\sT_\textup{aux}$-admissible.
\end{example}
\begin{proof}
The two covering sets are nonempty and disjoint, while the ring of functions $K\langle T\rangle$ on the closed semi-affinoid unit disc is a domain. The presheaf of semi-affinoid functions being a sheaf for $\sT_\srig$, the considered covering cannot be $\sT_\srig$-admissible.
\end{proof}

If $X$ is a semi-affinoid $K$-space with ring of global functions $A$ and if $M$ is a finite $A$-module, the presheaf  $\O_X$-module $M\otimes\O_X$ sending a semi-affinoid subdomain $U$ in $X$ to $M\otimes_A\O_X(U)$ is called the $\O_X$-module \emph{associated} to $M$. A presheaf $\O_X$-module $\sF$ is called \emph{associated}\index{module!associated} if it is isomorphic to $M\otimes\O_X$ for some finite $A$-module $M$. We sometimes abbreviate $\tilde{M}\mathrel{\mathop:}=M\otimes\O_X$.

\begin{cor}\label{moduleacythmcor}
Let $X$ be a semi-affinoid $K$-space, and let $\sF$ be an associated $\O_X$-module. Then every $\sT_\aux$-admissible covering $(X_i)_{i\in I}$ of $X$ is $\sF$-acyclic.
\end{cor}
\begin{proof}
By \cite{BGR} 8.1.4/3, we may assume that $I$ is finite. Let $A$ denote the ring of functions on $X$, and let $M$ be a finite $A$-module such that $\sF$ is associated to $M$. If $M$ is free over $A$, the augmented \v{C}ech complex $C^\bullet((X_i)_{i\in I},\sF)$ is a finite direct sum of the complexes $C^\bullet((X_i)_{i\in I},\O_X)$, so in this case the statement is a direct consequence of Theorem \ref{acyclicitytheorem}. In the general case, let us choose a short exact sequence
\[
0\rightarrow M'\rightarrow F\rightarrow M\rightarrow 0
\]
where $F$ is a finite free $A$-module. We obtain a sequence of augmented \v{C}ech complexes
\[
0\rightarrow C_\textup{aug}^\bullet((X_i)_{i\in I},\tilde{M}')\rightarrow C_\textup{aug}^\bullet((X_i)_{i\in I},\tilde{F})\rightarrow C_\textup{aug}^\bullet((X_i)_{i\in I},\tilde{M})\rightarrow 0
\]
which is exact since the rings of functions on semi-affinoid subdomains in $X$ are flat $A$-algebras. By what we have shown so far, the associated long exact cohomology sequence contains isomorphisms 
\[
H^{q-1}_\textup{aug}((X_i)_{i\in I},\tilde{M})\cong H^q_\textup{aug}((X_i)_{i\in I},\tilde{M}')\quad.
\]
Since $I$ is finite, there exists an integer $n$ such that 
\[
H^q_\textup{aug}((X_i)_{i\in I},\tilde{N})=H^q((X_i)_{i\in I},\tilde{N})=0
\]
for all $q\geq n$ and all $A$-modules $N$, cf.\ [BGR] 8.1.3/2. By induction, we see that 
\[
H^q_\textup{aug}((X_i)_{i\in I},\tilde{M})=0
\]
for all $q\geq 0$.
\end{proof}

In particular, $M\otimes\O_X$ is a $\sT_\aux$-sheaf. By \cite{BGR} 9.2.3/1, $M\otimes\O_X$ extends uniquely to a $\sT_\srig$-sheaf on $X$ which we again denote by $M\otimes\O_X$. The $\sT_\srig$-sheaf $\O_X$ itself will be called the \emph{sheaf of semi-affinoid functions}\index{sheaf!of semi-affinoid functions} on $X$.

If $U\subseteq X$ is a semi-affinoid \emph{pre}-subdomain that is $\sT_\srig$-admissible, then $\O_X(U)=\O_U(U)$. Indeed, $U$ admits a $\sT_\srig|_X$-admissible covering by semi-affinoid subdomains in $X$; since morphisms of semi-affinoid spaces are continuous for $\sT_\srig$, this covering is also $\sT_\srig|_U$-admissible, so the statement follows from the fact that both $\O_X$ and $\O_U$ are $\sT_\srig$-sheaves. However, it is not clear whether $\sT_\srig|_X$ restricts to $\sT_\srig|_U$. For example, we do not know whether a semi-affinoid subdomain of $U$ is $\sT_\srig|_X$-admissible. Of course, this does not affect our theory since we do not deal with general semi-affinoid \emph{pre}-subdomains.

The category of abelian sheaves on $(X,\sT_\srig|_X)$ has enough injective objects, so the functor $\Gamma(X,\cdot)$ from the category of abelian sheaves on $X$ to the category of abelian groups has a right derived functor $H^\bullet(X,\cdot)$. In virtue of the Acyclicity Theorem and its Corollary \ref{moduleacythmcor}, this right derived functor can be calculated in terms of \v{C}ech cohomology:

\begin{cor}
Let $X$ be a semi-affinoid $K$-space, and let $\sF$ be an associated $\O_X$-module. Then the natural homomorphism
\[
\check{H}^q(U,\sF)\rightarrow H^q(U,\sF)
\]
is an isomorphism for all $\sT_\srig$-admissible subsets $U\subseteq X$. In particular, 
\[
H^q(U,\sF)\,=\,0
\]
for all $q>0$ and all semi-affinoid subdomains $U\subseteq X$.
\end{cor}
\begin{proof}
Let us consider the system $S$ of semi-affinoid subdomains in $X$. Then the following holds:
\begin{packed_enum}
\item $S$ is stable under the formation of intersections,
\item every $\sT_\srig$-admissible covering $(U_i)_{i\in I}$ of a $\sT_\srig$-admissible subset $U\subseteq X$ admits a $\sT_\srig$-admissible refinement by sets in $S$, and
\item  $\check{H}^q(U,\sF)$ vanishes for all $q>0$ and all $U\in S$.
\end{packed_enum}
The statement now follows via the standard \v{C}ech spectral sequence argument.
\end{proof}

			\subsection{Semi-affinoid spaces as locally G-ringed spaces}\label{locgringedsec}
				Let $X$ be a semi-affinoid $K$-space, and let $A$ denote its ring of global functions. It follows from Lemma \ref{retrolocallem} and from Lemma \ref{rigurigadmsubsetcomparelem} that the stalk $\O_{X,x}$ of $\O_X$ in a point $x\in X$ coincides with the corresponding stalk on the rigid-analytic space obtained via Berthelot's construction. 

In particular, $\O_{X,x}$ is an excellent noetherian local ring. Moreover, it follows that $X$ together with its uniformly rigid $G$-topology and its sheaf of semi-affinoid functions $\O_X$ is a locally $G$-ringed $K$-space. Transcribing the discussion at the beginning of \cite{BGR} Section 9.3.1, we see that it depends on $X$ in a functorial way.


Let $\m\subseteq A$ be the maximal ideal corresponding to $x$. In complete analogy with \cite{BGR} 7.3.2/1, one can directly prove that $\O_{X,x}$ is a local ring with maximal ideal $\{h\in\O_{X,x}\,;\,h(x)=0\}$ and that the maximal ideal of $\O_{X,x}$ coincides with $\m\O_{X,x}$. Moreover, one proves the following:

\begin{lem}\label{stalkscomplem}
The natural homomorphism $A_\m\rightarrow\O_{X,x}$ is injective, and for each integer $n\in\N$, the natural homomorphism
\[
A/\m^{n+1}\rightarrow\O_{X,x}/\m^{n+1}\O_{X,x}
\]
is an isomorphism. In particular, we obtain an isomorphism $\hat{A}_\m\overset{\sim}{\rightarrow}\hat{\O}_{X,x}$ of maximal-adic completions.
\end{lem}


This is the semi-affinoid analogon of \cite{BGR} 7.3.2/3, and the \emph{proof} carries over \emph{verbatim}. 

It follows that a semi-affinoid function on $X$ is zero if and only if it vanishes in all stalks. In particular, the vanishing of semi-affinoid functions can be checked with respect to an arbitrary, possibly not $\sT_\srig$-admissible, covering of $X$ by semi-affinoid subdomains. 



\begin{prop}\label{ffprop}
The category of semi-affinoid $K$-spaces is a full subcategory of the category of locally G-ringed $K$-spaces.
\end{prop}

If suffices to literally transcribe the \emph{proof} of \cite{BGR} 9.3.1/2, using the fact that the semi-affinoid subdomains form a basis for $\sT_\srig$.

Let us recall that a morphism $\phi\colon Y\rightarrow X$ of locally G-ringed $K$-spaces is called an \emph{open immersion} if it is an isomorphism onto an admissible open subset of $Y$. It is called a \emph{local isomorphism} if every point $y$ of $Y$ admits an admissible open neighborhood $U$ such that the restriction $U\rightarrow X$ of $\phi$ is an open immersion. 
Example \ref{nonadmdisccovex} shows that the analog of \cite{BGR} 8.2.1/4 does not hold: There exist injective local isomorphisms of semi-affinoid $K$-spaces which do not even represent semi-affinoid pre-subdomains.



			\subsection{Uniformly rigid spaces}\label{srigspacessec}
				We can now easily \emph{globalize} the category of semi-affinoid $K$-spaces, thereby obtaining the category of \emph{uniformly rigid} $K$-spaces. 

\begin{defi}\label{srigspacedefi}\index{uniformly rigid space}\index{category!of uniformly rigid $K$-spaces}
A \emph{uniformly rigid} $K$-space is a \emph{saturated} locally G-ringed $K$-space $(X,\O_X)$ admitting an admissible  \emph{semi-affinoid covering}\index{covering!semi-affinoid}, that is, an admissible covering $(X_i)_{i\in I}$ such that for each $i\in I$, $(X_i,\O_X|_{X_i})$ is isomorphic to a semi-affinoid $K$-space. 
\end{defi}

A morphism of uniformly rigid $K$-spaces is a morphism of locally G-ringed $K$-spaces. We let $\sRig_K$ denote the category of uniformly rigid $K$-spaces. By Proposition \ref{saturatedprop} and Proposition \ref{ffprop}, $\sAff_K$ is a full subcategory of $\sRig_K$.
An admissible open subset $U$ of a uniformly rigid $K$-space $X$ is called a \emph{semi-affinoid subspace}\index{subspace!semi-affinoid} of $X$ if $(U,\O_X|_U)$ is isomorphic to a semi-affinoid $K$-space. A admissible open subset $U$ of a semi-affinoid $K$-space $X$ is called \emph{weakly retrocompact}\index{subspace!weakly retrocompact} if it is a finite union of retrocompact semi-affinoid subdomains in $X$.

We do not know whether a semi-affinoid subspace $U$ of a semi-affinoid $K$-space $X$ is necessarily a semi-affinoid subdomain in $X$.  One verifies that $U$ is a semi-affinoid pre-subdomain in $X$.

\begin{lem}
Let $X=\sSp A$ be a semi-affinoid $K$-space, and let $U=\sSp B$ be a semi-affinoid subspace of $X$. Then the restriction homomorphism $\phi\colon A\rightarrow B$ is flat.
\end{lem}
\begin{proof}
For every maximal ideal $\n\subseteq B$ with corresponding point $x\in U$ and $\phi$-preimage $\m\subseteq A$, the induced homomorphism $\phi_\n\colon A_\m\rightarrow B_\n$ induces an isomorphism of maximal-adic completions, the stalks $\O_{X,x}$ and $\O_{U,x}$ being naturally isomorphic. By the Flatness Criterion \cite{Bourbaki} III.5.2 Theorem 1, we conclude that $A\rightarrow B_\n$ is flat for all maximal ideals $\n$ in $B$, which means that $\phi$ is flat. 
\end{proof}






\begin{lem}\label{saffbasislem}
The semi-affinoid subspaces of a uniformly rigid $K$-space $X$ form a \emph{basis} for the G-topology on $X$.
\end{lem}
\begin{proof}
Let $(X_i)_{i\in I}$ be an admissible semi-affinoid covering of $X$, and let $U\subseteq X$ be an admissible open subset. Then $(X_i\cap U)_{i\in I}$ is an admissible covering of $U$. For each $i\in I$, $X_i\cap U$ is admissible open in $X_i$ and, hence, admits an admissible covering by semi-affinoid subdomains of $X_i$. Hence, $U$ has an admissible semi-affinoid covering.
\end{proof}

It follows that if $X$ is a uniformly rigid $K$-space and if $U\subseteq X$ is an admissible open subset, then $(U,\O_X|_U)$ is a uniformly rigid $K$-space, again.

It is now clear that the Glueing Theorem \cite{BGR} 9.3.2/1 and its proof carry over \emph{verbatim} to the uniformly rigid setting: $\sT_\srig$-descent is effective in the category of uniformly rigid $K$-spaces. Similarly, a morphism of uniformly rigid spaces can be defined locally on the domain; this is the uniformly rigid version of \cite{BGR} 9.3.3/1, and again the proof is obtained by literal transcription. The semi-affinoid subspaces of a uniformly rigid $K$-space forming a basis for the G-topology on that space, a uniformly rigid $K$-space is determined by its functorial points with values in \emph{semi-affinoid} $K$-spaces. 

We can also copy the proof of \cite{BGR} 9.3.3/2 to see that if $X$ is a semi-affinoid $K$-space and if $Y$ is a uniformly rigid $K$-space, then the set of morphisms from $Y$ to $X$ is naturally identified with the set of $K$-algebra homomorphisms from $\O_X(X)$ to $\O_Y(Y)$. 

Let $\fX$ be an affine formal $R$-scheme of ff type with semi-affinoid generic fiber $X$. The associated specialization map $\sp_\fX$ which we discussed in Section \ref{specmapsec} is naturally enhanced to a morphism of G-ringed $R$-spaces $\sp_\fX\colon X\rightarrow\fX$. Morphisms of uniformly rigid $K$-spaces being defined locally on the domain, we see that $\sp_\fX$ is \emph{universal} among all morphisms of G-ringed $R$-spaces from uniformly rigid $K$-spaces to $\fX$. Using this universal property, we can construct the \emph{uniformly rigid generic} fiber $\fX^\srig$\index{generic fiber!uniformly rigid} of a general formal $R$-scheme of locally ff type $\fX$, together with a functorial specialization map $\sp_\fX\colon\fX^\srig\rightarrow\fX$ which is universal among all morphisms of G-ringed $R$-spaces from uniformly rigid $K$-spaces to $\fX$; this process does \emph{not} involve Berthelot's construction. Arguing as in the proof of Corollary \ref{rigflatfaithfulcor}, we see that $\urig$ is \emph{faithful} on the category of \emph{flat} formal $R$-schemes of locally ff type. A formal \emph{$R$-model}\index{model!of a uniformly rigid space} of a uniformly rigid $K$-space $X$ is a formal $R$-scheme $\fX$ of locally ff type together with an isomorphism $X\cong\fX^\srig$. The map $\sp_\fX$ is surjective onto the closed points of $\fX$ whenever $\fX$ is flat over $R$. This follows from Lemma \ref{spsurjlem}, together with the remark that the underlying topological space of $\fX$ is a Jacobson space, cf.\ \cite{EGAIn} 0.2.8 and 6.4, so that the condition on a point in $\fX$ of being closed is local.

By Proposition \ref{amalgsumsprop}, the category of semi-affinoid $K$-spaces has \emph{fibered products}\index{fibered product!of uniformly rigid spaces}; following the method outlined in \cite{BGR} 9.3.5, we see that the category of uniformly rigid $K$-spaces has fibered products as well and that these are constructed by glueing semi-affinoid fibered products of semi-affinoid subspaces. It is clear from this description that the $\srig$-functor preserves fibered products.

Semi-affinoid subspaces of semi-affinoid spaces can be described in the style of the Gerritzen-Grauert Theorem \cite{BGR} 7.3.5/3:

\begin{lem}\label{semaffsubspacecharlem}
Let $X$ be a semi-affinoid $K$-space, and let $U\subseteq X$ be a semi-affinoid subspace. Then $U$ admits a leaflike covering $(U_i)_{i\in I}$ such that each $U_i$ is a semi-affinoid subdomain in $X$. 
\end{lem}
\begin{proof}
By Proposition \ref{explicitprop}, $U$ admits a covering $(V_j)_{j\in J}$ by semi-affinoid subdomains $V_j$ of $X$ which is refined by a leaflike covering $(U_i)_{i\in I}$ of $U$; we claim that the $U_i$ are semi-affinoid subdomains in $X$. Via pullback, the $V_j$ are semi-affinoid subdomains of $U$. Let $\phi\colon I\rightarrow J$ denote a refinement map.  By Corollary \ref{semaffsubdomgoodnotioncor}, for each $i\in I$ the set $U_i$ is a semi-affinoid subdomain in $V_{\phi(i)}$ and, hence, in $X$, as desired. 
\end{proof}

A morphism of uniformly rigid $K$-spaces is called \emph{flat}\index{morphism!flat} in a point of its domain if it induces a flat homomorphisms of stalks in this point, and it is called flat if it is flat in all points. Clearly a morphism of semi-affinoid $K$-spaces is flat in the sense of locally G-ringed $K$-spaces if and only if the underlying homomorphism of rings of global sections is flat.
				
			\subsection{Comparison functors}\label{compfuncsec}

				In this section, we compare the categories of rigid and uniformly rigid $K$-spaces. In particular, we define a functor 
\[
X\mapsto X^\r
\]
from the category of uniformly rigid $K$-spaces to the category of rigid $K$-spaces, together with functorial comparison morphisms 
\[
\comp_X\colon X^\r\rightarrow X
\]
of locally G-ringed $K$-spaces which are bijections on physical points and which induce isomorphisms of stalks.

				\subsubsection{The rigid space associated to a uniformly rigid space}\label{rigspaceofsemrigspacesec}
					In Section \ref{assocrigspacesubsec}, we have introduced the rigid $K$-space $X^\r$ associated to a semi-affinoid $K$-space $X=\sSp A$. We have seen that there exists a $K$-homomorphism $\phi\colon A\rightarrow\Gamma(X^\r,\O_{X_\r})$ such that the pair $(X^\r,\phi)$ is characterized by the universal property stated in Proposition \ref{semaffspaceassocprop}. 

In fact, $\phi$ extends to a morphism of locally G-ringed $K$-spaces
\[
\textup{comp}_X\colon X^\r\rightarrow X
\]
which is universal among all morphisms from rigid $K$-spaces to $X$; it is called the \emph{comparison morphism}\index{comparison morphism} associated to $X$. On the level of points, we define $\textup{comp}_X$ to be the bijection induced by $\phi$, cf.\ Section \ref{assocrigspacesubsec}; then $\textup{comp}_X$ is \emph{continuous} by Proposition \ref{srigrigcompprop}. If $U\subseteq X$ is a semi-affinoid subdomain and if $\tau$ denotes the corresponding open immersion, the remark at the end of Section \ref{assocrigspacesubsec} and the discussion following Definition \ref{semaffsubdomdefdef} show that the morphism $\tau^\r$ provided by Proposition \ref{semaffspaceassocprop} is an open immersion; considering functorial points with values in finite extensions of $K$, we see that $\im\tau^\r$ is the $\comp_X$-preimage of $U$. The semi-affinoid subdomains of $X$ forming a basis for the uniformly rigid $G$-topology on $X$, it follows from Proposition \ref{semaffspaceassocprop} that $\phi$ extends naturally to a sheaf homomorphism $\comp_X^\sharp$. By our discussion at the beginning of Section \ref{locgringedsec}, the resulting morphism $\comp_X$ of G-ringed $K$-spaces induces isomorphisms of stalks; in particular, it is local.

\begin{prop}\label{compmorunivprop}\index{rigid space!of a uniformly rigid space}
The morphism $\comp_X$ is universal among all morphisms of locally G-ringed $K$-spaces from rigid $K$-spaces to $X$.
\end{prop}
\begin{proof}
Let $Y$ be a rigid $K$-space, and let $\psi\colon Y\rightarrow X$ be a morphism of locally G-ringed $K$-spaces. By Proposition \ref{semaffspaceassocprop}, there is a unique morphism $\psi^\r\colon Y\rightarrow X^\r$ such that $\psi$ and $\comp_X\circ\psi^\r$ coincide on global sections. Since the points and the completed stalks of $X$ are recovered from the $K$-algebra of global sections of $X$, it follows that $\psi$ and $\comp_X\circ\psi^\r$ coincide.
\end{proof}

Let $X$ be any uniformly rigid $K$-space. Since the semi-affinoid subspaces of $X$ form a basis for the G-topology on $X$, we can invoke standard glueing arguments to show that the comparison morphisms attached to these semi-affinoid subspaces glue to a universal comparison morphism
\[
\comp_X\colon X^\r\rightarrow X
\]
from a rigid $K$-space to $X$.

The functor $X\mapsto X^\r$ is faithful, yet forgetful. For example, it is easily seen that an unbounded function on the rigid open unit disc induces a morphism to the rigid projective line over $K$ which is not induced by a morphism from the semi-affinoid open unit disc $\sSp (R[[S]]\otimes_RK)$ to the uniformly rigid projective line over $K$. Likewise, the functor $\r$ forgets the distinction between the semi-affinoid open unit disc just mentioned and the uniformly rigid open unit disc that is the generic fiber of a quasi-paracompact formal $R$-model of locally tf type for the rigid open unit disc $\D^1_K$. However, one can prove that $X\mapsto X^\r$ is fully faithful on the full subcategory of \emph{reduced semi-affinoid} $K$-spaces.

The functor $X\mapsto X^\r$ preserves fibered products. Indeed, this may be checked in the semi-affinoid situation, where it follows from the fact that fibered products of semi-affinoid spaces are uniformly rigid generic fibers of fibered products of affine flat formal $R$-models, together with the fact that Berthelot's $\rig$-functor preserves fibered products, cf.\ \cite{dJ} 7.2.4 (g). In particular, $X\mapsto X^\r$ preserves group structures.

A uniformly rigid $K$-space $X$ is called \emph{smooth}\index{uniformly rigid space!smooth} in a point $x\in X$ if its associated rigid space $X^\r$ is smooth in $x$. Of course, we could define and further analyze this smoothness property intrinsically within the category of uniformly rigid $K$-spaces. However, since a uniformly rigid $K$-space is indistinguishable from its associated rigid $K$-space as far as local properties are concerned, no real benefit is to be expected from such an effort, and so we refrain from pursuing this any further.
			
				\subsubsection{The uniformly rigid space associated to a rigid space}\label{ufspofrigspsec}

We let $\fC$ temporarily denote the category of quasi-paracompact flat formal $R$-schemes of locally tf type, and we let $\fC_\Bl$ denote its localization with respect to the class of admissible formal blowups. It follows easily from the definitions that the functor $\srig|_\fC\colon\fC\rightarrow\sRig_K$ factorizes over a functor $\sr'\colon\fC_\Bl\rightarrow\sRig_K$. By \cite{FRG} Theorem 2.8/3, the functor $\rig$ induces an equivalence $\rig_\Bl$ between $\fC_\Bl$ and the category of quasi-paracompact and quasi-separated rigid $K$-spaces. The functor $\rig_\Bl$ will be called the \emph{Raynaud equivalence}. Composing $\sr'$ with a quasi-inverse of $\rig_\Bl$, we obtain a functor $\sr\colon\Rig_K\rightarrow\sRig_K$; if $Y$ is a quasi-paracompact and quasi-separated rigid $K$-space, we say that $Y^\sr\mathrel{\mathop:}=\sr(Y)$ is the \emph{uniformly rigid $K$-space associated to $X$}\index{uniformly rigid space!of a rigid space} or the \emph{uniform rigidification} of $Y$. Of course, it depends on the choice of a quasi-inverse of the Raynaud equivalence.


\begin{prop}\label{srigspaceqisoidprop}
The composite functor $\r\circ\sr$ is quasi-isomorphic to the identity on the category of quasi-paracompact and quasi-separated rigid $K$-spaces.
\end{prop}
\begin{proof} Let $\rig_\Bl^{-1}$ denote the chosen inverse of the Raynaud equivalence. Let $Y$ be a quasi-paracompact and quasi-separated rigid $K$-space; then $\rig_\Bl^{-1}(Y)$ is a quasi-paracompact flat formal $R$-model of locally tf type for $Y$, and $Y^\sr=\rig_\Bl^{-1}(Y)^\srig$, which implies that $(Y^\sr)^\r=\rig_\Bl^{-1}(Y)^\rig$, functorially in $Y$. That is, $\r\circ\sr=\rig\circ \rig_\Bl^{-1}$, which is isomorphic to the identity functor.
\end{proof}

In particular, via the choice of an isomorphism $\r\circ\sr\cong\id$, the comparison morphisms $\comp_{Y^\sr}$ induce functorial comparison morphisms
\[
\comp_Y\colon Y\cong (Y^\sr)^r\rightarrow Y^\sr
\]
for all quasi-paracompact and quasi-separated rigid $K$-spaces $Y$. As before, these comparison morphisms reside in the category of locally G-ringed $K$-spaces.

\begin{cor}\label{assocsrigspacecompcor}
Let $Y$ be a quasi-paracompact and quasi-separated rigid $K$-space. Then $\comp_Y$ is the universal morphism from $Y$ to a uniformly rigid $K$-space.
\end{cor}
\begin{proof}
Let $X$ be a uniformly rigid $K$-space, and let $\psi\colon Y\rightarrow X$ be a morphism of locally G-ringed $K$-spaces. The morphism $\comp_Y$ is a bijection on points, and it induces isomorphisms of stalks; hence the morphism $Y^\sr\rightarrow X$ that we seek is unique if it exists. If $Y$ is affinoid and $X$ is semi-affinoid, there is nothing to show. Let $(X_i)_{i\in I}$ be an admissible semi-affinoid covering of $X$, and let $(Y_j)_{j\in J}$ be an admissible affinoid covering of $Y$ refining $(\psi^{-1}(X_i))_{i\in I}$. It suffices to see that $(Y_j^\sr)_{j\in J}$ is an admissible covering of $Y^\sr$. By \cite{FRG} Lemma 2.8/4, there exists a flat quasi-paracompact $R$-model of locally tf type $\fY$ for $Y$ such that $(Y_j)_{j\in J}$ is induced by an open covering of $\fY$. Since $\fY^\srig=Y^\sr$, it follows that $(Y_j^\sr)_{j\in J}$ is an admissible covering of $Y^\sr$, as desired. 
\end{proof}

\begin{cor}\label{ffcor}
The functor $\sr$ is fully faithful.
\end{cor}
\begin{proof}
Let $X$ and $Y$ be quasi-paracompact and quasi-separated rigid $K$-spaces. By Proposition \ref{srigspaceqisoidprop},  Proposition \ref{compmorunivprop} and Corollary \ref{assocsrigspacecompcor}, we have functorial bijections
\begin{eqnarray*}
\Hom(Y,X)&\cong&\Hom(Y,(X^\sr)^\r)\\
&\cong&\Hom(Y,X^\sr)\\
&\cong&\Hom(Y^\sr,X^\sr)\quad.
\end{eqnarray*}
\end{proof}

Of course, if $X$ is any uniformly rigid $K$-space, then the comparison morphism 
\[
\comp_X\colon X^\r\rightarrow X
\]
does \emph{not} represent all morphisms from $X^\r$ to uniformly rigid $K$-spaces. For example, if $X$ is the semi-affinoid open unit disc $\sSp (R[[S]]\otimes_RK)$, the natural morphism $\comp_{X^\r}$ from the rigid open unit disc $X^\r$ to its uniform rigidification $(X^\r)^\sr$ does \emph{not} extend to a morphism $X\rightarrow (X^\r)^\sr$. Indeed, such a morphism would need to be the identity on points, but $X$ is quasi-compact while $(X^\r)^\sr$ is \emph{not} quasi-compact.

The functor $Y\mapsto Y^\sr$ does \emph{not} respect arbitrary open immersions. For example, if $Y'\subseteq Y$ is the inclusion of the open rigid unit disc into the closed rigid unit disc, the morphism $(Y')^\sr\rightarrow Y^\sr$ is not an open immersion: its image is the semi-affinoid open unit disc, while $(Y')^\sr$ is not quasi-compact. However, it follows from \cite{BL2} 5.7 that $\sr$ preserves open immersions of \emph{quasi-compact} rigid $K$-spaces.

Quasi-separated rigid $K$-spaces being obtained from affinoid $K$-spaces by glueing along \emph{quasi-compact} admissible open subspaces, it thus follows that $\sr$ preserves fibered products. Indeed, this can now be checked in an affinoid situation, where the statement is clear from the construction of semi-affinoid fibered products in Section \ref{amalgsumssec}. In particular, $Y\mapsto Y^\sr$ preserves group structures.

%

			\subsection{Change of the base field}\label{basefieldchangesec}				
				Let $K'/K$ be a possibly non-finite discrete analytic extension field of $K$, and let $R'$ denote its valuation ring. Let $A$ be a semi-affinoid $K$-algebra, and let $\ul{A}\subseteq A$ be an $R$-lattice of ff type for $A$. Then $\ul{A}\hat{\otimes}_RR'$ is a flat $R'$-algebra of ff type, and we define
\[
A\hat{\otimes}_KK'\,\mathrel{\mathop:}=\,(\ul{A}\hat{\otimes}_RR')\otimes_RK\quad.
\]
The $K$-algebra $A\hat{\otimes}_KK'$ together with the $K$-homomorphism $A\rightarrow A\hat{\otimes}_KK'$ does not depend on the choice of $\ul{A}$ up to unique $K$-isomorphism. Indeed, let $\ul{A}'$ be another $R$-lattice of ff type for $A$. After replacing $\ul{A}'$ by $\ul{A}+\ul{A}'$, we may assume that $\ul{A}\subseteq\ul{A}'$, which defines a finite admissible formal blowup by Corollary \ref{freesemaffcor} ($ii$). By Proposition \ref{basicblowupprop} and \cite{EGAIII1} 4.8.3 ($iii$), finite formal blowups are preserved by flat base change in the category of formal $R$-scheme of locally ff type, so we conclude that 
\[
\ul{A}\hat{\otimes}_RR'\rightarrow\ul{A}'\hat{\otimes}_RR'
\]
is a finite formal blowup; it is admissible since the topological ring extension $R'/R$ is adic. We thereby obtain an isomorphism of $K'$-algebras after inverting a uniformizer $\pi'$ of $R'$. Uniqueness follows from the universal property of $\ul{A}\hat{\otimes}_RR'$. We thus obtain a \emph{base extension functor}
\[
\cdot\hat{\otimes}_KK'\colon\sAff_K\rightarrow\sAff_{K'}
\]
from the category of semi-affinoid $K$-spaces to the category of semi-affinoid $K'$-spaces. It clearly preserves semi-affinoid subdomains and retrocompact semi-affinoid subdomains. 


\begin{lem}\label{semaffsubspacecovbclem}
Let $X$ be a semi-affinoid $K$-space, and let $(X_i)_{i\in I}$ be an admissible semi-affinoid subdomain covering of $X$. Then $(X_i\hat{\otimes}_KK')_{i\in I}$ is an admissible semi-affinoid subdomain covering of $X\hat{\otimes}_KK'$.
\end{lem}
\begin{proof}
After passing to a refinement, we may assume that $(X_i)_{i\in I}$ is leaflike; it suffices to show that $(X_i\hat{\otimes}_KK')_{i\in I}$ is leaflike as well. Let $(X_j)_{j\in J}$ be a treelike extension of $(X_i)_{i\in I}$, together with a suitable rooted tree structure on $J$ such that $I\subseteq J$ is identified with the set of leaves. By Lemma \ref{simpadlem}, treelike coverings are induced by treelike formal coverings, and treelike formal coverings are clearly preserved under the base change $\cdot\hat{\otimes}_RR'$. It follows that $(X_j\hat{\otimes}_KK')_{j\in J}$ is a treelike covering for $X\hat{\otimes}_KK'$ such that the given rooted tree structure on $J$ is suitable for this covering. We conclude that $(X_i\hat{\otimes}_KK')_{i\in I}$ is a leaflike covering of $X\hat{\otimes}_KK'$, as desired.
\end{proof}

By Proposition \ref{retroprop}, it follows that $\cdot\hat{\otimes}_KK'$ preserves open immersions of weakly retrocompact semi-affinoid subspaces:

\begin{lem}\label{semaffsubspacebclem}
Let $X$ be a semi-affinoid $K$-space, and let $U\subseteq X$ be a weakly retrocompact semi-affinoid subspace; then the natural morphism from $U\hat{\otimes}_KK'$ to $X\hat{\otimes}_KK'$ is an open immersion of semi-affinoid $K'$-spaces.
\end{lem}
\begin{proof}
Let $(U_i)_{i\in I}$ be a finite covering of $U$ by retrocompact semi-affinoid subdomains of $X$. By Lemma \ref{semaffsubspacecovbclem}, the semi-affinoid subdomains $U_i\hat{\otimes}_KK'$ in $U\hat{\otimes}_KK'$ form an admissible semi-affinoid subdomain covering of $U\hat{\otimes}_KK'$. Since $\cdot\hat{\otimes}_KK'$ preserves semi-affinoid subdomains, it follows that the natural morphism
\[
\phi\colon U\hat{\otimes}_KK'\rightarrow X\hat{\otimes}_KK'
\]
is an injective local isomorphism whose image is the finite union of the retrocompact semi-affinoid subdomains $U_i\hat{\otimes}_KK'$. By Proposition \ref{retroprop}, this union is admissible; hence $\phi$ is an open immersion.
\end{proof}

By rather tedious glueing arguments, one can show that $\cdot\hat{\otimes}_KK'$ extends to a functor on the category of uniformly rigid $K$-spaces $X$ admitting an admissible semi-affinoid covering $(X_i)_{i\in I}$ such that for all $i,j\in I$, $X_i\cap X_j$ is a finite union of \emph{retrocompact} semi-affinoid subdomains both in $X_i$ and in $X_j$. For example, this condition is satisfied for all $\fX^\srig$, where $\fX$ varies among the formal $R$-schemes of locally ff type, and
\[
\fX^\urig\hat{\otimes}_KK'\,=\,(\fX\hat{\otimes}_RR')^\urig\quad.
\]
In our later applications, the functor $\cdot\hat{\otimes}_KK'$ will be applied only to semi-affinoid $K$-spaces and to quasi-separated rigid $K$-spaces; we therefore abstain from elaborating on the details of the globalization of $\cdot\hat{\otimes}_KK'$. Let us just mention one lemma which is useful in this context:

\begin{lem}
Let $X$ and $Y$ be semi-affinoid $K$-spaces, and let 
\begin{eqnarray*}
&\phi\colon U\hookrightarrow X&\quad,\\
&\psi\colon U\hookrightarrow Y&
\end{eqnarray*}
be open immersions of uniformly rigid $K$-spaces such that $U$ is a finite union of retrocompact semi-affinoid subdomains in $X$ and a finite union of retrocompact semi-affinoid subdomains in $Y$. Then $U$ admits a finite covering by semi-affinoid subspaces which are retrocompact semi-affinoind subdomains in both $X$ and $Y$.
\end{lem}
\begin{proof}
Let $(X_i)_{i\in I}$ and $(Y_j)_{j\in J}$ be finite coverings of $U$ by retrocompact semi-affinoid subdomains in $X$ and $Y$ respectively. Let us fix indices $i\in I$, $j\in J$. Since $Y_j\subseteq U\subseteq X$ is a semi-affinoid subspace of $X$, the intersection $X_i\cap Y_j=X_i\times_XY_j$ is semi-affinoid. We claim that $X_i\cap Y_j$ is a retrocompact semi-affinoid subdomain in $X$; then by symmetry $X_i\cap Y_j$ is a retrocompact semi-affinoid subdomain in $Y$. Since $X_i$ is a retrocompact semi-affinoid subdomain in $X$, it suffices to see that $X_i\cap Y_j$ is a retrocompact semi-affinoid subdomain in $X_i$. Since $Y_j$ is a retrocompact semi-affinoid subdomain in $Y$ and since $X_i\cap Y_j=X_i\times_YY_j$, the statement follows from Corollary \ref{semaffpreimcor} and the remark at the end of Section \ref{preimsemaffsec}.
 \end{proof}

	\section{Coherent modules on uniformly rigid spaces}\label{cohmodsec}
			
\subsection{Coherent and strictly coherent modules}


Let $X$ be a G-ringed $K$-space, and let $\sF$ be an $\O_X$-module. Let us recall some standard definitions concerning the coherence property, cf.\ \cite{FRG} 1.14/2:
\begin{packed_enum}\index{module!coherent}
\item $\sF$ is called of \emph{finite type} if there exists an admissible covering $(X_i)_{i\in I}$ of $X$ together with exact sequences
\[
\O_X^{s_i}|_{X_i}\rightarrow\sF|_{X_i}\rightarrow 0\quad.
\]
\item $\sF$ is called \emph{coherent} if $\sF$ is of finite type and if for any admissible open subspace $U\subseteq X$, the kernel of any morphism $\O_X^s|_U\rightarrow\sF|_U$ is of finite type.
\end{packed_enum}

If $X$ is a semi-affinoid $K$-space with ring of functions $A$, the functor $M\mapsto \tilde{M}$ is well-behaved, as we recall in the following lemma. The proof of Lemma \ref{assocbasicpropertieslem} is identical to the proof of \cite{FRG} 1.14/1; one uses the fact that the restriction homomorphisms associated to semi-affinoid subdomains are flat:

\begin{lem}\label{assocbasicpropertieslem}\index{module!associated}
The functor $M\mapsto\tilde{M}$ is fully faithful, and it commutes with the formation of kernels, images, cokernels and tensor products. Moreover, a sequence of $A$-modules
\[
0\rightarrow M'\rightarrow M\rightarrow M''\rightarrow 0
\]
is exact if and only if the associated sequence
\[
0\rightarrow \tilde{M}'\rightarrow \tilde{M}\rightarrow \tilde{M}''\rightarrow 0
\]
of $\O_X$-modules is exact.
\end{lem}

For a semi-affinoid $K$-space $X=\sSp A$, we have $\O_X^r=A^r\otimes\O_X$. Since $A$ is \emph{noetherian}, it follows from Lemma \ref{assocbasicpropertieslem} that kernels and cokernels of morphisms of type $\O_X^r\rightarrow\O_X^s$ are  associated. We thus conclude that an $\O_X$-module on a uniformly rigid $K$-space $X$ is coherent if and only if there exists an admissible semi-affinoid covering $(X_i)_{i\in I}$ of $X$ such that $\sF|_{X_i}$ is associated for all $i\in I$.



In particular, the structural sheaf $\O_X$ of any uniformly rigid $K$-space $X$ is coherent. Moreover, it follows from Lemma \ref{assocbasicpropertieslem} that kernels and cokernels of  morphisms of coherent $\O_X$-modules are coherent.

\begin{lem}
Let $\phi\colon Y\rightarrow X$ be a morphism of uniformly rigid $K$-spaces, and let $\sF$ be a coherent $\O_X$-module. Then $\phi^*\sF$ is a coherent $\O_Y$-module. 
\end{lem}
\begin{proof}
Indeed, we may assume that $X$ and $Y$ are semi-affinoid, $X=\sSp A$, $Y=\sSp B$, and that $\sF$ is associated to a finite $A$-module $M$. Then $\phi^*\sF$ is associated to $M\otimes_AB$, where $B$ is an $A$-algebra via $\phi^*$.
\end{proof}


\begin{defi}
Let $X$ be a uniformly rigid $K$-space. An $\O_X$-module $\sF$ is called \emph{strictly coherent}\index{module!strictly coherent} if for any semi-affinoid subspace $U\subseteq X$, the restriction $\sF|_U$ is an associated module.
\end{defi}

For example, the structural sheaf of a uniformly rigid $K$-space is strictly coherent. Since we do not know whether a semi-affinoid subspace of a semi-affinoid $K$-space is a semi-affinoid sub\emph{domain}, it is not \emph{a priori} clear whether any associated module on a semi-affinoid $K$-space is strictly coherent. In Corollary \ref{frameembassoccor}, however, we will show that this is indeed the case.

\subsection{Integral models of submodules}

Let $X$ be a uniformly rigid $K$-space. We introduce the following terminology: A \emph{frame} on $X$ is a finite direct sum $\O_X^r$ of $\O_X$, and a coherent $\O_X$-module $\sF$ is called \emph{frame-embeddable}\index{module!frame-embeddable} if there exists an injective $\O_X$-homomorphism $\sF\hookrightarrow\O_X^r$ for some $r\in\N$; such a homomorphism is called a \emph{frame embedding}. In particular, coherent $\O_X$-\emph{ideals} are frame-embeddable. In this section, we will prove the important fact that frame-embeddable coherent $\O_X$-modules are strictly coherent.

\begin{lem}\label{frameembpblem}
Let $\phi\colon Y\rightarrow X$ be a flat morphism of uniformly rigid $K$-spaces, and let $\sF$ be a frame-embeddable $\O_X$-module. Then $\phi^*\sF$ is frame-embeddable.
\end{lem}
\begin{proof}
Let us choose a frame embedding $\psi\colon\sF\hookrightarrow\O_X^r$; we obtain an induced morphism $\phi^*\psi\colon\phi^*\sF\rightarrow\O_Y^r$, and it suffices to see that it is injective. We may work locally on $Y$ and, hence, on $X$; thereby we may assume that $X$ and $Y$ are affine and that $\sF$ is associated. The statement is now clear.
\end{proof}

In order to study the frame-embeddable coherent modules on $X$, we consider integral models of coherent $\O_X$-modules. If $\fX$ is a formal $R$-scheme of locally ff type and if $\sF$ is a coherent $\O_\fX$-module, we obtain a coherent $\O_X$-module $\sF^\srig$ on $\fX^\srig$ which we call the \emph{uniformly rigid generic fiber} of $\sF$. If $X$ is a uniformly rigid $K$-space, if $\sF$ is a coherent $\O_X$-module and if $\fX$ is a flat formal $R$-model of locally ff type for $X$, then an $R$-\emph{model}\index{model!of a coherent module} of $\sF$ on $\fX$ is a coherent $\O_\fX$-module $\ul{\sF}$ together with an isomorphism $\ul{\sF}^\srig\cong\sF$ that is compatible with the given identification $\fX^\srig\cong X$. Sometimes we will not mention the isomorphism $\ul{\sF}^\srig\cong\sF$ explicitly. Let us recall that we have also defined the \emph{generic fiber} $\sF_K$ of $\sF$ as a sheaf on $\fX$, cf.\ Definition \ref{semafffuncmodeldef}. It is clear that 
\[
\sF_K\,=\,\sp_{\fX,*}(\sF^\srig)
\]
and that the functor $\srig$ from coherent $\O_\fX$-modules to coherent $\O_{\fX^\srig}$-modules factorizes naturally over the factor sending a coherent $\O_\fX$-module $\sF$ to its generic fiber $\sF_K$.

A direct power $\O_X^r$ of the structural sheaf has the pleasant property of admitting a natural model $\O_\fX^r$ on every flat formal $R$-model of locally ff type $\fX$ for $X$. We will show that frame-embeddings enjoy the same property. Let us first consider the affine situation. 


\begin{lem}\label{submodulecanmodellem}
Let $\ul{A}$ be an $R$-algebra, let $\ul{M}$ be an $\ul{A}$-module, and let $N\subseteq \ul{M}\otimes_RK$ be an $\ul{A}\otimes_RK$-submodule. Then there exists a unique $\ul{A}$-submodule $\ul{N}\subseteq\ul{M}$ such that the natural homomorphism $\ul{N}\otimes_RK\rightarrow\ul{M}\otimes_RK$ is an isomorphism onto $N$ and such that $\ul{M}/\ul{N}$ is $R$-flat.
\end{lem}
\begin{proof}
Let us abbreviate $M\mathrel{\mathop:}=\ul{M}\otimes_RK$, and let us set
\[
\ul{N}\,\mathrel{\mathop:}=\,\ker(\ul{M}\rightarrow M/N)\quad;
\]
then $\ul{N}$ is an $\ul{A}$-submodule of $\ul{M}$. If $n\in N$, there exists an $s\in\N$ such that $\pi^sn$ lies in the image of $\ul{M}$ in $M$; the natural $K$-homomorphism $\ul{N}\otimes_RK\rightarrow N$ is thus bijective. As an $\ul{A}$-submodule of $M/N$, the quotient $\ul{M}/\ul{N}$ is free of $\pi$-torsion and, hence, $R$-flat. 

If $\ul{N}'\subseteq\ul{M}$ is another $\ul{A}$-submodule whose image in $M$ generates $N$ as an $\ul{A}\otimes_RK$-module, then $\ul{N}'$ lies in the kernel $\ul{N}$ of $\ul{M}\rightarrow M/N$. If in addition $\ul{M}/\ul{N}'$ is flat over $R$, then the natural homomorphism $\ul{M}/\ul{N}'\rightarrow \ul{M}/\ul{N}'\otimes_RK=M/N$ is injective, which proves that $\ul{N}'$ coincides with this kernel.
\end{proof}

\begin{prop}\label{frameembprop}
Let $X$ be a uniformly rigid $K$-space, let $\sF'\subseteq\sF$ be an inclusion of coherent $\O_X$-modules, and let $\fX$ be an $R$-model of locally ff type for $X$ such that $\sF$ admits be an $R$-model  $\ul{\sF}$ on $\fX$. Then there exists a unique coherent $\O_\fX$-submodule $\ul{\sF}'\subseteq\ul{\sF}$ such that $\ul{\sF}/\ul{\sF}'$ is $R$-flat and such that the given isomorphism $\ul{\sF}^\srig\cong\sF$ identifies $(\ul{\sF}')^\srig$ with $\sF'$.
\end{prop}
\begin{proof}
We may work locally on $\fX$ and thereby assume that $\fX$ is affine. Uniqueness of $\ul{\sF}'$ is a consequence of Lemma \ref{submodulecanmodellem}. Since $\sF'$ is coherent, there exists a leaflike covering $(X_i)_{i\in I}$ of $X$ such that $\sF'|_{X_i}$ is associated. We are in a situation where we need to relate a given leaflike covering of $X$ to a given affine formal $R$-model $\fX$ of $X$. It is thus natural to extend $(X_i)_{i\in I}$ to a treelike covering of $X$ and to consider a formal presentation of this treelike covering with respect to $\fX$.

Let $(X_j)_{j\in J}$ be a treelike covering of $X$ extending $(X_i)_{i\in I}$, where $J$ has a suitable rooted tree structure such that $I$ is identified with the set of leaves of $J$. By Lemma \ref{simpadlem}, we may assume that $(X_j)_{j\in J}$ is induced by a treelike formal covering with respect to $\fX$ in the sense of Definition \ref{formalprecovdefi}. That is, there exists the following data:
\begin{packed_enum}
\item a family $(\fX_j)_{j\in J}$ of affine flat formal $R$-schemes of ff type,
\item for each inner vertex $j\in J$ an admissible formal blowup $\ul{\beta}_j\colon\fX_j'\rightarrow\fX_j$ and
\item for each inner vertex $j\in J$ and for each child $l$ of $j$ an open immersion 
\[
\ul{\phi}_l\colon\fX_l\hookrightarrow\fX_j'
\]
\end{packed_enum}
such that $\fX_r=\fX$, where $r$ is the root of $J$, and such that for each inner vertex $j\in J$ and every child $l$ of $j$,
\[
\ul{\beta}_j\circ\ul{\phi}_l\colon\fX_l\rightarrow\fX_j
\]
defines $X_l$ as a semi-affinoid subdomain in $X_j$. 

Via pullback, we may define the restrictions $\ul{\sF}|_{\fX_j}$ and $\ul{\sF}|_{\fX_j'}$ for all vertices $j\in J$ and for all inner vertices $j\in J$ respectively. Here our notation should not divert the reader from keeping in mind that the morphisms $\ul{\beta}_j$ need \emph{not} be flat. Let $j$ be an inner vertex of $J$. Let us \emph{assume} that for each child $l$ of $j$, we are given a coherent submodule
\[
\ul{\sF}'_l\,\subseteq\,\ul{\sF}|_{\fX_l}
\]
such that $\ul{\sF}|_{\fX_l}/\ul{\sF}_l'$ is $R$-flat and such that $(\ul{\sF}'_l)^\urig=\sF'|_{X_l}$. Then
\[
(\ul{\phi}_{l_1l_2}^*\ul{\sF}'_{l_2})_K\,=\,(\ul{\phi}_{l_2l_1}^*\ul{\sF}'_{l_1})_K
\]
as submodules of $\ul{\sF}|_{\fX_{l_1l_2},K}$ for all children $l_1$ and $l_2$ of $j$, where $\fX_{l_1l_2}$, $\ul{\phi}_{l_1l_2}$ and $\ul{\phi}_{l_2l_1}$ are defined in terms of the cartesian diagram of open immersions
\[
\begin{diagram}
\fX_{l_1l_2}&\rTo^{\ul{\phi}_{l_1l_2}}&\fX_{l_2}\\
\dTo<{\ul{\phi}_{l_2l_1}}&&\dTo>{\ul{\phi}_{l_2}}\\
\fX_{l_1}&\rTo^{\ul{\phi}_{l_1}}&\fX'_j
\end{diagram}
\]
and where $\ul{\sF}|_{\fX_{l_1l_2}}$ is the restriction of $\ul{\sF}$ to $\fX_{l_1l_2}$. For example, Lemma \ref{submodulecanmodellem} shows that this assumption is satisfied when all children of $j$ are leaves of $J$. Of course, $\ul{\phi}_{l_1l_2}^*\ul{\sF}'_{l_2}$ and $\ul{\phi}_{l_2l_1}^*\ul{\sF}'_{l_1}$ are naturally submodules of 
\[
\ul{\sF}|_{\fX_{l_1l_2}}\,=\,\ul{\phi}_{l_1l_2}^*\ul{\sF}|_{\fX_{l_2}}\,=\,\ul{\phi}_{l_2l_1}^*\ul{\sF}|_{\fX_{l_1}}\quad,
\]
the open immersions $\ul{\phi}_{l_1l_2}$ and $\ul{\phi}_{l_2l_1}$ being flat. We have natural isomorphisms
\begin{eqnarray*}
\ul{\sF}|_{\fX_{l_1l_2}}/\ul{\phi}^*_{l_1l_2}\ul{\sF}'_{l_2}&=&\ul{\phi}_{l_1l_2}^*(\ul{\sF}|_{\fX_{l_2}}/\ul{\sF}'_{l_2})\\
\ul{\sF}|_{\fX_{l_1l_2}}/\ul{\phi}^*_{l_2l_1}\ul{\sF}'_{l_1}&=&\ul{\phi}_{l_2l_1}^*(\ul{\sF}|_{\fX_{l_1}}/\ul{\sF}'_{l_1})\quad,\\
\end{eqnarray*}
showing that the quotients on the left hand side are $R$-flat. By Lemma \ref{submodulecanmodellem}, we conclude that
\[
\ul{\phi}_{l_1l_2}^*\ul{\sF}'_{l_2}\,=\,\ul{\phi}_{l_2l_1}^*\ul{\sF}'_{l_1}
\]
as submodules of $\ul{\sF}|_{\fX_{l_1l_2}}$. Hence, the $\ul{\sF}'_l$ with $l\in\children(j)$ glue to a unique coherent submodule $\ul{\sG}_j\subseteq\sF|_{\fX'_j}$. That is,
\[
\ul{\sF}'_l\,=\,\ul{\phi}_l^*\ul{\sG}_j
\]
as submodules of $\ul{\sF}|_{\fX_l}$ for each child $l$ of $j$. The quotient $\ul{\sF}|_{\fX'_j}/\ul{\sG}_j$ is $R$-flat. Indeed, this condition can be checked locally on $\fX'_j$ with respect to the Zariski covering $(\ul{\phi}_l)_{l\in\children(j)}$, and
\[
\ul{\phi}_l^*(\ul{\sF}|_{\fX'_j}/\ul{\sG}_j)\,=\,\ul{\sF}|_{\fX_l}/\ul{\sF}'_l
\]
for all $l\in\children(j)$, the $\ul{\phi}_l$ being flat. By \cite{EGAIII1} 3.4.2, the functor $\ul{\beta}_{j,*}$ sends coherent $\O_{\fX_j'}$-modules to coherent $\O_{\fX_j}$-modules; hence $\ul{\beta}_{j,*}(\ul{\sF}|_{\fX'_j}/\ul{\sG}_j)$ is a co\-he\-rent $R$-flat $\O_{\fX_j}$-module. By definition, $\ul{\sF}|_{\fX'_j}=\ul{\beta}_j^*\ul{\sF}|_{\fX_j}$, so we have a natural homomorphism of coherent $\O_{\fX_j}$-modules
\[
\ul{\sF}|_{\fX_j}\rightarrow \ul{\beta}_{j,*}\ul{\sF}|_{\fX'_j}\rightarrow\ul{\beta}_{j,*}(\ul{\sF}|_{\fX_j'}/\ul{\sG}_j)\quad.
\]
We define $\ul{\sF}'_j$ to be its kernel, thus obtaining an exact sequence of coherent $\O_{\fX_j}$-modules
\[
0\rightarrow \ul{\sF}'_j\rightarrow\ul{\sF}|_{\fX_j}\rightarrow\ul{\beta}_{j,*}(\ul{\sF}|_{\fX_j'}/\ul{\sG}_j)
\]
which shows that $\ul{\sF}|_{\fX_j}/\ul{\sF}'_j$ is $R$-flat. There is a natural morphism $\ul{\sF}'_j\rightarrow\ul{\beta}_{j,*}\ul{\sG}_j$, and we claim that the induced morphism
\[
(\ul{\beta}_j^*\ul{\sF}_j')_K\rightarrow(\ul{\beta}_j^*\ul{\beta}_{j,*}\ul{\sG}_j)_K\rightarrow\ul{\sG}_{j,K}
\]
is an isomorphism. By Proposition \ref{blowupisoonsemiafffuncprop}, the second morphism is an isomorphism, so we must show that the first morphism is an isomorphism as well. Let $\ul{X}_j$ be the spectrum of the ring of global functions on $\fX_j$, and let $\ul{b}_j\colon\ul{X}_j'\rightarrow\ul{X}_j$ be the admissible blowup such that $\ul{\beta}_j=\ul{b}_j^\wedge$, a wedge indicating formal completion with respect to an ideal of definition of $\fX$, cf.\ \cite{EGAIII1} 5.1.4 and 5.4.1. Let $\ul{F}_j$, $\ul{F}_j'$ and $\ul{G}_j$ denote the algebraizations of $\ul{\sF}|_{\fX_j}$, $\ul{\sF}_j'$ and $\ul{\sG}_j$ respectively; then 
\[
\ul{\sF}|_{\fX'_j}=(\ul{b}_j^*\ul{F}_j)^\wedge\quad. 
\]
By \cite{EGAIII1} 4.1.5,
\[
\ul{\beta}_{j,*}(\ul{\sF}|_{\fX_j'}/\ul{\sG}_j)\,=\,(\ul{b}_{j,*}((\ul{b}^*_j\ul{F}_j)/\ul{G}_j)))^\wedge\quad,
\]
so we have a short exact sequence
\[
0\rightarrow\ul{F}'_j\rightarrow\ul{F}_j\rightarrow\ul{b}_{j,*}((\ul{b}^*_j\ul{F}_j)/\ul{G}_j))\quad.
\]
The functor $\cdot\otimes_RK$ yields a short exact sequence
\[
0\rightarrow\ul{F}'_{j,K}\rightarrow\ul{F}_{j,K}\rightarrow\ul{b}_{j,K,*}((\ul{b}^*_{j,K}\ul{F}_{j,K})/\ul{G}_{j,K}))\quad.
\]
Since $b_{j,K}$ is an isomorphism and, hence, flat, we obtain an induced short exact sequence
\[
0\rightarrow\ul{b}^*_{j,K}\ul{F}'_{j,K}\rightarrow\ul{b}^*_{j,K}\ul{F}_{j,K}\rightarrow\ul{b}^*_{j,K}\ul{b}_{j,K,*}((\ul{b}^*_{j,K}\ul{F}_{j,K})/\ul{G}_{j,K}))
\]
which shows that $\ul{b}^*_{j,K}\ul{F}'_{j,K}=\ul{G}_{j,K}$, the functor $\ul{b}^*_{j,K}\ul{b}_{j,K,*}$ being naturally identified with the identity functor. Hence, the natural morphism
\[
\ul{b}_j^*\ul{F}'_j\rightarrow\ul{b}_j^*\ul{b}_{j,*}\ul{G}_j
\]
becomes an isomorphism under $\cdot\otimes_RK$. That is, its kernel and cokernel are $\pi$-torsion. It follows that kernel and cokernel of the completed morphism
\[
\ul{\beta}_j^*\ul{\sF}_j'\rightarrow\ul{\beta}_j^*\ul{\beta}_{j,*}\ul{\sG}_j
\]
are $\pi$-torsion as well, which yields our claim. It follows that the coherent $X_j$-module $\ul{\sG}_j^\srig$ is associated to $\ul{\sF}_j'$.

Let us now prove the statement of the proposition by induction on the volume $v(J)$ of $J$. We may assume that $J$ has more than one vertex. Let $i$ be a leaf of $J$ whose path to the root has maximal length, and let $j$ be the parent of $i$. Then all $l\in\children(j)$ are leaves of $J$. By assumption, $\sF'|_{X_l}$ is associated for all $l\in\children(j)$, so by Lemma \ref{submodulecanmodellem} there exists, for each $l\in\children(j)$, a unique coherent $\O_{\fX_l}$-submodule $\ul{\sF}'_j\subseteq\ul{\sF}|_{\fX_j}$ such that the isomorphism 
\[
(\ul{\sF}|_{\fX_l})^\srig\,\cong\,\sF|_{X_l}
\]
induced by the given isomorphism $\ul{\sF}^\srig\cong\sF$ identifies $(\ul{\sF}'_l)^\srig$ with $\sF'|_{X_l}$ and such that $\ul{\sF}|_{\fX_l}/\ul{\sF}'_l$ is $R$-flat. By what we have shown so far, $\sF'|_{X_j}$ is associated to a unique coherent $\O_{\fX_j}$-submodule $\ul{\sF}'_j\subseteq\ul{\sF}|_{\fX_j}$ such that $\ul{\sF}|_{\fX_j}/\ul{\sF}'_j$ is $R$-flat. We may thus replace $\subtree(j)$ by $\{j\}$. By our induction hypothesis, it follows that $\sF'$ is associated to a coherent submodule $\ul{\sF}'\subseteq\ul{\sF}$ with the property that $\ul{\sF}/\ul{\sF'}$ is $R$-flat, as desired.
\end{proof}

\begin{cor}\label{frameembassoccor}Let us note some consequences of Proposition \ref{frameembprop}:
\begin{packed_enum}
\item A coherent submodule of an associated module on a semi-affinoid $K$-space is associated.
\item Coherent submodules and coherent quotients of strictly coherent modules are strictly coherent.
\item An associated module on a semi-affinoid $K$-space is strictly coherent.
\end{packed_enum}
\end{cor}
\begin{proof}
Let us first show ($i$). Let $X=\sSp A$ be a semi-affinoid $K$-space, let $\ul{A}\subseteq A$ be an $R$-lattice of ff type, and let $\sF'$ be a coherent submodule of an associated module $\tilde{M}$. Since $\tilde{M}$ admits a model $\ul{M}$ over $\Spf \ul{A}$, Proposition \ref{frameembprop} implies that $\sF'\cong(\ul{\sF}')^\srig$ for a coherent module $\ul{\sF}'$ on $\Spf\ul{A}$. Since coherent modules on affine formal schemes are associated, it follows that $\sF'$ is associated.


Let us prove statement ($ii$). Let $X$ be a uniformly rigid $K$-space, let $\sF$ be a strictly coherent $\O_X$-module and let $\sF'\subseteq\sF$ be a coherent submodule. For every semi-affinoid subspace  $U\subseteq X$, the restriction $\sF'|_U$ is a coherent submodule of $\sF|_U$, and $\sF|_U$ is associated by assumption on $\sF$. It follows from ($i$) that $\sF'|_U$ is associated; hence $\sF'$ is strictly coherent. Let now $\sF''$ be a coherent quotient of $\sF$. Then the kernel $\sF'$ of the projection $\sF\rightarrow \sF''$ is a coherent submodule of $\sF$ and, hence, strictly coherent by what we have seen so far. Let $U\subseteq X$ be a semi-affinoid subspace; then we have a short exact sequence
\[
0\rightarrow\sF'|_U\rightarrow\sF|_U\rightarrow\sF''|_U\rightarrow 0
\]
where the first two modules are associated. It follows from Lemma \ref{assocbasicpropertieslem} that $\sF''|_U$ is associated as well.

Finally, statement ($iii$) follows from statement ($ii$) because by Lemma \ref{assocbasicpropertieslem}, an associated module is a quotient of a frame.
\end{proof}

We do \emph{not} know whether the analog of Kiehl's Theorem \cite{KiehlAB} 1.2 holds in general, that is to say whether every coherent module on a semi-affinoid $K$-space is associated. If $\fX$ is a flat formal $R$-scheme of locally ff type and if $\ul{\sF}$ is a coherent $\O_\fX$-module, we do \emph{not} even know in general whether $\ul{\sF}^\srig$ is strictly coherent.

\subsection{Another approach to strictly coherent modules}

So far we have studied frame-embeddable modules by means of formal $R$-models. There is yet another instructive way to approach the subject; 
it yields a new and maybe more elementary proof of Corollary \ref{frameembassoccor} ($iii$). The present section will not be used in the sequel of this thesis, and the reader may choose to skip it. 


If $A$ is a semi-affionid $K$-algebra and if $M$ is a finite $A$-module admitting a frame-embedding $M\subseteq A^n$, then $M$ together with this embedding can be reconstructed from a restriction to a sufficiently large semi-affinoid subdomain whose ring of functions is actually \emph{affinoid}. This will be a consequence of the following Proposition. As always, it is instructive to keep in mind the case where $n=1$, that is, where $M$ is an ideal in $A$.

\begin{prop}\label{modulreconstructprop}
Let $A$ be a noetherian ring, and let $(A_i)_{i\in \N}$ be a family of flat noetherian $A$-algebras. For each $i\in\N$, we let $\phi_i^*\colon A\rightarrow A_i$ denote the structural flat ring homomorphism, and we write $\phi_i$ for the associated morphism of ring spectra. Let us assume that there exist homomorphisms $A_{i+1}\rightarrow A_i$ spanning commutative triangles with the $\phi_i^*$. Moreover, let us assume that the following holds:
\begin{packed_enum}
\item For all $i\in \N$, $\phi_i(\Max A_i)\subseteq\phi_{i+1}(\Max A_{i+1})\subseteq\Max A$,
\item $\Max A$ is the union of the $\phi_i(\Max A_i)$, and
\item for all $i\in\N$, $\phi_i$ induces isomorphisms of completed stalks in closed points.
\end{packed_enum}
Let $n$ be a natural number, and let $M\subseteq A^n$ be an $A$-module. There exists an $i_0\in\N$ such that for all $i\geq i_0$,
\[
M\,=\,(M A_i)\cap A^n\quad,
\]
where $M A_i\mathrel{\mathop:}=M_i\otimes_AA_i$ which is a submodule of $A_i^n$ by flatness of $A_i$ over $A$, and where $(M A_i)\cap A^n$ is the preimage of $M A_i$ under $(\phi_i^*)^n\colon A^n\rightarrow A_i^n$.
\end{prop}
\begin{proof}
Let us write $X=\Spec A$, and let $Y$ be the disjoint union of the $Y_i\mathrel{\mathop:}=\Spec A_i$, viewed as an $X$-scheme via the coproduct $\phi$ of the $\phi_i$. Since $\phi$ is flat and since its image contains all closed points of $X$, it follows from [EGA IV$_2$] 2.2.2 that $\phi$ is in fact surjective. Since $\phi_i(Y_i)\subseteq\phi_{i+1}(Y_{i+1})$ for all $i\in\N$, there exists an index $i_0$ such that $\phi_i(Y_i)$ contains all \emph{associated points} of $A^n/M$ for all $i\geq i_0$. We claim that $M=(M A_i)\cap A^n$ for all $i\geq i_0$. Since $M$ is contained in $(M A_i)\cap A^n$, we must show that for all $i\geq i_0$, the natural homomorphism of $A$-modules
\[
A^n/M\rightarrow A_i^n/(MA_i)\quad(*)
\]
is injective. Let us idenfity $A_i^n/(MA_i)$ with $(A^n/M)\otimes_AA_i$. Let $\n\subseteq A_i$ be a maximal ideal with  $\phi_i^*$-pullback $\m\subseteq A$; then the completed stalk of $A^n_i/(MA_i)$ in $\n$ is
\begin{eqnarray*}
(A^n_i/(MA_i))_\n^\wedge&\cong&A^n_i/(MA_i)\otimes_{A_i}\hat{A}_{i,\n}\\
&\cong&(A^n/M)\otimes_AA_i\otimes_{A_i}\hat{A}_{i,\n}\\
&\cong&(A^n/M)\otimes_A\hat{A}_{\m}\\
&\cong&(A^n/M)_\m^\wedge\quad.
\end{eqnarray*}
Let $f\in A^n/M$ be an element whose image in $A_i^n/(MA_i)$ is zero. By Krull's Intersection Theorem, $f_\m=0$ in $(A^n/M)_\m$ for every maximal ideal $\m\subseteq A$ that lies in the image of $\phi_i$. Let $\p\subseteq A$ be an associated prime ideal of $A^n/M$, let $\q\subseteq A_i$ be a prime above $\p$, and let $\n\subseteq A_i$ be a maximal specialization of $\p$; then $\m\mathrel{\mathop:}=\phi^*\n$ is a maximal specialization of $\p$, and $f$ is zero in $(A^n/M)_\m$. Since $(A^n/M)_\p$ is a localization of $(A^n/M)_\m$, we conclude that $f$ vanishes in $(A^n/M)_\p$. Thus $f\in A^n/M$ vanishes in every associated prime ideal of $A^n/M$, so it follows that $f$ must be identically zero. It follows that $(*)$ is injective, as desired.
\end{proof}

\begin{cor}\label{modulreconstructcor}
Let $A$ be a semi-affinoid $K$-algebra, let $\ul{A}$ be an $R$-lattice of ff type for $A$ with an ideal of definition $\a$, and for each $i\in\N$ let 
\[
\phi_i^*\colon A\rightarrow A_i\,\mathrel{\mathop:}=\,\ul{A}\langle\a^{i+1}/\pi\rangle\otimes_RK
\]
be the natural $K$-homomorphism. Let $n$ be a positive natural number, and let $M\subseteq A^n$ be a submodule. There exists an $i_0\in\N$ such that for all $i\geq i_0$,
\[
M\,=\,(M\cdot A_i)\cap A^n\quad.
\]
This applies in particular to the case where $n=1$, that is, where $M$ is an ideal in $A$. If $M$ is the zero ideal, we see in particular that $\phi^*_i$ is injective for $i\gg 0$.
\end{cor}

Indeed, the homomorphisms $\psi_i$ satisfy the conditions of Proposition \ref{modulreconstructprop}, cf. our discussion of Berthelot's construction in Section \ref{berthsec}.

We can now give an alternative proof of Corollary \ref{frameembassoccor}:

\begin{cor}
Let $X$ be a semi-affinoid $K$-space, and let $\sF$ be a frame-em\-bed\-da\-ble $\O_X$-module. Then $\sF$ is associated.
\end{cor}
\begin{proof}
Let us choose a frame embedding $\sF\subseteq\O_X^n$ of the coherent $\O_X$-module $\sF$. Let $(X_i)_{i\in I}$ be a leaflike covering of $X$ such that $\sF|_{X_i}$ is associated. Let us write $X_i=\sSp A_i$, and let $M_i$ be a finite $A_i$-module such that $\sF|_{X_i}$ is associated to $M_i$, for $i\in I$. By Corollary \ref{modulreconstructcor}, there exists  a semi-affinoid subdomain  $U\subseteq X$, $U=\sSp B$, where $B$ is an affinoid $K$-algebra and such that for each $i\in I$ the frame-embedding $M_i\subseteq A_i^n$ can be recovered from its restriction to $U\cap X_i$. The constituents of a leaflike covering are retrocompact, and by Lemma \ref{veryspecialinafflem} the retrocompact semi-affinoid subdomains of $U$ affinoid. Hence Kiehl's Theorem \cite{KiehlAB} 1.2 for affinoid spaces shows that the frame-embeddings $(M_i\subseteq A_i^n)|_{U\cap X_i}$ glue to a frame embedding $N\subseteq B^n$. We set $M\mathrel{\mathop:}=N\cap A^n$, and we claim that $\sF\subseteq\O_X^n$ is associated to $M\subseteq A^n$. This can be checked locally on the $X_i$ and now follows from our choice of $U$.
\end{proof}

Arguing as in the proof of Corollary \ref{frameembassoccor}, we see that quotients of frames are associated, which implies that associated modules are strictly coherent. In particular, coherent quotients of associated modules are associated. By Lemma \ref{assocbasicpropertieslem}, it follows that coherent submodules of associated modules are associated, so we recover the full statement of Corollary \ref{frameembassoccor} without using Proposition \ref{frameembprop}.

		\subsection{Closed uniformly rigid subspaces}
			\begin{defi}\label{closedimdefi}
A morphism of uniformly rigid $K$-spaces $\phi\colon Y\rightarrow X$ is called a \emph{closed immersion}\index{closed immersion!of uniformly rigid spaces} if there exists an admissible semi-affinoid covering $(X_i)_{i\in I}$ of $X$ such that for each $i\in I$, the restriction $\phi^{-1}(X_i)\rightarrow X_i$ of $\phi$ is a closed immersion of semi-affinoid $K$-spaces in the sense of Definition \ref{closedimsemaffdefin}.
\end{defi}

We easily see that closed immersions are injective on the level of physical points. 

\begin{lem}\label{closedimstrictcohlem}
Let $\phi\colon Y\rightarrow X$ be a closed immersion of uniformly rigid $K$-spaces. Then $\phi^\sharp\colon \O_X\rightarrow\phi_*\O_Y$ is an epimorphism of sheaves. Moreover, the $\O_X$-modules $\phi_*\O_Y$ and $\ker\phi^\sharp$ are strictly coherent. 
\end{lem}
\begin{proof}
Since $\O_X$ is strictly coherent, it suffices by Corollary \ref{frameembassoccor} ($ii$) to show that $\phi^\sharp$ is an epimorphism and that both $\ker\phi^\sharp$ and $\phi^*\O_Y$ are coherent. Considering an admissible semi-affinoid covering $(X_i)_{i\in I}$ of $X$ such that for all $i\in I$, the restriction $\phi^{-1}(X_i)\rightarrow X_i$ of $\phi$ is a closed immersion of semi-affinoid $K$-spaces, we reduce to the case where both $X$ and $Y$ are semi-affinoid and where $\phi$ corresponds to a surjective homomorphism of semi-affinoid $K$-algebras. Now the desired statements follow from Lemma \ref{assocbasicpropertieslem}.
\end{proof}

\begin{prop}\label{closedimprop}
Let $\phi\colon Y\rightarrow X$ be a morphism of uniformly rigid $K$-spaces. Then the following are equivalent:
\begin{packed_enum}
\item $\phi$ is a closed immersion.
\item For each semi-affinoid subspace $U\subseteq X$, the restriction $\phi^{-1}(U)\rightarrow U$ is a closed immersion of semi-affinoid $K$-spaces in the sense of Definition \ref{closedimsemaffdefin}.
\end{packed_enum}
\end{prop}
\begin{proof}
The implication ($ii$)$\Rightarrow$($i$) is trivial, the semi-affinoid subspaces forming a basis for the G-topology on $X$. Let us assume that ($i$) holds, let $\sI$ denote the kernel of $\phi^\sharp$, and let $U\subseteq X$ be a semi-affinoid subspace; then $\phi$ induces a short exact sequence
\[
0\rightarrow\sI|_U\rightarrow\O_U\rightarrow\phi^*\O_Y|_U\rightarrow 0\quad.
\]
Let $A$ denote the ring of functions on $U$. By Lemma \ref{closedimstrictcohlem}, $\sI$ and $\phi_*\O_Y$ are strictly coherent; hence the above short exact sequence is associated to a short exact sequence of $A$-modules
\[
0\rightarrow I\rightarrow A\rightarrow B\rightarrow 0\quad.
\]
Since morphisms from uniformly rigid $K$-spaces to semi-affinoid $K$-spaces correspond to $K$-homomorphisms of rings of global functions, we can now copy the proof of \cite{BGR} 9.4.4/1 to see that the restriction $\phi^{-1}(U)\rightarrow U$ of $\phi$ is associated to the projection $A\rightarrow B$.
\end{proof}

In particular, a morphism of semi-affinoid $K$-spaces is a closed immersion in the sense of Definition \ref{closedimdefi} if and only if it is a closed immersion of semi-affinoid $K$-spaces in the sense of Definition \ref{closedimsemaffdefin}. We can now define a \emph{closed uniformly rigid subspace}\index{subspace!closed uniformly rigid} as an equivalence class of closed immersions, in the usual way. By standard glueing arguments, we see that the closed uniformly rigid subspaces of a uniformly rigid $K$-space $X$ correspond to the coherent $\O_X$-ideals. We easily see that closed immersions of uniformly rigid $K$-spaces are preserved under base change.




It is clear that closed immersions of formal $R$-schemes of locally ff type induce closed immersions on uniformly rigid generic fibers. Conversely, given a uniformly rigid $K$-space $X$ together with an $R$-model of locally ff type $\fX$ and a closed uniformly rigid subspace $V\subseteq X$, there exists a unique $R$-flat closed formal subscheme $\fV\subseteq\fX$ such that the given isomorphism $\fX^\srig\cong X$ identifies $\fV^\srig$ with $V$. Indeed, this is an immediate consequence of Proposition \ref{frameembprop}. We say that $\fV$ is the \emph{schematic closure}\index{schematic closure!of a closed subspace} of $V$ in $\fX$.

%



The comparison functors studied in Section \ref{compfuncsec} preserve closed immersions. Indeed, this can be verified in the semi-affinoid and affinoid situations respectively. In the case of the functor $\sr$, there is nothing to show. In the case of the forgetful functor $\r$, the statement follows by looking at schematic closures and using the fact that Berthelot's construction preserves closed immersions, cf.\ \cite{dJ} 7.2.4 (e).

\subsubsection{Separated uniformly rigid spaces}\label{separatedsection}

As usual, a morphism $\phi\colon Y\rightarrow X$ of uniformly rigid $K$-spaces is called \emph{separated} if its \emph{diagonal morphism} 
\[
\Delta_\phi\colon Y\rightarrow Y\times_XY
\]
is a closed immersion. A uniformly rigid $K$-space $X$ is called separated\index{uniformly rigid space!separated} if its structural morphism $X\rightarrow\sSp K$ is separated. If $X$ is a uniformly rigid $K$-space, we let $\Delta_X$ denote the diagonal of its structural morphism.

Semi-affinoid $K$-spaces are visibly separated. Moreover, uniformly rigid generic fibers of separated morphisms of formal $R$-schemes of locally ff type are separated, the functor $\srig$ preserving fibered products and closed immersions.

\begin{lem}\label{sepintlem}
Let $X$ be a separated uniformly rigid $K$-space. The intersection of two semi-affinoid subspaces in $X$ is a semi-affinoid subspace in $X$.
\end{lem}
\begin{proof}
Let $U$ and $V$ be semi-affinoid subspaces in $X$. We easily see, using points with values in finite field extensions of $K$, that $U\cap V$ is the $\Delta_X$-preimage of $U\times_{\sSp K}V$ which is a semi-affinoid subspace of $X\times_{\sSp K}X$. Since $\Delta_X$ is a closed immersion by assumption on $X$, it follows from Proposition \ref{closedimprop} that $U\cap V$ is a semi-affinoid subspace of $X$.
\end{proof}

\begin{cor}\label{calcthecohcor}
Let $X$ be a separated uniformly rigid $K$-space, and let $\sF$ be a coherent $\O_X$-module. Then the natural morphism
\[
\check{H}^q(X,\sF)\overset{\sim}{\rightarrow}H^q(X,\sF)
\]
is an isomorphism for all $q\geq 0$.
\end{cor}
\begin{proof}
Let $S$ denote the set of semi-affinoid subspaces $U$ in $X$ with the property that $\sF|_U$ is associated. By Lemma \ref{sepintlem}, this set is stable under the formation of intersections. It is clearly a basis for the G-topology on $X$, and $\check{H}^q(U,\sF)=0$ for any $U$ in $S$ and any $q\geq 0$ by Corollary \ref{moduleacythmcor}. We conclude by the usual \v{C}ech spectral sequence argument.
\end{proof}

Now and in the following, fibered products without indication of the base are understood over $\sSp K$ or $\Spf R$, depending on the context.

If $X$ is a separated uniformly rigid $K$-space and if $\phi\colon Y\rightarrow X$ is a morphism of uniformly rigid $K$-spaces, then the graph $\Gamma_\phi\colon Y\rightarrow Y\times X$ of $\phi$ is a closed immersion since it is obtained from $\Delta_X$ via pullback. In particular, if $\fX$ and $\fY$ are $R$-models of locally ff type for $X$ and $Y$ respectively, the schematic closure of $\Gamma_\phi$ in $\fY\times\fX$ is well-defined.

Since the comparison functors studied in Section \ref{compfuncsec} preserve fibered products and closed immersions, they preserve the separatedness property.
			
	\section{Galois descent for uniformly rigid spaces}\label{galdescsec}
			We outline some aspects of Galois descent theory for uniformly rigid $K$-spaces. Let $K'/K$ be a finite Galois field extension with Galois group $\Gamma$, and let $R'$ denote the discrete valuation ring of $K'$.

\begin{prop}\label{semaffdescprop}
Let $A'$ be a semi-affinoid $K'$-algebra together with a $\Gamma$-action extending the tautological action on $K'$. Then $A\mathrel{\mathop:}=(A')^\Gamma$ is a semi-affinoid $K$-algebra, and the natural $K'$-homomorphism
\[
A\otimes_KK'\rightarrow A'
\]
is a $\Gamma$-equivariant isomorphism.
\end{prop}
\begin{proof}
The last statement follows by means of standard Galois descent arguments; it thus suffices to see that $A$ is semi-affinoid over $K$. Let $\ul{A}'$ be an $R'$-lattice of ff type for $A'$. After replacing $\ul{A}'$ by the sum of its $\Gamma$-translates, we may assume that $\ul{A}'$ is $\Gamma$-stable. Since $R'/R$ is finite, we may consider $\ul{A}'$ as an $R$-algebra of ff type. Let 
\[
((s_1,\ldots,s_m),(t_1,\ldots,t_n))
\]
be a formal generating system of $\ul{A}'$ over $R'$, and let $r$ denote the cardinality of $\Gamma$. For each $1\leq i\leq m$, let $u_{i,1},\ldots,u_{i,r}$ be the elementary symmetric polynomials in the $\gamma(s_i)$, where $\gamma$ varies in $\Gamma$. Similarly, for each $1\leq j\leq n$ let $v_{j,1},\ldots,v_{j,r}$ be the elementary symmetric polynomials in the $\gamma(t_j)$, where again $\gamma$ varies in $\Gamma$. Then the $u_{i,k}\in\ul{A}'$ are topologically nilpotent. Let $\ul{B}\subseteq\ul{A}'$ denote the $R$-algebra of ff type generated by the $u_{i,k}$ and the $v_{j,k}$. That is, $\ul{B}$ is the image of the $R$-homomorphism
\[
R[[(X_{ik})_{1\leq i\leq m\,,\,1\leq k\leq r}]]\langle (Y_{jk})_{1\leq j\leq n\,,\,1\leq k\leq r}\rangle\rightarrow\ul{A}'
\]
sending $X_{ik}$ to $u_{ik}$ and $Y_{jk}$ to $v_{jk}$. Then $\ul{B}$ is $\Gamma$-stable, and hence $B\mathrel{\mathop:}=\ul{B}\otimes_RK$ is a $K$-subalgebra of $A$. Since $B$ is semi-affinoid, it suffices to see that the inclusion $B\subseteq A$ is finite. This may be verified after the base change $\cdot\otimes_KK'$; hence it suffices to see that $\ul{A}'$ is finite over $\ul{B}'\mathrel{\mathop:}=\ul{B}\otimes_RR'$, which is the $R'$-subalgebra of ff type of $\ul{A}'$ generated by the $u_{ik}$ and the $v_{jk}$. Let $\b'$ denote the ideal of definition of $\ul{B}'$ induced by the formal presentation above. By the formal version of Nakayama's Lemma, cf.\ \cite{Eis} Ex. 7.2, it suffices to see that $\ul{A}'/\b'\ul{A}'$ is finite over $\ul{B}'/\b$. The continuous homomorphism $\ul{B}'\subseteq\ul{A}'$ is \emph{adic} since each $s_i\in\ul{A}'$ satisfies an integral equation over $\ul{B}'$ with coefficients in $\b'$. Hence we have reduced to an algebraic situation, and it suffices to remark that the $t_j$ are likewise integral over $\ul{B}'$. 
\end{proof}

\begin{prop}\label{descffprop}
The functor $X\mapsto X\otimes_KK'$ from the category of uniformly rigid $K$-spaces to the category of uniformly rigid $K'$-spaces together with a semi-linear $\Gamma$-action is faithful; that is, the natural map
\[
\Hom_K(Y,X)\rightarrow\Hom_{K'}(Y_{K'},X_{K'})^\Gamma
\]
is injective for all uniformly rigid $K$-spaces $X$ and $Y$. If $X$ is semi-affinoid, this map is bijective.
\end{prop}
\begin{proof}
We use a subscript $K'$ to indicate the image under the functor $\cdot\otimes_KK'$. Let $X$ and $Y$ be uniformly rigid $K$-spaces, and let $\phi,\psi\colon Y\rightarrow X$ be morphisms such that $\phi_{K'}=\psi_{K'}$. The projection $Y_{K'}\rightarrow Y$ being surjective, $\phi$ and $\psi$ must coincide on physical points; we may thus assume that $X$ and $Y$ are semi-affinoid. Let $V\subseteq Y$ be the closed uniformly rigid coincidence subspace of $\phi$ and $\psi$, that is, the pullback of the diagonal of $X\times X$ under the morphism $(\phi,\psi)$. Its formation commutes with base change, which implies that $V_{K'}=Y_{K'}$; we conclude that $V=Y$, which means that $\phi=\psi$. 

Conversely, let us assume that $X$ is semi-affinoid, and let $\phi'\colon Y_{K'}\rightarrow X_{K'}$ be a morphism that is equivariant with respect to the $\Gamma$-actions. By what we have shown so far, we may assume that $Y$ is semi-affinoid as well. Since morphisms of semi-affinoid $K$-spaces correspond to homomorphisms of rings of global functions, finite Galois descent for $K$-algebras shows that $\phi'$ descends to a morphism $\phi\colon Y\rightarrow X$.
\end{proof}

\begin{cor}\label{dessaffeffcor}
Let $X'$ be a semi-affinoid $K'$-space with a semi-linear $\Gamma$-action. Then there exists a semi-affinoid $K$-space $X$ together with a $\Gamma$-equivariant isomorphism
\[
\phi\colon X'\overset{\sim}{\rightarrow}X\otimes_KK'\quad,
\]
and the pair $(X,\phi)$ is unique up to unique isomorphism.
\end{cor}
\begin{proof}
By Proposition \ref{semaffdescprop}, there exists a semi-affinoid $K$-space $X$ such that $X'$ is $\Gamma$-equivariantly isomorphic to $X\otimes_KK'$. If $(X_1,\phi_1)$ and $(X_2,\phi_2)$ are two pairs as in the statement of the Proposition, then Proposition \ref{descffprop} shows that there exists a unique $K$-isomorphism $X_1\cong X_2$ such that the induced $K'$-isomorphism $X_1\otimes_KK'\cong X_2\otimes_KK'$ spans a commutative triangle together with $\phi_1$ and $\phi_2$.
\end{proof}

We do not know whether a semi-affinoid $K'$-space with a semi-linear $\Gamma$-action might possibly descend to a uniformly rigid $K$-space which is not semi-affinoid.

The category of semi-affinoid $K$-spaces being closed under the formation of fibered products, a semi-affinoid $K'$-group equipped with a semi-linear $\Gamma$-action descends uniquely to a semi-affinoid $K$-group.

Let $X'$ be a semi-affinoid $K'$-space equipped with a semi-linear $\Gamma$-action, and let $U'\subseteq X'$ be a $\Gamma$-stable semi-affinoid subdomain. By what we have seen so far, $U'$ descends to a semi-affinoid $K$-space $U$ together with a morphism $U\rightarrow X$ defining $U$ as a semi-affinoid \emph{pre}-subdomain in $X$. However, we do not know whether $U\subseteq X$ is admissible open, since we are lacking in a semi-affinoid version of the Gerritzen-Grauert Theorem. 


Similarly, if $X$ and $Y$ are global uniformly rigid $K$-spaces and if $\phi'\colon Y\otimes_KK'\rightarrow X\otimes_KK'$ is a $\Gamma$-equivariant morphism, we do not know whether the induced map on physical points $|\phi|\colon|Y|\rightarrow|X|$ is continuous: If $U\subseteq X$ is admissible open, then $\phi^{-1}(U)\otimes_KK'$ is admissible open in $Y\otimes_KK'$, but we do not know whether this implies the admissibility of $\phi^{-1}(U)$ in $X$. 


%



\chapter{Formal Néron models}\label{fnmchap}

	\section{Complements on formal geometry, Part II}\label{compl2sec}
	
		\subsection{Faithful flatness}\label{faithflatsec}
			

We begin by gathering some standard facts on flatness and faithful flatness in locally noetherian formal geometry. A morphism $\phi\colon\fY\rightarrow\fX$ of locally noetherian formal schemes is called \emph{flat} at a point $y$ of $\fY$ if the induced local homomorphism of stalks is flat. It is called flat if it is flat in every point of $\fY$ and it is called \emph{faithfully flat}\index{morphism!flat}\index{morphism!faithfully flat} if it is flat and surjective.

If $\fY$ and $\fX$ are affine, then $\phi$ is flat if and only if its underlying homomorphism of rings of global sections is flat. This is easily seen by considering completed stalks and by using faithfully flat descent of flatness, cf.\ \cite{AJL} 7.1.1. The corresponding statement for \emph{faithful} flatness is \emph{not} true unless $\phi$ is adic. For example, if $\fX=\Spf A$ is affine and if $\phi$ is the completion of $\fX$ along an ideal $I\subseteq A$ such that $A$ is $I$-adically complete and such that $I$ strictly contains the biggest ideal of definition of $A$, then $\phi$ induces an isomorphic on rings of global sections, but $\phi$ is not surjective.

Flatness and faithful flatness are preserved under base change of locally ff type:

\begin{prop}\label{flatnessbasechangeprop}
Let us consider a cartesian diagram
\[
\begin{diagram}
\fY\times_\fX\fX'&\rTo^{\phi'}&\fX'\\
\dTo<{\psi'}&&\dTo>\psi\\
\fY&\rTo^\phi&\fX
\end{diagram}
\]
of locally noetherian formal schemes, and let us assume that $\phi$ is (faithfully) flat. If  $\phi$ is adic or if $\psi$ is of locally ff type, then $\phi'$ is (faithfully) flat as well.
\end{prop}

\begin{proof}
If $\phi$ is flat, then $\phi'$ is flat by the Bourbaki Flatness Criterion; cf.\ \cite{AJL} 7.1. By \cite{EGAIn} 2.6.1 (ii), surjectivity of morphisms of schemes is stable under base change; hence we obtain the desired statement on faithful flatness.
\end{proof}

An adic morphism of locally noetherian formal schemes is an isomorphism if and only if it is an isomorphism after faithfully flat base change of locally ff type:

\begin{prop}\label{isodescprop}
Let us consider a cartesian diagram
\[
\begin{diagram}
\fY\times_\fX\fX'&\rTo^{\phi'}&\fX'\\
\dTo<{\psi'}&&\dTo>\psi\\
\fY&\rTo^\phi&\fX
\end{diagram}
\]
of locally noetherian formal schemes, where $\fX$ is quasi-compact, where $\phi$ is adic and separated and where $\psi$ is faithfully flat and of ff type. If $\phi'$ is an isomorphism, then $\phi$ is an isomorphism as well.
\end{prop}

\begin{proof}
Since $|\phi'|$ is bijective and since the natural map $|\fY\times_\fX\fX'|\rightarrow|\fY|\times_{|\fX|}|\fX'|$ is surjective, the projection
\[
|\fY|\times_{|\fX|}|\fX'|\rightarrow|\fX'|
\]
must be bijective. Since $|\psi|$ is surjective, it follows that $|\phi|$ is bijective. By Zariski's Main Theorem \cite{EGAIV4} 18.12.13, it follows that $\phi$ restricts to a finite morphism on subschemes of definition. By \cite{EGAIII1} 4.8.1, it follows that $\phi$ is finite. To prove that $\phi$ is an isomorphism, we may thus assume that $\fX$ and $\fY$ are affine and that $\phi$ is associated to a finite ring homomorphism. In the finite case, complete tensor products agree with ordinary tensor products, so we conclude by means of descent theory for schemes that $\phi$ is indeed an isomorphism.
\end{proof}
	
		\subsection{Higher direct images and flat base change}\label{hbcsec}			
			In the case of schemes, it is well known -- and easy to prove -- that the formation of higher direct images of a quasi-coherent module under a quasi-compact and quasi-separated morphism commutes with flat base change. The corresponding statement for the higher direct images of a coherent $\O_\fY$-module $\sG$ under a possibly non-algebraic proper morphism $\phi\colon\fY\rightarrow\fX$ of locally noetherian formal schemes is less obvious. Since we were unable to find a reference, we decided to provide the proof of the formal flat base change theorem. We derive it from \cite{EGAIII1} 3.4.4, which gives very precise information on how -- in the case where $\fX$ is affine, say -- the cohomology groups $H^q(\fY,\sG)$ relate to the cohomology groups $H^q(\fY,\sG_n)$, where $\sG_n$ denotes the reduction of $\sG$ modulo the $(n+1)$-st power of a fixed ideal of definition of $\fX$.

We begin by collecting some elementary facts on projective systems.

Let $A$ be a ring, let $(M_\alpha)_{\alpha\in\N}$ be a system of $A$-modules, and let $\phi_\alpha\colon M_{\alpha+1}\rightarrow M_\alpha$ be a system of $A$-homomorphisms. For natural numbers $\beta\geq \alpha$, we write
\[
\phi_{\alpha\beta}\,\mathrel{\mathop:}=\,\phi_\alpha\circ\cdots\circ\phi_{\beta-1}\colon M_{\beta}\rightarrow M_\alpha\quad.
\]
Let us recall from [EGAIII]$_1$ 0.13.1.1 that the projective system $(M_\alpha,\phi_\alpha)_\alpha$ is said to satisfy the Mittag-Leffler (ML) condition if for each $\alpha\in\N$, there exists some $\beta\geq \alpha$ such that for all $\gamma\geq\beta$,
\[
\phi_{\alpha\gamma}(M_\gamma)=\phi_{\alpha\beta}(M_\beta)\quad.
\]
We say that the system $(M_\alpha,\phi_\alpha)_\alpha$ is Artin-Rees (AR) null if there exists some $\beta\in\N$ such that $\phi_{\alpha,\alpha+\beta}=0$ for all $\alpha\in\N$. This property is preserved under any base change $\cdot\otimes_AA'$. Moreover, if $(M_\alpha,\phi_\alpha)_\alpha$ is (AR) null, then $(M_\alpha,\phi_\alpha)_\alpha$ clearly satisfies (ML), and $\varprojlim (M_\alpha,\phi_\alpha)_\alpha=0$.

\begin{lem}\label{ml1lem}
If $(M_\alpha,\phi_\alpha)_\alpha$ satisfies \textup{(ML)} and has trivial limit, then for each $\alpha\in\N$, there exists some $\beta\geq\alpha$ such that $\phi_{\alpha\beta}=0$.
\end{lem}
\begin{proof}
For each $\alpha\in\N$, we set 
\[
M'_\alpha\,\mathrel{\mathop:}=\,\bigcap_{\beta\geq\alpha}\phi_{\alpha\beta}(M_\beta)\quad;
\]
the submodule $M'_\alpha\subseteq M_\alpha$ is called the submodule of universal images. The homomorphisms $\phi_\alpha$ restrict to a projective system $(M'_\alpha,\phi'_\alpha)_\alpha$, and $\varprojlim (M'_\alpha,\phi'_\alpha)_\alpha=\varprojlim (M_\alpha,\phi_\alpha)_\alpha=0$. Since $(M_\alpha,\phi_\alpha)_\alpha$ satisfies (ML), the morphisms $\phi'_\alpha$ are surjective, cf.\ [EGAIII]$_1$ 0.13.1.2. It follows that all $M'_\alpha$ must be zero. Moreover, since $(M_\alpha,\phi_\alpha)_\alpha$ satisfies (ML), we see that for each $\alpha$ there exists some $\beta\geq\alpha$ such that $M'_\alpha=\phi_{\alpha\beta}(M_\beta)$. Hence, for each $\alpha$ there exists some $\beta\geq\alpha$ such that $\phi_{\alpha\beta}=0$.
\end{proof}

In particular, in this case
\[
\plim(M_\alpha\otimes_AA',\phi_\alpha\otimes_AA')_\alpha=0
\]
for any $A$-algebra $A'$.


Let $\a\subseteq A$ be an ideal. We recall from \cite{EGAIII1} 0.13.7.7 that a filtration $(N_\alpha)_{\alpha\in\N}$ of an $A$-module $N$ is called $\a$-\emph{good} if $\a N_n\subseteq N_{n+1}$ for all $n$ and if equality holds for all $n$ greater than some $n_0$.

\begin{lem}\label{ml2lem}
Let $A$ be a noetherian adic ring, let $\a$ be an ideal of definition of $A$, let $M$ be a finite $A$-module, let $(M_\alpha,\phi_\alpha)_\alpha$ be a projective system of $A$-modules, and let
\[
(\rho_\alpha\colon M\rightarrow M_\alpha)_{\alpha\in\N}
\]
be a system of $A$-module homomorphisms that is compatible with the transition homomorphisms $\phi_\alpha$. For $\alpha\in\N$, we let $N_\alpha$ denote the kernel of $\rho_\alpha$. Let us assume that the following holds:
\begin{packed_enum}
\item For each $\alpha\in\N$, multiplication by elements in $\a^{\alpha+1}$ is trivial on $M_\alpha$.
\item The filtration $(N_\alpha)_{\alpha\in\N}$ of $M$ is $\a$-good.
\item The system $(M_\alpha,\phi_\alpha)_{\alpha\in \N}$ satisfies the \textup{(ML)} condition.
\item The morphisms $\rho_\alpha$ induce an isomorphism $M\overset{\sim}{\rightarrow}\varprojlim (M_\alpha,\phi_\alpha)_\alpha$.
\end{packed_enum}
Then for any noetherian adic flat topological $A$-algebra $A'$, the morphisms $\rho_\alpha\otimes_AA'$ induce an isomorphism
\[
M\otimes_AA'\rightarrow\varprojlim(M_\alpha\otimes_AA')\quad.
\]
\end{lem}

We emphasize that the topology on $A'$ is \emph{not} required to be induced by the topology on $A$.

\begin{proof}
For each $\alpha\in\N$, the $A$-homomorphism $\rho_\alpha$ induces an $A$-homomorphism
\[
\tilde{\rho}_\alpha\colon M/\a^{\alpha+1}M\rightarrow M_\alpha\quad.
\]
Let $\tilde{N}_\alpha$ and $Q_\alpha$ denote its kernel and its cokernel respectively; we obtain a projective system of exact sequences
\[
0\rightarrow\tilde{N}_\alpha\rightarrow M/\a^{\alpha+1}M\rightarrow M_\alpha\rightarrow Q_\alpha\rightarrow 0\quad.
\]
Since $(N_\alpha)_{\alpha\in\N}$ is $\a$-good, there exists some $\beta\in\N$ such that $\a N_\gamma=N_{\gamma+1}$ for all $\gamma\geq\beta$, which implies that $N_{\alpha+\beta}=\a^\alpha N_{\beta}$ for all $\alpha\in\N$. It follows that all transition maps
\[
\tilde{N}_{\alpha+1+\beta}\rightarrow\tilde{N}_\alpha
\]
are trivial; that is, $(\tilde{N}_\alpha)_\alpha$ is (AR) null. The above exact sequences decompose into two projective systems of short exact sequences
\[
0\rightarrow\tilde{N}_\alpha\rightarrow M/\a^{\alpha+1}M\rightarrow Q_\alpha'\rightarrow 0
\]
and
\[
0\rightarrow Q_\alpha'\rightarrow M_\alpha\rightarrow Q_\alpha\rightarrow 0\quad,
\]
where $(Q_\alpha')_\alpha$ visibly satisfies (ML) and where the system $(\tilde{N}_\alpha)_\alpha$ satisfies (ML) since it is (AR) null. Hence, we may pass to the limit without losing exactness; reassembling the resulting two short exact sequences, we obtain an exact sequence
\[
0\rightarrow\varprojlim\tilde{N}_\alpha\rightarrow \varprojlim (M/\a^{\alpha+1}M)\rightarrow \varprojlim M_\alpha\rightarrow \varprojlim Q_\alpha\rightarrow 0\quad.
\]
Since $M$ is $\a$-adically complete and since $\varprojlim\rho_\alpha\colon M\rightarrow\varprojlim M_\alpha$ is an isomorphism by assumption, it follows that $\varprojlim Q_\alpha=0$. Since $(M_\alpha)_\alpha$ satisfies (ML), the same holds for $(Q_\alpha)_\alpha$, cf.\ \cite{EGAIII1} 0.13.2.1 (i); it follows that $(Q_\alpha)_\alpha$ satisfies the conditions of Lemma \ref{ml1lem}. Since $A'$ is flat over $A$, we have exact sequences
\[
0\rightarrow\tilde{N}_\alpha\otimes_AA'\rightarrow M/\a^{\alpha+1}M\otimes_AA'\rightarrow M_\alpha\otimes_AA'\rightarrow Q_\alpha\otimes_AA'\rightarrow 0\quad.
\]
Now $(\tilde{N}_\alpha\otimes_AA')_\alpha$ is (AR)-null, and $\varprojlim (Q_\alpha\otimes_AA')=0$ by Lemma \ref{ml1lem} and the subsequent remark. We may argue as above to see that we may pass to the limit without losing exactness. We thus obtain an isomorphism
\[
\varprojlim(M/\a^{\alpha+1}M\otimes_AA')\overset{\sim}{\rightarrow}\varprojlim(M_\alpha\otimes_AA')\quad.
\]
Since $M\otimes_AA'$ is finite over the noetherian ring $A'$ and since $A'$ is $\a A'$-adically complete and separated, cf. \cite{EGAIn} 0.7.2.4, we conclude that
\[
M\otimes_AA'\rightarrow\varprojlim(M_\alpha\otimes_AA')
\]
is an isomorphism, as desired.
\end{proof}

\begin{theorem}\label{mainflatbasechangethm}
Let $\phi\colon\fY\rightarrow\fX$ be a proper morphism of locally noetherian formal schemes, let $\psi\colon\fX'\rightarrow\fX$ be a flat morphism of locally noetherian formal schemes, and let
\[
\begin{diagram}
\fY'&\rTo^{\phi'}&\fX'\\
\dTo<{\psi'}&&\dTo>\psi\\
\fY&\rTo^\phi&\fX
\end{diagram}
\]
be the induced cartesian diagram. Then $\fY'$ is locally noetherian, $\phi'$ is proper, $\psi'$ is flat, and for any coherent $\O_\fY$-module $\sG$, the natural morphism of coherent $\O_{\fX'}$-modules
\[
\psi^*R^q\phi_*\sG\rightarrow R^q\phi'_*(\psi')^*\sG
\]
is an isomorphism.
\end{theorem}
\begin{proof}
By \cite{EGAI} 10.13.5 (ii), $\fY'$ is locally noetherian, and $\phi'$ is of tf type. We let a subscript $0$ denote reduction modulo the biggest ideal of definition of $\fX$. Since $\phi'_0$ is obtained from $\phi_0$ via base change, it follows that $\phi'$ is proper. By Proposition \ref{flatnessbasechangeprop}, $\psi'$ is flat. To prove the statement on the higher direct image sheaves, it suffices to consider the case where $\fX$ and $\fX'$ are affine, $\fX=\Spf A$, $\fX'=\Spf A'$. We must show that for all $q\geq 0$, the natural homomorphism of finite $A'$-modules
\[
H^q(\fY,\sG)\otimes_AA'\rightarrow H^q(\fY',(\psi')^*\sG)
\]
is an isomorphism. Let $\a$ be the biggest ideal of definition of $A$. A subscript $n$ denotes reduction modulo $\a^{n+1}$. By \cite{EGAIII1} 3.4.4, there is a natural system of $A$-homomorphisms
\[
H^q(\fY,\sG)\rightarrow H^q(\fY,\sG_n)
\]
satisfying the conditions of Lemma \ref{ml2lem}. By Lemma \ref{ml2lem}, we obtain an induced isomorphism
\[
H^q(\fY,\sG)\otimes_AA'\rightarrow \varprojlim(H^q(\fY,\sG_n)\otimes_AA')\quad.
\]
Now $H^q(\fY,\sG_n)=H^q(\fY_n,\sG_n)$, and $\fY_n$ is a scheme over $\Spec A_n$. It follows that
\[
H^q(\fY_n,\sG_n)\otimes_AA'\,=\,H^q(\fY_n\otimes_AA',\sG_n\otimes_AA')
\]
by the flat base change theorem for proper morphisms of schemes. Since $\fY'_n$ is obtained from $\fY_n\otimes_AA'$ via formal completion, the Comparison Theorem \cite{EGAIII1} 4.1.5 shows that 
\[
H^q(\fY_n\otimes_AA',\sG_n\otimes_AA')\cong H^q(\fY'_n,(\psi')^*\sG_n)\quad,
\]
so we obtain a natural isomorphism
\[
H^q(\fY,\sG)\otimes_AA'\,\cong\,\varprojlim_n(H^q(\fY'_n,(\psi')^*\sG_n)\quad.
\]
Let $\a'$ be the biggest ideal of definition of $A'$, and let a subscript $m$ denote reduction modulo $(\a')^{m+1}$. We note that $\a'$ induces an ideal of definition of $\fY'_n$ for all $n$. By \cite{EGAIII1} 3.4.4, 
\[
H^q(\fY'_n,(\psi')^*\sG_n)\cong\varprojlim_m H^q(\fY'_n,((\psi')^*\sG_n)_m)\quad,
\]
so we conclude that 
\begin{eqnarray*}
H^q(\fY,\sG)\otimes_AA'&\cong&\varprojlim_n(H^q(\fY'_n,(\psi')^*\sG_n)\\
&\cong&\varprojlim_n \varprojlim_m H^q(\fY'_n,((\psi')^*\sG_n)_m)\\
&\cong&\varprojlim_m H^q(\fY',((\psi')^*\sG)_m)\\
&\cong&H^q(\fY',(\psi')^*\sG)\quad,
\end{eqnarray*}
where the last equality is again \cite{EGAIII1} 3.4.4, and where we have used that $((\psi')^*\sG)_{n,m}=((\psi')^*\sG)_m$ for all $m\geq n$. One verifies that the natural isomorphism thus obtained agrees with the natural homomorphism in the statement of the proposition.
\end{proof}

		\subsection{Higher direct images and generic fibers}
			In this section, we study the behavior of higher direct images under proper morphisms of formal $R$-schemes of locally ff type with respect to the uniformly rigid generic fiber functor. 

\begin{lem}\label{cohkacycliclem}
Let $\fX$ be an affine formal $R$-scheme of ff type, and let $\sF$ be a coherent $\O_\fX$-module. Then 
\[
H^q(\fX,\sF_K)\,=\,\check{H}^q(\fX,\sF_K)\,=\,0
\]
for all $q\geq 1$.
\end{lem}
\begin{proof}
By \cite{EGAIII1} 3.4.4, we know that $H^q(\fU,\sF)=0$ for any affine open part $\fU$ of $\fX$ and for all $q\geq 1$. By the usual \v{C}ech spectral sequence argument, it follows that $\check{H}^q(\fU,\sF)=0$ for any affine open part $\fU$ of $\fX$ and for all $q\geq 1$. Since the functor $\cdot\otimes_RK$ is exact and since it commutes with direct limits, it follows that $\check{H}^q(\fU,\sF_K)=0$ for all affine opens $\fU$ in $\fX$ and all $q\geq 1$. Again, a \v{C}ech spectral sequence argument shows that $H^q(\fU,\sF_K)=0$ for all affine opens $\fU$ in $\fX$ and all $q\geq 1$.
\end{proof}

\begin{lem}\label{genfibbclem}
Let $\phi\colon\fY\rightarrow\fX$ be a proper morphism of formal $R$-schemes of locally ff type, and let $\sG$ be a coherent $\O_\fY$-module. Then for any $q\geq 0$, the natural homomorphism of $\O_{\fX,K}$-modules
\[
(R^q\phi_*\sG)_K\overset{\sim}{\rightarrow} (R^q\phi_*)(\sG_K)
\]
is an isomorphism.
\end{lem}
\begin{proof}
It suffices to show that in the case where $\fX$ is affine, the natural homomorphism
\[
H^q(\fY,\sG)\otimes_RK\rightarrow H^q(\fY,\sG_K)
\]
is an isomorphism. By Lemma \ref{cohkacycliclem}, we may work with \v{C}ech cohomology instead of derived functor cohomology, and hence the statement follows from the fact that the functor $\cdot\otimes_RK$ is exact and that it commutes with direct limits, together with the fact that $\fY$ is quasi-compact.
\end{proof}

\begin{prop}\label{blowupdirimcompprop}
Let $\phi\colon\fY\rightarrow\fX$ be a proper morphism of formal $R$-schemes of locally ff type, let $\psi\colon\fX'\rightarrow\fX$ be an admissible formal blowup, and let
\[
\begin{diagram}
\fY'&\rTo^{\phi'}&\fX'\\
\dTo<{\psi'}&&\dTo>{\psi}\\
\fY&\rTo^{\phi}&\fX
\end{diagram}
\]
be the induced cartesian diagram in the category of \emph{flat} formal $R$-schemes of locally ff type. Let $\sG$ be a coherent $\O_Y$-module. For each $q\geq 0$, the natural comparison morphism
\[
\psi^*R^q\phi_*\sG\rightarrow R^q\phi'_*(\psi')^*\sG
\]
induces an isomorphism
\[
(\psi^*R^q\phi_*\sG)_K\rightarrow (R^q\phi'_*(\psi')^*\sG)_K\quad.
\]
\end{prop}
\begin{proof}
By Proposition \ref{blowupisoonsemiafffuncprop}, it suffices to show that the induced homomorphism
\[
R^q\phi_*\sG_K\rightarrow \psi_*R^q\phi'_*(\psi')^*\sG_K
\]
is an isomorphism; by Lemma \ref{genfibbclem}, our notation is unambiguous. By Lemma \ref{genfibbclem} and Proposition \ref{blowupisoonsemiafffuncprop}, the functors $R^q\psi_*$ and $R^q\psi'_*$ vanish on sheaves obtained from coherent sheaves via $\cdot\otimes_RK$. Hence, the spectral sequences
\begin{eqnarray*}
R^p\psi_*R^q\phi'_*&\Rightarrow&R^{p+q}(\psi_*\phi'_*)\\
R^p\phi_*R^q\psi'_*&\Rightarrow&R^{p+q}(\phi_*\psi'_*)
\end{eqnarray*}
yield, for all $q\geq 0$, edge morphisms
\begin{eqnarray*}
R^q(\psi_*\phi'_*)&\rightarrow&\psi_*R^q\phi'_*\\
(R^q\phi_*)\psi'_*&\rightarrow&R^q(\phi_*\psi'_*)
\end{eqnarray*}
which are isomorphisms when evaluated in a sheaf that is obtained from a coherent sheaf via $\cdot\otimes_RK$. Since $\psi_*\phi'_*=\phi_*\psi'_*$ and since higher direct images of coherent sheaves under proper morphisms are coherent by \cite{EGAIII1} 3.4.2, we obtain natural identifications
\begin{eqnarray*}
\psi_*R^q\phi'_*(\psi')^*\sG_K&\cong&R^q\psi_*\phi'_*(\psi')^*\sG_K\\
&\cong&R^q(\phi_*\psi'_*)(\psi')^*\sG_K\\
&\cong&R^q(\phi_*)\psi'_*(\psi')^*\sG_K\\
&\cong&R^q\phi_*\sG_K\quad,
\end{eqnarray*}
where the last equality is again Proposition \ref{blowupisoonsemiafffuncprop}. One checks that the resulting isomorphism agrees with the homomorphism that is induced by the comparison homomorphism.
\end{proof}

\begin{cor}\label{highdirimsamaffpbcor}
The statement of Proposition \ref{blowupdirimcompprop} still holds when $\psi$ is a composition of simple morphisms.
\end{cor}
\begin{proof}
We may assume that $\psi$ is an open immersion, a completion morphism or an admissible formal blowup. If $\psi$ is an open immersion, there is nothing to show, and if $\psi$ is a completion morphism, the statement follows from Theorem \ref{mainflatbasechangethm}, completion morphisms of locally noetherian formal schemes being flat. Hence it suffices to consider the case where $\psi$ is an admissible formal blowup; in this case, the claim follows from Proposition \ref{blowupdirimcompprop}.
\end{proof}

\begin{cor}\label{propdirimgenfibercor}
Let $\ul{\phi}\colon\fY\rightarrow\fX$ be a proper morphism of flat formal $R$-schemes of locally ff type, let $\phi\colon Y\rightarrow X$ denote its uniformly rigid generic fiber, and let $\sG$ be a coherent $\O_\fY$-module. Then $R^q\phi_*(\sG^\srig)$ is a coherent $\O_X$-module for each $q\geq 0$, and the morphism 
\[
\sp_\fY^\sharp\colon\O_\fY\rightarrow\sp_{\fY,*}\O_Y
\]
induces natural isomorphisms
\[
(R^q\ul{\phi}_*\sG)^\srig\overset{\sim}{\rightarrow}R^q\phi_*(\sG^\srig)
\]
for $q\geq 0$.
\end{cor}
\begin{proof}
Let $q\geq 0$ be fixed. By \cite{Tohoku} 3.7.2, $R^q\phi_*(\sG^\srig)$ is associated to presheaf 
\[
V\mapsto H^q(\phi^{-1}(V),\sG^\srig)\quad,
\]
where $V$ varies among the admissible subsets of $X$. Let $(\fX_i)_{i\in I}$ be an affine open covering of $\fX$; then the semi-affinoid subdomains in the semi-affinoid subspaces $X_i\mathrel{\mathop:}=\fX_i^\srig$ of $X$ form a basis of the uniformly rigid G-topology of $X$. Let $U$ be a semi-affinoid subdomain of $X_i$ for some $i\in I$, and let $\fU\rightarrow\fX_i$ be a defining morphism; by Corollary \ref{highdirimsamaffpbcor}, we have a natural identification
\[
(R^q\ul{\phi}_*\sG)^\srig|_U\,\cong\,R^q\ul{\phi}_{\fU}^*(\sG|_{\fY_\fU})^\srig\quad,
\]
where $\ul{\phi}_{\fU}$ and $\sG|_{\fY_\fU}$ are defined via pullback to $\fU$ and $\fY_\fU\mathrel{\mathop:}=\fY\times'_\fX\fU$ respectively. It thus suffices to define, in the case where $\fX$ is affine, a natural isomorphism
\[
H^q(\fY,\sG)\otimes_RK\overset{\sim}{\rightarrow}H^q(Y,\sG^\srig)\quad.
\]
We may use the \v{C}ech cohomology groups attached to an affine covering $(\fY_i)_{i\in I}$ of $\fY$ to calculate the derived functor cohomology groups. Indeed, $(\fY_i)_{i\in I}$ is a Leray covering of $\fY$ for $\sG$, and $(\fY_i^\srig)_{i\in I}$ is a Leray covering for $Y$ by the Acyclicity Theorem \ref{acyclicitytheorem}. Since $\cdot\otimes_RK$ is exact, we thus obtain the desired natural isomorphisms.
\end{proof}

		\subsection{Schematic images of proper morphisms}\label{schemimsec}
			\begin{defi}\label{schemimdefi}
Let $\phi\colon\fY\rightarrow\fX$ be a morphism of locally noetherian formal schemes. A closed formal subscheme $\fV\subseteq\fX$ is called the \emph{schematic image} of $\phi$ if $\phi$ factorizes through $\fV$ and if $\fV$ is minimal among all closed formal subschemes of $\fX$ with this property.
\end{defi}

The schematic image of $\phi$ is clearly unique if it exists. Indeed, if $\fV_1$ and $\fV_2$ are closed formal subschemes of $\fX$ such that $\phi$ factorizes through $\fV_1$ and $\fV_2$, then $\phi$ also factorizes through the closed formal subscheme $\fV_1\times_\fX\fV_2$ of $\fX$; hence if both $\fV_1$ and $\fV_2$ are minimal, it follows that $\fV_1=\fV_1\times_\fX\fV_2=\fV_2$. If $\fX$ is defined over $R$ and if $\fY$ is $R$-flat, then the schematic image of $\phi$, if it exists, is $R$-flat, the ideal of $\pi$-torsion of a locally noetherian formal $R$-scheme being coherent.

Moreover, it is clear from Definition \ref{schemimdefi} that the schematic image of $\phi$ exists if and only if the kernel of $\phi^\sharp$ admits a biggest coherent $\O_\fX$-submodule. For example, this is the case when $\ker\phi^\sharp$ is itself coherent. In this case, the formation of the schematic image of $\phi$ is local on $\fX$.

Let us consider an example where the schematic image of $\phi$ exists without $\ker\phi^\sharp$ being coherent:

\begin{example}
Let $\phi\colon\D^1_k\rightarrow\A^1_k$ be the completion of the affine line over the field $k$ in the origin. Then the schematic image of $\phi$ exists and in fact coincides with $\A^1_k$, while $\ker\phi^\sharp$ is not coherent.
\end{example}
\begin{proof}
Any closed subscheme of the affine scheme $\A^1_k$ is associated to an ideal in the ring of global sections of $\A^1_k$. Since $\phi$ induces an injection of rings of global sections, $\phi$ cannot factor through a proper closed subscheme of $\A^1_k$; hence the schematic image of $\phi$ exists and coincides with $\A^1_k$. If the ideal $\ker\phi^\sharp$ was coherent, it would be trivial since it has no nonzero global section. However, for any nonempty open subscheme $U\subseteq\A^1_k$ not containing the origin, $(\ker\phi^\sharp)(U)=\O_{\A^1_k}(U)$ since $\phi^{-1}(U)=\emptyset$.
\end{proof}


We are only interested in the case where $\ker\phi^\sharp$ is coherent. For example, this is the case when $\ker\phi^\sharp$ is trivial or when $\phi_*\O_\fY$ is coherent. If $\phi$ is proper, then $\phi_*\O_\fY$ is coherent by \cite{EGAIII1} 3.4.2.

\subsubsection{Schematic images and flat base change}

\begin{prop}\label{flatbasechangeprop}
Let $\phi\colon\fY\rightarrow\fX$ be a proper morphism of locally noetherian formal schemes, let $\psi\colon\fX'\rightarrow\fX$ be a flat morphism of locally noetherian formal schemes, and let
\[
\begin{diagram}
\fY'&\rTo^{\phi'}&\fX'\\
\dTo<{\psi'}&&\dTo>\psi\\
\fY&\rTo^\phi&\fX
\end{diagram}
\]
be the induced cartesian diagram. If $\fV\subseteq\fX$ and $\fV'\subseteq\fX'$ are the schematic images of $\phi$ and $\phi'$ respectively, then
\[
\fV'\,=\,\psi^{-1}(\fV)\,\mathrel{\mathop:}=\,\fV\times_\fX\fX'\quad.
\]
\end{prop}
\begin{proof}
Let $\sI$ denote the coherent $\O_\fX$-ideal defining $\fV$. Then $\sI$ is defined by the natural exact sequence of coherent $\O_\fX$-modules
\[
0\rightarrow\sI\rightarrow\O_\fX\rightarrow\phi_*\O_\fY\quad.
\]
Since $\psi$ is flat, we obtain an induced exact sequence
\[
0\rightarrow\psi^*\sI\rightarrow\O_{\fX'}\rightarrow\psi^*\phi_*\O_\fY\quad,
\]
and $\psi^*\sI$ is the coherent $\O_{\fX'}$-ideal defining $\psi^{-1}(\fV)$. By Theorem \ref{mainflatbasechangethm}, $\psi^*\phi_*\O_\fY$ is naturally identified with $(\phi')_*\O_{\fY'}$, so we see that $\psi^*\sI$ also defines the schematic image $\fV'$ of $\phi'$, as desired.
\end{proof}

\subsubsection{Schematic images and schematic closures}

\begin{lem}\label{schemcldirimlem}
Let $\ul{\phi}\colon\fY\rightarrow\fX$ be a proper morphism of formal $R$-schemes of locally ff type whose uniformly rigid generic fiber $\phi\colon Y\hookrightarrow X$ is a closed immersion. Then the schematic closure of $Y$ in $\fX$ coincides with the schematic image of $\ul{\phi}$. Moreover, the natural morphism from $\fY$ to the schematic image of $\ul{\phi}$ is surjective on physical points.
\end{lem}
\begin{proof}
Since $\fY$ is $R$-flat, the schematic image $\fV$ of $\ul{\phi}$ is $R$-flat as well. To show that $\fV$ is the schematic closure of $Y$ in $\fX$, it remains to see that the unique factorization $\ul{\phi}'\colon\fY\rightarrow\fV$ of $\ul{\phi}$ induces an isomorphism of uniformly rigid generic fibers. Let $\sI\subseteq\O_\fX$ be the coherent ideal defining $\fV$ as a closed formal subscheme of $\fX$. The exact sequence of coherent $\O_\fX$-modules
\[
0\rightarrow\sI\rightarrow\O_\fX\rightarrow\ul{\phi}_*\O_\fY
\]
induces an exact sequence of coherent $\O_X$-modules
\[
0\rightarrow\sI^\srig\rightarrow\O_X\rightarrow(\ul{\phi}_*\O_\fY)^\srig\quad,
\]
the functor $\srig$ being exact. By Corollary \ref{propdirimgenfibercor}, it follows that
\[
0\rightarrow\sI^\srig\rightarrow\O_X\rightarrow \phi_*\O_Y\rightarrow 0
\]
is exact. Hence, $Y$ and the uniformly rigid generic fiber of $\fV$ coincide as closed uniformly rigid subspaces of $X$, as stated.

It remains to see that $\ul{\phi}'$ is surjective. Since $\fV$ is $R$-flat, the specalization map $\sp_{\fV}$ is surjective onto the set of closed points of $\fV$;  since $\ul{\phi}'$ is an isomorphism on uniformly rigid generic fibers, it follows that the image of $\ul{\phi}'$ contains every closed point of $\fV$. Since $\ul{\phi}'$ is of locally ff type, the image of $\ul{\phi}'$ is locally constructible by \cite{EGAIV1} 1.8.5. Hence, its complement its locally constructible as well.
Since the closed points of $\fV$ lie very dense, we conclude from \cite{EGAIn} 0.2.6.2 (b) that $\ul{\phi}'$ is surjective, as desired.
\end{proof}

		\subsection{Schematically dominant morphisms}\label{schemunivsec}
			Let us recall the notion of a schematically dominant morphism of schemes as it is defined in \cite{EGAIV3} 11.10: A morphism $\phi\colon Y\rightarrow X$ of schemes is called \emph{schematically dominant} if $\phi^\sharp$ is a monomorphism. Let $S$ be any scheme such that $X$ is defined over $S$. By \cite{EGAIV3} 11.10.1, the following are equivalent: 
\begin{packed_enum}
\item $\phi$ is schematically dominant.
\item For any open subscheme $U\subseteq X$, the restriction $\phi^{-1}(U)\rightarrow U$ of $\phi$ does not factor though a proper closed subscheme of $U$.
\item For any open subscheme $U\subseteq X$ and any $S$-scheme $X'$, two $U$-valued points $x_1,x_2\in X'(U)$ of $X'$ coincide if and only if the $\phi$-induced $\phi^{-1}(U)$-valued points of $X'$ coincide.
\end{packed_enum}
If $\phi\colon Y\rightarrow X$ is a morphism of $S$-schemes, then according to \cite{EGAIV3} 11.10.8, $\phi$ is called \emph{universally schematically dominant over} $S$ if for all $S$-schemes $T$, the base change $\phi_T\colon Y_T\rightarrow X_T$ of $\phi$ from $S$ to $T$ is schematically dominant.

We generalize this concept to the setting of morphisms of locally noetherian formal schemes that are \emph{adic} over a given base. Let $\fS$ be a locally noetherian formal base scheme. We recall that a locally noetherian formal $\fS$-scheme $\fX$ is called \emph{adic} if the ring topology on $\O_\fX$ is induced by the ring topology on $\O_\fS$.

\begin{defi}\label{univschemdom}\index{morphism!schematically dominant}
A morphism $\phi\colon\fY\rightarrow\fX$ of \emph{adic} locally noetherian formal $\fS$-schemes is called \emph{universally schematically dominant above} $\fS$ if for any \emph{scheme} $T$ and any morphism $T\rightarrow\fS$, the induced morphism of schemes
\[
\phi_T\colon\fY\times_\fS T\rightarrow\fY\times_\fS T
\]
is schematically dominant. 
\end{defi}

Here we do not assume $T$ to be locally noetherian. 

In analogy with \cite{EGAIV3} 11.10.9, universal schematic dominance is, in the flat case, equivalent to schematic dominance in the fibers:

\begin{lem}\label{univschemdomfiberlem}
Let $\phi\colon\fY\rightarrow\fX$ be a morphism of locally noetherian adic $\fS$-schemes, and let us assume that $\fY$ is flat over $\fS$. Then the following are equivalent:
\begin{packed_enum}
\item $\phi$ is universally schematically dominant above $\fS$.
\item For any point $s\in\fS$, the fiber 
\[
\phi_s\colon\fY\times_\fS k(s)\rightarrow\fX\times_\fS k(s)
\]
is a schematically dominant morphism of schemes.
\end{packed_enum}
\end{lem}
\begin{proof}
The implication ($i$)$\Rightarrow$($ii$) is trivial, so let us assume that ($ii$) is satisfied. We show that ($i$) holds. Let $T$ be any $\fS$-scheme; we must show that $\phi_T$ is schematically dominant. Since schematic dominance is local on the target, we may assume that $T$ is affine and, hence, quasi-compact. Then the morphism $T\rightarrow \fS$ factorizes through an infinitesimal neighborhood $\fS_n$ of the smallest subscheme of definition $\fS_0$ of $\fS$ in $\fS$. It thus suffices to show that $\phi_{\fS_n}$ is universally schematically dominant. By Proposition \ref{flatnessbasechangeprop}, the scheme $\fY\times_\fS\fS_n$ is flat over $\fS_n$, so by \cite{EGAIV3} 11.10.9, $\phi_{\fS_n}$ is universally schematically dominant if and only if $\phi_{\fS_n,s}$ is schematically dominant for each point $s\in\fS_n$. However, $|\fS|=|\fS_n|$, and for any point $s$ in this set, $\phi_{\fS_n,s}=\phi_s$.
\end{proof}

In particular, an open immersion $\fU\hookrightarrow\fX$ of flat adic locally noetherian formal $R$-schemes is universally schematically dominant if and only if $\fU_k\subseteq\fX_k$ is schematically dense, in which case we say that $\fU$ is \emph{$R$-dense}\index{subspace!R@$R$-dense} in $\fX$. If in addition $\fX_k$ is reduced, this condition is equivalent to $\fU_k$ being dense in $\fX_k$.

Let us emphasize that in the present section, we did not impose any relative finiteness conditions on morphisms.

		\subsection{Smoothness}\label{smoothnesssec}
			Smooth morphisms of locally noetherian formal schemes have been studied in \cite{AJR1}, \cite{AJR3} and, essentially, in \cite{EGAIV1} Chapter 0. Let us briefly recall the main definitions and results. According to \cite{AJR1} 2.6, a morphism of locally ff type of locally noetherian formal schemes $\phi\colon\fY\rightarrow\fX$ is called \emph{smooth}\index{morphism!smooth}, \emph{unramified} or \emph{étale} if it satisfies the corresponding infinitesimal lifting property, cf.\ \cite{AJR1} 2.1, 2.3. These properties are stable under composition as well as under base change with respect to morphisms of locally ff type, cf.\ \cite{AJR1} 2.9. Moreover, they are local on the domain, cf.\ \cite{AJR1} 4.1. The morphism $\phi$ is called smooth, unramified or étale in a point $x\in\fX$ if it satisfies the corresponding property in an open neighborhood of $\fX$, cf.\ \cite{AJR1} 4.3. Open immersions are étale, and closed immersions are unramified, cf.\ \cite{AJR1} 2.12. By \cite{AJR1} 4.4, completion morphisms are étale. Smoothness is conveniently characterized in terms of sheaves of relative \emph{continuous} differentials: Let $\phi\colon\fY\rightarrow\fX$ be of locally ff type, where $\fX$ and $\fY$ are locally noetherian. Then the sheaf of relative continuous differential forms $\Omega^1_{\fY/\fX}$ is coherent, and the morphism $\phi\colon\fY\rightarrow\fX$ is étale if and only if $\Omega^1_{\fY/\fX}$ is trivial, cf.\ \cite{AJR1} 4.6. If $\phi$ is smooth, then $\Omega^1_{\fY/\fX}$ is locally free by \cite{AJR1} 4.8; its rank is called the \emph{relative dimension} of $\fY$ over $\fX$. Let $\fS$ be a locally noetherian formal base scheme, and let $\fX$ be smooth of locally ff type over $\fS$. Then $\phi$ is smooth if and only if $\fY$ is smooth over $\fS$ and the natural sequence
\[
0\rightarrow \phi^*\Omega^1_{\fX/\fS}\rightarrow\Omega^1_{\fY/\fS}\rightarrow\Omega^1_{\fY/\fX}\rightarrow 0
\]
is exact and locally split, cf.\ \cite{AJR1} 4.12. One derives the \emph{Jacobian criterion} for smoothness in locally noetherian formal geometry: Let us assume that $\fX$ is smooth of locally ff type over $\fS$ and that $\phi\colon\fY\hookrightarrow\fX$ is a closed immersion. Let $y\in \fY$ be a point, and let us consider local sections of $\fX$ defining $\phi$ near $y$. Then $\fY$ is smooth in $y$ over $\fS$ if and only if the associated Jacobian matrix has an appropriate minor that does not vanish in $y$, cf.\ \cite{AJR1} 4.15 and \cite{AJR3} 5.11, 5.12. 

By \cite{AJR3} 5.4, smoothness of a morphism $\phi\colon\fY\rightarrow\fX$ of locally ff type in a point $y$ is equivalent to formal smoothness of the corresponding homomorphism of (completed) stalks with respect to the maximal-adic topologies. Moreover, by \cite{AJR3} 5.4, $\phi$ is smooth in $y$ if and only if it is flat in $y$ and if its fiber over the image of $y$ in $\fX$ is smooth in $y$. Let us point us that this fiber is in general not a scheme but a formal scheme, unless $\phi$ is adic.


A smooth morphism of locally ff type is locally an étale morphism of ff type to a relative closed formal unit ball, cf.\ \cite{AJR3} 5.9. Moreover, a smooth morphism of locally ff type is locally a completion of a smooth morphism of tf type, cf.\ \cite{AJR3} 7.12.

\begin{prop}\label{regprop}\index{morphism!regular}
Let $\phi\colon\fY\rightarrow\fX$ be a smooth morphism of ff type of affine formal $R$-schemes of ff type. Then $\phi$ induces a \emph{regular} homomorphism of rings of global sections.
\end{prop}
\begin{proof}
By \cite{EGAIV2} 6.8.3 ($ii$), we may work locally on $\fY$. By \cite{EGAIV2} 7.8.3 ($v$), completion homomorphisms of \emph{excellent} noetherian rings are regular; by \cite{AJR3} 7.12, we may thus assume that $\phi$ is adic. By \cite{Elk} Théorème 7, $\fY$ is obtained from a smooth algebra over the ring of global functions on $\fX$ via formal completion. By \cite{EGAIV2} 6.8.6, smooth morphisms of locally noetherian schemes are regular. Since algebras of finite type over excellent noetherian rings are excellent, the statement now follows from a second application of \cite{EGAIV2} 7.8.3 ($v$).
\end{proof}

In particular, it follows from Proposition \ref{regprop} and from \cite{EGAIV2} 6.5.4 ($ii$) that if $\phi\colon\fY\rightarrow\fX$ is a smooth morphism of locally ff type of formal $R$-schemes of locally ff type and if $\fX$ is \emph{normal}, then $\fY$ is normal as well.

\subsubsection{Local smooth formal $R$-schemes}

In our subsequent arguments, we will recurrently complete smooth formal $R$-schemes of locally ff type in closed points in order to reduce to a local situation. Completion morphisms being étale, the local formal $R$-schemes of ff type thus obtained are smooth. We are therefore interested in the structure of local smooth formal $R$-schemes of ff type. In the remainder of this section, we show that under a suitable finite extension of $R'/R$ which we allow to ramify in the residue fields, a smooth local formal $R$-scheme of ff type splits as a finite union of open formal unit balls.

\begin{lem}\label{localstructurerationalpointlem}
Let $\fX$ be a local smooth formal $R$-scheme of ff type whose residue field naturally coincides with $k$. Then there exists an $R$-isomorphism 
\[
\fX\cong\Spf R[[T_1,\ldots,T_d]]\quad,
\]
where $d$ denotes the relative dimension of $\fX$ over $R$.
\end{lem}
\begin{proof}
Let us write $\fX=\Spf A$, $A_k=A/\pi A$. By \cite{EGAIV1} 0.19.6.4, there exists a $k$-isomorphism $\phi_k\colon A_k\overset{\sim}{\rightarrow} k[[T_1,\ldots,T_d]]$. Since $A$ is $\pi$-adically complete and formally $R$-smooth for the maximal-adic topologies,  \cite{EGAIV1} 0.19.7.1.5 shows that $\phi_k$ extends to an $R$-isomorphism $A\overset{\sim}{\rightarrow} R[[T_1,\ldots,T_d]]$.
\end{proof}

\begin{prop}\label{smoothlocalstructureprop}
Let $\fX$ be a local smooth formal $R$-scheme of ff type. There exists a finite extension $R'/R$ of discrete valuation rings with ramification index $e_{R'/R}=1$ such that the base extension $\fX\otimes_RR'$ is $R'$-isomorphic to a disjoint union of finitely many copies of $\Spf R'[[T_1,\ldots,T_d]]$, where $d$ denotes the relative dimension of $\fX$ over $R$.
\end{prop}
\begin{proof}
Let us write $\fX=\Spf A$, let $k_A$ denote the residue field of $A$, and let $k'$ be the normal envelope of $k_A$ over $k$. Then $k_A\otimes_kk'$ is, up to nilpotents, a finite direct sum of copies of $k'$. By \cite{Liu} Lemma 10.3.32, the extension $k'/k$ lifts to a finite extension of discrete valuation rings $R'/R$ such that $e_{R'/R}=1$. By Hensel's Lemma, $\fX\otimes_RR'$ is a finite disjoint union of local smooth formal $R'$-schemes of ff type whose residue fields all coincide with $k'$. The statement thus follows from Lemma \ref{localstructurerationalpointlem}.
\end{proof}
			
		\subsection{A complement on Berthelot's construction}\label{berthcompsec}				
			In \cite{RZ} Proposition 5.12, the following result is attributed to R. Huber. Since no reference is given, we provide a proof:

\begin{lem}\label{qctftlem}
Let $\fX$ be a flat formal $R$-scheme of ff type such that $\fX^\rig$ is quasi-compact. Then $\fX$ is of tf type over $R$.
\end{lem}


\begin{proof}
The property of being of locally tf type over $R$ can be checked locally on $\fX$; we may thus assume that $\fX$ is affine. We claim that there exists a flat $R$-model of tf type $\fY$ of $\fX^\rig$ together with a morphism $\phi\colon\fY \rightarrow \fX$ inducing an isomorphism of rigid generic fibers. Indeed: Since $\fX^\rig$ is quasi-compact and quasi-separated, there exists a flat $R$-model of tf type for $\fX^\rig$, cf.\ \cite{BL1} 4.1. Let $\ul{A}$ denote the ring of global functions on $\fX$, and let $(\ul{s},\ul{t})$ be a formal generating system for $\ul{A}$. By Proposition \ref{affineblowupprop}, there exists a finite admissible formal blowup $\fY'$ of $\fY$ such that the components of $(\ul{s},\ul{t})$ give rise to global functions on $\fY'$. After replacing $\fY$ by $\fY'$, the existence of $\phi$ follows from the fact that $R$-morphisms to $\fX$ correspond to continuous $R$-homomorphisms of rings of global sections, cf.\ \cite{EGAIn} 10.4.6. 

Since $\fY$ is adic over $R$, there exists a natural number $n\geq 1$ such that the special fiber $\fY_k$ of $\fY$ factorizes through the $n$-th infinitesimal neighborhood of the smallest subscheme of definition $\fX_0$ of $\fX$ in $\fX$. Let $ \a\subseteq\ul{A}$ denote the ideal corresponding to $\fX_0$, let $\fX''\rightarrow \fX$ be the admissible formal blowup of $\fX$ in $(\pi, \a^n)$, and let $\fX'\subseteq \fX''$ be the corresponding formal dilatation. By the universal property of formal dilatations, $\phi$ extends uniquely to a morphism $\phi'\colon\fY \rightarrow \fX'$. Since $\fX$ is $R$-flat by assumption, the specialization map $\sp_{\fX''}$ is surjective. Since $\phi^\rig$ is an isomorphism, it follows that the inclusion $\fX'\subseteq \fX''$ is in fact an equality. Let $\ul{A}'$ denote the flat $R$-algebra of global functions on the affine formal scheme $\fX'$; then $\ul{A}'$ is of tf type over $R$, cf.\ Lemma \ref{berthconstrbasiclem}. By Corollary \ref{afflattfintypecor}, it follows that $\ul{A}$ is of tf type over $R$ as well.
\end{proof}

We will generalize Lemma \ref{qctftlem} to a relative setting. To do so, we first generalize a construction in \cite{dJ} 7.1.13 and \cite{dJErr}:

\begin{lem}\label{dJlem}
Let $\fX$ be an affine flat formal $R$-scheme of ff type, and let $\fV\subseteq\fX$ be a closed formal subscheme. For each $n\in\N$, we let $\fV_{(n)}\subseteq\fX$ denote the $(n+1)$-st infinitesimal neighborhood of $\fV$ in $\fX$, and we let $\fX_{(n)}$ denote the formal dilatation of $\fX$ in $\fV_{(n)}\cap\fX_k$. Then there exists an integer $n_0$ such that for all $n\geq n_0$, there exists a natural closed immersion $\phi_n\colon\fV_{(n-n_0)}\hookrightarrow \fX_{(n)}$ such that the natural diagram
\[
\begin{diagram}
&&\fX_{(n)}\\
&\ruInto<{\phi_n}&\dTo\\
\fV_{(n-n_0)}&\rInto&\fX
\end{diagram}
\]
commutes. Moreover, for varying $n$ the $\phi_n$ are compatible.
\end{lem}
\begin{proof}
Let $A$ denote the $R$-algebra of global functions on $\fX$, and let $I\subseteq A$ be the ideal corresponding to $\fV$. By the Artin-Rees Lemma, there exists an integer $n_0$ such that for all $n\geq n_0$,
\[
\pi A\cap I^n\,=\, I^{n-n_0}(\pi A\cap  I^{n_0})\,\subseteq\,\pi I^{n-n_0}\quad.
\]
We claim that for any integer $t\geq 1$ and any $n\geq n_0$,
\[
\pi^t A\cap(\sum_{i=1}^t \pi^{t-i} I^{in})\,\subseteq\,\pi^t I^{n-n_0}\quad.
\]
Indeed, let us argue by induction on $t$. We have chosen $n_0$ such that the statement holds for $t=1$, so we may assume that $t>1$. Let $a$ be an element on the left hand side. Since $a$ lies in $\sum_{i=1}^t \pi^{t-i} I^{in}$, we can write
\[
a\,=\,\sum_{i=1}^t\pi^{t-i}a_i
\]
with elements $a_i\in I^{in}$. Since $a$ is divisible by $\pi$, the same must be true for $a_t$; hence $a_t\in\pi A\cap I^{tn}$, and $\pi A\cap I^{tn}\subseteq\pi I^{tn-n_0}$ by our choice of $n_0$. Since $n\geq n_0$, it follows that $a_t\in\pi I^{(t-1)n}$. We set $a_i'\mathrel{\mathop:}=a_i$ for $1\leq i \leq t-2$, $a_{t-1}'\mathrel{\mathop:}=a_{t-1}+\pi^{-1}a_t\in I^{(t-1)n}$; then $a_i'\in I^{in}$ for $1\leq i\leq t-1$, and
\[
a\,=\,\pi\sum_{i=1}^{t-1}\pi^{t-1-i}a'_i\quad.
\]
Since $\pi^{-1}a\in\pi^{t-1}A$, the induction hypothesis implies that $a\in\pi^t I^{n-n_0}$, as desired; thus our claim has been shown. 
Let us now consider an element $a\in A[\frac{ I^n}{\pi}]\subseteq A\otimes_RK$. We choose a representation
\[
a\,=\,\sum_{i=0}^t\frac{a_i}{\pi^i}\quad,\quad a_i\in I^{in}
\]
of $a$. The class of $a_0$ modulo $ I^{n-n_0}$ does not depend on the choice of the representation. Indeed, if $a=0$, then $\pi^ta_0\in\pi^t I^{n-n_0}$ by the above claim. It is now clear that $a\mapsto a_0$ defines a surjective $A$-algebra homomorphism $A[I^n/\pi]\rightarrow A/ I^{n-n_0}$. The statement thus follows via formal completion with respect to an ideal of definition of $A$. 
\end{proof}

We can now provide the announced relative version of Lemma \ref{qctftlem}:

\begin{prop}\label{rqcadicprop}
Let $\phi\colon\fY\rightarrow \fX$ be a morphism of flat formal $R$-schemes of locally ff type. If $\phi^\rig$ is quasi-compact, then $\phi$ is adic.
\end{prop}
\begin{proof}
We may assume that $\fX$ and $\fY$ are affine, $\fX=\Spf A$, $\fY=\Spf B$. Using $\phi$, we consider $B$ as an $A$-algebra. Let $ \a\subseteq A$ be an ideal of definition; we must show that $\a B$ is an ideal of definition of $B$. Let $n\in\N$ be a natural number, and let $\fX'$ and $\fY'$ denote the formal dilatations of $\fX$ and $\fY$ in the respective ideals generated by $\pi$ and $ \a^{n+1}$; then $(\fY')^\rig$ is the $\phi^\rig$-preimage of $(\fX')^\rig$. By assumption on $\phi$, $(\fY')^\rig$ is quasi-compact, and by Lemma \ref{qctftlem}, it follows that $\fY'$ is of tf type over $R$. By Lemma \ref{dJlem}, there exists an $n_0\in\N$ such that for all $n\geq n_0$, there exists a closed immersion of $\Spf B/\a^{n-n_0}B$ into $\fY'$. It follows that $\a B$ is an ideal of definition for $B$.
\end{proof}

	\section{Morphisms to uniform rigidifications}\label{vpsec}
	
There may exist smooth affine formal $R$-schemes of ff type $\fX$ such that $\fX(R')=\emptyset$ for all unramified local extensions of discrete valuation rings $R'/R$, where these extensions are allowed to be infinite. For example, if the residue field $k$ of $R$ is not perfect, such an example is obtained by completing the finite type $R$-scheme $\A^1_R$ in a closed point of its special fiber where the residue field of this point is inseparable over $k$. Even for general $R$ and in the case of the open formal unit disc $\fX=\Spf R[[S]]$, there are \emph{not enough étale points}: If $R'/R$ is any unramified local extension and if $x\in \fX(R')$ is any $R'$-valued point, then $x$ factorizes over the closed subdisc of radius $\pi$; that is, the $x^*$-image of the coordinate $S$ in $R'$ is $\pi$-divisible. We can thus say that the unramified points of $\fX$ do not see the boundary of $\fX$ where $|S|\rightarrow 1$. 

If we replace the Jacobson-adic topology on $R[[S]]$ by the $\pi$-adic topology, new unramified points appear: Let $R[[S]]^\pi$ denote the ring $R[[S]]$ equipped with its $\pi$-adic topology, and let $\fX_\pi$ denote its formal spectrum. Then $\fX_\pi$ is an \emph{adic} noetherian formal $R$-scheme. Its special fiber $\fX_{\pi,k}$ is the \emph{scheme} $\Spec k[[S]]$, and the natural morphism $\fX\rightarrow\fX_\pi$ is a completion morphism. However, $\fX_\pi$ is \emph{not} of ff type over $R$. The localization $R[[S]]_{(\pi)}$ of $R[[S]]$ in the generic point $\eta$ of $\fX_{\pi,k}$ is a discrete valuation ring which is an infinite unramified local extension of $R$. Its completion $R'$ defines a point $|\cdot|_\Gauss\in \fX_\pi(R')$ which is not induced from a point in $\fX(R')$ and which corresponds to the Gauss point in $M(R[[S]]\otimes_RK)$, cf.\ Proposition \ref{integralredshilovprop}. 

Of course, we can define $\fX_\pi$ when $R[[S]]$ is replaced by a general $R$-algebra of ff type. We say that $\fX_\pi$ is the \emph{envelope}\index{envelope!affine} of $\fX$. The notion of envelope also generalizes to a slightly more global setting, cf.\ Section \ref{envsec}. The points of $\fX_\pi$ with values in local extensions $R'/R$ of complete discrete valuation rings may be called the \emph{virtual points} of $\fX$. It is one of our main objectives to make these virtual points more accessible and to make use of \emph{formally unramified} virtual points in the context of formal Néron models.

Envelopes of affine formal $R$-schemes of ff type have already been used in \cite{S} and \cite{Hubercohpaper}. Their formation does not commute with complete localization; in particular, we do not know how to define the envelope of a global formal $R$-scheme of locally ff type. Moreover, since envelopes do not satisfy any relative finiteness condition, fibered products like $\fX_\pi\times_R\fX_\pi$ need not be locally noetherian. Finally, generic fibers of envelopes are not defined in the categories of rigid or uniformly rigid $K$-spaces, but in the category of analytic adic $K$-spaces, cf.\ \cite{Huberbuch}.

In this section, we will prove the following: Let $\fY$ be an \emph{affine smooth} formal $R$-scheme of ff type, let $\fX$ be a separated formal $R$-scheme of locally \emph{tf} type, let $\phi\colon Y\rightarrow X$ be a morphism of uniformly rigid generic fibers, and let $\ul{\Gamma}_\phi$ denote the schematic closure of the graph of $\phi$ in $\fY\times\fX$, cf.\ the end of Section \ref{separatedsection}. Then the projection $p_\fY|_{\ul{\Gamma}_\phi}$ from $\ul{\Gamma}_\phi$ to $\fY$ is proper. Moreover, it admits a unique envelope $\ul{\Gamma}_\pi\rightarrow\fY_\pi$, and this morphism is an isomorphism above the generic points of $\fY_{\pi,k}$. Hence, if $\eta$ is a virtual unramified point of $\fY$ which specializes to a generic point of $\fY_{\pi,k}$, we can canonically define the image of $\eta$ in $\fX$. In particular, this applies to the case where $\fY$ is the formal open unit disc $\Spf R[[S]]$ and where $\eta$ is the Gauss point.
		\subsection{Envelopes}\label{envsec}	
			As we have just explained, if $\fX=\Spf A$ is an affine formal $R$-scheme of ff type, then its \emph{envelope} $\fX_\pi$ is defined to be the formal spectrum of $A^\pi$, where $A^\pi$ is the ring $A$ equipped with the $\pi$-adic topology. The notion of envelopes can be generalized from the affine case to the case of formal $R$-schemes that are \emph{proper} over an affine formal $R$-scheme of ff type:

\begin{defi}\label{envdefi}\index{envelope!proper}
Let $\fS$ be an affine formal $R$-scheme of ff type, and let $\fX$ be a proper formal $\fS$-scheme. An envelope of $\fX$ over $\fS$ is a proper formal $\fS_\pi$-scheme $\fX_\pi$ such that $\fX\rightarrow\fS$ is obtained from $\fX_\pi\rightarrow\fS_\pi$ by completion along a subscheme of definition of $\fS$.
\end{defi}

In the above situation, the envelope of $\fX$ over $\fS$ is unique if it exists. Indeed, after reducing modulo powers of $\pi$ this is an immediate consequence of Grothendieck's Formal Existence Theorem, cf.\ \cite{EGAIII1} 5.4.1. Moreover, it is clear that $\fX_\pi$ exists when $\fX\rightarrow\fS$ is the completion of a proper morphism of schemes $X\rightarrow S$, where $S$ denotes the spectrum of the ring of global functions on $\fS$. Indeed, in this case $\fX_\pi\rightarrow\fS_\pi$ is simply obtained from $X\rightarrow S$ via $\pi$-adic formal completion.

As we have already pointed out, it is possible -- and will be necessary -- to define the generic fiber\index{generic fiber!of an envelope} of an envelope, at the expense of passing to Huber's category of adic spaces. For an introduction to the theory of adic spaces, we refer to \cite{Hubermainadspacepaper} and \cite{Huberbuch} Chapter 1. The adic space $t(\fX)$ associated to a locally noetherian formal scheme $\fX$ is defined in \cite{Hubermainadspacepaper} 4.1. By \cite{Hubermainadspacepaper}, the functor  $t:\fX\mapsto t(\fX)$ is fully faithful and commutes with open immersions. If $\fX=\Spf A$ is affine, then $t(\fX)\mathrel{\mathop:}=\Spa(A,A)$ is affinoid. By \cite{Huberbuch} 1.2.2, the fibered product 
\[
t(\fX)\otimes_RK\,\mathrel{\mathop:}=\,t(\fX)\times_{\Spa(R,R)}\Spa(K,R)
\]
exists for any locally noetherian formal $R$-scheme $\fX$. Moreover, it is seen from the explicit construction in the proof of \cite{Huberbuch} 1.2.2 that $t(\fX)\otimes_RK$ is \emph{affinoid} whenever $\fX$ is affine and \emph{adic} over $R$. Let us point out that $t(\fX)\otimes_RK$ needs not even be quasi-compact when $\fX$ is not adic over $R$. Indeed, if $\fX$ is of locally ff type over $R$, it is seen from the construction in the proof of \cite{Huberbuch} 1.2.2 that $t(\fX)\otimes_RK$ is the adic space $(\fX^\rig)^\ad$ associated to the Berthelot generic fiber  of $\fX$, cf.\ also \cite{RZ} Remark 5.11. Here $\ad$ denotes the natural fully faithful functor from the category of rigid $K$-spaces to the category of adic $\Spa(K,R)$-spaces. We refer to \cite{Hubermainadspacepaper} 4.3 for the definition of $\ad$.

If $\fX$ is a locally noetherian formal $R$-scheme, not necessarily of locally ff type, the adic space $t(\fX)\otimes_RK$ is called the \emph{generic fiber} of $\fX$. We will briefly write
\[
\fX_K\,\mathrel{\mathop:}=\,t(\fX)\otimes_RK\quad.
\]

		\subsection{The schematic closure of a graph is proper}\label{graphpropersec}
				Let $\fX$ and $\fY$ be flat formal $R$-schemes of locally ff type with uniformly rigid generic fibers $X$ and $Y$. Let us assume that $\fX$ is separated; then $X$ is separated as well. Let $\phi\colon Y\rightarrow X$ be a uniformly rigid morphism, and let $\ul{\Gamma}_\phi$ denote the schematic closure of the graph of $\phi$ in $\fY\times\fX$. In this section, we prove that the projection $p_\fY|_{\ul{\Gamma}_\phi}$ from $\ul{\Gamma}_\phi$ to $\fY$ is \emph{proper}.

\begin{lem}\label{graphclosureadiclem}
The projection $p_\fY|_{\ul{\Gamma}_\phi}$ is \emph{adic}.
\end{lem}
\begin{proof}
Indeed, $p_\fY|_{\ul{\Gamma}_\phi}$ induces an isomorphism of uniformly rigid generic fibers and, hence, an isomorphism of rigid-analytic generic fibers. Since isomorphisms are quasi-compact, it follows from Proposition \ref{rqcadicprop} that $p_\fY|_{\ul{\Gamma}_\phi}$ is adic.
\end{proof}



\begin{lem}\label{graphproperlem}
Let us assume that there exists a proper morphism $\ul{\psi}\colon\fY'\rightarrow\fY$ inducing an isomorphism of uniformly rigid generic fibers such that $\phi$ extends to a morphism $\ul{\phi}\colon\fY'\rightarrow\fX$. Then $p_\fY|_{\ul{\Gamma}_\phi}$ is proper.
\end{lem}
\begin{proof}
By Lemma \ref{schemcldirimlem}, $\ul{\Gamma}_\phi$ coincides with the schematic image of the proper morphism
\[
(\ul{\psi}\times\id_\fX)\circ\Gamma_{\ul{\phi}}\quad\colon\quad\fY'\hookrightarrow\fY'\times\fX\rightarrow\fY\times\fX\quad.
\]
Let $\ul{\phi}'\colon\fY'\rightarrow\ul{\Gamma}_\phi$ denote the unique factorization of this morphism; then $\ul{\phi}'$ is an isomorphism on uniformly rigid generic fibers. Let us consider the commutative diagram
\[
\begin{diagram}
\fY'&\rTo^{\ul{\phi}'}&\ul{\Gamma}_\phi\\
&\rdTo<{\ul{\psi}}&\dTo>{p_\fY|_{\ul{\Gamma}_\phi}}\\
&&\fY
\end{diagram}
\]
of morphisms inducing isomorphisms of uniformly rigid generic fibers. By Proposition \ref{rqcadicprop}, the morphisms $\ul{\phi}'$ and $p_\fY|_{\ul{\Gamma}_\phi}$ are of tf type, cf. the proof of Lemma \ref{graphclosureadiclem}. By Lemma \ref{schemcldirimlem}, $\ul{\phi}'$ is surjective onto the physical points of $\ul{\Gamma}_\phi$. Since $\fX$ is separated, $p_\fY$ is separated, and hence $p_\fY|_{\ul{\Gamma}_\phi}$ is separated as well. After reducing the above diagram of adic morphisms of locally noetherian formal schemes modulo an ideal of definition, we conclude from \cite{EGAII} 5.4.3 ($ii$) that $p_\fY|_{\ul{\Gamma}_\phi}$ is proper, as desired.
\end{proof}


The conclusion of Lemma \ref{graphproperlem} holds without the assumption on the existence of the morphism $\ul{\psi}$:

\begin{theorem}\label{graphproperthm}
The morphism $p_\fY|_{\ul{\Gamma}_\phi}\colon\ul{\Gamma}_\phi\rightarrow\fY$ is proper. 
\end{theorem}
\begin{proof}
By Lemma \ref{graphclosureadiclem}, $\ul{\Gamma}_\phi$ is adic over $\fY$. To show that $p_\fY|_{\ul{\Gamma}_\phi}$ is proper, we may assume that $\fY$ is affine such that $Y$ is semi-affinoid. Let $(\fX_i)_{i\in I'}$ be an affine open covering of $\fX$, and let $(X_i)_{i\in I'}$ be the induced admissible covering of $X$.  Let $(Y_i)_{i\in I}$ be a leaflike refinement of $(\phi^{-1}(X_i))_{i\in I'}$, and let $(Y_j)_{j\in J}$ be a treelike covering of $Y$ enlarging $(Y_i)_{i\in I}$ together with a suitable rooted tree structure on $J$ such that $I$ is identified with the set of leaves of $J$. Let us choose a formal presentation 
\[
(\fY_j,\ul{\phi}_j,\ul{\beta}_j)_j
\]
of $(Y_j)_{j\in J}$ with respect to $\fY$ according to Lemma \ref{simpadlem}, where we may assume that for each leaf $j$ of $J$ the restriction 
\[
\phi|_{Y_j}\colon Y_j\rightarrow X
\]
extends uniquely to a morphism $\ul{\phi}_j\colon\fY_j\rightarrow\fX$. We forget that the morphisms 
\[
\ul{\beta}_j\colon\ul{\fY}_j'\rightarrow\ul{\fY}_j
\]
are admissible formal blowups; we merely retain the information that they are proper and that they induce isomorphisms of uniformly rigid generic fibers. Let us argue by induction on the volume of $J$, the statement being trivially true in the leaves of $J$. By our induction hypothesis, for any child $j$ of the root $r$ of $J$, the schematic closure $\ul{\Gamma}_{\phi_j}$ of the graph of $\phi_j\mathrel{\mathop:}=\phi|_{Y_j}$ in $\fY_j\times\fX$ is proper over $\fY_j$. By uniqueness, the formation of schematic closures commutes with localization; hence the $\ul{\Gamma}_{\phi_j}$ glue to a closed formal subscheme $\ul{\Gamma}'_{\phi_r}\subseteq\fY_r'\times\fX$ that is proper over $\fY_r'$ and, hence, proper over $\fY_r$. We are now in the situation of Lemma \ref{graphproperlem}, with $\fY'=\ul{\Gamma}'_{\phi_r}$, so we may conclude that $\ul{\Gamma}_\phi$ is indeed proper over $\fY=\fY_r$.
\end{proof}

\begin{cor}\label{graphpropercor}
Under the assumptions stated at the beginning of this section, there exists a proper morphism $\ul{\psi}\colon\fY'\rightarrow\fY$ inducing an isomorphism of uniformly rigid generic fibers such that $\phi$ extends to a morphism $\ul{\phi}\colon\fY'\rightarrow\fX$.
\end{cor}

Indeed, by Theorem \ref{graphproperthm} it suffices to take 
\[
\fY'=\ul{\Gamma}_\phi\quad,\quad\ul{\psi}=p_\fY|_{\ul{\Gamma}_\phi}\quad.
\]
Let us point out that the proper morphism $p_\fY|_{\ul{\Gamma}_\phi}$ needs \emph{not} be induced, locally on $\fY$, by proper morphisms of schemes via formal completion. 
				
		\subsection{The envelope of a graph}\label{graphenvsec}
				In the situation studied in the previous section, let us additionally assume that the formal $R$-scheme $\fX$ is  of locally \emph{tf} type and that $\fY$ is \emph{affine}. It is still true that $p_\fY|_{\ul{\Gamma}_\phi}$ needs not be algebraic. However, in virtue of Grothendieck's Formal Existence Theorem, the \emph{envelope}\index{envelope!of a graph} of $p_\fY|_{\ul{\Gamma}_\phi}$ exists:

\begin{prop}\label{graphenvexprop}
Under the above assumptions, the envelope 
\[
(p_\fY|_{\ul{\Gamma}_\phi})_\pi\colon(\ul{\Gamma}_{\phi})_\pi\rightarrow\fY_\pi
\]
of $p_\fY|_{\ul{\Gamma}_\phi}$ exists. More precisely speaking, $\ul{\Gamma}_\phi$ is induced by a closed formal subscheme 
\[
(\ul{\Gamma}_\phi)_\pi\,\subseteq\,\fY_\pi\times\fX\quad,
\]
and $(p_\fY|_{\ul{\Gamma}_\phi})_\pi$ is given by $p_{\fY_\pi}|_{(\ul{\Gamma}_\phi)_\pi}$.
\end{prop}
\begin{proof}
For $n\in\N$, we let a subscript $n$ indicate reduction modulo $\pi^{n+1}$. The affine formal scheme $\fY_n$ is the completion of the affine scheme $\fY_{\pi,n}$ along a subscheme of definition of $\fY$, and their rings of global sections coincide. The closed formal subscheme $\ul{\Gamma}_\phi\subseteq\fY\times\fX$ induces closed formal subschemes $(\ul{\Gamma}_\phi)_n\subseteq\fY_n\times\fX_n$, and $\fY_n\times\fX_n$ is the completion of $\fY_{\pi,n}\times\fX_n$. By Corollary \ref{graphpropercor}, $(\ul{\Gamma}_\phi)_n$ is proper over $\fY_n$; hence by \cite{EGAIII1} 5.1.8, $(\ul{\Gamma}_\phi)_n$ algebrizes to a uniquely determined closed subscheme 
\[
(\ul{\Gamma}_\phi)_{\pi,n}\subseteq\fY_{\pi,n}\times\fX_n\quad.
\]
By uniqueness, $(\ul{\Gamma}_\phi)_{\pi,n}=(\ul{\Gamma}_\phi)_{\pi,n+1}\cap(\fY_{\pi,n}\times\fX_n)$. By \cite{EGAIn} 10.6.3, 
\[
(\ul{\Gamma}_\phi)_\pi\,\mathrel{\mathop:}=\,\indlim_n\,(\ul{\Gamma}_\phi)_{\pi,n}
\]
is a formal scheme, and we easily see that it is a closed formal subscheme of $\fY_\pi\times\fX$ which is proper over $\fY_\pi$ via $p_{\fY_\pi}$ and whose completion is $\ul{\Gamma}_\phi$.
\end{proof}

		\subsection{Localization of envelopes}\label{envlocsec}
				Before proceeding any further, we need to discuss how envelopes behave with respect to localization. We have already mentioned the fact that the formation of envelopes does not commute with localization: If $\fX$ is an affine formal $R$-scheme of ff type and if $\fU\subseteq\fX$ is an affine open formal subscheme, then $\fU_\pi$ needs not be an open formal subscheme of $\fX_\pi$. However, if $(\fU_i)_{i\in I}$ is an affine open \emph{covering} of $\fX$, it is easily seen that the induced family of morphisms $(\fU_{i,\pi}\rightarrow\fX_\pi)_{i\in I}$ is faithfully flat. Indeed, since these morphisms are adic, faithful flatness is verified modulo powers of $\pi$. In this subsection, we prove a corresponding statement for envelopes of algebraic \emph{proper} morphisms. We recall that algebraizations of proper formal morphisms are unique in virtue of Grothendieck's Formal Existence Theorem, cf.\ \cite{EGAIII1} 5.4.1.
 
\begin{lem}\label{algcovrefinelem}
Let $\fX$ be an affine flat formal $R$-scheme of ff type, let $X$ denote the spectrum of its ring of global functions, let 
\[
\ul{\phi}\colon\fX'\rightarrow\fX
\]
be a proper morphism induced from a proper morphism of schemes $\phi\colon X'\rightarrow X$ via formal completion along an ideal of definition of $\fX$, and let $(\fU_i)_{i\in I}$ be any affine open covering of $\fX'$. Then there exists a finite affine open covering $(V_j)_{j\in J}$ of $X'$ whose completion $(\fV_j)_{j\in J}$ along an ideal of definition of $\fX$ refines $(\fU_i)_{i\in I}$.
\end{lem}
\begin{proof}
Let $X_0\subseteq X$ be a closed subscheme corresponding to a subscheme of definition of $\fX$, and let $X'_0$ denote its $\phi$-preimage; then $\fX'$ is identified with the completion of $X'$ along $X'_0$, and $X'_0$ is identified with a subscheme of definition of $\fX'$. A locally noetherian formal scheme is affine if and only if it admits an affine subscheme of definition, cf.\ \cite{EGAIn} 2.3.5; hence we conclude that the affine open coverings of $\fX'$ correspond to the affine open coverings of $X'_0$. We may assume that $I$ is finite, $X'_0$ being quasi-compact. For each $i\in I$, we extend $\fU_i\cap X'_0$ to an open subset $U_i$ of $X'$. Since $X'$ is noetherian, $U_i$  is quasi-compact and, hence, admits a finite affine open covering $(V_{ij})_{j\in J_i}$. Let $(V_j)_{j\in J}$ be the disjoint union of these families; we claim that it covers $X'$. Indeed, the closed complement of the union of the $V_j$ in $X'$ does not meet $X'_0$; hence its closed image in $X$ does not meet $X_0$. However, since the ring of functions on $X$ is complete for the topology that is defined by the ideal corresponding to $X_0$, every closed point of $X$ lies in $X_0$. It follows that the complement of the union of the $V_j$ in $X'$ is empty, as desired.
\end{proof}

\begin{prop}\label{thefaithflatprop}
Let $\fX$ be an affine flat formal $R$-scheme of ff type, and let 
\[
\ul{\phi}\colon\fX'\rightarrow\fX
\]
be a proper morphism that induces an isomorphism of rigid generic fibers and that is obtained from a proper morphism of schemes $\phi\colon X'\rightarrow X$, where $X$ denotes the spectrum of the ring of global functions on $\fX$. Let $(\fU_i)_{i\in I}$ be a finite affine open covering of $\fX'$ that is induced by a finite affine open covering $(U_i)_{i\in I}$ of $X'$. For each $i\in I$, we let $U_i'$ denote the spectrum of the ring of functions on $\fU_i$; then the canonical family of morphisms of schemes
\[
\psi_i\colon U_i'\rightarrow U_i\subseteq X'
\]
is faithfully flat.
\end{prop}
\begin{proof}
The $\psi_i$ are flat, adic completions of noetherian rings being flat. Since the $U_i$ cover $X'$, it suffices to show that for any $i\in I$, the restricted family is faithfully flat over $U_i$; hence we must show that for every $i\in I$ and every closed point $x\in U_i$, there exists an index $j\in I$ such that $x\in\im\psi_j$. We must beware of the fact that the property of a point in $U_i$ of being closed is \emph{not} local on $U_i$. Let $y$ denote the $\phi$-image of $x$ in $X$, let $A$ be the ring of functions on $\fX$, let $\a$ be an ideal of definition for $A$, and let $\p\subseteq A$ be the prime ideal corresponding to $y$. If $y$ is a \emph{closed} point of $X$, then $y\in X_0\mathrel{\mathop:}=V(\a)$, and hence $x\in\phi^{-1}(V(\a))$ which is identified with a subscheme of definition of $\fX'$. Since $\fX'$ is covered by the $\fU_j$, there exists an index $j\in J$ such that $x\in\fU_j$, and we obtain $x\in\im\psi_j$. Let us now consider the case where $y\in X$ is \emph{not} closed.

Let $B$ denote the ring of functions on $U_i$, and let $\m\subseteq B$ be the maximal ideal corresponding to $x$. Let $M$ and $L$ denote the residue fields of $B$ and $A$ in $x$ and $y$ respectively. Since $M=B/\m$ is of finite type over $A/\p$, it is of finite type over $L=Q(A/\p)$ and, hence, finite over $L$. By \cite{AT} Theorem 1, it follows that $L$ is of finite type over $A/\p$. By \cite{AT} Theorem 4, we conclude that $A/\p$ has only a finite number of prime ideals and that each nonzero prime ideal of $A/\p$ is maximal. Since $\p$ is not maximal in $A$ by our assumption on $y$, $A/\p$ has Krull dimension $1$. The maximal ideals of $A/\p$ correspond to the maximal ideals of $A$ which contain $\p$; hence they are open in the separated and complete $\a(A/\p)$-adic quotient topology on $A/\p$. We know from Lemma \ref{jacobsontoplem} that this topology coincides with the Jacobson-adic topology. Since $A/\p$ is integral, we conclude from \cite{Bourbaki} III.2.13 Corollary of Proposition 19 that the one-dimensional domain $A/\p$ is local and complete in its maximal-adic topology.

Let $R'$ denote the integral closure of $A/\p$ in $L$. Since $A$ is excellent, $R'$ is finite over $A/\p$, and we easily see that $R'$ is a complete discrete valuation ring whose valuation topology coincides with the $\a$-adic topology. Since $R'$ is finite over $A$, the rings of functions on the affine $R'$-schemes $U_i'\otimes_AR'$ are $\a$-adically complete. Hence, we may perform the base change $\cdot\otimes_AR'$ and thereby assume that $A=R'$ is a complete discrete valuation ring such that $y\in\Spec R'$ is the generic point. 

We are now in a classical rigid-analytic situation over $R'$. The closed points in the generic fiber $X'_y$ of $X'$ correspond to the points of the rigid analytification ${X'_y}^\an$ of $X'_y$, which is the generic fiber of the completion $\fX'$ of $X'$ since $X'$ is proper over $R'$. Since $x$ is closed in $X'$ and since ${X'_y}^\an$ is covered by the maximal spectra of the rings of functions on the generic fibers of the $U'_j$, we see that there exists a point in some $U_j'$ mapping to $x$, as desired.
\end{proof}

\begin{cor}\label{thefaithflatcor}
Let $\fX$ be an affine flat formal $R$-scheme of ff type, let $X$ denote the spectrum of the ring of functions on $\fX$, and let 
\[
\ul{\phi}\colon\fX'\rightarrow\fX
\]
be a proper morphism that induces an isomorphism of rigid generic fibers and that is induced by a proper morphism of schemes $\phi\colon X'\rightarrow X$ via formal completion along an ideal of definition of $\fX$. Any affine open covering of $\fX'$ admits a finite affine open refinement $(\fU_i)_{i\in I}$ such that the canonical family 
\[
(\fU_{i,\pi}\rightarrow\fX'_\pi)_{i\in I}
\]
is faithfully flat.
\end{cor}

Indeed, this follows from Lemma \ref{algcovrefinelem} and Proposition \ref{thefaithflatprop} via $\pi$-adic formal completion. The statement applies in particular in the case where $\ul{\phi}$ is an admissible formal blowup.

Let us introduce the following terminology: A point in a locally noetherian $R$-scheme is called $R$-\emph{regular}\index{point!R@$R$-regular} if it admits an affine open neighborhood whose ring of functions is $R$-regular in the sense of \cite{EGAIV2} 6.8.1.

\begin{prop}\label{liftrregpointsprop}
Let $\fX$ be an affine flat formal $R$-scheme of ff type, let $X$ denote the spectrum of the ring of functions on $\fX$, let 
\[
\ul{\phi}\colon\fX'\rightarrow\fX
\]
be an admissible formal blowup, and let $\phi\colon X'\rightarrow X$ be the corresponding admissible blowup of noetherian $R$-schemes. Let $(\fU_i)_{i\in I}$ be a finite affine open covering of $\fX'$ that is induced by a finite affine open covering $(U_i)_{i\in I}$ of $X'$. For each $i\in I$, we let $U_i'$ denote the spectrum of the ring of functions on $\fU_i$. Then for every generic point $\eta\in X_k$ that is $R$-regular in $X$, there exist 
an index $i\in I$ and a generic point $\eta'\in U'_{i,k}$ above $\eta$ that is $R$-regular in $U_i'$.
\end{prop}
\begin{proof}
Since $X_k$ is regular and, hence, reduced in $\eta$, Lemma \ref{genisolem} shows that there exists an open subscheme $U\subseteq X'$ that maps, via $\phi$, isomorphically onto an $R$-regular open neighborhood of $\eta$. By Proposition \ref{thefaithflatprop}, the family of morphisms
\[
\psi_i\colon U_i'\rightarrow U_i\subseteq X'\quad.
\]
is faithfully flat. Since the schemes $U_i$ are of finite type over $X$, they are excellent; hence by \cite{EGAIV2} 7.8.3 (v), the morphisms $\psi_i$ are regular. For each $i\in I$ we let $V_i\subseteq U_i'$ denote the $\psi_i$-preimage of $U$; then the disjoint union of the $V_i$ is $R$-regular and faithfully flat over $U$. Let $i\in I$ be an index such that $\eta$ lies in the image of $V_i$. Let $\xi\in V_{i,k}$ be a point above $\eta$, and let $\eta'$ be a maximal generalization of $\xi$ in $V_{i,k}$. Then $\eta'$ is a generic point of $U'_{i,k}$ that is $R$-regular in $U_i$ and that lies above $\eta$.
\end{proof}

In the situation of the above Proposition \ref{liftrregpointsprop}, let us assume that both $\fX$ and $\fX'$ are normal. Then by \cite{EGAIV2} 7.8.3 (v) the $U_i'$ are normal as well, and we easily see that the completed stalks of $\fU_{i,\pi}$ and $\fX_\pi$ in $\eta'$ and $\eta$ respectively are complete discrete valuation rings.

		\subsection{Generically defined morphisms}\label{gendefmorph}
				Let us recall the situation that we considered in Section \ref{graphenvsec}: Let $\fY$ be an affine flat formal $R$-scheme of ff type with uniformly rigid generic fiber $Y$, let $\fX$ be a separated flat formal $R$-scheme of locally tf type with uniformly rigid generic fiber $X$, let $\phi\colon Y\rightarrow X$ be a uniformly rigid morphism, and let $\ul{\Gamma}_\phi$ denote the schematic closure of the graph of $\phi$ in $\fY\times\fX$. In Section \ref{graphpropersec} and Section \ref{graphenvsec}, we have seen that the restricted projection $p_\fY|_{\ul{\Gamma}_\phi}$ is a proper morphism which has an envelope $(p_\fY|_{\ul{\Gamma}_\phi})_\pi$. In this section, we prove that if $\fY$ is smooth, then $(p_\fY|_{\ul{\Gamma}_\phi})_\pi$ is an isomorphism above the generic points of $\fY_\pi$.

\begin{theorem}\label{complfiberisothm}
Let us assume that $\fY$ is generically normal, and let $\eta$ be a generic point in $\fY_{\pi,k}$ that is $R$-regular in the spectrum of the ring of functions on $\fY$. Then the restricted projection 
\[
(p_\fY|_{\ul{\Gamma}_\phi})_\pi\colon(\ul{\Gamma}_\phi)_\pi\rightarrow\fY_\pi
\]
induces an isomorphism under the base change $\cdot\times_{\fY_\pi}\Spf\hat{\O}_{\fY_\pi,\eta}$.
\end{theorem}
\begin{proof}
If $\phi$ extends to a morphism from $\fY$ to $\fX$, then $p_\fY|_{\ul{\Gamma}_\phi}$ is an isomorphism, and there is nothing to show. Let $(Y_i)_{i\in I}$ be a treelike covering of $Y$ together with a treelike formal covering $(\fY_i)_{i\in I}$ with respect to $\fY$ inducing $(Y_i)_{i\in I}$ such that for each leaf $i\in I$, the restriction of $\phi$ to $Y_i$ extends to a morphism from $\fY_i$ to $\fX$; such a covering exists in virtue of Proposition \ref{ffprop}, Lemma \ref{simpadlem} and Corollary \ref{freesemaffcor} ($iv$). Since $\fY$ is generically normal, Corollary \ref{normisblowupcor} shows that we may assume all $\fY_i$, $i\in I$, to be \emph{normal}. After possibly passing to a refinement, Lemma \ref{algcovrefinelem} allows us to assume that all affine open coverings in this treelike formal covering are induced by affine open coverings on the respective admissible algebraic blowups. By Proposition \ref{liftrregpointsprop}, there exists a child $i$ of the root of $I$ together with a point $\eta'$ in the spectrum of the ring of functions on $\fY_i$ such that $\eta'$ is $R$-regular in this spectrum, such that $\eta'$ is generic in the special fiber $\fY_{i,\pi,k}$ and such that $\eta'$ lies above $\eta$. Let $R_{\eta'}$ and $R_\eta$ denote the completed stalks of $\fY_{i,\pi}$ and $\fY_\pi$ in $\eta'$ and $\eta$ respectively; according to the remark following Proposition \ref{liftrregpointsprop}, $R_{\eta'}/R_\eta$ is a local extension of complete discrete valuation rings. Let $\ul{\Gamma}_{\phi,i}$ denote the schematic closure of the graph of $\phi|_{Y_i}$ in $\fY_i\times\fX$, and let $(\ul{\Gamma}_{\phi,i})_\pi\subseteq\fY_{i,\pi}\times\fX$ denote its envelope.
Using induction on the volume of $I$, we may assume that $p_{\fY_{i,\pi}}|_{(\ul{\Gamma}_{\phi,i})_\pi}$ becomes an isomorphism under the base change
\[
\cdot\times_{\fY_{i,\pi}}\Spf R_{\eta'}\quad.
\]
The pullback $\ul{\Gamma}_\phi\times_\fY\fY_i$ is a closed formal subscheme of $\fY_i\times\fX$ which is proper over $\fY_i$, and its generic fiber is the graph of $\phi|_{Y_i}$; hence 
\[
\ul{\Gamma}_{\phi,i}\,\subseteq\,\ul{\Gamma}_\phi\times_\fY\fY_i
\]
is the closed formal subscheme defined by the $\pi$-torsion ideal. Since completions of locally noetherian schemes are flat and since envelopes are unique, it follows that
\[
(\ul{\Gamma}_{\phi,i})_\pi\,\subseteq\,(\ul{\Gamma}_\phi\times_\fY\fY_i)_\pi
\]
is likewise the closed formal subscheme defined by the $\pi$-torsion ideal. Since completion commutes with fibered products and since envelopes are unique, we have a natural identification 
\[
(\ul{\Gamma}_\phi\times_\fY\fY_i)_\pi\,=\,(\ul{\Gamma}_{\phi})_\pi\times_{\fY_\pi}\fY_{i,\pi}\quad.
\]
We have to show that the morphism $(\ul{\Gamma}_\phi)_\pi\times_{\fY_\pi}\Spf R_\eta\rightarrow \Spf R_\eta$ is an isomorphism. Since it is of tf type and since $R_{\eta'}$ is faithfully flat over $R_\eta$, descent theory reduces the problem to proving that
\[
(\ul{\Gamma}_\phi)_\pi\times_{\fY_\pi}\Spf R_{\eta'}\rightarrow\Spf R_{\eta'}
\]
is an isomorphism. We observe that $R_{\eta'}$ is flat over $R_\eta$ and, hence, over $\fY_\pi$; in particular, the domain of the above morphism is \emph{flat} over $R$. Moreover, since $(\ul{\Gamma}_{\phi,i})_\pi\subseteq(\ul{\Gamma}_\phi)_\pi\times_{\fY_\pi}\fY_{i,\pi}$ is defined by $\pi$-torsion and since $R_{\eta'}$ is flat over $\fY_{i,\pi}$, the closed formal subscheme
\begin{eqnarray*}
(\ul{\Gamma}_{\phi,i})_\pi\times_{\fY_{i,\pi}}\Spf R_{\eta'}&\subseteq&((\ul{\Gamma}_\phi)_\pi\times_{\fY_\pi}\fY_{i,\pi})\times_{\fY_{i,\pi}}\Spf R_{\eta'}\\
&\cong&(\ul{\Gamma}_\phi)_\pi\times_{\fY_\pi}\Spf R_{\eta'}
\end{eqnarray*}
is defined by $\pi$-torsion as well. We conclude that the above closed immersion is in fact an isomorphism.  By our induction hypothesis, the projection to $\Spf R_{\eta'}$ is an isomorphism, as desired.
\end{proof}

In the above proof, the $\pi$-torsion in $(\ul{\Gamma}_\phi)_\pi\times_{\fY_\pi}\fY_{i,\pi}$ is due to the fact that admissible formal blowups are in general not flat. We have used the fact that under the above assumptions, they are flat locally on their envelopes.

\begin{cor}\label{genisocor}
In the situation of Theorem \ref{complfiberisothm}, the morphism $(p_\fY|_{\ul{\Gamma}_\phi})_\pi$ is an isomorphism over an open neighborhood of $\eta$.
\end{cor}
\begin{proof}
Let us abbreviate $\phi\mathrel{\mathop:}=(p_\fY|_{\ul{\Gamma}_\phi})_\pi$; we have to show that there exists an open neighborhood $\fU$ of $\eta$ in $\fY_\pi$ such that the restriction $\phi^{-1}(\fU)\rightarrow \fU$ of $\phi$ is an isomorphism. Let $k(\eta)$ denote the residue field of $\fY_\pi$ in $\eta$. By Theorem \ref{complfiberisothm}, we know that the special fiber $\phi_k$ of $\phi$ becomes an isomorphism under the base change $\cdot\times_{\fY_{\pi,k}}\Spec k(\eta)$. Since $\phi_k$ is proper, it follows from \cite{EGAIV4} 18.12.7 that there exists an open neighborhood $U_k$ of $\eta$ in $\fY_{\pi,k}$ such that the restriction $\phi_k^{-1}(U_k)\rightarrow U_k$ of $\phi_k$ is a closed immersion. Let $\fU\subseteq\fY_\pi$ be the open formal subscheme whose special fiber is $U_k$. By \cite{EGAIII1} 4.8.1, the restriction $\phi^{-1}(\fU)\rightarrow \fU$ of $\phi$ is a finite morphism. By \cite{Eis} Ex. 7.2, a finite homomorphism of $\pi$-adically complete and separated $R$-algebras that is surjective modulo $\pi$ is surjective; hence $\phi^{-1}(\fU)\rightarrow \fU$ is a closed immersion as well. We may shrink $\fU$ and thereby assume that $\fU$ as affine, $\fU=\Spf A$ for some $\pi$-adically complete $R$-algebra $A$. Then the above restriction of $\phi$ corresponds to an ideal $I$ in $A$. Let $\p\subseteq A$ be the prime ideal corresponding to $\eta\in U_k$. By Theorem \ref{complfiberisothm} and using faithfully flat descent, we see that $A_\p/IA_\p=A_\p$, that is, $I A_\p=0$. Since $A$ is noetherian, $I$ is finitely generated. After shrinking $\fU$ further, we may thus assume that the restriction of $\phi$ is an isomorphism, as desired.
\end{proof}

	\section{The Weil Extension Theorem}\label{weilextsec}

		In this section, we will investigate conditions under which a morphism of uniformly rigid $K$-spaces extends to given formal $R$-models of locally ff type of given uniformly rigid subspaces, particularly in the case when we are dealing with group objects. Let us note that extensions of uniformly rigid morphisms to \emph{flat} formal models are uniquely determined, the functor $\urig$ being faithful on the category of flat formal $R$-schemes of locally ff type.

		\subsection{Reduction to the formal fibers}\label{redformfibsec}
			Let $\fX$ and $\fY$ be flat formal $R$-schemes of locally ff type, and let $\phi\colon Y\rightarrow X$ be a morphism of uniformly rigid generic fibers. We begin by proving, in the case where $\fY$ is normal or where $\fY_k$ is reduced, that $\phi$ extends to an open neighborhood of a closed point $y\in\fY$ if and only if $\phi$ extends to the completion $\fY|_y$ of $\fY$ in $y$. This fact will later allow us to reduce to the case where $\fY$ is local. 

\begin{lem}\label{affineextlem}
Let us assume that $\fX$ is \emph{affine} and that $\fY$ is normal or that $\fY_k$ is reduced. Then $\phi$ extends uniquely to a morphism
$\ul{\phi}\colon\fY\rightarrow\fX$.
\end{lem}
\begin{proof}
We may work locally on $\fY$ and thereby assume that $\fY$ is affine as well; then both $X$ and $Y$ are semi-affinoid. Let us write $\fX=\Spf \ul{A}$, $\fY=\Spf \ul{B}$, $X=\sSp A$ and $Y=\sSp B$. By Proposition \ref{ffprop} and Corollary \ref{freesemaffcor} (iv), $\phi$ extends to a morphism $\ul{\phi}'\colon\fY'\rightarrow\fX$, where $\fY'=\Spf \ul{B}'$ is a finite admissble formal blowup of $\fY$. We obtain inclusions 
\[
\ul{B}\subseteq\ul{B}'\subseteq \mathring{B}\quad,
\]
so it suffices to show that $\ul{B}=\mathring{B}$. If $\fY$ is normal, this follows from Proposition \ref{normallatticeprop}. If $\fY_k$ is reduced, we argue as follows: Let $f\in B$ be power-bounded, and let $n$ be the natural number with the property that $\pi^nf\in \ul{B}-\pi \ul{B}$. If $n\geq 1$, then $|\pi^n f|_\sup<1$, and hence the class of $\pi^n f$ modulo $\pi$ is nilpotent. This cannot be, since $\ul{B}_k$ is reduced and since $\pi^nf\notin\pi \ul{B}$. Hence $n=0$, which means that $f\in \ul{B}$.
\end{proof}




If $V\subseteq\fY$ is a closed or an open subset, we define
\[
]V[_\fY\,\mathrel{\mathop:}=\,\sp_\fY^{-1}(V)\quad,
\]
adopting the notation introduced by Berthelot in the rigid-analytic setting, cf.\ \cite{Berth} 1.1.2. The subset $]V[_\fY$ is called the \emph{formal fiber}\index{formal fiber} of $\fY$ over $V$. It is an open uniformly rigid subspace of $Y$. Indeed, if $V$ is closed, then 
\[
]V[_\fY\,=\,(\fY|_V)^\srig\quad,
\]
and if $V$ is open, then 
\[
]V[_\fY\,=\,V^\srig\quad,
\]
where $V$ is equipped with the open formal subscheme structure. We see that if $V$ is affine, then $]V[_\fY$ is semi-affinoid. 

\begin{prop}\label{redtoffprop}
Let $y\in \fY$ be a closed point where $\fY$ is normal or where $\fY_k$ is reduced. Then the following are equivalent:
\begin{packed_enum}
\item $\phi|_{]y[_\fY}$ extends to a morphism $\ul{\phi}_y\colon\fY|_y\rightarrow\fX$.
\item $\phi$ extends to an open neighborhood of $y$ in $\fY$.
\end{packed_enum}
\end{prop}
\begin{proof}
The implication ($ii$)$\Rightarrow$($i$) is trivial; we need to prove its converse. Let $\fU\subseteq\fX$ be an \emph{affine} open neighborhood of the $\ul{\phi}_y$-image of the unique physical point of $\fY|_y$; then 
\[
]y[_\fY\,\subseteq\,\phi^{-1}(]\fU[_\fX)\quad.
\] 
By Lemma \ref{affineextlem}, it suffices to show that there exists an open neighborhood $\fV\subseteq\fY$ of $y$ such that $]\fV[_\fY\subseteq\phi^{-1}(]\fU[_\fX)$. Since the underlying topological space of $\fY$ is a Jacobson space, it thus suffices to see that the set of closed points $y'\in \fY$ satisfying
\[
]y'[_\fY\subseteq\phi^{-1}(]\fU[_\fX)
\]
 is \emph{open} in the set $\fY_\cl$ of closed points of $\fY$. This set is identified with
\[
\fY_\cl\setminus\sp_\fY(Y\setminus \phi^{-1}(]\fU[_\fX))\quad;
\]
hence, we may conclude our argument by Proposition \ref{mainweilresopenprop} below.
\end{proof}

\begin{prop}\label{mainweilresopenprop}
Let $\fX$ and $\fY$ be flat formal $R$-schemes of locally ff type with uniformly rigid generic fibers $X$ and $Y$ respectively, let $\phi\colon Y\rightarrow X$ be a uniformly rigid morphism, and let $\fU\subseteq\fX$ be an open formal subscheme. Then
\[
\fY_\cl\setminus\sp_\fY(Y\setminus\phi^{-1}(\fU^\srig))
\]
is \emph{open} in the set $\fY_\cl$ of closed points in $\fY$.
\end{prop}
\begin{proof}
Let us first assume that $\phi$ extends to a morphism $\ul{\phi}\colon\fY\rightarrow\fX$. In this case, the set in the statement of the proposition is identified with $\ul{\phi}^{-1}(\fU)_\cl$, and there is nothing to show. Indeed: If $y$ is a closed point in $\fY$ such that $\ul{\phi}(y)\in\fU$, then $\sp_\fY^{-1}(y)$ maps to $\fU^\srig$ under $\phi$, the specialization map being functorial. Conversely, if $\sp_\fY^{-1}(y)$ is mapped to $\fU^\srig$ via $\phi$, then $y$ maps to $\fU$ under $\ul{\phi}$, again by functoriality of the specialization map and since $\sp_\fY^{-1}(y)$ is non-empty, $\fY$ being $R$-flat.

In the general case, we choose a treelike formal covering $(\fY_i)_{i\in I}$ of $\fY^\srig$ with respect to $\fY$ such that for each leaf $i\in \leaves(I)$, the restriction of $\phi$ to $Y_i\mathrel{\mathop:}=\fY_i^\srig$ extends to a morphism 
\[
\ul{\phi}_i\colon\fY_i\rightarrow\fX\quad;
\]
such a treelike formal covering exists in virtue of Proposition \ref{ffprop}, Lemma \ref{simpadlem} and Corollary \ref{latticeexcor} ($iv$). By what we have seen so far, the set
\[
S_i\,\mathrel{\mathop:}=\,(\fY_i)_\cl\setminus \sp_{\fY_i}(Y_i\setminus(Y_i\cap\phi^{-1}(\fU^\srig)))
\]
is open in $(\fY_i)_\cl$ for each \emph{leaf} $i\in I$. We claim that the same statement holds in fact for all $i\in I$. Using induction on the volume of $I$, we may assume that it holds for all children $i$ of the root $r$ of $I$. It remains to show that the statement holds for the root $r$ itself. Let $\ul{\psi}\colon\fY'\rightarrow\fY$ denote the admissible formal blowup of $\fY$ that is part of the chosen treelike formal covering; for each $i\in\children(r)$, we consider $\fY_i$ as an open formal subscheme of $\fY'$. It is clear that
\[
S'\,\mathrel{\mathop:}=\,\fY'_\cl\setminus\sp_{\fY'}(Y\setminus\phi^{-1}(\fU^\srig))
\]
is the union of the sets $S_i$ considered above, where $i$ varies in $\children(r)$. Let us write
\[
S\,\mathrel{\mathop:}=\,\fY_\cl\setminus\sp_\fY(Y\setminus\phi^{-1}(\fU^\srig))
\]
to denote the subset of $\fY_\cl$ that we are interested in. We claim that 
\[
S=\fY_\cl\setminus\ul{\psi}(\fY'_\cl\setminus S')\quad;
\]
then the openness of $S$ follows from the fact that $\ul{\psi}$ is proper and, hence, a closed map of underlying topological spaces. It remains to justify the above claim. Let $y$ be a closed point in $\fY$. Since the specialization map is functorial and since $\ul{\psi}^\srig$ is bijective, the formal fiber of $y$ is identified with union of the formal fibers of the closed points $y'$ in $\ul{\psi}^{-1}(y)$. Hence, $y\in S$ is equivalent to $\ul{\psi}^{-1}(y)\subseteq S'$, where by abuse of notation we consider $\ul{\psi}$ as a map on sets of closed points. The claim is now obvious.
\end{proof}

		\subsection{Tubes around closed uniformly rigid subspaces}\label{tubessec}
			Before proceeding further, we need some more preparations. In this section, which may be considered as an interlude on uniformly rigid geometry, we prove that if $X$ is a semi-affinoid $K$-space and if $Z\subseteq X$ is a closed uniformly rigid subspace that is contained in a weakly retrocompact admissible open subspace $U\subseteq X$, then $U$ contains a so-called \emph{tube} around $Z$ in $X$. For the corresponding statement in rigid geometry, cf.\ \cite{KiehldeRh} Satz 1.6 and \cite{Kisin} Lemma 2.3. In our proof, we follows Kisin's approach. In doing so, we have to overcome the problem that his arguments involve a nontrivial existence result concerning formal models of tf type for certain admissible open subspaces, cf.\ \cite{BL1} 4.4, which is not literally available in the uniformly rigid setting.

The following Proposition generalizes Proposition \ref{zarsrigprop}. The corresponding statement for affinoid spaces is proven in \cite{Kisin} Lemma 2.1. 

\begin{prop}\label{zarretromixadmprop}
Let $X$ be a semi-affinoid $K$-space, and let $(X_i)_{i\in I}$ be a finite covering of $X$, where for each $i\in I$, $X_i\subseteq X$ is
\begin{packed_enum}
\item Zariski-open in the sense of Section \ref{gtopsec} or
\item a retrocompact semi-affinoid subdomain.
\end{packed_enum}
Then the covering $(X_i)_{i\in I}$ of $X$ has a retrocompact refinement; in particular, it is $\sT_\srig$-ad\-mis\-si\-ble.
\end{prop}
\begin{proof}
Let us write $I$ as a disjoint union of subsets $I_z$ and $I_r$ such that $X_i$ is Zariski-open in $X$ for every $i\in I_z$ and such that $X_i$ is a retrocompact semi-affinoid subdomain in $X$ for every $i\in I_r$. After passing to a refinement, we may assume that for each $i\in I_z$, $X_i=D(f_i)$ for some semi-affinoid function $f_i$ on $X$. We proceed by induction on the cardinality of $I_z$. If $I_z$ is empty, there is nothing to show; so let us assume that $J_z$ is nonempty. For each $j\in I_z$, we let $Z_j=V(f_j)\subseteq X$ denote the Zariski-closed complement of $X_j$; we consider it as a reduced closed semi-affinoid subspace of $X$. By our induction hypothesis, the covering $(X_i\cap{Z_j})_{i\in I\setminus\{j\}}$ of $Z_j$ has a retrocompact refinement $(Z_{jl})_{l\in L_j}$. Let us choose a refinement map $\alpha_j\colon L_j\rightarrow I\setminus\{j\}$. Let $l$ be an element of $\alpha_j^{-1}(I_z\setminus\{j\})$; then the retrocompact subdomain $Z_{jl}\subseteq Z_j$ is contained in the Zariski-open subset $Z_j\cap X_{\alpha_j(l)}$ of $Z_j$. We claim that there exists a retrocompact semi-affinoid subdomain $X_{jl}'$ in $X$ that is contained in $X_{\alpha_j(l)}$ and that contains $Z_{jl}$. Indeed, by Proposition \ref{zarsrigprop} and the subsequent remark, there exists a $\sT_\srig$-admissible covering of $X_{\alpha_j(i)}$ by retrocompact semi-affinoid subdomains $X_{\alpha_j(l),\geq\varepsilon}\subseteq X$, with $\varepsilon\in\sqrt{|K^*|}$ tending to zero. It follows that $(X_{\alpha_j(l),\geq\varepsilon}\cap Z_j)_{\varepsilon\rightarrow 0}$ is an admissible covering of $Z_j\cap X_{\alpha_j(l)}$. Since $Z_{jl}$ is quasi-compact, it suffices to take $X'_{jl}=X_{\alpha_j(l),\geq\varepsilon}$ with $\varepsilon>0$ small enough. We may enlarge the covering $(X_i)_{i\in I}$ by adding the finitely many retrocompact subdomains $X'_{jl}$; indeed, the new finite covering that we obtain is a refinement of $(X_i)_{i\in I}$ since $X'_{jl}\subseteq X_{\alpha_j(l)}$. Hence, we may assume that
\[
\bigcup_{j\in I_z}Z_j\subseteq\bigcup_{i\in I_r}X_i\quad.
\]
Let $\fX$ be a flat affine formal $R$-model of ff type for $X$. For each $i\in I_r$, we consider a formal presentation $\fX_i\rightarrow\fX$ of $X_i$ as a retrocompact semi-affinoid subdomain of $X$ with respect to $\fX$, together with a factorization into open immersions and admissible formal blowups. Let $J$ denote the almost linear rooted tree obtained by glueing the resulting linear rooted trees along their roots, and let $(\fX_j)_{j\in J}$ be the associated formal data. We construct a rooted tree $J'$, as follows: Let $J'$ first be a copy of $J$, and let $r\in J$, $r'\in J'$ denote the respective roots. For each $j\in\children(r)$ and for each $j'\in\children(r)$ with $j'\neq j$, we consider a copy of $\subtree(j')$ and glue it to the vertex in $J'$ corresponding to $j\in J$. Now $J'$ is a rooted tree such that $\children(r)$ is naturally identified with $\children(r')$ and such that for each $j\in\children(r')$, $\subtree(j)$ is a linear rooted tree with volume $v(J)-1$. We can now recursively apply this construction to each such $\subtree(j)$, for $j\in\children(r')$; since the volumes strictly decrease, this process must terminate. Let $J'$ denote the rooted tree that we obtain in the end. Using pullback in the category of flat formal $R$-schemes of ff type, we can equip $J'$ naturally with formal structure. In particular, for each $j\in J'$ we have an attached flat affine formal $R$-scheme of ff type $\fX_j$ with a uniformly rigid generic fiber $X_j$. If $j\in\leaves(J')$, then $X_j$ is a finite intersection of the $X_i$, with $i$ varying in $I_r$.

For each $i\in I_z$ and each $j\in J'$, we let $\fZ_{ij}$ denote the schematic closure of $Z_i\cap X_j$ in $\fX_j$. Then the uniformly rigid generic fiber of $\fX_j\setminus\fZ_{ij}$ is a finite union of retrocompact semi-affinoid subdomains in $X$ that does not meet $Z_i$ and, hence, is contained in the complement $X_i$ of $Z_i$. It thus suffices to see that the $X_j$, $j\in \leaves(J')$, together with the $(\fX_j\setminus\fZ_{ij})^\srig$ with $i\in I_z$, $j\in J'$ cover $X$. Let $x\in X$ be a point; we must show that
\[
x\in \bigcup_{j\in\leaves(J')}X_j\quad\cup\quad\bigcup_{j\in J',i\in I_z}(\fX_j\setminus\fZ_{ij})^\srig\quad.
\]
Let us consider the affine $R$-model $\fX_{r'}=\fX$ of $X$ that is attached to the root $r'$ of $J'$. Since $\bigcup_{i\in I_z}Z_i\subseteq\bigcup_{j\in \leaves(J')} X_j$, we have in particular an inclusion $\bigcup_{i\in I_z} Z_i\subseteq\bigcup_{j\in \children(r')} X_j$. Since the $\fZ_{ij}$ are $R$-flat, their specialization maps are surjective onto their closed points; since the closed points lie very dense, we thus obtain an inclusion $\bigcup_{i\in I_z}\fZ_{ir'}\subseteq\bigcup_{j\in\children(r')}\fX_j$. Hence, $(\fX_j)_{j\in\children(r')}$ together with $(\fX_{r'}\setminus\fZ_{ir})_{i\in I_z}$ covers $\fX_{r'}$, and it follows that the corresponding uniformly rigid generic fibers cover $X$. If $x\in(\fX_{r'}\setminus\fZ_{ir})^\srig$ for some $i\in I_z$, there is nothing more to show. If this is not the case, then $x\in X_j$ for some $j\in\children(r')$, and we may replace $J'$ by $\subtree(j)$. By induction on the volume of $J'$, the statement follows.
\end{proof}

We can now prove the uniformly rigid analog of \cite{Kisin} Lemma 2.3:

\begin{cor}\label{tubecor}
Let $X=\sSp A$ be a semi-affinoid $K$-space, and let $Z\subseteq X$ be a closed uniformly rigid subspace that is defined by functions $f_1,\ldots,f_n\in A$. Let $U\subseteq X$ be weakly retrocompact uniformly rigid open subspace containing $Z$. There exists some $\varepsilon\in\sqrt{|K^*|}$ such that the \emph{tube}\index{tube} 
\[
X_{\leq\varepsilon}\,\mathrel{\mathop:}=\,X(\varepsilon^{-1}f_1,\ldots,\varepsilon^{-1}f_n)\,\mathrel{\mathop:}=\,\{x\in X\,;\,|f_i(x)|\leq\varepsilon\,\,\forall\, i=1,\ldots,n\}
\]
is contained in $U$.
\end{cor}
\begin{proof}
Let $X_f\subseteq X$ denote the complement of $Z$ in $X$; then Proposition \ref{zarretromixadmprop} shows that $(X_f,U)$ is an admissible covering of $X$. By Proposition \ref{zarsrigprop} and the subsequent remark, $(X_{\geq\delta})_\delta$ is an admissible covering of $X_f$, where the $X_{\geq\delta}$ are defined in terms of the $f_i$ in the obvious way. By transitivity of admissibility for coverings, it follows that $(X_{\geq\delta})_\delta\cup U$ is an admissible covering of $X$. Since $X$ is quasi-compact, this covering admits a finite subcovering, which means that $X=U\cup X_{\geq\delta}$ for $\delta$ small enough; now $X_{\leq\varepsilon}\subseteq U$ for any $\varepsilon<\delta$ in $\sqrt{|K^*|}$.
\end{proof}

The following easy lemma shows that Corollary \ref{tubecor} can be applied in interesting situations:

\begin{lem}\label{tubeauxlem}
Let $Y$ be a semi-affinoid $K$-space, let $\fX$ be a flat formal $R$-scheme of locally tf type with uniformly rigid generic fiber $X$, let $\fU\subseteq\fX$ be an affine open subscheme, and let $\phi\colon Y\rightarrow X$ be a uniformly rigid morphism. Then $\phi^{-1}(\fU^\srig)$ is weakly retrocompact in $Y$.
\end{lem}
\begin{proof}
Let $(\fX_i)_{i\in I}$ be an affine open covering of $\fX$, and let $(X_i)_{i\in I}$ be the induced admissible semi-affinoid covering of $X$. Let $(Y_j)_{j\in J}$ be a retrocompact refinement of $(\phi^{-1}(X_i))_{i\in I}$, and let $\psi\colon J\rightarrow I$ be a suitable refinement map. Let us write $U\mathrel{\mathop:}=\fU^\srig$, and for each $i\in I$, let us write $U_i\mathrel{\mathop:}=U\cap X_i$. Then
\[
\phi^{-1}(U)\,=\,\bigcup_{j\in J}(\phi|_{Y_j})^{-1}(U_{\psi(j)})\quad.
\]
Since $Y_j$ is retrocompact in $Y$ for all $j$, it suffices to show that $(\phi|_{Y_j})^{-1}(U_{\psi(j)})$ is weakly retrocompact in $Y_j$ for all $j$. Hence we may assume that $\fX$ is affine, which implies that $X$ is semi-affinoid. Clearly $U\subseteq X$ is weakly retrocompact. By the retrocompact version of Corollary \ref{semaffpreimcor}, $\phi^{-1}(U)$ is weakly retrocompact in $Y$, as desired.
\end{proof}

		\subsection{Analytic continuation}\label{ancontsec}
			We need some more preparations. Let $\fX=\Spf R[[S_1,\ldots,S_n,Z]]$ denote the $n+1$-dimensional open formal unit disc over $R$, and let $f$ be a $Z$-distinguished global function on $\fX$; that is, a formal power series whose reduction in $k[[Z]]$ is a nonzero non-unit. While the non-vanishing locus $\fX_f$ of $f$ on $\fX$ is empty, the non-vanishing locus $\fX_{\pi,f}$ of $f$ on the \emph{envelope} of $\fX$ is nonempty. Let $r$ be a natural number, and let $\fX'$ denote the affine part of the admissible formal blowup of $\fX$ in the ideal $(\pi^{r+1},S_1,\ldots,S_n)$ where the pullback ideal is generated by $\pi^{r+1}$. The natural morphism of affine flat formal $R$-schemes of ff type $\fX'\rightarrow\fX$ is induced by a unique morphism of envelopes $(\fX')_\pi\rightarrow\fX_\pi$ which induces a morphism $(\fX')_{\pi,f}\rightarrow\fX_{\pi,f}$. The resulting cartesian diagram
\[
\begin{diagram}
(\fX')_{\pi,f}&\rInto&(\fX')_\pi\\
\dTo&&\dTo\\
\fX_{\pi,f}&\rInto&\fX_\pi
\end{diagram}
\]
may be considered as a \emph{formal Hartogs figure}\index{Hartogs figure}; cf.\ \cite{Ltke} Section 3 for the notion of Hartogs figures in rigid geometry. We will prove the following analytic continuation result, in analogy with \cite{Ltke} Satz 6:

\begin{theorem}\label{ancontthm}
Let $g$ and $h$ be global functions on $\fX_{\pi,f}$ and $(\fX')_\pi$ respectively whose pullbacks to $(\fX')_{\pi,f}$ coincide. Then $g$ and $h$ are the pullbacks of a unique global function on $\fX_\pi$.
\end{theorem}

We first recall an elementary and well-known lemma. Let $A$ be a complete noetherian local ring with residue field $\kappa$. A formal power series $f\in A[[Z]]$ is called \emph{distinguished} if its reduction in $\kappa[[Z]]$ is a nonzero non-unit. 

\begin{lem}\label{distfinlem}
If $f\in A[[Z]]$ is a distinguished power series, then the continuous $A$-endomorphism 
\[
\phi\colon A[[Z]]\rightarrow A[[Z]]\quad,\quad Z\mapsto f
\]
is finite and free.
\end{lem}
\begin{proof}
We first show that $\phi$ is finite. Let $\m$ denote the maximal ideal of $A$. Since $\phi$ is $A$-linear, $\phi(\m)A[[Z]]=\m A[[Z]]$, and it follows that $(\m,f)$ is an ideal of definition of $A[[Z]]$. Since $A[[Z]]$ is $(\m,Z)$-adically separated and complete, it follows from \cite{Eis} Ex. 7.2 that $\phi$ is finite if and only if $A/(\m,f)$ is finite over the residue field $\kappa=A/(\m,Z)$ of $A$. By the Weierstrass Division Theorem for formal power series in $\kappa[[Z]]$, $A/(\m,f)$ is finite free over $\kappa$ of rank $n$, as desired, where $n$ is the $Z$-adic valuation of the reduction of $f$ modulo $\m$. We have thus shown that $\phi$ is finite.

To show that $\phi$ is free, it suffices to prove that $\phi$ is flat, the noetherian ring $A[[Z]]$ being local. Let $\tilde{\phi}$ denote the reduction of $\phi$ modulo $\m$. We easily see that $\tilde{\phi}$ is injective. Since $\kappa[[Z]]$ is a principal ideal domain and since $\tilde{\phi}$ is injective, it follows that $\tilde{\phi}$ is flat. By the Flatness Criterion (\cite{Bourbaki} Chap. III \S 5.2 Theorem 1 ($i$)$\Leftrightarrow$($iii$)), it suffices to show that $\m A[[Z]]\otimes_{A[[Z]],\phi}A[[Z]]\rightarrow \phi(\m)A[[Z]]$ is bijective, which follows from the fact that $\phi$ is an $A$-homomorphism and that $A[[Z]]$ is $A$-flat.
\end{proof}

Let us now give the \emph{proof} of Theorem \ref{ancontthm}, returning to the situation that we have introduced at the beginning of this section. Let us write $\fY=\fX$. By Lemma \ref{distfinlem}, the $R$-morphism $\phi\colon\fY\rightarrow\fX$ given by $\phi^*S_i=S_i$ for $1\leq i\leq n$ and $\phi^*Z=f$ is finite and free, and it induces a finite free morphism of envelopes $\fY_\pi\rightarrow\fX_\pi$ that restricts to a finite and free morphism $\fY_{\pi,f}\rightarrow\fX_{\pi,Z}$. Let $\fX''$ denote the admissible formal blowup of $\fX$ in $(\pi^{r+1},S_1,\ldots,S_n)$; since $\phi^*S_i=S_i$, the strict transform of $\fX''$ under the flat morphism $\phi$ is the admissible formal blowup of $\fY$ in the ideal generated by $\pi^{r+1}$ and the $S_i$ on $\fY$. By the explicit description of admissible formal blowups that we recalled in Section \ref{blowupssec}, the affine locus $\fY'$ in $\fY''$ where the pullback ideal is generated by $\pi^{r+1}$ is obtained from the corresponding locus $\fX'\subseteq\fX''$ via $\phi$-pullback. Hence, the cartesian diagram of affine noetherian formal $R$-schemes
\[
\begin{diagram}
(\fY')_{\pi,f}&\rInto&(\fY')_\pi\\
\dTo&&\dTo\\
\fY_{\pi,f}&\rInto&\fY_\pi
\end{diagram}
\]
is obtained from the cartesian diagram
\[
\begin{diagram}
(\fX')_{\pi,Z}&\rInto&(\fX')_\pi\\
\dTo&&\dTo\\
\fX_{\pi,Z}&\rInto&\fX_\pi
\end{diagram}
\]
via $\phi$-pullback. After choosing a basis for the ring of functions on $\fY$ over the ring of functions on $\fX$, the problem is thus reduced to the case where $f=Z$.

The ring of functions on $\fX_{\pi,Z}$ is easily identified with
\[
R[[S_1,\ldots,S_n,Z]]^\pi\langle Z^{-1}\rangle\,=\,\{\sum_{i=-\infty}^\infty a_i Z^i\,;\,a_i\in R[[S_1,\ldots,S_n]]\,,\,|a_i|\overset{i\rightarrow-\infty}{\longrightarrow} 0\}
\]
with uniquely determined coefficients $a_i$ where, as in Section \ref{envsec}, the superscript $\pi$ indicates that a ring is being equipped with the $\pi$-adic topology, and where $|\cdot|$ denotes the Gauss norm. The ring of functions on $(\fX')_{\pi,Z}$ is given by
\[
R[[S_1,\ldots,S_n]]\langle\frac{S_1}{\pi^{r+1}},\ldots,\frac{S_n}{\pi^{r+1}}\rangle[[Z]]^\pi\langle Z^{-1}\rangle\quad,
\]
which admits the analogous explicit description. Let us write the function $g$ on $\fX_{\pi,Z}$ uniquely as $g=\sum_i a_i Z^i$; we must show that $a_i=0$ for $i<0$. Since the natural pullback of $g$ to $(\fX')_{\pi,f}$ coincides with the natural pullback of
\[
h\,\in\,R[[S_1,\ldots,S_n]]\langle\frac{S_1}{\pi^{r+1}},\ldots,\frac{S_n}{\pi^{r+1}}\rangle[[Z]]\quad,
\]
it suffices to note that the natural homomorphism
\[
R[[S_1,\ldots,S_n]]\rightarrow R[[S_1,\ldots,S_n]]\langle \frac{S_1}{\pi^{r+1}},\ldots,\frac{S_n}{\pi^{r+1}}\rangle
\]
is injective, which is clear from the fact that after inverting $\pi$, it is a flat homomorphism of domains. \qed

		\subsection{The structure of formal tubes}\label{tubesstrsec}
			Let $\fX$ be a smooth affine formal $R$-scheme of ff type, let $f_1,\ldots,f_r$ be an $R$-regular sequence in the ring of global functions on $\fX$ such that the closed formal subscheme $\fV$ defined by the $f_i$ is $R$-smooth, and let $n\in\N$ be a natural number. 
\begin{prop}\label{formaltubestructureprop}
Let
\[
\fX\langle\pi^{-(n+1)}f\rangle\,\mathrel{\mathop:}=\,\fX\left\langle\frac{f_1}{\pi^{n+1}},\ldots,\frac{f_r}{\pi^{n+1}}\right\rangle
\]
denote the affine open part of the admissible formal blowup of $\fX$ in the ideal $(\pi^{n+1},f_1,\ldots,f_r)$ where the pullback ideal is generated by $\pi^{n+1}$. Then $\fX\langle\pi^{-(n+1)}f\rangle$ is $R$-smooth, and if $\fV$ is connected, then $\fX\langle\pi^{-(n+1)}f\rangle$ is connected as well.
\end{prop}
\begin{proof}
Let us write $f_0\mathrel{\mathop:}=\pi^{n+1}$, and let $A$ denote the ring of global functions on $\fX$. Since $(f_i)_{1\leq i\leq r}$ is an $R$-regular sequence, the sequence $(f_i)_{0\leq i\leq r}$ is regular. By \cite{Eis} Exercise 17.14, it follows that the blowup
of $\Spec A$ in the ideal $(f_0,\ldots,f_r)$ is covered by the affine open subschemes 
\[
\Spec A[T_0,\ldots,\hat{T}_i,\ldots,T_r]/((f_iT_j-f_j)_{0\leq j\leq r\,,\,j\neq i})\quad;
\]
hence
\[
\fX\langle\pi^{-(n+1)}f\rangle\,=\,\Spf A\langle T_1,\ldots,T_r\rangle/((\pi^{n+1}T_j-f_j)_{1\leq j\leq r})\quad.
\]
Since $\fX\langle\pi^{-(n+1)}f\rangle$ is $R$-flat, the fiber criterion for smoothness \cite{AJR3} 5.4 shows that it thus suffices to prove that 
\[
\fX\langle\pi^{-(n+1)}f\rangle_k\,=\,\Spf ((A/(f_1,\ldots,f_r)\otimes_Rk)\langle T_1,\ldots,T_r\rangle)
\]
is smooth over $k$. Since $\fV$ is $R$-smooth by assumption, its special fiber 
\[
\fV_k=\Spf (A/(f_1,\ldots,f_r)\otimes_Rk)
\]
is $k$-smooth, so the statement is clear. The statement concerning connectedness is now equally obvious.
\end{proof}

		\subsection{Weil Extension Theorems}\label{mainweilextsec}
		In this section, we prove uniformly rigid analogs of extension theorems for rational maps to smooth group objects as they are found in \cite{Weil} \S II n. 15 Prop. 1, \cite{BLR} 4.4 and \cite{BS} 2.6. Our general setup is the following: Let $G$ be a smooth uniformly rigid $K$-group, let $\fG$ be a smooth formal $R$-group scheme of locally ff type together with an open immersion $\fG^\srig\hookrightarrow G$ preserving group structures, let $\fX$ be a smooth formal $R$-scheme of locally ff type, and let $\phi\colon\fX^\srig\rightarrow G$ be a morphism of uniformly rigid $K$-spaces. We specify conditions under which $\phi$ extends, necessarily uniquely, to a morphism $\ul{\phi}\colon\fX\rightarrow\fG$. 

We observe that rigid and uniformly rigid $K$-groups and formal $R$-groups of ff type are automatically separated.

Let us first mention a useful descent lemma which is an analog of \cite{BLR} 2.5/5:

\begin{lem}\label{extensiondesclem}
Let $\fY$ and $\fX$ be flat formal $R$-schemes of locally ff type, where $\fX$ is separated, and let $\phi\colon Y\rightarrow X$ be a morphism of uniformly rigid generic fibers. Let us assume that there exists a flat formal $R$-scheme of locally ff type $\fY'$ together with a faithfully flat quasi-compact morphism $\ul{\psi}\colon\fY'\rightarrow\fY$ and a morphism $\ul{\phi}'\colon\fY'\rightarrow\fX$ such that the induced diagram of uniformly rigid generic fibers
\[
\begin{diagram}
Y'&&\\
\dTo>\psi&\rdTo>{\phi'}&\\
Y&\rTo^\phi&X
\end{diagram}
\]
commutes. Then $\phi$ extends uniquely to a morphism $\ul{\phi}\colon\fY\rightarrow\fX$.
\end{lem}
\begin{proof}
We may work locally on $\fY$ and, hence, assume that $\fY$ is quasi-compact. Let $\ul{\Gamma}_\phi\subseteq\fY\times\fX$ denote the schematic closure of the graph of $\phi$ in $\fY\times\fX$. Then $\ul{\psi}^*\ul{\Gamma}_\phi$ coincides with the graph $\Gamma_{\ul{\phi}'}$ of $\ul{\phi}'$, and hence the restricted projection 
\[
p_\fY|_{\ul{\Gamma}_\phi}\colon\ul{\Gamma}_\phi\rightarrow\fY
\]
becomes an isomorphism after base change with respect to $\ul{\psi}$. By Theorem \ref{graphproperthm}, $p_\fY|_{\ul{\Gamma}_\phi}$ is proper, hence adic and separated, so we deduce from Proposition \ref{isodescprop} that $p_\fY|_{\ul{\Gamma}_\phi}$ is an isomorphism, as desired.
\end{proof}

\subsubsection{Uniformly rigid morphisms to smooth formal groups}

We us first consider the case where $\fG^\srig$ coincides with $G$ and where $\fX$ is of locally \emph{tf} type. Let us recall that whenever we write fibered products without indicating the base, the latter is understood to be $\Spf R$ or $\sSp K$, depending on the context.

\begin{prop}\label{tftypedomainprop}
Let $\fG$ be a smooth formal $R$-group scheme of locally ff type, let $\fX$ be a smooth formal $R$-scheme of locally tf type, and let 
\[
\phi\colon\fX^\srig\rightarrow\fG^\srig
\]
be a morphism of uniformly rigid generic fibers. Then $\phi$ extends uniquely to an $R$-morphism 
\[
\ul{\phi}\colon\fX\rightarrow\fG\quad.
\]
\end{prop}
\begin{proof}
We may assume that $\fX$ is affine and connected; let $X$ denote the semi-affinoid generic fiber of $\fX$, and let $G$ denote the uniformly rigid generic fiber of $G$. By \cite{BL1} 4.4, every admissible covering of $X$ is refined by an affine open covering of a suitable admissible formal blowup of $\fX$. By Proposition \ref{retroprop}, Proposition \ref{ffprop} and \cite{BS} Lemma 2.2 (i), it follows that $\phi$ extends to a morphism 
\[
\ul{\phi}|_\fW\colon\fW\rightarrow\fG
\]
for some nonempty open formal subscheme $\fW\subseteq\fX$. Let 
\[
\psi\colon X\times X\rightarrow G
\]
denote the morphism sending $(x_1,x_2)$ to $\phi(x_1)\cdot \phi(x_2)^{-1}$. The product $\fW\times\fW$ is an open formal subscheme of $\fX\times\fX$, and the morphism 
\[
\ul{\psi}|_{\fW\times\fW}\colon\fW\times\fW\rightarrow\fG
\]
sending $(x_1,x_2)$ to $\ul{\phi}|_\fW(x_1)\cdot\ul{\phi}|_\fW(x_2)^{-1}$ extends $\psi$. We claim that $\psi$ extends to an open formal subscheme of $\fX\times\fX$ containing the diagonal $\Delta_\fX$. Let $\fH\subseteq\fG$ be an affine open neighborhood of the identify section, let $H$ denote its semi-affinoid generic fiber, and let $\fU\subseteq\fW\times\fW$ be the $\ul{\psi}|_{\fW\times\fW}$-preimage of $\fH$. Let $(x,x)$ be any closed point of $\Delta_\fX$; we must show that $\psi$ extends to an open neighborhood of $(x,x)$ in $\fX\times\fX$. If $(x,x)\in\fW\times\fW$, there is nothing to show, so we may assume that $(x,x)\notin\fW\times\fW$. Let $\fV$ be any connected affine open neighborhood of $(x,x)$ in $\fX\times\fX$. Then $\fU\cap\fV$ is $R$-dense in $\fV$. Indeed, since $\fV_k$ is an integral scheme, it suffices to observe that $\fU\cap\fV$ is nonempty, and already $\fU\cap\fV\cap\Delta_\fX$ is nonempty since $(\Delta_\fX)_k\cong\fX_k$ is an integral scheme and since both $\fU\cap\Delta_\fX$ and $\fV\cap\Delta_\fX$ are nonempty open formal subschemes of $\Delta_\fX$, the former containing $\Delta_\fW$ and the latter containing $(x,x)$. 

Let $d$ denote the relative dimension of $\fX$ in $x$. After shrinking $\fV$, we may assume that conditions ($i$)--($iii$) in the proof of \cite{BS} Theorem 2.6 are satisfied on $\fV$. That is, $\fV$ does not meet any component of $(\fX\times\fX)\setminus\fU$ not containing $(x,x)$, and there exists an $R$-regular sequence of global functions $f_1,\ldots,f_{2d-1}$ on $\fV$ such that $\Delta_\fX\cap\fV$ is defined by $f_1,\ldots,f_d$, such that the vanishing locus of the $f_i$ is a formal subscheme $\fC\subseteq\fV$ of relative dimension one and such that $\fC\cap(\fV\setminus\fU)=\{(x,x)\}$. Indeed, this follows easily from the fact that $\fX$ is $R$-smooth and that $\fU\cap\fV\cap\Delta_\fX$ is $R$-dense in $\Delta_\fX$, by choosing a suitable regular parameter system for the stalk of $\fX\times\fX$ in $(x,x)$.

Let $V$ denote the semi-affinoid generic fiber of $\fV$, and let us now pass to associated \emph{rigid} spaces. Since $\psi$ maps $\Delta_X$ to the unit section of $G$, \cite{Kisin} Lemma 2.3 shows that there exists an $\varepsilon>0$ in $\sqrt{|K^*|}$ such that the tube
\[
V^\r(\varepsilon^{-1}f_1,\ldots,\varepsilon^{-1}f_{2d-1})
\]
in the affinoid $K$-space $V^\r$ maps to $H^\r$ via $\psi^\r$. Hence, $V(\varepsilon^{-1}f_1,\ldots,\varepsilon^{-1}f_{2d-1})$ maps to $H$ via $\psi$, and so the same holds for the Hartogs figure
\[
V(\varepsilon^{-1}f_1,\ldots,\varepsilon^{-1}f_{2d-1})\cup(\fV\cap\fU)^\srig\quad.
\]
By \cite{Ltke} Theorem 7, the $\psi$-pullback of any global function on $H$ extends to a global function on $V$. By Proposition \ref{ffprop}, we thus obtain a morphism of semi-affinoid $K$-spaces $\psi'_V\colon V\rightarrow H$ coinciding with $\psi|_V$ on the above Hartogs figure. The coincidence subspace $(\psi'_V,\psi|_V)^{-1}(\Delta_G)$ is a closed semi-affinoid subspace of $V$ which contains the above Hartogs figure and, hence, coincides with $V$. It follows that $\psi'_V$ and $\psi|_V$ coincide, which implies that $\psi|_V$ factorizes through $H$. Since the affine formal $R$-schemes $\fV$ and $\fH$ are smooth, they are normal, cf.\ the remark following Proposition \ref{regprop}, and by Proposition \ref{normallatticeprop} it follows that $\psi|_V$ extends uniquely to a morphism $\ul{\phi}|_\fV\colon\fV\rightarrow\fH$, as desired. Our claim has been shown.

Let $\fV$ now denote any open neighborhood of $\Delta_\fX$ in $\fX\times\fX$ such that the restriction of $\psi$ to the uniformly rigid generic fiber $V$ of $\fV$ extends to a morphism $\ul{\psi}|_\fV$ from $\fV$ to $\fG$. Let us consider the morphism $\fV\cap(\fX\times\fW)\rightarrow\fG$ sending $(x_1,x_2)$ to $\ul{\psi}|_\fV(x_1,x_2)\cdot\ul{\phi}|_\fW(x_2)$. On uniformly rigid generic fibers, it spans a commutative triangle together with $u$ and the generic fiber of the projection $p_1\colon\fX\times\fX\rightarrow\fX$ to the first factor. Moreover, the restriction of $p_1$ to $\fV\cap(\fX\times\fW)$ is faithfully flat. Indeed, flatness is clear, and surjectivity follows from the fact that $\fV$ contains $\Delta_\fX$ and that $\fW\subseteq\fX$ is $R$-dense. It now follows from Lemma \ref{extensiondesclem} that $\phi$ extends uniquely to a morphism $\ul{\phi}\colon\fX\rightarrow\fG$, as desired.
\end{proof}

Without much further effort, we can pass to the more general situation where $\fX$ is of locally \emph{ff} type over $R$:

\begin{theorem}\label{ffteverywheredefinedthm}
Let $\fG$ be a smooth formal $R$-group scheme of locally ff type, let $\fX$ be a smooth formal $R$-scheme of locally ff type, and let 
\[
\phi\colon\fX^\srig\rightarrow\fG^\srig
\]
be a morphism of uniformly rigid generic fibers. Then $\phi$ extends uniquely to an $R$-morphism 
\[
\ul{\phi}\colon\fX\rightarrow\fG\quad.
\]
\end{theorem}
\begin{proof}
Let $X$ and $G$ denote the uniformly rigid generic fibers of $\fX$ and $\fG$ respectively, and let 
\[
\psi\colon X\times X\rightarrow G
\]
be the morphism sending $(x_1,x_2)$ to $\phi(x_1)\cdot \phi(x_2)^{-1}$.  We claim that $\psi$ extends to an open neighborhood of $\Delta_\fX$ in $\fX\times\fX$. 

Let $\fH\subseteq\fG$ be an affine open neighborhood of the identity section, and let $H\subseteq G$ denote the uniformly rigid generic fiber of $\fH$. The formal $R$-scheme $\fX\times\fX$ is $R$-smooth and, hence, normal; cf.\ the remark following Proposition \ref{regprop}. By Proposition \ref{redtoffprop}, it thus suffices to show that for every closed point $(x,x)$ in $\Delta_\fX$, the formal fiber $](x,x)[_{\fX\times\fX}$ maps to $H$ under $\psi$. To establish our claim, we may thus replace $\fX$ by its completion along a closed point $x$ and thereby assume that $\fX$ is local. Moreover, after some finite possibly ramified base extension we may assume that the residue field of $\fX$ naturally coincides with $k$, cf.\ Proposition \ref{smoothlocalstructureprop} and its proof. By Lemma \ref{localstructurerationalpointlem}, $\fX$ is then isomorphic to $\D^d_R$, where $d$ denotes the relative dimension of $\fX$ over $R$. Then $\fX\times\fX=\D^{2d}_R$, and we must show that $\psi$ factorizes through $H$. Since $X$ is the (non-admissible) union of the admissible open subspaces 
\[
\B^d_{\leq\varepsilon,K}\times\B^d_{\leq\varepsilon,K}\quad\textup{for}\quad\varepsilon\rightarrow 1\quad\textup{in}\quad\sqrt{|K^*|}\quad,
\]
we may, after some finite ramified base extension, assume that $\fX=\B^d_R$ and that $\phi$ extends to a morphism $\phi'$ on a strictly \emph{larger} rational closed polydisc $X'$. Let $\fX'=\B^d_R$ be the natural smooth $R$-model of $X'$, and let $\ul{\tau}\colon\fX\rightarrow\fX'$ be the $R$-morphism corresponding to the inclusion of $X$ into $X'$. Since the physical image of $\ul{\tau}\times\ul{\tau}$ is a single point, it suffices to show that $\phi'$ extends to a morphism $\fX'\rightarrow\fG$. This follows from Proposition \ref{tftypedomainprop}, so our claim has been shown.

Let us now return to the situation where $\fX$ is a general smooth formal $R$-scheme of locally ff type. Let $\fV\subseteq\fX\times\fX$ be an open neighborhood of $\Delta_\fX$ with uniformly rigid generic fiber $V$ such that $\psi|_V$ extends to a morphism $\ul{\psi}|_\fV\colon\fV\rightarrow\fG$, let $\ul{\Gamma}_\phi$ denote the schematic closure of the graph of $\phi$ in $\fX\times\fG$, and let $\ul{\phi}'\colon\ul{\Gamma}_\phi\rightarrow\fG$ denote the restriction of the projection $p_\fG\colon\fX\times\fG\rightarrow\fG$. Let $\fV'$ denote the preimage of $\fV$ under the morphism $(\id_\fX\times p_\fX|_{\ul{\Gamma}_\phi})\colon\fX\times\ul{\Gamma}_\phi\rightarrow\fX\times\fX$; we claim that $\fV'$ is faithfully flat over $\fX$ via the first projection $p_1$. Flatness follows from the fact that $\ul{\Gamma}_\phi$ is $R$-flat, from the fact that flatness is preserved under base change of ff type, cf.\ Proposition \ref{flatnessbasechangeprop}, and from the fact that open immersions of locally noetherian formal schemes are flat. To establish surjectivity, we argue as follows: Since images of morphisms of $k$-schemes of finite type are constructible, cf.\ \cite{EGAIV1} 1.8.5, and since the closed points in a $k$-scheme of finite type lie very dense, if suffices to see that $p_1|_{\fV'}$ hits all closed points in $\fX$. Let $x\in\fX$ be a closed point; since $p_\fX|_{\ul{\Gamma}_\phi}$ induces an isomorphism of uniformly rigid generic fibers and since $\sp_\fX$ is surjective, there exists a closed point $x'\in\ul{\Gamma}_\phi$ above $x$; then $(x,x')$ is a point in $\fV'$ projecting to $x$. Surjectivity has thus been shown. Let us now consider the morphism $\fV'\rightarrow\fG$ sending a point $(x,x')$ to $\ul{\psi}|_\fV(x,p_\fX|_{\ul{\Gamma}_\phi}(x'))\cdot\ul{\psi}'(x')$. Together with $p_1$ and $\phi$, we obtain a commutative triangle on uniformly rigid generic fibers. It follows from Lemma \ref{extensiondesclem} that $\phi$ extends to a morphism $\ul{\phi}\colon\fX\rightarrow\fG$, as desired.
\end{proof}

\subsubsection{Formally unramified points}

Let us finally pass to the more challenging and more important situation where $\fG^\srig$ is an admissible open subgroup of a larger uniformly rigid group $G$. Here we have to confine ourselves to the case where $\fG$ is of tf type over $R$, where $G$ is the uniform rigidification of a quasi-paracompact rigid $K$-group and where the open immersion $\fG^\rig\subseteq G^\r$ is retrocompact.

A not necessarily finite discrete analytic extension field $K'$ of $K$ is called \emph{formally unramified} over $K$ if the induced local extension of complete discrete valuation rings is formally unramified. If $X$ is a quasi-separated rigid space and if $U\subseteq X$ is a rigid subspace, we say that $U$ contains all formally unramified points\index{point!formally unramified} of $X$ if for every formally unramified discrete analytic extension field $K'$ over $K$, every $K'$-valued point of $X\hat{\otimes}_KK'$ lies in the admissible open subspace $U\hat{\otimes}_KK'$. 

\begin{theorem}\label{mainweilextthm}
Let $G$ be a quasi-paracompact smooth rigid $K$-group, and let $\fG$ be a quasi-paracompact smooth formal $R$-group scheme of locally tf type together with a retrocompact open immersion $\fG^\rig\hookrightarrow G$ that is compatible with the group structures such that $\fG^\rig$ contains all formally unramified points of $G$. Let $\fX$ be a smooth formal $R$-scheme of locally ff type. Then any uniformly rigid morphism 
\[
\phi\colon\fX^\srig\rightarrow G^\sr
\]
extends uniquely to a morphism 
\[
\ul{\phi}\colon\fX\rightarrow\fG\quad.
\]
\end{theorem}
\begin{proof}
By \cite{FRG} 2.8/3, there exist a quasi-paracompact $R$-model of locally tf type $\fG''$ for $G$ and an admissible formal blowup $\fG'\rightarrow\fG$ such that the inclusion of $\fG^\rig$ into $G$ is induced by a morphism $\fG'\rightarrow\fG''$. The argument in the proof of \cite{BL1} 4.7 shows that $\fG''$ is separated. Since $\fG^\rig\subseteq G$ is retrocompact, \cite{BL2} Corollary 5.4 and \cite{FRG} 2.6/13-14 show that after replacing $\fG''$ and $\fG'$ by suitably chosen admissible formal blowups, we may assume that $\fG'\subseteq\fG''$ is an open formal subscheme. 

By Proposition \ref{redtoffprop}, we may assume that $\fX$ is local. Let $\ul{\Gamma}_\phi$ denote the schematic closure of the graph of $\phi$ in $\fX\times\fG''$. By Corollary \ref{genisocor}, the envelope $(p_\fX|_{\ul{\Gamma}_\phi})_\pi$ is an isomorphism over an open neighborhood $\fW\subseteq\fX_\pi$ of the generic point $\eta$ in $\fX_{\pi,k}$. Let us abbreviate
\[
\ul{\phi}'_\pi\,\mathrel{\mathop:}=\,p_{\fG''}^\pi|_{(\ul{\Gamma}_\phi)_\pi}\quad,
\]
where 
\[
p_{\fG''}^\pi\colon\fX_\pi\times\fG''\rightarrow\fG''
\]
denotes the projection to the second factor.
The completed stalk $R_\eta$ of $\fX_\pi$ in $\eta$ is a complete discrete valuation ring, cf.\ the remark following Proposition \ref{liftrregpointsprop}; by Proposition \ref{regprop}, it is a formally unramified local extension of $R$. Since $\fG^\rig$ contains all formally étale points of $G$, $\eta\in(\ul{\phi}'_\pi)^{-1}(\fG')$. After shrinking $\fW$, we may therefore assume that 
\[
\fW\,\subseteq\,(\ul{\phi}'_\pi)^{-1}(\fG')
\]
such that we obtain a morphism $\ul{\phi}''_\pi|_\fW\colon\fW\rightarrow\fG'\rightarrow\fG$. 

Let $\delta\colon G^\sr\times G^\sr\rightarrow G^\sr$ denote the twisted multiplication morphism sending $(x,y)$ to $x\cdot y^{-1}$, and let 
\[
\psi\colon X\times X\rightarrow G^\sr
\]
denote the morphism $\delta\circ(\phi\times\phi)$ sending $(x_1,x_2)$ to $\phi(x_1)\cdot \phi(x_2)^{-1}$. Let $\fH\subseteq\fG$ be an affine open neighborhood of the identity section, and let $H$ denote its uniformly rigid generic fiber.
We claim that $\psi$ extends to a morphism $\ul{\psi}|_\fV\colon\fV\rightarrow\fH$ for some open neighborhood $\fV$ of $\Delta_\fX$ in $\fX\times\fX$. 
As in the proof of Theorem \ref{ffteverywheredefinedthm}, we reduce the proof of this claim to the case where $\fX=\Spf R[[T_1,\ldots,T_d]]$, using the fact that for any finite local extension of discrete valuation rings $R'/R$, $(\fX\otimes_RR')_{\pi'}=\fX_\pi\otimes_RR'$, with $\pi'$ denoting a uniformizer of $R'$, together with the fact that $\fW\otimes_RR'\subseteq\fX_\pi\otimes_RR'$ is $R'$-dense. Let us write 
\[
\fU\mathrel{\mathop:}=(p_1^\pi)^{-1}(\fW)\cap (p_2^\pi)^{-1}(\fW)\quad,
\]
where for $i=1,2$, $p_i^\pi\colon(\fX\times\fX)_\pi\rightarrow\fX_\pi$, denotes the morphism induced by the $i$-th projection $p_i\colon\fX\times\fX\rightarrow\fX$. We observe that the closed diagonal embedding
\[
\Delta_\fX\colon\fX\hookrightarrow\fX\times\fX
\]
is induced by a closed immersion
\[
\Delta_{\fX_\pi}\colon\fX_\pi\hookrightarrow(\fX\times\fX)_\pi
\]
of envelopes such that $p_i^\pi\circ\Delta_{\fX_\pi}=\id_{\fX_\pi}$ for $i=1,2$. Indeed, this follows immediately by looking at rings of global sections. In particular, $\Delta_{\fX_\pi}^{-1}(\fU)=\fW$ is $R$-dense in $\fX_\pi$. By abuse of notation, we will not distinguish between the closed immersion $\Delta_{\fX_\pi}$ and the closed formal subscheme that it defines.

Let $\ul{\delta}\colon\fG\times\fG\rightarrow\fG$ denote the twisted multiplication morphism sending $(x,y)$ to $x\cdot y^{-1}$; then $\ul{\delta}$ is an $R$-model of the restriction of $\delta$ to $\fG^\srig$. We define a morphism $\ul{\psi}_\pi|_\fU\colon\fU\rightarrow\fG$ by setting 
\[
\ul{\psi}_\pi|_\fU\,\mathrel{\mathop:}=\,\ul{\delta}\circ(\ul{\phi}''_\pi\circ p_1^\pi|_\fU,\ul{\phi}''_\pi\circ p_2^\pi|_\fU)\quad,
\]
where for $i=1,2$, we let $p_i^\pi|_\fU\colon\fU\rightarrow\fW$ denote the restriction of $p_i^\pi$. The restriction of $\ul{\psi}_\pi|_\fU$ to $\fU\cap\Delta_{\fX_\pi}$ factorizes through the unit section of $\fG$; hence in particular
\[
\fU\cap\Delta_{\fX_\pi}\,\subseteq\,\ul{\psi}_\pi|_\fU^{-1}(\fH)\quad.
\]
We may thus shrink $\fU$ such that $\ul{\psi}_\pi|_\fU$ factorizes through $\fH$, without losing the property that $\fU\cap\Delta_{\fX_\pi}$ is $R$-dense in $\Delta_{\fX_\pi}$. The open formal subscheme $\fU\subseteq(\fX\times\fX)_\pi$ is a finite union of basic open formal subschemes; the special fiber of $\Delta_{\fX_\pi}\cong\fX_\pi$ being integral, we may shrink $\fU$ further such that $\fU$ is a basic open subset in $(\fX\times\fX)_\pi$, again without losing the property that $\fU\cap\Delta_{\fX_\pi}$ is $R$-dense in $\Delta_{\fX_\pi}$. Let $f$ be a global function on $(\fX\times\fX)_\pi$ such that $\fU$ is the  locus where $f$ does not vanish. We may assume that $f$ is nonzero, for if $\fU$ was empty, then $\Delta_{\fX_\pi}$ and, hence, $X$ would be empty, and there would be nothing to show. Again, since $\fU\cap\Delta_{\fX_\pi}$ is $R$-dense in $\Delta_{\fX_\pi}$, we may choose an $R$-isomorphism 
\[
\fX\times\fX\cong\Spf R[[S_1,\ldots,S_{2d-1},Z]]
\]
such that the diagonal $\Delta_\fX$ is defined by $S_1,\ldots,S_d$ and such that $f$ is either a unit or $Z$-distinguished. We will briefly write $\ul{S}$ to denote the system $S_1,\ldots,S_{2d-1}$. 

Since $\psi$ maps $\Delta_X$ to the identity of $G^\sr$, Corollary \ref{tubecor} and Lemma \ref{tubeauxlem} show that $\psi^{-1}(H)$ contains some tube around the diagonal and, hence, some tube around the vanishing locus of the $S_i$. That is, there exists an $r\in\N$ such that 
\[
(X\times X)(\pi^{-r}\ul{S})\subseteq \psi^{-1}(H)\quad.
\]
Let $\ul{\tau}\colon(\fX\times\fX)\langle\pi^{-r}\ul{S}\rangle\rightarrow\fX\times\fX$ denote the natural formal dilatation inducing this tube, and let
\[
\ul{\tau}_\pi\colon(\fX\times\fX)\langle\pi^{-r}\ul{S}\rangle_\pi\rightarrow(\fX\times\fX)_\pi
\]
denote its envelope. By Proposition \ref{formaltubestructureprop}, $(\fX\times\fX)\langle\pi^{-r}\ul{S}\rangle$ is $R$-smooth and, hence, normal. It thus follows from Proposition \ref{normallatticeprop} and Proposition \ref{ffprop} that the restriction of $\psi$ to the above tube is induced by a morphism of affine formal $R$-schemes $\ul{\psi}'\colon(\fX\times\fX)\langle\pi^{-r}\ul{S}\rangle\rightarrow\fH$; let
\[
\ul{\psi}'_\pi\colon(\fX\times\fX)\langle\pi^{-r}\ul{S}\rangle_\pi\rightarrow\fH
\]
denote its envelope. Let us set
\[
\fU'\,\mathrel{\mathop:}=\,\ul{\tau}_\pi^{-1}(\fU)\quad;
\]
we claim that
\[
\ul{\psi}'_\pi|_{\fU'}\,=\,\ul{\psi}_\pi|_\fU\circ\ul{\tau}_\pi|_{\fU'}\quad.
\]

Let us write 
\[
H'\mathrel{\mathop:}=\delta^{-1}(H)\subseteq G^\sr\times G^\sr\quad.
\]
By \cite{FRG} 2.8/3, there exists an admissible formal blowup $\fG'''\rightarrow\fG''\times\fG''$ such that $\delta$ extends to a morphism $\fG'''\rightarrow\fG''$. Let $\fH'\subseteq\fG'''$ denote the open preimage of $\fH$ under the induced morphism $\fG'''\rightarrow\fG$; then $\fH'$ is a model of $H'$. Let $\ul{\Gamma}_{\phi\times\phi}$ denote the schematic closure of the graph of $\phi\times\phi$ in $\fX\times\fX\times\fG'''$; then $\ul{\Gamma}_{\phi\times\phi}$ is the strict transform of $\ul{\Gamma}_\phi\times\ul{\Gamma}_\phi$, and the strict transform $\ul{\Gamma}_{\phi\times\phi}\langle\pi^{-r}\ul{S}\rangle$ of $\ul{\Gamma}_{\phi\times\phi}$ under $\ul{\tau}$ is the schematic closure of the graph of $\phi\times\phi$ restricted to $(X\times X)(\pi^{-r}\ul{S})$ in $(\fX\times\fX)\langle\pi^{-r}\ul{S}\rangle\times\fG'''$. Since the specialization map
\[
\sp_{\ul{\Gamma}_{\phi\times\phi}\langle\pi^{-r}\ul{S}\rangle}
\]
is surjective onto closed points, $\ul{\Gamma}_{\phi\times\phi}\langle\pi^{-r}\ul{S}\rangle$ lies in $(\fX\times\fX)\langle\pi^{-r}\ul{S}\rangle\times\fH'$. Since $\fX\times\fX$ and $(\fX\times\fX)\langle\pi^{-r}\ul{S}\rangle$ are affine, the envelopes $(\ul{\Gamma}_{\phi\times\phi})_\pi$ and $(\ul{\Gamma}_{\phi\times\phi}\langle\pi^{-r}\ul{S}\rangle)_\pi$ exist. 

For each $i=1,2$, there is a natural commutative diagram
\[
\begin{diagram}
(\ul{\Gamma}_{\phi\times\phi}\langle\pi^{-r}\ul{S}\rangle)_\pi&\rTo&(\ul{\Gamma}_{\phi\times\phi})_\pi&\rTo&(\ul{\Gamma}_\phi\times\ul{\Gamma}_\phi)_\pi&\rTo^{q_i^\pi}&(\ul{\Gamma}_\phi)_\pi\\
\dInto&&\dInto&&\dInto&&\dInto\\
(\fX\times\fX)\langle\pi^{-r}\ul{S}\rangle_\pi\times\fG'''&\rTo&(\fX\times\fX)_\pi\times\fG'''&\rTo&(\fX\times\fX)_\pi\times\fG''\times\fG''&\rTo&\fX_\pi\times\fG''\\
\end{diagram}
\]
where the lower right horizontal morphism is induced by $p_i^\pi$.
Indeed, the morphism $(\ul{\Gamma}_{\phi\times\phi})_\pi\rightarrow(\ul{\Gamma}_\phi\times\ul{\Gamma}_\phi)_\pi$ is the strict transform of the admissible formal blowup $\fG'''\rightarrow\fG''\times\fG''$, and since completion morphisms of noetherian formal schemes are flat, we see that the left square is cartesian in the category of flat formal $R$-schemes. Finally, for $i=2$ let us consider the morphism $X\times X\rightarrow G$ sending $(x,y)$ to $\phi(y)$. The closure of its graph in $\fX\times\fX\times\fG''$ is $\fX\times\fX\times_{\fX}\ul{\Gamma}_\phi=\fX\times\ul{\Gamma}_\phi$, and its envelope over $(\fX\times\fX)_\pi$ is 
\[
(\fX\times\fX)_\pi\times_{\fX_\pi}(\ul{\Gamma}_\phi)_\pi\subseteq(\fX\times\fX)_\pi\times\fG''\quad.
\]
By \cite{EGAIII1} 5.4.1 the natural $\fX\times\fX$-morphism $\ul{\Gamma}_\phi\times\ul{\Gamma}_\phi\rightarrow\fX\times\ul{\Gamma}_\phi$ extends uniquely to a $(\fX\times\fX)_\pi$-morphism $(\ul{\Gamma}_\phi\times\ul{\Gamma}_\phi)_\pi\rightarrow(\fX\times\fX)_\pi\times_{\fX_\pi}(\ul{\Gamma}_\phi)_\pi$, and we obtain $q_2^\pi$ via projection to $(\ul{\Gamma}_\phi)_\pi$. The case $i=1$ is dealt with in the same way.

Let us note that the morphism $(\ul{\Gamma}_{\phi\times\phi})\langle\pi^{-r}\ul{S}\rangle_\pi\rightarrow(\fX\times\fX)\langle\pi^{-r}\ul{S}\rangle_\pi$ induces an isomorphism of rings of global sections. Indeed, by the Comparison Theorem \cite{EGAIII1}  4.1.5 this may be verified after completion. Since the ring of global functions on $(\ul{\Gamma}_\phi\times\ul{\Gamma}_\phi)\langle\pi^{-r}\ul{S}\rangle$ is naturally identified with a subring of the ring of power-bounded semi-affinoid functions on $(X\times X)(\pi^{-r}\ul{S})$, the statement follows from Proposition \ref{normallatticeprop} together with the fact that $(\fX\times\fX)\langle\pi^{-r}\ul{S}\rangle$ is normal, cf.\ Proposition \ref{formaltubestructureprop}.

Let us write $\fH''\mathrel{\mathop:}=\ul{\delta}^{-1}(\fH)$, let 
\[
\fH'''\subseteq\fG'\times\fG'\subseteq\fG''\times\fG''
\]
denote the strict transform of $\fH''$, and let $\tilde{\fH}'''\subseteq\fG'''$ denote the strict transform of $\fH'''$. Moreover, let $\fH''''$ denote the preimage of $\fH'''$ in $\fW\times\fW\subseteq\fX_\pi\times\fX_\pi$; then $\fU$ is contained in the preimage of $\fH''''$ under the morphism $(\fX\times\fX)_\pi\rightarrow\fX_\pi\times\fX_\pi$. Let us now consider the resulting commutative diagram
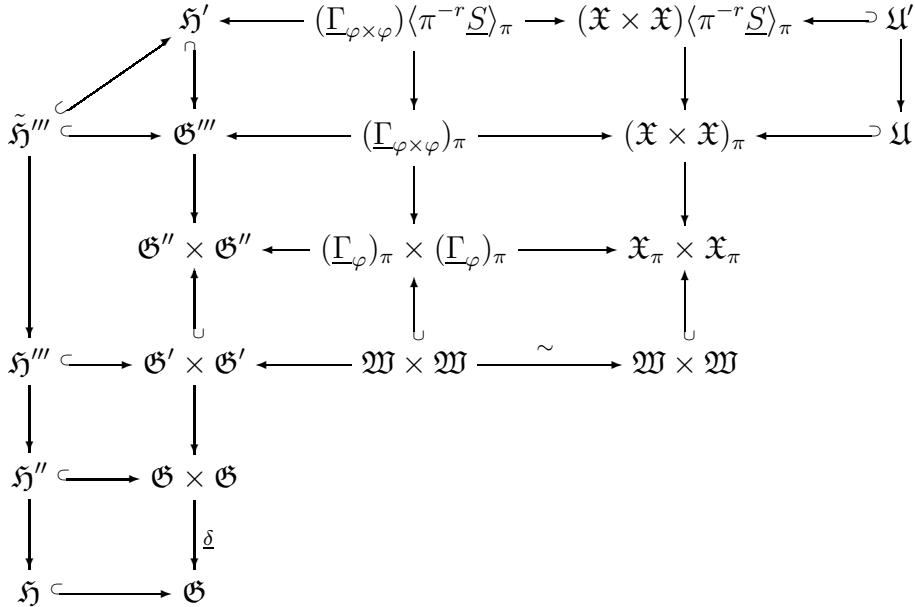
\[
\begin{diagram}
&&\fH'&\lTo&(\ul{\Gamma}_{\phi\times\phi})\langle\pi^{-r}\ul{S}\rangle_\pi&\rTo&(\fX\times\fX)\langle\pi^{-r}\ul{S}\rangle_\pi&\lInto&\fU'&\\
&\ruInto&\dInto&&\dTo&&\dTo&&\dTo&\\
\tilde{\fH}'''&\rInto&\fG'''&\lTo&(\ul{\Gamma}_{\phi\times\phi})_\pi&\rTo&(\fX\times\fX)_\pi&\lInto&\fU&\\
&&\dTo&&\dTo&&\dTo&&&\\
\dTo&&\fG''\times\fG''&\lTo&(\ul{\Gamma}_\phi)_\pi\times(\ul{\Gamma}_\phi)_\pi&\rTo&\fX_\pi\times\fX_\pi&&&\\
&&\uInto&&\uInto&&\uInto&&&\\
\fH'''&\rInto&\fG'\times\fG'&\lTo&\fW\times\fW&\rTo^{\sim}&\fW\times\fW&&&\\
\dTo&&\dTo&&&&&&&\\
\fH''&\rInto&\fG\times\fG&&&&&&&\\
\dTo&&\dTo>{\ul{\delta}}&&&&&&&\\
\fH&\rInto&\fG&&&&&&&\quad.
\end{diagram}
\]

Let $g$ be a global function on $\fH$. Then $g$ induces a function on $\fH''''$ and, hence, a function on $\fU$, which pulls back to a function $g_1$ on $\fU'$. On the other hand, $g$ pulls back to a global function on $(\ul{\Gamma}_{\phi\times\phi})\langle\pi^{-r}\ul{S}\rangle_\pi$ and thus induces a global function on $(\fX\times\fX)\langle\pi^{-r}\ul{S}\rangle_\pi$ which restricts to a function $g_2$ on $\fU'$. We claim that $g_1$ and $g_2$ coincide. To verify this, we may work locally on the locally integral noetherian formal scheme $\fU'$. By Proposition \ref{formaltubestructureprop}, the generic point in the special fiber of $(\fX\times\fX)\langle\pi^{-r}\ul{S}\rangle_\pi$ is $R$-regular in the ring of global functions on $(\fX\times\fX)\langle\pi^{-r}\ul{S}\rangle_\pi$; hence Corollary \ref{genisocor} shows that $(\ul{\Gamma}_{\phi\times\phi})\langle\pi^{-r}\ul{S}\rangle_\pi\rightarrow(\fX\times\fX)\langle\pi^{-r}\ul{S}\rangle_\pi$ is an isomorphism above $\fU'$. The statement now follows from a straightforward diagram chase, using the fact that the function induced by $g$ on the strict transform of $\fH'''$ in $\fG'''$ coincides with the restriction of the function induced by $g$ on $\fH'$, which is due to the fact that the uniformly rigid generic fiber of $\ul{\delta}$ is the restriction of $\delta$ to $\fG^\srig$.

Applying the Continuation Theorem \ref{ancontthm} in the case where the global function $f$ on $(\fX\times\fX)_\pi$ defining $\fU$ is not a unit, it now follows that the $\ul{\psi}'_\pi$-pullback of $g$ extends to all of $\fX\times\fX$. Varying $g$, we thus obtain a morphism 
\[
\ul{\psi}\colon\fX\times\fX\rightarrow\fH
\]
whose uniformly rigid generic fiber coincides with $\psi$ on $(X\times X)(\pi^{-r}\ul{S})$; that is, the closed uniformly rigid coincidence subspace of $\psi$ and $\ul{\psi}^\srig$ in $X\times X$ contains $(X\times X)(\pi^{-r}\ul{S})$. By Proposition \ref{closedimprop}, closed uniformly rigid subspaces of a semi-affinoid $K$-space correspond to ideals in the ring of global functions on that space. Since the restriction homomorphism from global functions on $X\times X$ to global functions on $(X\times X)(\pi^{-r}\ul{S})$ is injective, it follows that the coincidence subspace of $\psi$ and $\ul{\psi}^\srig$ coincides with $X\times X$, and we conclude that $\psi=\ul{\psi}^\srig$, as desired.



Let us return to the case where $\fX$ is a general smooth local formal $R$-scheme of ff type. By what we have seen so far, there exists an open formal neighborhood $\fV$ of $\Delta_\fX$ in $\fX\times\fX$ such that the restriction $\psi|_V$ of $\psi$ to the uniformly rigid generic fiber $V$ of $\fV$ extends to a morphism $\ul{\psi}|_\fV\colon\fV\rightarrow\fH$. The underlying topological space of $\fX\times\fX$ is a finite collection of closed points; hence $\fX\times\fX$ is a disjoint union of \emph{local} smooth formal $R$-schemes of ff type. Since the diagonal $\Delta_\fX\cong\fX$ is connected, it is contained in precisely one connected component of $\fX\times\fX$ which we call the diagonal component. After shrinking $\fV$, we may assume that $\fV$ coincides with the diagonal component. Clearly $(\fX\times\fX)_\pi$ is the disjoint union of the envelopes of the connected components of $\fX\times\fX$; hence $\fV_\pi$ is the diagonal component of $(\fX\times\fX)_\pi$, by which we mean the connected component containing the connected closed formal subscheme $\Delta_{\fX_\pi}$. 


Following the notation introduced in Section \ref{envsec}, we use a subscript $K$ to indicate the formation of generic fibers in the category of Huber's \emph{adic spaces}. The open formal subschemes 
\begin{eqnarray*}
\fV_\pi&\subseteq&(\fX\times\fX)_\pi\\
(p_2^\pi)^{-1}(\fW)&\subseteq&(\fX\times\fX)_\pi
\end{eqnarray*}
induce open adic subspaces
\begin{eqnarray*}
\fV_{\pi,K}&\subseteq&(\fX\times\fX)_{\pi,K}\\
(p_{2,K}^\pi)^{-1}(\fW_K)&\subseteq&(\fX\times\fX)_{\pi,K}\quad.
\end{eqnarray*}
By \cite{Hubercohpaper} 3.1, the adic space $\fX_K$ associated to the rigid space $\fX^\rig$ is an open adic subspace of $\fX_{\pi,K}$; hence we obtain an open adic subspace
\[
(p_{1,K}^\pi)^{-1}(\fX_K)\subseteq(\fX\times\fX)_{\pi,K}\quad.
\]
On the triple intersection
\[
\fT\,\mathrel{\mathop:}=\,\fV_{\pi,K}\cap(\pi_{2,K}^\pi)^{-1}(\fW_K)\cap(p_{1,K}^\pi)^{-1}(\fX_K)
\]
we can now define two morphisms to $G^\ad$: On the one hand, we may consider the morphism
\[
\alpha_\fT\,\mathrel{\mathop:}=\,\phi^{\r,\ad}\circ p_{1,K}^\pi
\]
which first projects to $\fX_K$ via $p_{1,K}^\pi$ and then applies the morphism of adic spaces associated to $\phi$. On the other hand, we may consider the morphism
\[
\beta_\fT\,\mathrel{\mathop:}=\,\mu^\ad\circ\left((\ul{\psi}|_\fV)_{\pi,K},(\ul{\phi}''_{\pi}|_\fW)_K\circ p_{2,K}^\pi\right),
\]
where $\mu$ denotes the multiplication on $G$. We claim that $\alpha_\fT$ and $\beta_\fT$ coincide. Indeed, it suffices to show that $(\ul{\psi}|_\fV)_{\pi,K}|_\fT\colon\fT\rightarrow\fH_K$ coincides with 
\[
\delta^\ad\circ(\phi^\ad\circ p_{1,K}^\pi,\ul{\phi}'_{\pi,K}\circ p_{2,K}^\pi)\quad.
\]
This follows from the fact that morphisms to affinoid adic spaces are determined by their global sections, cf.\ \cite{Huberhabil} 3.2.9 ($ii$), together with the commutative diagram above.

We now claim that the morphism from the triple intersection above to $\fX_K$ that is induced by $p_{1,K}^\pi$ hits every classical point in $\fX_K$. An affinoid adic space admitting a noetherian ring of definition comes equipped with a surjective support map onto the spectrum of its ring of global sections, cf.\ \cite{Huberhabil} 3.6.11. Moreover, if a point in the affinoid adic space associated to a semi-affinoid $K$-space has support is a maximal ideal, then this point must correspond to a classical point, the residue fields of semi-affinoid $K$-algebras in maximal ideals being finite extensions of $K$. It thus suffices to see that the morphism 
\[
\fV_\pi\cap (p_2^\pi)^{-1}(\fW)\rightarrow\fX_\pi
\]
induced by $p_1^\pi$ is faithfully flat. Flatness is clear from the fact that $p_1^\pi$ itself is flat. Since $\fX_\pi$ is affine, it suffices to show that the considered morphism hits every closed point in $\fX_\pi$. Let $x\in \fX_\pi$ be a closed point; then $x\in \fX$, and the residue field $k(x)$ of $\fX_\pi$ in $x$ is $k$-finite. Let us write $T=\Spec k(x)$, and let $T\rightarrow \fX_\pi$ be the morphism corresponding to $x$. Since $T\rightarrow\fX_\pi$ is a closed immersion, we see that $T\times_{\fX_\pi}(\fX\times\fX)_\pi$ is naturally identified with $T\times\fX_\pi$. Since $\fW\subseteq\fX_\pi$ is $R$-dense, $T\times\fW\subseteq T\times\fX_\pi$ is schematically dense. Since $\fV_\pi$ contains the diagonal, $T\times_{\fX_\pi}\fV_\pi\subseteq T\times\fX_\pi$ is nonempty open and, hence, has nonempty intersection with $T\times\fW$, as desired.

Since $\beta_\fT$ factorizes through $\fG^{\rig,\ad}$, we finally conclude that $\phi$ factorizes through $\fG^\srig$. Hence, we are in the situation of Theorem \ref{ffteverywheredefinedthm}, so we see that $\phi$ indeed extends to a morphism $\ul{\phi}\colon\fX\rightarrow\fG$, as desired.
\end{proof}

	\section{Néron models of uniformly rigid spaces}\label{nmsec}
		The formal Néron model of a uniformly rigid $K$-space is defined as follows:

\begin{defi}\label{ffnmdefi}\index{Néron model!of a uniformly rigid space}
Let $Z$ be a smooth uniformly rigid $K$-space. A \emph{formal Néron model} for $Z$ is a smooth formal $R$-scheme $\fZ$ of locally ff type together with a morphism $\fZ^\srig\rightarrow Z$ such that for every smooth formal $R$-scheme of locally ff type $\fX$, every morphism 
\[
\phi\colon\fX^\srig\rightarrow Z
\]
admits a unique extension 
\[
\ul{\phi}\colon\fX\rightarrow\fZ\quad.
\]
\end{defi}

Of course, the formal Néron model of $Z$ is unique up to unique isomorphism if it exists. Moreover, it is clear that in order to verify the \emph{universal Néron property} in Definition \ref{ffnmdefi}, we may and will tacitly assume that the so-called \emph{test object} $\fX$ is affine. 
Let us point out that we do not require the \emph{universal} morphism $\fZ^\srig\rightarrow Z$ to be an open immersion.

		\subsection{General existence results}\label{genexressec}
			
We begin with an immediate observation:

\begin{prop}\label{affinenmprop}
A smooth affine formal $R$-scheme of ff type is the formal Néron model of its uniformly rigid generic fiber.
\end{prop}
\begin{proof}
Morphisms to affine formal schemes correspond to continuous homomorphisms of rings of global sections. By Lemma \ref{jacobsonequivlem}, homomorphisms of $R$-algebras of ff type are automatically continuous, and by the remark following Proposition \ref{regprop}, smooth formal $R$-schemes are normal. Hence the statement follows from Proposition \ref{normallatticeprop}. 
\end{proof}

As an immediate consequence of Theorem \ref{ffteverywheredefinedthm}, we obtain the following result:

\begin{prop}\label{genfibneronprop}
Let $\fG$ be a smooth formal $R$-group scheme of locally ff type; then $\fG$ is the formal Néron model of its uniformly rigid generic fiber.
\end{prop}

In \cite{BS} 1.1, Bosch and Schlöter have defined \emph{formal Néron models}\index{Néron model!of a rigid space} of rigid $K$-spaces. Their definition is analogous to our definition \ref{ffnmdefi} above, except that rigid spaces and smooth formal $R$-schemes of locally \emph{tf} type are considered and except for the fact that the universal morphism is required to be an open immersion. They also define the more general notion of \emph{formal Néron quasi-models} of rigid $K$-spaces, cf.\ \cite{BS} 6.4, where the universal morphism is merely required to be a quasi-open immersion, that is, an injective local isomorphism that is an open immersion on every \emph{quasi-compact} admissible open subspace of its domain. In our setting, the distinction between formal Néron models and formal Néron quasi-models disappears because we impose no conditions on the universal morphism whatsoever.

The following observation is immediate from the definitions:

\begin{prop}\label{neronmodinverseprop}
Let $G$ be a quasi-paracompact rigid $K$-group such that the formal Néron model $\fG$ of $G^\ur$ exists and such that $\fG$ is of locally tf type. Then $\fG$ is the formal Néron (quasi-)model of $G$ in the sense of \textup{\cite{BS}}, provided that the natural morphism $\fG^\rig\rightarrow G$ obtained from the universal morphism $\fG^\urig\rightarrow G^\ur$ via $\r$ is a (quasi-)open immersion.
\end{prop}
\begin{proof}
Let $\fZ$ be an affine formal $R$-scheme of tf type. It suffices to establish a natural bijection $\Hom(\fZ^\rig,G)\cong\Hom(\fZ^\urig, G^\ur)$. 
Since $\fZ$ is of tf type over $R$, we have an identification $(\fZ^\rig)^\ur\cong\fZ^\urig$; hence the desired natural bijection is provided by Corollary \ref{ffcor}.
\end{proof}

In \cite{BS} 1.4, the following important criterion is stated: Let $G$ be a smooth rigid $K$-group, and let $\fG$ be a smooth formal $R$-group scheme of locally tf type together with a retrocompact open immersion $\fG^\rig\hookrightarrow G$ respecting the group structures. Then $\fG$ is a formal Néron model of $G$ if and only if $\fG^\rig$ contains all unramified points of $G$, that is, all points with values in finite unramified field extensions of $K$. Using Theorem \ref{mainweilextthm}, we immediately obtain the uniformly rigid analogon of this criterion:

\begin{theorem}\label{maincompthm}
Let $G$ be a quasi-paracompact smooth rigid $K$-group, and let $\fG$ be a quasi-pa\-ra\-com\-pact smooth formal $R$-group scheme of locally tf type together with a retrocompact open immersion $\fG^\rig\hookrightarrow G$ respecting the group structures. If $\fG^\rig$ contains all \emph{formally} unramified points of $G$, then $\fG$ is a formal Néron model of $G^\sr$.
\end{theorem}

\begin{cor}\label{neronex1cor}
Let $G$ be a smooth and quasi-paracompact rigid $K$-group whose group $G(K^\sh)$ of unramified points is bounded. Then the formal Néron model $\fG$ of $G$ in the sense of \textup{\cite{BS}} exists. It is of tf type, and it is the formal Néron model of the uniformly rigid $K$-group $G^\sr$ associated to $G$.
\end{cor}
\begin{proof}
By \cite{BS} 1.2, $G$ admits a quasi-compact formal Néron $\fG$ which is a formal $R$-group scheme of tf type. Since $G$ is quasi-separated and since $\fG^\rig$ is quasi-compact, the open immersion $\fG^\rig\hookrightarrow G$ is automatically retrocompact. By \cite{Wegel} Theorem 4, the formation of $\fG$ commutes with base change with respect to formally unramified local extensions of $R$. Hence, $\fG^\rig$ contains all formally unramified points of $G$, and so the statement follows from Theorem \ref{maincompthm}. 
\end{proof}

Let us recall from \cite{BLR} 10.1.1 that an \emph{Néron lft-model}\index{Néron model!algebraic} of a smooth $K$-scheme $\mathscr{X}$ is a smooth $R$-model $\ul{\mathscr{X}}$ of $\mathscr{X}$, not necessarily quasi-compact, satisfying the universal property that for every smooth $R$-scheme $\ul{\mathscr{Z}}$ with generic fiber $\mathscr{Z}$, every morphism $\mathscr{Z}\rightarrow \mathscr{X}$ extends uniquely to a morphism $\ul{\mathscr{Z}}\rightarrow\ul{\mathscr{X}}$.

\begin{cor}\label{neronex2cor}
Let $\mathscr{G}$ be a quasi-compact smooth $K$-group scheme that has a Néron lft-model $\ul{\mathscr{G}}$, and let $G$ denote the rigid analytification of $\mathscr{G}$. If $\ul{\mathscr{G}}$ is quasi-compact or if $\mathscr{G}$ is commutative,  then the completion $\fG$ of $\ul{\mathscr{G}}$ along its special fiber is the formal Néron model of $G^\sr$.
\end{cor}
\begin{proof}
Analytifications are quasi-paracompact, so $G^\sr$ is well-defined. By \cite{BS} 6.2, the natural morphism $\fG^\rig\rightarrow G$ is a retrocompact open immersion, and $\fG$ is the formal Néron model of $G$ in the sense of \cite{BS}. By \cite{BLR} 10.1/3, the formation of $\ul{\mathscr{G}}$ commutes with formally unramified base change; it follows that the same holds for $\fG$. We thus conclude that $\fG^\rig$ contains all formally unramified points of $G$. The statement is now a consequence of Theorem \ref{maincompthm}.
\end{proof}

In particular, we see that formal Néron models of uniformly rigid $K$-groups exist in abundance.

		\subsection{Construction techniques}\label{constrtechressec}
			We can now construct interesting new examples of uniformly rigid $K$-spaces admitting formal Néron models, implementing ideas of Chai in the uniformly rigid setting, cf.\ \cite{ChaiPre} 5.4. We thereby obtain formal Néron models that are not of locally tf type and whose generic fibers do not coincide with the considered uniformly rigid $K$-space.

\subsubsection{Completion and descent}

Let $K'/K$ be a finite Galois extension with Galois group $\Gamma$, let $R'/R$ denote the associated extension of discrete valuation rings, and let $k'/k$ be the induced residue field extension. Let $Z$ be a uniformly rigid $K$-space that is associated to a quasi-paracompact and quasi-separated \emph{rigid} $K$-space, and let $Z'\mathrel{\mathop:}=Z\otimes_KK'$ be the uniformly rigid $K'$-space induced via base extension. Let us assume that the formal Néron models $\fZ$ and $\fZ'$ of $Z$ and $Z'$ exist, that they are separated and of locally \emph{tf} type and that their universal morphisms are injective morphisms inducing isomorphisms of completed stalks. Let
\[
\psi\colon\fZ\otimes_RR'\rightarrow\fZ'
\]
denote the \emph{base change morphism}\index{base change!morphism} induced by the universal property of $\fZ'$. Its generic fiber is an injective morphism inducing isomorphisms of completed stalks. In general, $\psi^\rig$ will \emph{not} be surjective: The 'volume' of the formal Néron model may grow under base change due to the appearance of new unramified points.

Let $V_{k'}$ be a $\Gamma$-stable affine closed subscheme of $\fZ'_{k'}$ that is contained in an \emph{affine} open formal subscheme of $\fZ'$. Since $\fZ'$ is separated, $V_{k'}$ then has a $\Gamma$-\emph{stable} affine open neighborhood $\fU'$ in $\fZ'$. Let us write
\[
Y'\,\mathrel{\mathop:}=\,(\fZ'|_{V_{k'}})^\srig\quad,
\]
and let us set $U'\mathrel{\mathop:}=(\fU')^\srig$; then $Y'\subseteq U'$ is a $\Gamma$-equivariant open immersion of semi-affinoid $K'$-spaces. By Proposition \ref{semaffdescprop} and Proposition \ref{descffprop}, it descends uniquely to a morphism of semi-affinoid $K$-spaces $Y\rightarrow U$ which is injective and which induces isomorphisms of completed stalks; these latter properties follow from Proposition \ref{descffprop} by looking at points with values in finite Artinian $k$-algebras. By finite Galois descent for \emph{rigid} spaces and by Corollary \ref{ffcor}, the morphism $U'\rightarrow Z'$ descends uniquely to a morphism $U\rightarrow Z$. We thus see that $Y'\rightarrow Z'$ descends uniquely to a morphism $Y\rightarrow Z$ which is injective and which induces isomorphisms of completed stalks. 

Let $W_{k'}$ denote the $\psi_{k'}$-preimage of $V_{k'}$. There exists a unique reduced closed subscheme $W_k\subseteq\fZ_k$ such that the underlying reduced subschemes of $W_{k'}$ and $W_k\otimes_kk'$ coincide. Indeed, since $k'/k$ is normal and since base change with respect to finite purely inseparable field extensions induces homeomorphisms, this follows from finite Galois descent for schemes in the special fibers. Let
\[
\fY\,\mathrel{\mathop:}=\,\fZ|_{W_k}
\]
denote the completion of $\fZ$ along $W_k$; then $\fY\otimes_RR'$ is the completion of $\fZ\otimes_RR'$ along $W_{k'}$. By Proposition \ref{descffprop}, the morphism $\fY^\srig\otimes_KK'\rightarrow Y'$ descends uniquely to a morphism $\fY^\srig\rightarrow Y$ which is injective and which induces isomorphisms of completed stalks.

\begin{prop}\label{constrnmexprop}
The morphism
\[
\fY^\srig\rightarrow Y
\]
identifies $\fY$ with the formal Néron model of $Y$.
\end{prop}
\begin{proof}
Let $\fX$ be a smooth formal $R$-scheme of locally ff type, let $X$ denote its uniformly rigid generic fiber, and let $\phi\colon X\rightarrow Y$ be a morphism; we must show that $\phi$ extends uniquely to a morphism $\ul{\phi}\colon\fX\rightarrow\fY$. Since $\fY^\srig\rightarrow Y$ is a monomorphism, uniqueness of $\ul{\phi}$ follows from the fact that $\fX$ is $R$-flat. Since $\fZ$ is the formal Néron model of $Z$, the composition $\tilde{\phi}$ of $\phi$ with the monomorphism $Y\rightarrow Z$ extends uniquely to a morphism 
\[
\ul{\tilde{\phi}}\colon\fX\rightarrow\fZ\quad;
\]
in particular the image of $\tilde{\phi}\otimes_KK'$ lies in the image of $(\fZ\otimes_RR')^\srig\rightarrow Z'$. Since $\tilde{\phi}\otimes_KK'$ factorizes through $Y'$, all points in its image specialize to $V_{k'}$ in $\fZ'$. It follows that all points in the image of $\tilde{\phi}$ specialize to $W_k\subseteq\fZ$. Since $\sp_\fX$ is surjective onto the closed points of $\fX$, we conclude that $\ul{\tilde{\phi}}$ factorizes uniquely as a morphism $\ul{\phi}\colon\fX\rightarrow\fY$ composed with the completion morphism $\fY\rightarrow\fZ$. The generic fiber of $\ul{\phi}$ must coincide with $\phi$ since $Y\rightarrow Z$ is a monomorphism.
\end{proof}

By Corollary \ref{neronex2cor}, this construction applies in particular to the case where $Z$ is the analytification of a quasi-compact commutative $K$-group scheme admitting Néron lft-models both over $K$ and over $K'$.

\subsubsection{Product decomposition}

\begin{lem}\label{proddecomplem}
Let $G_1$ and $G_2$ be uniformly rigid $K$-groups such that $G_1\times G_2$ has a formal Néron model. Then $G_1$ and $G_2$ have formal Néron models as well.
\end{lem}
\begin{proof}
Let $\fG$ denote the formal Néron model of $G\mathrel{\mathop:}=G_1\times G_2$. The product group structure on $G$ induces group endomorphisms $h_1$ and $h_2$ of $G$ such that $G_i$ is identified with the kernel of $h_i$ via the natural morphism $G_i\hookrightarrow G_1\times G_2$, such that $h_i^2=h_i$, $h_1\circ h_2=h_2\circ h_1=e_G$ and $h_1\cdot h_2=h_2\cdot h_1=\id_G$, where $e_G$ denotes the endomorphism of $G$ that is induced by the unit section. By the universal property of $\fG$, we obtain group endomorphisms $\ul{h}_i$ of $\fG$ satisfying the same identities. Let us set $\fG_i\mathrel{\mathop:}=\ker\ul{h}_i$; then the universal morphism $\fG^\srig\rightarrow G$ induces morphisms $\fG_i^\srig\rightarrow G_i$. Clearly the natural morphism $\fG_1\times\fG_2\rightarrow\fG$ is an isomorphism, and $\fG_i$ is the formal Néron model of $G_i$.
\end{proof}




		\subsection{Strongly concordant groups}\label{concordantsec}
			Let us recall that the \emph{formal multiplicative group} is the formal open unit disc 
\[
\hat{\G}_{m,R}\,\mathrel{\mathop:}=\,\Spf R[[T]]
\]
with its multiplicative group structure, that is, with the group structure whose comultiplication $R[[T]]\rightarrow R[[X,Y]]$ is described by 
\[
T\mapsto (1+X)(1+Y)-1\quad.
\]
Let $d$ be a natural number. A formal $R$-group scheme $\fG$ of ff type is called a \emph{split formal torus}\index{split formal torus} of rank $d$ over $R$ if it is isomorphic to the $d$-fold fibered product of $\hat{\G}_{m,R}$ over $R$ with itself. 

Let us assume that $K$ has positive residue characteristic $\textup{char}\,k=p>0$. If $\fG$ is a split formal $R$-torus of rank $d$, its \emph{character group} \index{character group!of a split formal torus}
\[
X^*(\fG)\,\mathrel{\mathop:}=\,\Hom_{R-\textup{groups}}(\fG,\hat{\G}_{m,R})
\]
is a free $\Z_p$ module of rank $d$. Indeed, it suffices to remark that the endomorphism ring of $\hat{\G}_{m,R}$ is isomorphic to $\Z_p$, where an element $r\in\Z_p$ is identified with the endomorphism sending the coordinate $T$ to $(1+T)^r-1$, cf.\ \cite{Lubin} 2.3.0. If $\fG$ and $\fH$ are split formal $R$-tori, there is an obvious natural identification
\[
\Hom_{R-\textup{groups}}(\fG,\fH)\,\cong\,\Hom_{\Z_p-\textup{modules}}(X^*(\fH),X^*(\fG))\quad.
\]

In \cite{ChaiPre} 4.3 ($ii$), Chai has defined strongly concordant rigid $K$-groups. We transfer his definition and his terminology to the uniformly rigid setting:

\begin{defi}\index{group!strongly concordant}
A semi-affinoid $K$-group $G$ is called \emph{strongly concordant} if there exist a natural number $d$ and a finite Galois extension $K'/K$ together with a $K'$-group isomorphism
\[
G\otimes_KK'\,\cong\,(\hat{\G}_{m,R'}^d)^\srig\quad,
\]
where $R'$ denotes the valuation ring of $K'$. We say that $G$ \emph{splits} over $K'$.
\end{defi}

In other words, strongly concordant $K$-groups are finite Galois twists of open uniformly rigid polydiscs with their multiplicative group structures. If $G$ is a strongly concordant $K$-group, then $G$ is smooth and connected, and the integer $d$ in the above definition coincides with the dimension of $G^\r$.

If $G$ is a \emph{split} strongly concordant $K$-group, Proposition \ref{affinenmprop} allows us to identify its underlying split formal torus with its formal Néron model. In particular, split strongly concordant $K$-groups admit formal Néron models, and their universal morphisms are isomorphisms. 

If $G$ is a not necessarily split strongly concordant $K$-group and if $K'/K$ is a finite Galois extension splitting $G$, then
\[
X^*(G)\,\mathrel{\mathop:}=\,\Hom_{K'-\textup{groups}}(G\otimes_KK',(\hat{\G}_{m,R})^\srig)
\]
together with its natural $\Gal(K^\sep/K)$-action is called the \emph{character group}\index{character group!of a concordant group} of $G$. If $\fG'$ is the formal Néron model of $G\otimes_KK'$, the universal Néron property of $\hat{\G}_{m,R}$ yields an identification of Galois-modules
\[
X^*(G)\,\cong\,X^*(\fG')\quad.
\]
In particular, $X^*(G)$ is a finite free $\Z_p$-module, and the Galois module $X^*(G)$ does not depend on the choice of $K'$. Let $\sCon_K$ denote the category whose objects are pairs $(M,\varrho)$, where $M$ is a finite free $\Z_p$-module and where $\rho$ is a action of $\Gal(K^\sep/K)$ on $M$ that is trivial on a finite index subgroup; morphisms in $\sCon_K$ are given by $\Z_p$-module homomorphisms commuting with the Galois actions. We have seen that $X^*(\cdot)$ defines a contravariant functor from the category of strongly concordant $K$-groups to $\sCon_K$. We show that it is in fact an anti-equivalence; cf.\ also \cite{ChaiPre} 4.4.2:

\begin{prop}
The functor $X^*(\cdot)$ is an anti-equivalence.
\end{prop}
\begin{proof}
To show that $X^*(\cdot)$ is fully faithful, let us consider two strongly concordant $K$-groups $G$ and $H$. Let $K'/K$ be a finite Galois extension splitting both $G$ and $H$, let $\Gamma$ denote its Galois group, and let $\fG'$ and $\fH'$ denote the formal Néron models of $G\otimes_KK'$ and $H\otimes_KK'$ respectively. By Proposition \ref{descffprop}, the natural morphism 
\[
\Hom_{K-\textup{groups}}(G,H)\rightarrow\Hom_{K'-\textup{groups}}(G\otimes_KK',H\otimes_KK')^\Gamma
\]
is bijective. By the universal property of $\fH'$, the set on the right hand side is naturally identified with $\Hom_{R'-\textup{groups}}(\fG',\fH')^\Gamma=\Hom_{\Z_p\textup{-mod}}(X^*(\fH'),X^*(\fG'))^\Gamma$, which coincides naturally with $\Hom_{\sCon_K}(X^*(G),X^*(H))$, as desired.

To show that $X^*(\cdot)$ is essentially surjective, let us consider a finite Galois extension $K'/K$ with Galois group $\Gamma$ and a finite free $\Z_p$-module $M$ with an action of $\Gamma$. Let $d$ denote the rank of $M$, and let us identify $M$ with $\Z_p^d$. By what we have seen so far, the action of $\Gamma$ on $M$ corresponds contravariantly to an action of $\Gamma$ on the canonical split formal torus of rank $d$. The existence of a strongly concordant $K$-group $G$ such that $X^*(G)\cong M$ now follows from Proposition \ref{semaffdescprop}.
\end{proof}


Let $\mathscr{T}$ be an algebraic $K$-torus. We define its associated strongly concordant $K$-group $T$ by setting
\[
X^*(T)\,\mathrel{\mathop:}=\,X^*(\mathscr{T})\otimes_\Z\Z_p\quad,
\]
where $X^*(\mathscr{T})$ is the character group of $\mathscr{T}$. Explicitly, $T$ is constructed as follows: Let $K'/K$ be a splitting field for $\mathscr{T}$, and let $\ul{\mathscr{T}}'$ be the Néron lft-model of $\mathscr{T}\otimes_KK'$. Let $T'$ be the semi-affinoid generic fiber of the completion of $\ul{\mathscr{T}}'$ along the identity section of its special fiber; then $T$ is the semi-affinoid $K$-group obtained from $T'$ via descent. Indeed, this is clear from the fact that the identity component of $\ul{\mathscr{T}}'$ is isomorphic to a split $R'$-torus.

Similarly, we can define and construct the strongly concordant $K$-group associated to an abelian $K$-variety with potentially split multiplicative reduction. It follows from Proposition \ref{constrnmexprop} and the subsequent remark that the strongly concordant $K$-groups associated to algebraic tori or to abelian varieties with potentially split multiplicative reduction admit formal Néron models and that these formal Néron models are obtained from the respective algebraic Néron models over $R$ by formal completion along closed subschemes of their special fibers.

In the remainder of this section, we prove that all strongly concordant $K$-groups admit formal Néron models, even if they are not induced from \emph{algebraic} groups. We use the existence result for strongly concordant $K$-groups induced by algebraic tori, together with the fact that every strongly concordant $K$-group is a \emph{direct summand} of the strongly concordant $K$-group associated to an algebraic torus.

Let $\Tori_K$ denote the category of algebraic $K$-tori; then $X^*(\cdot)$ induces an anti-equivalence from $\Tori_K$ to the category of finite free $\Z$-modules with continuous $\Gal(K^\sep/K)$-action. In particular, $\Tori_K$ is additive. Let $\Tori_{K,\Z_p}$ denote the $\Z_p$-linear category whose objects are the algebraic $K$-tori and whose morphisms are defined by 
\[
\Hom_{\Tori_{K,\Z_p}}(T,T')\,\mathrel{\mathop:}=\,\Hom_{\Tori_K}(T,T')\otimes_\Z\Z_p\quad.
\]
Let us write $T\mapsto T\otimes_\Z\Z_p$ to denote the natural functor from $\Tori_K$ to $\Tori_{K,\Z_p}$.
\begin{prop}\label{karoubiprop}
The natural functor 
\[
\Tori_{K,\Z_p}\rightarrow\sCon_K
\]
identifies $\sCon_K$ with the Karoubian envelope of $\Tori_{K,\Z_p}$, that is, with the $\Z_p$-linear category of pairs $(T,i)$ where $T$ is an algebraic $K$-torus and where $p$ is an idempotent endomorphism of $T\otimes_\Z\Z_p$.
\end{prop}

This will be a consequence of statement ($ii$) in the following elementary lemma, which is \cite{ChaiPre} 5.8 and which we reprove here for the convenience of the reader:

\begin{lem}\label{ellatticelem}
Let $\Gamma$ be a finite group, and let $V$ be a finite-dimensional linear representation of $\Gamma$ over $\Q_p$. 
\begin{enumerate}
\item There exists a $\Q$-rational finite-dimensional representation $W$ of $\Gamma$ such that $V$ is $\Gamma$-equivariantly isomorphic to a $\Gamma$-stable direct summand of $W\otimes_\Q\Q_p$.
\item Let $V'$ be a $\Gamma$-invariant linear complement of $V$ in $W\otimes_\Q\Q_p$, and let 
$\ul{V}\subseteq V$,
$\ul{V}'\subseteq V'$ be $\Gamma$-stable $\Z_p$-lattices. Then there exists 
a $\Gamma$-stable $\Z$-lattice $\ul{W}\subseteq W$ such that the $\Gamma$-equivariant decomposition $W\otimes_\Q\Q_p\cong V\oplus V'$ restricts to a decomposition
\[
\ul{W}\otimes_\Z\Z_p\cong\ul{V}\oplus\ul{V}'\quad.
\]
\item If $\ul{V}_1\subseteq\ul{V}_2$ are $\Gamma$-stable $\Z_p$-lattices of $V$, there exist $\Gamma$-stable $\Z$-lattices $\ul{W}_1$, $\ul{W}_2$ of $W$ such that
\[
\ul{W}_i\otimes_\Z\Z_p\cong\ul{V}_i\oplus\ul{V}'\quad\textup{for $i=1,2$}
\]
and such that $\ul{W}_1\subseteq\ul{W}_2$.
\end{enumerate}
\end{lem}
\begin{proof}
To prove the first statement, we may assume that $V$ is irreducible. Then $V$ is a direct summand of the regular $\Q_p$-rational representation $\Q_p[\Gamma]$, and so it suffices to take $W=\Q[\Gamma]$. 
Let us prove the second statement. Clearly $\ul{V}\oplus\ul{V}'$ is a $\Gamma$-stable $\Z_p$-lattice in $W\otimes_\Q\Q_p$; we must descend it to a $\Gamma$-stable $\Z$-lattice in $W$. To do so, we choose finite systems of generators $(v_i)_{i\in I}$ and $(v'_i)_{i\in I'}$ of the $\Z_p[\Gamma]$-modules $\ul{V}$ and $\ul{V}'$ respectively. Since $W$ is $p$-adically dense in $W\otimes_\Q\Q_p$, there exist systems $(w_i)_{i\in I}$ and $(w_i')_{i\in I'}$ in $W$ such that 
\begin{eqnarray*}
v_i-w_i\,\in\,p(\ul{V}\oplus\ul{V}')&\forall&i\in I\quad\textup{and}\\
v'_i-w'_i\,\in\,p(\ul{V}\oplus\ul{V}')&\forall&i\in I'\quad.
\end{eqnarray*}
In particular, the $w_i$ and the $w_i'$ are contained in $\ul{V}\oplus\ul{V}'$. Let $\ul{W}\subseteq W$ be the $\Z[\Gamma]$-submodule generated by the $w_i$ and the $w_i'$. Since the inclusion $\ul{W}\otimes_\Z\Z_p\subseteq\ul{V}\oplus\ul{V}'$ reduces to an isomorphism of $\F_p$-vector spaces modulo $p$, Nakayama's Lemma implies that $\ul{W}\otimes_\Z\Z_p=\ul{V}\oplus\ul{V}'$ as $\Z_p$-modules; hence $\ul{W}$ is a $\Z$-lattice of $W$ with the desired property. 

Let us show statement ($iii$). We first use the above construction in order to find a $\Z$-lattice $\ul{W}_2$ of $W$ adapted to $\ul{V}_2$ as in ($ii$). We then use the same strategy to construct $\ul{W}_1$, adapted to $\ul{V}_1$, where we choose the $w_i$ and the $w_i'$ within $\ul{W}_2$; this is possible since $\ul{W}_2\cap (\ul{V}_1\oplus\ul{V}')$ is dense for the $p$-adic topology on $\ul{V}_1\oplus\ul{V}'$. Indeed, by the Artin-Rees Lemma the $p$-adic topology on $\ul{V}_1\oplus\ul{V}'$ coincides with the subspace topology relative to the inclusion $\ul{V}_1\oplus\ul{V}'\subseteq \ul{V}_2\oplus\ul{V}'$, and $\ul{V}_1\oplus\ul{V}'$ is open in $\ul{V}_2\oplus\ul{V}'$ since it contains $p^n(\ul{V}_2\oplus\ul{V}')$ for $n$ big enough. It is now clear that the resulting $\Z$-lattice $\ul{W}_1$ is contained in $\ul{W}_2$.
\end{proof}

The \emph{proof} of Proposition \ref{karoubiprop} is now straightforward: We argue in the dual categories of respective character groups, using statements ($i$) and ($ii$) of Lemma \ref{ellatticelem}. Let $\Gamma$ be a finite group, let $\fC$ denote the additive category of $\Z[\Gamma]$-modules that are finite and free over $\Z$, and let $\fC'$ denote the corresponding category over $\Z_p$. For $\ul{W}\in\fC$, the direct summands of $\ul{W}\otimes_\Z\Z_p$ in the category of $\Z_p[\Gamma]$-modules correspond to the idempotent elements in 
\[
\End_{\Z_p[\Gamma]}(\ul{W}\otimes_\Z\Z_p)\,=\,\End_{\Z[\Gamma]}(\ul{W})\otimes_\Z\Z_p\quad.
\]
The statement is now obvious.\qed

Combining the results of this section with Lemma \ref{proddecomplem}, we obtain the following statement:

\begin{theorem}
Every strongly concordant uniformly rigid $K$-group has a formal Néron model.
\end{theorem}

In his preprint \cite{ChaiPre}, Chai also defines and studies rigid $K$-groups which he calls \emph{concordant}. In terms of their character groups, they correspond to finite free $\Z_p$-modules equipped with a \emph{continuous} $\Gal(K^\sep/K)$-action that is trivial on a finite index subgroup of the inertia subgroup; here continuity is understood with respect to the $p$-adic topology on $\Z_p$. We will not pursue the study of formal Néron models for semi-affinoid concordant groups. Instead, we will explain how formal Néron models for strongly concordant semi-affinoid $K$-groups can be used to derive a formula for the base change conductor of an abelian variety with potentially multiplicative reduction. Being insensitive to unramified base field extensions, the base change conductor can be computed after base change to the maximal unramified extension of $K$, where the notions of concordant and strongly concordant groups coincide.
			
	\section{The base change conductor}\label{chaismethod}
			Let $\sG$ be a smooth $K$-group scheme, and let us assume that the Néron lft-model of $\sG\otimes_KK'$ exists for all finite extensions $K'/K$. Moreover, let us assume that there exists a finite field extension $K'/K$ such that for all finite extensions $K''/K'$, the identity component of the Néron lft-model $\ul{\sG}''$ of $\sG\otimes_KK''$ is obtained from the identity component of the Néron lft-model $\ul{\sG}'$ of $\sG\otimes_KK'$ via base change:
\[
(\ul{\sG}'')^\circ\,=\,(\ul{\sG}')^\circ\otimes_{R'}R''\quad,
\]
where $R'$ and $R''$ are the valuation rings of $K'$ and $K''$ respectively and where $(\cdot)^\circ$ indicates the formation of the identity component, cf. \cite{EGAIV3} 15.6.5. We say that the identity component of the Néron model of $\sG$ \emph{stabilizes} and that $K'/K$ is a \emph{stabilizing field extension} for $\sG$. Let $K'/K$ be a stabilizing field extension for $\sG$, let $e_{K'/K}$ denote its ramification index, and let $\ul{\sG}$, $\ul{\sG}'$ denote the Néron models of $\sG$ and $\sG\otimes_KK'$ over the valuation rings $R$ and $R'$ of $K$ and $K'$ respectively. The natural \emph{base change morphism}
\[
\phi\colon\ul{\sG}\otimes_RR'\rightarrow\ul{\sG}'
\]
induces covariantly a homomorphism 
\[
\phi_*\colon\Lie(\ul{\sG}\otimes_RR')\rightarrow\Lie(\ul{\sG}')
\]
of relative tangent spaces in unit sections. 

\begin{defi}\label{bcconddefi}\index{base change!conductor}
The non-negative rational number
\[
c(\sG,K)\,\mathrel{\mathop:}=\,\frac{1}{e_{K'/K}}\,\length_{R'}\,\coker\phi_*
\]
is called the \emph{base change conductor} of $\sG$.
\end{defi}

Clearly the number $c(\sG,K)$ does not depend on the choice of the stabilizing extension $K'/K$. Let us assume that $K$ has positive residue field characteristic $\textup{char}\,k=p>0$ and that $\sG$ is an abelian $K$-variety with potentially multiplicative reduction. In this section, we sketch Chai's strategy to derive a formula for $c(\sG,K)$ in terms of the Galois action on the character group of the formal torus that is obtained from a stable Néron model $\ul{\sG}'$ via completion along the identity section, in the case where $k$ is perfect. We implement Chai's ideas within the framework of uniformly rigid geometry, thereby illustrating how the theory of formal Néron models for uniformly rigid spaces can profitably be applied to problems in arithmetic algebraic geometry. 

Using the existence of formal Néron models for strongly concordant $K$-groups, we can define the \emph{base change conductor} of a strongly concordant $K$-group $G$ by simply transcribing Definition \ref{bcconddefi}. If $\sG$ is an algebraic $K$-torus or an algebraic $K$-variety with potentially split multiplicative reduction and if $G$ is the strongly concordant $K$-group associated to $\sG$, then 
\[
c(G,K)=c(\sG,K)\quad;
\]
this is clear from the fact that the formal Néron models appearing in the definition of $c(G,K)$ are obtained from the corresponding algebraic Néron models via formal completion. We thus see that the category of strongly concordant groups provides a means of relating the base change conductor of an abelian variety with potentially multiplicative reduction to the base change conductor of an algebraic torus.

\begin{defi}
A morphism $\phi\colon H\rightarrow G$ of strongly concordant $K$-groups is called an \emph{isogeny}\index{isogeny!of concordant groups} if the corresponding homomorphism $\phi^*\colon X^*(G)\rightarrow X^*(H)$ of character groups induces an isomorphism via the functor $\cdot\otimes_{\Z_p}\Q_p$.
\end{defi}

There exists an isogeny $H\rightarrow G$ if and only if there exists an isogeny $G\rightarrow H$. Moreover, isogenies of algebraic $K$-tori and abelian $K$-varieties with potentially split multiplicative reduction induce isogenies of their associated strongly concordant $K$-groups. 

From now on, let us moreover assume that the residue field $k$ of $K$ is \emph{perfect}. Using the highly nontrivial fact (cf.\ \cite{ChaiYudeSh} 11 and 12.1) that the base change conductor of an algebraic torus is an \emph{isogeny invariant}, one proves, cf.\ \cite{ChaiPre} 6.3:

\begin{prop}\label{isoginvprop}
The base change conductor of a strongly concordant $K$-group is an isogeny invariant.
\end{prop}

Indeed, using the corresponding result for algebraic tori, the statement follows from Lemma \ref{ellatticelem} ($iii$) and from the obvious fact that the base change conductor is additive with respect to products. \qed

In particular, it follows that the base change conductor for abelian varieties with potentially multiplicative reduction is an isogeny invariant.

In order to compute base change conductors, we may and do assume that the perfect residue field $k$ of $K$ is algebraically closed; this is due to the fact that Néron models and formal Néron models are preserved under unramified base change. Let us now discuss Chai's explicit formula for the base change conductor of a strongly concordant group. 

For any finite Galois extension $K'/K$ Galois group $\Gamma$, Chai defines a class function $\bA_\Gamma$ on $\Gamma$ with values in a cyclotomic extension of $\Q_p$ such that:

\begin{theorem}\label{chaismainthm}
For every strongly concordant $K$-group $G$ split by $K'$, the base change conductor of $G$ is given by the formula
\[
c(G,K)\,=\,(\bA_\Gamma,\chi_G)\quad,
\]
where $\chi_G$ is the character of $X^*(G)\otimes_{\Z_p}\Q_p$.
\end{theorem}

Here $(\cdot,\cdot)$ is the natural pairing of class functions defined in \cite{Serre} I.2.2 Remark (1).
The above formula applies in particular to the case where $G$ is associated to an abelian $K$-variety with potentially multiplicative reduction. Let us recall the definition and some properties of the class function $\bA_\Gamma$ and afterwards sketch the proof of Theorem \ref{chaismainthm}. Let $(\Gamma_i)_{i\in\Z_{\geq -1}}$ denote the lower numbering ramification filtration of $\Gamma$, and let $i_\Gamma$ be the integer-valued function on $\Gamma$ that is defined by $s\in\Gamma_i\Leftrightarrow i_\Gamma(s)\geq i+1$, cf.\ \cite{SerreLF} IV.1. Since we are in a totally ramified situation, $\Gamma_0=\Gamma$. 

Let $n$ denote the order of the tame Galois group $\Gamma/\Gamma_1$; then $n$ is prime to $p$, and one can canonically identify $\Gamma/\Gamma_1$ with $\mu_n$ in the residue field $k'$ of $K'$. We set $F\mathrel{\mathop:}=K'$ when $K$ is characteristic zero and $F\mathrel{\mathop:}=\Frac(W(k'))$ otherwise, where $W(k')$ denotes the ring of Witt vectors with coefficients in the residue field $k'$ of $K'$. Let $\omega\colon\Gamma/\Gamma_1\rightarrow\mu_n(F)$ be the natural isomorphism. In \cite{ChaiPre} 3.2.1, Chai defines an $F$-valued class function $\bA_\Gamma$ on $\Gamma$ by setting
\[
\bA_\Gamma(s)\,\mathrel{\mathop:}=\,
	\begin{cases}
		\frac{1}{\omega(s)-1}&\textup{if}\,\,s\in\Gamma\setminus\Gamma_1\\
		-\frac{1}{2}i_\Gamma(s)&\textup{if}\,\,s\in\Gamma_1\setminus\{\id_\Gamma\}\\
		\frac{1}{2}	\sum_{t\in\Gamma\setminus\{\id_\Gamma\}}i_\Gamma(t)&\textup{if}\,\,s=\id_\Gamma\quad.
	\end{cases}
\]
In \cite{ChaiPre} 3.2.2, Chai proves that the sum of $\bA_\Gamma$ and its conjugate equals the Artin character for $\Gamma$; he thus calls $\bA_\Gamma$ a \emph{bisection of the Artin character}. In \cite{ChaiPre} 3.3.2, Chai shows that $\bA_\Gamma$ is an element of the rational virtual character ring $R(\Gamma)_\Q$; that is, $\bA_\Gamma$ is a $\Q$-linear combination of characters of $\bar{F}$-rational representations, where $\bar{F}$ is an algebraic closure of $F$. The system $(\bA_\Gamma)_\Gamma$ defines an $\bar{F}$-valued distribution on the space of locally constant $\bar{F}$-valued functions on $\Gal(K^\sep/K)$, cf.\ \cite{ChaiPre} 3.4.1. 

An elementary yet important step in the proof of Theorem \ref{chaismainthm} is the following re\-pre\-sen\-ta\-tion-theoretic argument whose proof relies on Artin's theorem \cite{Serre} Theorem 17; we refer to \cite{ChaiPre} 7.1 for the details:

\begin{lem}\label{inductionlem}
Let $\chi$ be the character of a finite-dimensional $\Q_p$-rational representation of $\Gamma$. Then $\chi$ can be written as a finite sum
\[
\chi\,=\,\sum_{i\in I} a_i\chi_i\quad,
\]
where $a_i\in\Q $ and where the $\chi_i$ are \emph{induced} from finite-dimensional $\Q_p$-rational representations of cyclic subgroups $\Gamma_i\subseteq\Gamma$ such that for each $i\in I$, $\Gamma_i$ has order prime to $p$ or is a $p$-group.
\end{lem}

In \cite{ChaiPre} 3.5.4, the pairing of $\bA_\Gamma$ with an induced character is explicitly computed. To approach the proof of Theorem \ref{chaismainthm} via Lemma \ref{inductionlem}, one must compute how the base change conductor behaves with respect to \emph{induction} on the level of character groups:

Let $L$ be a subfield of $K'/K$, and let $\Gamma_L\subseteq\Gamma$ denote the corresponding subgroup of $\Gamma$. Let $M_L$ be a $\Z_p[\Gamma_L]$-module that is finite free over $\Z_p$, and let $G_L$ denote the corresponding strongly concordant $L$-group which splits over $K'$. Then the induction 
\[
\Ind^\Gamma_{\Gamma_L}M_L\,=\,M_L\otimes_{\Z_p[\Gamma_L]}\Z_p[\Gamma]
\]
is finite free over $\Z_p$; let $\Res_{L/K}G_L$ denote the corresponding strongly concordant $K$-group which splits over $K'$. Let us consider a $\Z[\Gamma_L]$-module $T_L$ that is finite and free over $\Z$ together with a direct sum decomposition
\[
T_{L,\Z_p}\,=\,M_L\oplus M_L'
\]
of $\Z_p[\Gamma_L]$-modules; such data exists by Lemma \ref{ellatticelem}. Then
\[
\Ind^\Gamma_{\Gamma_L}T_L\,=\,T\otimes_{\Z[\Gamma_L]}\Z[\Gamma]
\]
is finite free over $\Z$, and we obtain an induced decomposition of $\Z_p[\Gamma]$-modules
\[
(\Ind^\Gamma_{\Gamma_L}T_L)_{\Z_p}\,=\,\Ind^\Gamma_{\Gamma_L}M_L\oplus \Ind^\Gamma_{\Gamma_L}M'_L\quad.\quad(*)
\]
Let $\mathscr{T}_L$ denote the algebraic $L$-torus corresponding to $T_L$; then the algebraic $K$-torus corresponding to $\Ind^\Gamma_{\Gamma_L}T_L$ is given by the \emph{Weil restriction}
\[
\Res_{L/K}\mathscr{T}_L
\]
of $\mathscr{T}_L$ under $L/K$. Let $\ul{\mathscr{T}_L}$ denote the Néron model of $\mathscr{T}_L$; then the formal Néron model of $(\Ind^\Gamma_{\Gamma_L}T_L)_{\Z_p}$ is obtained from the Néron model
\[
\Res_{R_L/R}\ul{\mathscr{T}_L}
\]
of $\Res_{L/K}\mathscr{T}_L$ via formal completion, where $R_L$ denotes the valuation ring of $L$. The formal Néron model of $\Res_{L/K}G_L$ is thus obtained by splitting off, in this formal completion, the direct factor corresponding to the decomposition $(*)$ above. Using this description, one proves the formula
\[
c(\Res_{L/K}G_L,K)\,=\,c(G_L,L)+\frac{1}{2}v_K(\textup{disc}(L/K))\dim(G_L)\quad,
\]
where $v_K$ is the valuation on $K$, where $\textup{disc}(L/K)$ denotes the discriminant of $L/K$ and where $\dim(G_L)$ is the dimension of $G_L$, cf.\ \cite{ChaiPre} 7.3. By \cite{ChaiPre} 3.5.4, one has the analogous formula for the pairing of $\bA_\Gamma$ with induced characters. By Lemma \ref{inductionlem} and by isogeny invariance of $c$, one can thus reduce to the case where $\Gamma$ is a cyclic $p$-group or a cyclic group of order prime to $p$.

Let $n$ be some natural number. The decomposition of the regular $\Q_p$-rational representation $\Q_p[\Z/p^n\Z]$ is already defined over $\Q$ because the cyclotomic polynomials are defined over $\Q$ and irreducible over $\Q_p$. Hence, if $\Gamma$ is a cyclic $p$-group, one may use isogeny invariance of the base change conductor to reduce to the case where $G$ corresponds to an algebraic torus, and thus the statement of Theorem \ref{chaismainthm} follows from \cite{ChaiYudeSh} 11 and 12.1. Finally, it remains to consider the situation where $\Gamma$ is cyclic of order prime to $p$. Using isogeny invariance of the base change conductor and the fact that $p$ does not divide the order of $\Gamma$, one reduces to the case where $X^*(G)$ is a direct summand of $\Z_p[\Gamma]$. Then $G$ is a direct summand of the strongly concordant group associated to the induced torus $\Res_{K'/K}\G_m$, where $K'/K$ is cyclic and totally \emph{tamely} ramified. This case is settled by means of an explicit computation, cf.\ \cite{Chai} 7.4.1. 

We thus conclude our sketch of the proof of Theorem \ref{chaismainthm}. \qed

\def\MakeUppercase#1{#1}

\renewcommand\bibname{References}

%
%
%

\newpage
\addcontentsline{toc}{chapter}{Index}
\printindex




\addcontentsline{toc}{chapter}{\bibname}

{\small
\begin{thebibliography}{\hspace{0.6in}}
                      \setlength{\itemsep}{0.1ex plus0.05ex minus0.1ex}
                      \setlength{\parsep}{0ex}
                      \setlength{\itemsep}{0ex}


\bibitem[AJL]{AJL} Alonso Tarrío, L.; Jeremías 	López, A.; Lipman, J.: \textit{Duality and flat base change on formal schemes}. Proc.\ Amer.\ Math.\ Soc. 131 (2003), no.\ 2, 351--357.

\bibitem[AJR1]{AJR1} Alonso Tarrío, L.; Jeremías López, A.; Rodríguez, M. P.: \textit{Infinitesimal lifting and Jacobi criterion for smoothness on formal schemes}.\ Comm.\ Algebra 35 (2007), no\. 4, 1341--1367.


\bibitem[AJR2]{AJR3} Alonso Tarrío, L.; Jeremías López, A.; Rodríguez, M. P.: \textit{Local structure theorems for smooth maps of formal schemes}. J.\ Pure Appl.\ Algebra 213 (2009), no\. 7, 1373--1398.


\bibitem[AT]{AT} Artin, E.; Tate, J.: \textit{A note on finite ring extensions}. J.\ Math.\ Soc.\ Japan 3, (1951). 74--77.

\bibitem[Bar]{Bar} Bartenwerfer, W.: \textit{Der Kontinuitätssatz für reindimensionale $k$-affinoide Räume}. Math.\ Ann.\ 193 1971 139--170.

\bibitem[Berk1]{Berk1} Berkovich, V.: \textit{Spectral theory and analytic geometry over non-Archimedean fields.} 
Mathematical Surveys and Monographs, 33. American Mathematical Society, Providence, RI, 1990. x+169 pp.

\bibitem[Berk]{Berk2} Berkovich, V.: \textit{Vanishing cycles for formal schemes II}. Invent.\ Math.\ 125 (1996), no.\ 2, 367--390.

\bibitem[Bert]{Berth} Berthelot, P.: \textit{Cohomologie rigide et cohomologie rigide à supports propres}. Prépublication de l'université de Rennes 1, 1991.

\bibitem[Bo]{FRG} Bosch, S.: \textit{Lectures on formal and rigid geometry}. SFB-Preprint 378, Münster, 2008.

\bibitem[BGR]{BGR} Bosch, S.; Güntzer, U.; Remmert, R.: \textit{Non-Archimedean analysis. A systematic approach to rigid analytic geometry}. Grundlehren der Mathematischen Wissenschaften, 261. Springer-Verlag, Berlin, 1984. xii+436 pp.

\bibitem[BLR]{BLR} Bosch, S.; Lütkebohmert, W.; Raynaud, M.: \textit{Néron models}. Ergebnisse der Mathematik und ihrer Grenzgebiete (3), 21. Springer-Verlag, Berlin, 1990. x+325 pp. 

\bibitem[BL1]{BL1} Bosch, S.; Lütkebohmert, W.: \textit{Formal and rigid geometry. I. Rigid spaces.} Math.\ Ann.\ 295 (1993), no.\ 2, 291--317.

\bibitem[BL2]{BL2} Bosch, S.; Lütkebohmert, W.: \textit{Formal and rigid geometry. II. Flattening techniques.} Math.\ Ann.\ 296 (1993), no.\ 3, 403--429.

\bibitem[BS]{BS} Bosch, S.; Schlöter, K.: \textit{Néron models in the setting of formal and rigid geometry.} Math.\ Ann.\ 301 (1995), no.\ 2, 339--362.

\bibitem[Bou]{Bourbaki} Bourbaki, N.; \textit{Commutative algebra.} Chapters 1--7. Translated from the French. Reprint of the 1989 English translation. Elements of Mathematics (Berlin). Springer-Verlag, Berlin, 1998. xxiv+625 pp. 


\bibitem[Ch1]{Chai} Chai, C.-L.: \textit{Néron models for semiabelian varieties: congruence and change of base field.} Asian J.\ Math.\ 4 (2000), no.\ 4, 715--736.

\bibitem[Ch2]{ChaiPre} Chai, C.-L.: \textit{A bisection of the Artin conductor}, preprint 2003.

\bibitem[ChYu]{ChaiYudeSh} Chai, C.-L.; Yu, J.-K.: \textit{Congruences of Néron models for tori and the Artin conductor.} With an appendix by Ehud de Shalit. Ann.\ of Math.\ (2) 154 (2001), no.\ 2, 347--382.

\bibitem[Co]{ConradIrr} Conrad, B.: \textit{Irreducible components of rigid spaces.} Ann.\ Inst.\ Fourier (Grenoble) 49 (1999), no.\ 2, 473--541.

\bibitem[dJ]{dJ} de Jong, J.: \textit{Crystalline Dieudonné module theory via formal and rigid geometry.} Inst.\ Hautes Études Sci.\ Publ.\ Math.\ No.\ 82 (1995), 5--96 (1996).

\bibitem[dJErr]{dJErr} de Jong, J.: \textit{Erratum to: ``Crystalline Dieudonné module theory via formal and rigid geometry''.} Inst.\ Hautes Études Sci.\ Publ.\ Math.\ No.\ 87 (1998), p.\ 175-175.

\bibitem[Eis]{Eis} Eisenbud, D.: \textit{Commutative algebra. With a view toward algebraic geometry.} Graduate Texts in Mathematics, 150. Springer-Verlag, New York, 1995. xvi+785 pp

\bibitem[Elk]{Elk} Elkik, R.: \textit{Solutions d'équations à coefficients dans un anneau hensélien.} Annales scientifiques de l'É.N.S. 4ième série, tome 6, no 4 (1973), p.\ 553-603

\bibitem[Fa]{FaltingsMac} Faltings, G.: \textit{Über Macaulayfizierung.} Math.\ Ann.\ 238 (1978), no.\ 2, 175--192.

\bibitem[Gr]{Tohoku} Grothendieck, A.: \textit{Sur quelques points d'algèbre homologique.} Tôhoku Math.\ J.\ (2) 9 1957 119--221.

\bibitem[EGA I]{EGAI} Grothendieck, A.: \textit{Éléments de géométrie algébrique. I. Le langage des schémas.} Inst.\ Hautes Études Sci.\ Publ.\ Math.\ No.\ 4 1960 228 pp.

\bibitem[EGA I$_\textup{n}$]{EGAIn} Grothendieck, A.; Dieudonné, J.: \textit{Éléments de géométrie algébrique. I.} Die Grundlehren der mathematischen Wissenschaften in Einzeldarstellungen, 166. Springer-Verlag, Berlin, 1971. ix+466 pp.

\bibitem[EGA II]{EGAII} Grothendieck, A.: \textit{Éléments de géométrie algébrique. II. Étude globale élémentaire de quelques classes de morphismes.} Inst.\ Hautes Études Sci.\ Publ.\ Math.\ No.\ 8 1961 222 pp.

\bibitem[EGA III$_1$]{EGAIII1} Grothendieck, A.: \textit{Éléments de géométrie algébrique. III. Étude cohomologique des faisceaux cohérents. I.} Inst.\ Hautes Études Sci.\ Publ.\ Math.\ No.\ 11 1961 167 pp.


\bibitem[EGA IV$_1$]{EGAIV1} Grothendieck, A.: \textit{Éléments de géométrie algébrique. IV. Étude locale des schémas et des morphismes de schémas. I.} Inst.\ Hautes Études Sci.\ Publ.\ Math.\ No.\ 20 1964 259 pp.

\bibitem[EGA IV$_2$]{EGAIV2} Grothendieck, A.: \textit{Éléments de géométrie algébrique. IV. Étude locale des schémas et des morphismes de schémas. II.} Inst.\ Hautes Études Sci.\ Publ.\ Math.\ No.\ 24 1965 231 pp.


\bibitem[EGA IV$_3$]{EGAIV3} Grothendieck, A.: \textit{Éléments de géométrie algébrique. IV. Étude locale des schémas et des morphismes de schémas. III.} Inst.\ Hautes Études Sci.\ Publ.\ Math.\ No.\ 28 1966 255 pp.

\bibitem[EGA IV$_4$]{EGAIV4} Grothendieck, A.: \textit{Éléments de géométrie algébrique. IV. Étude locale des schémas et des morphismes de schémas IV.} Inst.\ Hautes Études Sci.\ Publ.\ Math.\ No.\ 32 1967 361 pp.



\bibitem[H1]{Huberhabil} Huber, R.: \textit{Bewertungsspektrum und rigide Geometrie.} Regensburger Mathematische Schriften, 23. Universität Regensburg, Fachbereich Mathematik, Regensburg, 1993. xii+309 pp.

\bibitem[H2]{Hubermainadspacepaper} Huber, R.: \textit{A generalization of formal schemes and rigid analytic varieties.} 
Math.\ Z.\ 217 (1994), no.\ 4, 513--551. 

\bibitem[H3]{Huberbuch} Huber, R.: \textit{Étale cohomology of rigid analytic varieties and adic spaces.} Aspects of Mathematics, E30. Friedr.\ Vieweg \& Sohn, Braunschweig, 1996. x+450 pp.

\bibitem[H4]{Hubercohpaper} Huber, R.: \textit{A finiteness result for the compactly supported cohomology of rigid analytic varieties.} II. Ann.\ Inst.\ Fourier (Grenoble) 57 (2007), no.\ 3, 973--1017.

\bibitem[Kie1]{KiehlAB} Kiehl, R.: \textit{Theorem A und Theorem B in der nichtarchimedischen Funktionentheorie.}
Invent.\ Math.\ 2 1967 256--273. 

\bibitem[Kie2]{KiehldeRh} Kiehl, R.: \textit{Die de Rham Kohomologie algebraischer Mannigfaltigkeiten über einem bewerteten Körper.} Inst.\ Hautes Études Sci.\ Publ.\ Math.\ No.\ 33 1967 5--20. 


\bibitem[Kis]{Kisin} Kisin, M.: \textit{Local constancy in $p$-adic families of Galois representations.} Math.\ Z.\ 230 (1999), no.\ 3, 569--593.

\bibitem[Li]{Liu} Liu, Q.; \textit{Algebraic geometry and arithmetic curves.} Translated from the French by Reinie Erné. Oxford Graduate Texts in Mathematics, 6. Oxford Science Publications. Oxford University Press, Oxford, 2002. xvi+576 pp.

\bibitem[LR]{LiRo} Lipshitz, L.; Robinson, Z.: \textit{Rings of separated power series and quasi-affinoid geometry.} Astérisque No.\ 264 (2000), vi+171 pp. 


\bibitem[Lü1]{Ltke} Lütkebohmert, W.: \textit{Fortsetzbarkeit $k$-meromorpher Funktionen.} Math.\ Ann.\ 220 (1976), no.\ 3, 273--284. 

\bibitem[Lü2]{LtkeFRG} Lütkebohmert, W.: \textit{Formal-algebraic and rigid-analytic geometry.} Math.\ Ann.\ 286 (1990), no.\ 1-3, 341--371.

\bibitem[Lu]{Lubin} Lubin, J.: \textit{One-parameter formal Lie groups over $p$-adic integer rings.} Ann.\ of Math.\ (2) 80 1964 464--484.


\bibitem[Ma]{MatsumuraCRT} Matsumura, H.: \textit{Commutative ring theory.} Translated from the Japanese by M.\ Reid. Second edition. Cambridge Studies in Advanced Mathematics, 8. Cambridge University Press, Cambridge, 1989. xiv+320 pp.

\bibitem[Ni]{N} Nicaise, J.; \textit{A trace formula for rigid varieties, and motivic Weil generating series for formal schemes.} 
Math.\ Ann.\ 343 (2009), no.\ 2.

\bibitem[RZ]{RZ} Rapoport, M.; Zink, T.: \textit{Period spaces for $p$-divisible groups.} Annals of Mathematics Studies, 141. Princeton University Press, Princeton, NJ, 1996. xxii+324 pp.

\bibitem[Se1]{SerreLF} Serre, J.-P.: \textit{Corps locaux.} Deuxième édition. Publications de l'Université de Nancago, No.\ VIII. Hermann, Paris, 1968. 245 pp.

\bibitem[Se2]{Serre} Serre, J.-P.: \textit{Linear representations of finite groups.} Translated from the second French edition by Leonard L.\ Scott. Graduate Texts in Mathematics, Vol.\ 42. Springer-Verlag, New York-Heidelberg, 1977.


\bibitem[St]{S} Strauch, M.: \textit{Deformation spaces of one-dimensional formal modules and their cohomology.} Adv.\ Math.\ 217 (2008), no.\ 3, 889--951.

\bibitem[Sw]{Sw} Swan, R. G.: \textit{Néron-Popescu desingularization.} Algebra and geometry (Taipei, 1995), 135--192, 
Lect.\ Algebra Geom., 2, Int.\ Press, Cambridge, MA, 1998. 

\bibitem[Te]{T} Temkin, M.: \textit{Desingularization of quasi-excellent schemes in characteristic zero.}
Adv.\ Math.\ 219 (2008), no.\ 2, 488--522. 

\bibitem[V1]{V1} Valabrega, P.: \textit{On the excellent property for power series rings over polynomial rings.} J.\ Math.\ Kyoto Univ.\ 15 (1975), no.\ 2, 387--395.

\bibitem[V2]{V2} Valabrega, P.: \textit{A few theorems on completion of excellent rings. Nagoya Math.} J.\ 61 (1976), 127--133.

\bibitem[Wei]{Weil} Weil, A.: \textit{Variétés abéliennes et courbes algébriques.} Hermann, Paris (1948), republished in Courbes Algébriques et Variétés Abéliennes. Hermann, Paris (1971).

\bibitem[Weg]{Wegel} Wegel, O.: \textit{Das Verhalten formeller Néron-Modelle unter Basiswechsel.} Dissertation, Münster 1997.

\bibitem[Ye]{Yek} Yekutieli, A.: \textit{Smooth formal embeddings and the residue complex.} Canad.\ J.\ Math.\ 50 (1998), no.\ 4, 863--896.

\end{thebibliography}
}
\end{document}